\documentclass[11pt]{article}

\usepackage{geometry}
\usepackage{comment} 
\excludecomment{allproofs} 

\usepackage[utf8]{inputenc}
\usepackage{amsthm}
\usepackage{graphicx}
\usepackage{amssymb}
\usepackage{color}
\usepackage[dvipsnames]{xcolor}
\usepackage{tikz}
\usepackage{tikz-cd}
\usetikzlibrary{shapes.geometric}
\usepackage[many]{tcolorbox} 
\tcbuselibrary{theorems}
\usepackage{amsmath}
\usepackage[hidelinks]{hyperref}
\usepackage{float}
\usepackage{mathtools}
\usepackage{enumitem}
\usepackage{subcaption}

\usepackage{bbm} 

\theoremstyle{plain}
\newtheorem{theorem}{Theorem}[section]
\newtheorem{lemma}[theorem]{Lemma}

\newtheorem{proposition}[theorem]{Proposition}
\newtheorem{corollary}[theorem]{Corollary}
\newtheorem{claim}[theorem]{Claim}

\newtheorem{problem}[theorem]{Problem}

\theoremstyle{definition}
\newtheorem{definition}[theorem]{Definition}
\newtheorem{example}[theorem]{Example}

\theoremstyle{remark}
\newtheorem{remark}[theorem]{Remark}

\theoremstyle{plain}
\newtheorem*{theorem*}{Theorem}
\newtheorem*{lemma*}{Lemma}

\counterwithin{equation}{section}

\newcommand{\RR}{\mathbb{R}}
\newcommand{\CC}{{\mathbb C}}
\newcommand{\HH}{{\mathbb H}}
\newcommand{\DD}{{\mathbb D}}

\newcommand{\fcal}{\mathcal{F}}

\newcommand{\ical}{\mathcal{I}}

\newcommand{\pcal}{\mathcal{P}}

\def\id{{\operatorname{id}}}
\def\tr{\hbox{\textrm{tr}}}

\newcommand{\diag}{{\operatorname{diag}}}
\def\Sym{\operatorname{Sym}}
\def\Herm{\operatorname{Herm}}
\def\QHerm{\operatorname{QHerm}}
\def\SHS{\operatorname{SHS}}

\def\intr{\operatorname{int}}
\def\del{\partial}
\def\dbar{\overline{\partial}}
\def\ddbar{\del \dbar}
\def\dbard{\dbar \del}
\def\II{\mathrm{II}}

\def\signdist{\sigma}

\usepackage[style=numeric, giveninits=true, sortcites=true, backend=biber, maxbibnames=99]{biblatex} 
\AtEveryBibitem{\clearfield{isbn}}
\AtEveryBibitem{\clearfield{issn}}
\AtEveryBibitem{\clearfield{doi}}
\AtEveryBibitem{\clearfield{url}}
\AtEveryBibitem{\clearfield{number}}
\AtEveryBibitem{\clearfield{series}}

\DeclareFieldFormat{edition}{#1~ed\adddot}

\DeclareNameAlias{default}{given-family}

\DefineBibliographyExtras{english}{%
}

\renewbibmacro{in:}{
}

\DeclareFieldFormat*{title}{\mkbibemph{#1}}

\renewbibmacro*{title}{%
  \printfield{title}%
  \ifentrytype{misc}{\addcomma\addspace preprint}{}%
}

\DeclareFieldFormat{pages}{#1}

\DeclareFieldFormat*{journaltitle}{\textnormal{#1}}

\DeclareFieldFormat*{volume}{\mkbibbold{#1}}

\title{
Regularity via subequations
}
\author{Vasanth Pidaparthy, Yanir A. Rubinstein}
\date{\today}

\begin{document}

\maketitle

\font\itsml=cmti7




\noindent
\hfill {\itsml ``Today's science is addicted to more and more bursts of energy.

\vglue-4pt
\hfill
But what really makes humanity smarter and which we as society 

\vglue-4pt
\noindent
\hfill
should value is, instead, opportune drops in entropy.'' --- T. Tokieda}


\begin{abstract}
Is there a unified theory of a priori estimates for
purely second-order constant coefficient elliptic PDE?
A strong indication that such a theory exists is provided by
the foundational a priori estimates of Caffarelli--Nirenberg--Spruck (1984) and Trudinger (1995) that treated many equations, yet left out many others.


This article develops such a unified theory based on Harvey--Lawson's theory of subequations, the aforementioned a priori estimates, and some new ingredients. 
Serendipitously, this also solves a problem posed by Nirenberg (2007).
In addition to unifying all previously known such estimates, we find a vast number of new a priori estimates (e.g., spanning real, complex, quaternionic, $G_2$, among other natural equations, and even differential inclusions that generalize PDEs). 


\end{abstract}

\tableofcontents

\section{Introduction}

This article is concerned with regularity for Dirichlet boundary value problems (BVP) for purely second-order constant coefficient elliptic PDE,
\begin{equation}
    \label{eq:Subeq-Reg-intro-purely-second-order-PDE}
    \begin{aligned}
        \fcal(D^2u)=0& \qquad \text{on } \Omega,\\
        u = \varphi& \qquad \text{on } \partial \Omega,
    \end{aligned}
\end{equation}
for a {bounded domain} $\Omega \subset \RR^n$, a boundary condition $\varphi \in C^\infty(\partial \Omega)$, and a smooth partial differential operator $\fcal \in  C^\infty(\Sym^2(\RR^n))$ that is elliptic,
\begin{equation}
    \label{eq:Subeq-Reg-intro-purely-second-order-PDE-elliptic}
    \fcal(A)\ge \fcal(B), \qquad \forall\ A\ge B.
\end{equation}

\begin{problem}
    \label{problem:Subeq-Reg-main-problem-about-general-fcal-pde}
    Under what assumptions on $\fcal$ does \eqref{eq:Subeq-Reg-intro-purely-second-order-PDE} admit:

    \vglue4pt
    \noindent (i) weak solutions $u\in C^0(\overline{\Omega})$? (ii) classical solutions $u\in C^2(\overline{\Omega})$?
\end{problem}

As we describe in \S\ref{subsection:Subeq-Reg-intro-weak-solutions}, Problem~\ref{problem:Subeq-Reg-main-problem-about-general-fcal-pde}~(i) is well-understood, with the most satisfactory result due to Harvey--Lawson in 2009 (Theorem~\ref{thm:Subeq-Reg-Har-Law-exist-thm}). Problem~\ref{problem:Subeq-Reg-main-problem-about-general-fcal-pde}~(ii) has been the topic of intense research for the past 50 years, with foundational results by Caffarelli--Nirenberg--Spruck in the 1980's \cite{caffarelli-nirenberg-spruck-I-MR0739925, CKNS-II-MR780073, caffarelli-nirenberg-spruck-dirichlet-problem-III-MR0806416}, Guan--Spruck and Trudinger in 1990's \cite{Guan-Spruck-adm-subsol-MR1247995, guan-subsol-MR1265805, Trudinger1995-DHE}, among many others. As explained in \S\ref{StrSubSec} (and Remark \ref{remark:Subeq-Reg-subequations-and-operators}), a more natural (and  general) version of 
Problem~\ref{problem:Subeq-Reg-main-problem-about-general-fcal-pde}~(ii) has been posed by Nirenberg (Problem~\ref{problem:Subeq-Reg-for-which-subeq-regularity}), so we, by slight abuse, refer to both
Problem~\ref{problem:Subeq-Reg-main-problem-about-general-fcal-pde}~(ii)
and Problem~\ref{problem:Subeq-Reg-for-which-subeq-regularity} as
 {\it Nirenberg's Problems}.

The main result of the present article is, as we will try to convince the reader, a fairly satisfactory resolution of Nirenberg's Problems,  
stated in Theorem~\ref{thm:Subeq-Reg-a-priori-estimates} and Corollary~\ref{corollary:Subeq-Reg-existence-thm}. Among other things, it unifies essentially all known a priori estimates for equations of type \eqref{eq:Subeq-Reg-intro-purely-second-order-PDE} as well as establishes new a priori estimates for many equations. Furthermore, the class of equations it covers is essentially optimal.

\subsection{Background: the story for weak solutions}
\label{subsection:Subeq-Reg-intro-weak-solutions}

Weak solutions to non-linear elliptic PDE have a long history. One of the earliest notions of a weak solution to \eqref{eq:Subeq-Reg-intro-purely-second-order-PDE} was the Alexandrov solution to the Monge--Amp\`ere equation, when $\fcal = \det$ \cite{Alexandrov-book-MR71041, Rauch-Taylor-MR454331}. Weak solutions to the BVP~\eqref{eq:Subeq-Reg-intro-purely-second-order-PDE} for general elliptic operators $\fcal$ were pioneered in the work of Crandall, Evans, Ishii, Lions and others~\cite{CIL-MR1118699}.

Harvey--Lawson~\cite{HarLaw-dirichletduality} (and others before, e.g., Krylov~\cite{krylov-geneal-notion-of-elliptic-PDE-MR1284912}) observed that the PDE \eqref{eq:Subeq-Reg-intro-purely-second-order-PDE} does not depend on the full operator $\fcal$, but only on the level set $\{\fcal =0 \}$, a hypersurface in $\Sym^2(\RR^n)$, and that, furthermore, $\{\fcal=0\}$ and $\{\fcal\ge 0\}$ determine each other.

\begin{definition}[Subequation]
\label{def:Subeq-Reg-subequations}
    {\textrm{(Harvey--Lawson~\cite[{\S}3]{HarLaw-dirichletduality})}} A \textit{subequation} is a closed, non-empty, proper subset of all symmetric matrices $F\subsetneq \Sym^2(\RR^n)$, satisfying
    \begin{equation}
    \label{eq:Subeq-Reg-subequation-def-prop-basic}
    F+ \pcal \subset F,
    \end{equation}
    where
\begin{equation}
    \label{eq:Subeq-Reg-positive-semi-def-subeq}
    \pcal := \left\{P \in \Sym^2(\RR^n): P \ge 0 \right\}.
\end{equation}
\end{definition}

Subequations generalize elliptic operators in that $\fcal$ determines a subequation $\{\fcal\ge 0\}$
(cf. Remark~\ref{remark:Subeq-Reg-subequations-and-operators}). Thus the BVP \eqref{eq:Subeq-Reg-intro-purely-second-order-PDE} is a special case of the \textit{$F$-BVP},
\begin{equation}
    \label{eq:Subeq-Reg-const-coeff-subequation-PDE}
    \begin{aligned}
        D^2u \in \partial F& \qquad \text{on } \Omega,\\
        u = \varphi& \qquad \text{on } \partial \Omega.
    \end{aligned}
\end{equation}
A key idea in their approach is identifying a \textit{dual subequation}
\begin{equation}
    \label{eq:Subeq-Reg-dual-to-a-subequation}
    \widetilde{F} := -\overline{\Sym^2(\RR^n) \setminus F}.
\end{equation}
The dual is also a subequation \cite[p.~13, (6)]{HarLaw-dirichletduality}, and since $\partial F = F\cap -\widetilde{F}$, \eqref{eq:Subeq-Reg-const-coeff-subequation-PDE} can equivalently be phrased in terms of the dual.

\begin{theorem}
    \label{thm:Subeq-Reg-Har-Law-exist-thm}
    {\rm{(Harvey--Lawson~\cite[Theorem~6.2]{HarLaw-dirichletduality})}} For a bounded domain $\Omega \subset \RR^n$ with $C^2$ boundary $\partial \Omega$ that is both strictly $\vec{F}$- and $\vec{\widetilde{F}}$-convex, and $\varphi \in C^0(\partial \Omega)$, there exists a unique weak solution $u\in C^0(\overline{\Omega})$ to \eqref{eq:Subeq-Reg-const-coeff-subequation-PDE}.
\end{theorem}

Theorem~\ref{thm:Subeq-Reg-Har-Law-exist-thm} answers Problem~\ref{problem:Subeq-Reg-main-problem-about-general-fcal-pde}~(i). The hypotheses on 
$\Omega$ (in terms of $F$) are mild, and essentially optimal. This captures a vast class of equations including, e.g., all branches of the Monge--Amp\`ere, $k$-Hessian, and special Lagrangian equations. The conditions on the boundary $\partial \Omega$ are naturally associated to a subequation $F$ \cite[\S5]{HarLaw-dirichletduality}, and in particular are satisfied for uniformly convex domains (see \S\ref{subsection:Subeq-Reg-unif-convex-domains}). An in-depth treatment can be found in the textbooks \cite{Payne-Redaelli-MR4952961, Cirant-Payne-Harvey-Lawson-MR4671835, Subequation-book-I}.

\subsection{Results: the story for strong solutions}
\label{StrSubSec}

Our goal is to solve Problem~\ref{problem:Subeq-Reg-main-problem-about-general-fcal-pde}~(ii). A key observation is that a more general (and natural) problem is Nirenberg's problem~\cite[p. 97]{Nirenberg-note-MR2641120}:

\begin{problem}[Nirenberg]
\label{problem:Subeq-Reg-for-which-subeq-regularity}
   When is the weak solution from Theorem~\ref{thm:Subeq-Reg-Har-Law-exist-thm} smooth?
\end{problem}

Note that Theorem~\ref{thm:Subeq-Reg-Har-Law-exist-thm} places no restriction on the subequation $F$. By contrast, the gist of Problem~\ref{problem:Subeq-Reg-for-which-subeq-regularity} is to find the `right' class of subequations as there are many examples of subequations $F$ for which the BVP \eqref{eq:Subeq-Reg-const-coeff-subequation-PDE} does \textit{not} admit classical solutions (Example~\ref{example:Subeq-Reg-examples-for-optimality}).

\begin{definition}
    \label{def:Subeq-Reg-intro-non-deg-subequation}
    A subequation $F\subset \Sym^2(\RR^n)$ is \textit{non-degenerate} if
    \begin{equation}
        \label{eq:Subeq-Reg-non-degenerate-subequation}
        F+\pcal \setminus \{0\} \subset \intr F.
    \end{equation}
\end{definition}
Somewhat surprisingly, this very natural class of subequation does  not seem to have been studied  before systematically. When $F$ is convex, the condition \eqref{eq:Subeq-Reg-non-degenerate-subequation} is analogous to strict ellipticity of a partial differential operator $\fcal$ (Lemma~\ref{lemma:Subeq-Reg-subequation-is-elliptic}).

This article provides an answer to Problem~\ref{problem:Subeq-Reg-for-which-subeq-regularity}.
In essence, the `right' conditions on $F$ are smoothness, convexity, non-degeneracy, and, finally, invariance under a Lie group $G\subset O(n)$ acting transitively on $S^{n-1}$:

\begin{theorem}[A priori estimates]
\label{thm:Subeq-Reg-a-priori-estimates}
        For a Lie group $G\subset O(n)$ acting transitively on $S^{n-1}$, a convex, non-degenerate, $G$-invariant subequation $F \subset \Sym^2(\RR^n)$ with $C^2$ boundary $\partial F$, a bounded domain $\Omega \subset \RR^n$ with $C^4$ boundary $\partial \Omega$ that is strictly $\vec{F}$-convex, and a boundary condition $\varphi \in C^4(\overline{\Omega})$, there exists a constant $C = C\left( n, \Omega, F, \|\varphi \|_{C^4} \right) >0$ such that any solution $u\in C^4( \overline{ \Omega}) $ to the Dirichlet problem \eqref{eq:Subeq-Reg-const-coeff-subequation-PDE} satisfies
        \[\|u\|_{C^2(\overline{\Omega})} \le C.\]
\end{theorem}
Here $G$-invariance of $F$ means invariance under the conjugation action of $G\subset O(n)$ on $\Sym^2(\RR^n)$
\begin{equation}
    \label{eq:Subeq-Reg-group-invariant-subequation}
    O^T FO = F, \qquad \forall\ O\in G.
\end{equation}
Theorem~\ref{thm:Subeq-Reg-a-priori-estimates} is proved in \S\ref{section:Subeq-Reg-proof-of-main-a-priori-theorem}. The existence of a classical solution then follows from Theorem~\ref{thm:Subeq-Reg-a-priori-estimates} by standard arguments, i.e., the Evans--Krylov theory~\cite{Evans-C2alpha-est-1982-MR649348, Krylov-C2alpha-domain-bdry-MR688919} (see \S\S\ref{subsection:Subeq-Reg-evans-krylov-theorem}--\ref{subsection:Subeq-Reg-schauder-estimates}) and the continuity method adapted to this setting (see \S\ref{section:continuity-method}):

\begin{corollary}[Solution to Problem~\ref{problem:Subeq-Reg-for-which-subeq-regularity}]
\label{corollary:Subeq-Reg-existence-thm}
    Under the assumptions of Theorem~\ref{thm:Subeq-Reg-a-priori-estimates}, if $\partial F$ is a $C^3$, $\partial \Omega$ is $C^5$, and $\varphi \in C^5(\overline{\Omega})$, then the BVP \eqref{eq:Subeq-Reg-const-coeff-subequation-PDE} has a unique solution $u \in C^{4,\alpha}(\overline{\Omega})$ for some $\alpha = \alpha\left(n, F, \Omega, \|\varphi\|_{C^4}  \right) \in (0,1)$.
\end{corollary}

\begin{remark}[What equations does this result cover?]~
\label{remark:Subeq-Reg-which-equations-are-covered}
    \begin{enumerate}[label={(\roman*)}]
        \item Theorem~\ref{thm:Subeq-Reg-a-priori-estimates} (and Corollary~\ref{corollary:Subeq-Reg-existence-thm}) simultaneously cover essentially all known a priori $C^2$ estimates and existence theorems for \eqref{eq:Subeq-Reg-const-coeff-subequation-PDE}, including, e.g., the top branch of the real Monge--Amp\`ere equation, the $k$-Hessian equation, the top two branches of the special Lagrangian equation \eqref{eq:Subeq-Reg-subequations-for-SL}.

        \item Theorem~\ref{thm:Subeq-Reg-a-priori-estimates} can also be applied to so called generalized to equations determined by pairs of subequations $(F_1,F_2)$ where either existence or uniqueness might fail (see \S\ref{section:Subeq-Reg-equations-determined-by-pairs-of-subequations}). Interesting examples of such generalized equations are the twisted Monge--Amp\`ere equations studied by Streets--Warren~\cite{Streets-Warren-MR3498927, Har-Law-PDE-Krylov-MR4127162}.
        
        \item Theorem~\ref{thm:Subeq-Reg-a-priori-estimates} applies to equations like those discussed in Example~\ref{example:Subeq-Reg-fails-CNS-for-all-translations} which do not satisfy some assumptions of Caffarelli--Nirenberg--Spruck~\cite[Theorem~2]{caffarelli-nirenberg-spruck-dirichlet-problem-III-MR0806416}, and Trudinger~\cite[Theorem~1.2]{Trudinger1995-DHE}.

        \item The Lie groups $G\subset O(n)$ that act transitively on the unit sphere $S^{n-1}$ include $SO(n)$, $U(m)$ and $SU(m)$ when $n=2m$, $Sp(m)$ when $n=4m$, and the exceptional groups $G_2$ when $n=7$, $Spin(7)$ when $n=8$, and $Spin (9)$ when $n=16$ \cite[Theorem~59]{Petersen-Riem-Geom-MR3469435}. To the best of our knowledge, Theorem~\ref{thm:Subeq-Reg-a-priori-estimates} is the first a priori $C^2$ estimate on domains for the exceptional groups $ G_2$, $Spin(7)$ and $Spin(9)$, though, e.g., Alesker--Gordon~\cite{Alesker-Gordon-MR4771382} have studied a $Spin(9)$-invariant octonionic Monge--Amp\`ere equation on compact $GL_2(\mathbb{O})$-affine manifolds.
        
        The only known a priori $C^2$ estimates for $Sp(m)$-invariant operators are those for the quaternionic Monge--Amp\`ere equation, due to Alesker and Zhu \cite{Alesker-quaternionic-MA-MR1967025, Zhu-quaternionic-MA-MR3544695}. Theorem~\ref{thm:Subeq-Reg-a-priori-estimates} generalizes far beyond this operator to the BVP \eqref{eq:Subeq-Reg-intro-purely-second-order-PDE} for concave strictly elliptic quaternionic PDE:
        \[\begin{aligned}
            \fcal(D^2u) = f(\lambda([D^2u]^{\QHerm})) =1 \qquad& \text{on } \Omega,\\ 
            u = \varphi \qquad& \text{on } \partial \Omega,
        \end{aligned}\]
        where $A^{\QHerm}$ is the quaternionic Hermitian part of any $A\in \Sym^2(\RR^{4m})$ \eqref{eq:Subeq-Reg-quaternionic-hermitian-part-of-a-matrix}, $\lambda(A^{\QHerm}) \in \RR^m$ is the vector of eigenvalues of $A^{\QHerm}$ (note each eigenvalue has multiplicity $4$ since $A^{\QHerm}$ is quaternionic Hermitian), and $f$ is a $C^2$, permutation invariant, concave function that is strictly elliptic, $\partial f/\partial\lambda_i >0$ for $i=1,\dots, m$.
        
        Even in the $U(m)$-invariant case, the typically studied equations have the form $\fcal(\ddbar u) = 0$ where $\ddbar u$ denotes the complex Hessian. Theorem~\ref{thm:Subeq-Reg-a-priori-estimates} generalizes a priori estimates to operators $\fcal(D^2u)$ that can depend on the non-Hermitian part of $D^2u$, as long as they are $U(m)$-invariant, e.g., $\fcal \in C^2(\intr \pcal)$ (see \eqref{eq:Subeq-Reg-positive-semi-def-subeq}),
        \[\fcal(A) := \log \det \left( A^{\Herm} + \frac{1}{2} A^\perp \right) , \qquad \forall\ A\in \intr \pcal,\]
        where $A^{\Herm} = (A - JAJ)/2 \in \Sym^2(\RR^{2m})$ is the Hermitian part and $A^{\perp} = (A+JAJ)/2 \in \Sym^2(\RR^{2m})$ is the skew-Hermitian part \eqref{eq:Subeq-Reg-hermitian-part-of-a-matrix}.

        \item Theorem~\ref{thm:Subeq-Reg-a-priori-estimates} also applies to the $G$-invariant G\r{a}rding--Dirichlet (GD) polynomials (Theorem~\ref{thm:Subeq-Reg-existence-for-GD-polynomials}). $C^0$ estimates were studied by Harvey--Lawson~\cite{HarLaw-dirichletduality, Har-Law-det-majorization-MR4589710, HL-duke-paper-MR4965192}. While $C^2$ estimates for the associated BVP \eqref{eq:Subeq-Reg-BVP-for-GD-polynomials} follow for $O(n)$- and $U(m)$-invariant GD polynomials  by Caffarelli--Nirenberg--Spruck~\cite[Theorems~2]{caffarelli-nirenberg-spruck-dirichlet-problem-III-MR0806416}, Trudinger~\cite[Theorem~1.2]{Trudinger1995-DHE}, and S.-Y.~Li~\cite[Theorem~1.1]{Li-song-ying-MR2128299}, to the best of our knowledge, Theorem~\ref{thm:Subeq-Reg-a-priori-estimates} provides the first $C^2$ estimate for the $Sp(m)$-invariant quaternionic GD polynomials (see \S\ref{subsec:Subeq-Reg-GD-thm-applied-to-CMA-and_QMA}).
    \end{enumerate}
\end{remark}

\begin{remark}[Why are the results (essentially) optimal?]~
\label{remark:Subeq-Reg-why-are-results-optimal}
    \begin{enumerate}[label={(\roman*)}]
        \item (Convexity) The subequations $F_c$ \eqref{eq:Subeq-Reg-subequations-for-SL} for the special Lagrangian equation are smooth, non-degenerate, and $O(n)$-invariant, but \textit{non}-convex for $|c|<(n-2)\pi/2$. And the $F_c$-BVP \eqref{eq:Subeq-Reg-const-coeff-subequation-PDE} admits solutions $u \in C^{1,\alpha}(D) \setminus C^{1,1}(D)$ (see Example~\ref{example:Subeq-Reg-examples-for-optimality}~(i)).

        \item (Non-degeneracy) The $\pcal$-BVP \eqref{eq:Subeq-Reg-homogeneous-RMA-intro} with smooth boundary condition $\varphi$ admits a solution $u \notin C^{1,1} ({\overline{\Omega}})$. In earnest, $\pcal$ \eqref{eq:Subeq-Reg-positive-semi-def-subeq} is convex and $O(n)$-invariant, but neither non-degenerate nor smooth. However, in dimension $n=2$, it is possible to smooth out  $\pcal$ in a neighborhood of the origin (see Figure~\ref{fig:Subeq-Reg-smoothing-of-pcal-near-the-origin}), to construct a convex, $O(n)$-invariant and smooth (but degenerate) subequation $\pcal'$ admitting a non-classical solution (see Example~\ref{example:Subeq-Reg-examples-for-optimality}~(ii)).

        \item ($G$-invariance): For $B:=B(0,1)\subset \RR^3$ denoting the unit ball, Mooney--Savin~\cite{Mooney-Savin-u11u22-MR4026172} constructed  $u\in C^{1,1/3}(B)$ solving the $F$-BVP \eqref{eq:Subeq-Reg-const-coeff-subequation-PDE} for the subequation
        \[F := \left\{A \in \Sym^2(\RR^3): a_{11} \ge0,\quad a_{22} \ge0,\quad a_{33} \ge 0,\quad a_{11} a_{22} a_{33} \ge 1\right\}.\]
        $F$ is smooth, non-degenerate and convex, but is neither $O(3)$- nor $SO(3)$-invariant. Note that the constructed solution has the same regularity on the boundary, $u|_{\partial B}\in C^{1,1/3}(\partial B)$. It would be interesting to find a counter-example, with better boundary regularity.
        
        \item The degenerate special Lagrangian (DSL) equation is a degenerate elliptic PDE on a product domain $\Omega = (0,1)\times D$, representing the geodesic equation in the space of positive Lagrangians~\cite{Rub-Sol-DSL-MR3620699}. It corresponds to a subequation $\fcal_c$~\cite[(35)]{Rub-Sol-DSL-MR3620699} that is convex (for $c\ge (n-1)\pi/2$)~\cite{VP-YAR-convexity} but not smooth, $G$-invariant, or non-degenerate. It admits weak solutions $u\in C^0(\overline{\Omega})$~\cite[Theorem~1.1]{Rub-Sol-DSL-MR3620699}, with $u\notin C^2(\overline{\Omega}) $ \cite{Kapota-P-R-Solomon-in-prep}.
    \end{enumerate}
\end{remark}

\medskip

\begin{remark}[Technical aspects]~
    \begin{enumerate}[label={(\roman*)}]
        \item Comparing the assumptions in Theorem~\ref{thm:Subeq-Reg-Har-Law-exist-thm} and Theorem~\ref{thm:Subeq-Reg-a-priori-estimates}, $\vec{\widetilde{F}}$-convexity of $\partial \Omega$ has been dropped. This is because when $F$ is convex, any bounded domain $\Omega \subset \RR^n$ with $C^2$ boundary is automatically strictly $\vec{\widetilde{F}}$-convex.

        \item Strict $\vec{F}$-convexity is only needed to invoke Proposition~\ref{prop:Subeq-Reg-existence-of-delta-subsol-and-supersol}~(i) and construct an admissible $\delta$-subsolution $\underline{u}$ \eqref{eq:Subeq-Reg-strongly-admissible-subsolution}. If such a subsolution exists, then strict $\vec{F}$-convexity is needed neither for Theorem~\ref{thm:Subeq-Reg-a-priori-estimates} nor for Corollary~\ref{corollary:Subeq-Reg-existence-thm} (cf. Guan~\cite[Theorem~1.1]{guan-subsol-MR1265805}).

        \item Each of the Lie groups listed in Remark~\ref{remark:Subeq-Reg-which-equations-are-covered} can arise as the holonomy of a Riemannian manifold, and such a Riemannian manifold is called a $G$-manifold. For example an $U(m)$-manifold is K\"ahler, and an $SU(m)$-manifold is Ricci flat K\"ahler, etc. The complex case $U(m)$, when $n=2m$, is very well-studied, especially the complex Monge--Amp\`ere (CMA) equation in particular has a long history ~\cite{Bedford-Taylor-1976-MR445006, CKNS-II-MR780073}.
        
        For a subequation $F\subset \Sym^2(\RR^n)$ invariant under any one of these groups $G \subset O(m)$, the PDE \eqref{eq:Subeq-Reg-const-coeff-subequation-PDE} is also well-defined on the corresponding Riemannian $G$-manifold \cite[Lemma~5.2]{Har-Law-riem-mflds-MR2844439}. The strategy of this article can also be employed to obtain existence in this general setting by following Guan~\cite{Guan-dirichlet-prob-on-riem-mfld-subsol-MR4545123}.
    \end{enumerate}
\end{remark}

\begin{remark}[Applying Theorem~\ref{thm:Subeq-Reg-a-priori-estimates} to the operator formulation \eqref{eq:Subeq-Reg-intro-purely-second-order-PDE}]~
\label{remark:Subeq-Reg-subequations-and-operators}

In practice, PDEs are traditionally phrased in terms of {\it operators},
as in \eqref{eq:Subeq-Reg-intro-purely-second-order-PDE}. 
Why then do we phrase our main results in terms of subequations?
The point is that this is the strongest way to phrase the a priori estimate, even if one is purely interested to apply it to an operator!
(Even so, we will next explain how to apply  Theorem~\ref{thm:Subeq-Reg-a-priori-estimates} to operators in practice.)

To convince the reader, 
    consider an operator $\fcal:\Gamma \rightarrow \RR$ defined on an open set $\Gamma \subset \Sym^2(\RR^n)$, with superlevel set
        \begin{equation}
            \label{eq:Subeq-Reg-intro-superlevel-set-of-operator}
            F_{\fcal} :=\{ \fcal \ge 0 \}.
        \end{equation}
    Under the compatibility condition
    \begin{equation}
        \label{eq:Subeq-Reg-intro-Ffcal-compatibility-condition}
        \partial F_{\fcal} = \{\fcal = 0\},
    \end{equation}
    the $\fcal$-BVP~\eqref{eq:Subeq-Reg-intro-purely-second-order-PDE} is equivalent to the $F_{\fcal}$-BVP. Furthermore, when $\fcal$ is degenerate elliptic \eqref{eq:Subeq-Reg-intro-purely-second-order-PDE-elliptic}, $F_{\fcal}$ is a subequation and is determined by the level set 
    \cite[(1.7)]{Har-Law-PDE-Krylov-MR4127162}
    \begin{equation}
        \label{eq:Subeq-Reg-Ffcal-determined-by-bdry}
        F_{\fcal} = \partial F_{\fcal} + \pcal = \{\fcal = 0\} + \pcal.
    \end{equation}

    \noindent
        Thus, Theorem~\ref{thm:Subeq-Reg-a-priori-estimates} and Corollary~\ref{corollary:Subeq-Reg-existence-thm} tell us that the $\fcal$-BVP~\eqref{eq:Subeq-Reg-intro-purely-second-order-PDE} admits smooth solutions (assuming \eqref{eq:Subeq-Reg-intro-Ffcal-compatibility-condition}), if
    \begin{equation}
        \label{eq:Subeq-Reg-assumptions-on-Ffcal-superlevel-set}
        F_{\fcal} \quad \text{is a smooth, convex, non-degenerate and $G$-invariant subequation}.
    \end{equation}

    What hypotheses on $\fcal$ imply \eqref{eq:Subeq-Reg-assumptions-on-Ffcal-superlevel-set}? If $\fcal$ is itself smooth, concave, strictly elliptic $D\fcal >0$ and $G$-invariant, then \eqref{eq:Subeq-Reg-assumptions-on-Ffcal-superlevel-set} holds. Unfortunately, many natural operators $\fcal$ fail to be concave. For example, the Monge--Amp\`ere equation is the BVP~\eqref{eq:Subeq-Reg-intro-purely-second-order-PDE} for the operator $\fcal = \det -1$, which is $C^2$ and strictly elliptic, but non-concave on $\Gamma = \intr \pcal$ \eqref{eq:Subeq-Reg-positive-semi-def-subeq}. However, both $\fcal_1 = \sqrt[n]{\det} - 1$ and $\fcal_2 = \log \det$ are concave on $\intr \pcal$. And importantly, $ F_{\fcal} = F_{\fcal_1} = F_{\fcal_2}$ \eqref{eq:Subeq-Reg-intro-superlevel-set-of-operator}. Hence, \eqref{eq:Subeq-Reg-assumptions-on-Ffcal-superlevel-set} holds for the Monge--Amp\`ere operator $\fcal = \det -1$.

    \noindent
    In other words, a given operator $\fcal$ might not be the best choice of a compatible operator for the subequation $F_{\fcal}$ \eqref{eq:Subeq-Reg-assumptions-on-Ffcal-superlevel-set}. This is why we formulate assumptions on the more geometric object $F_{\fcal}$ (which is determined by the equation $\fcal(D^2u) = 0$ \eqref{eq:Subeq-Reg-Ffcal-determined-by-bdry}) instead of the operator $\fcal$.

    \noindent
    For a more complete discussion about compatible operator--subequation pairs $(\fcal, F)$, the reader is referred to Cirant--Harvey--Lawson--Payne~\cite[Chapter~11]{Cirant-Payne-Harvey-Lawson-MR4671835} and Harvey--Payne~\cite{Harvey--Payne-MR5056882}.
\end{remark}

\subsection{New ingredients}

Perhaps the most basic new idea in the present work is that {\it a basic problem about a priori estimates can be approached via subequations.} This is genuinely surprising since 
subequation theory (so far at least) is a theory of weak (viscosity) solutions. Essentially all previous papers using subequations concern such weak solutions and yield at best $C^0$ estimates, but not the $C^2$ and higher estimates needed for regularity.

The really simple, but key, new notion that hints that subequations could be relevant to regularity is that of a
{\it non-degenerate subequation} (Definition~\ref{def:Subeq-Reg-intro-non-deg-subequation}) with its geometric properties developed in \S\ref{section:Subeq-Reg-subequations-prelim}.

Another key new ingredient is the use of 
the {\it signed  distance to a subequation} \eqref{eq:Subeq-Reg-signed-dist-def}. This behaves particularly nicely with respect to linearization
(and allows to go back and forth between subequations and operators).

The {\it
conormal subequation} \eqref{eq:Subeq-Reg-F+-projected-subequation} is a natural new notion that helps zoom-in on the most difficult a priori Hessian estimate: the normal-normal estimate (\S\ref{sec:Subeq-Reg-normal-hessian-estimates}). It is the geometric construction underlying Trudinger's asymptotic operator (cf. \eqref{eq:Subeq-Reg-formula-for-sign-dist-to-F+}  and \eqref{eq:Subeq-Reg-Trudinger-asymp-operator}). 
As a bonus, it yields a conceptual geometric formulation of Trudinger's trick that, in passing,
allows us to generalize it considerably.

Finally, identifying that all of the estimates previously obtained separately for real, complex, or quaternionic matrices (as well as many new estimates), actually work for subequations invariant under {\it any $G\subset O(n)$ acting transitively on $S^{n-1}$}, is made possible by combining all of the simplifying notions above to help us see through the forest of a priori estimates.

\subsection{Comparison to previous works}

Traditionally, existence of smooth solutions to \eqref{eq:Subeq-Reg-const-coeff-subequation-PDE} have been studied for smooth, concave, $O(n)$-invariant partial differential \textit{operators} $\fcal$ \cite{caffarelli-nirenberg-spruck-dirichlet-problem-III-MR0806416, Trudinger1995-DHE,guan-subsol-MR1265805}. 

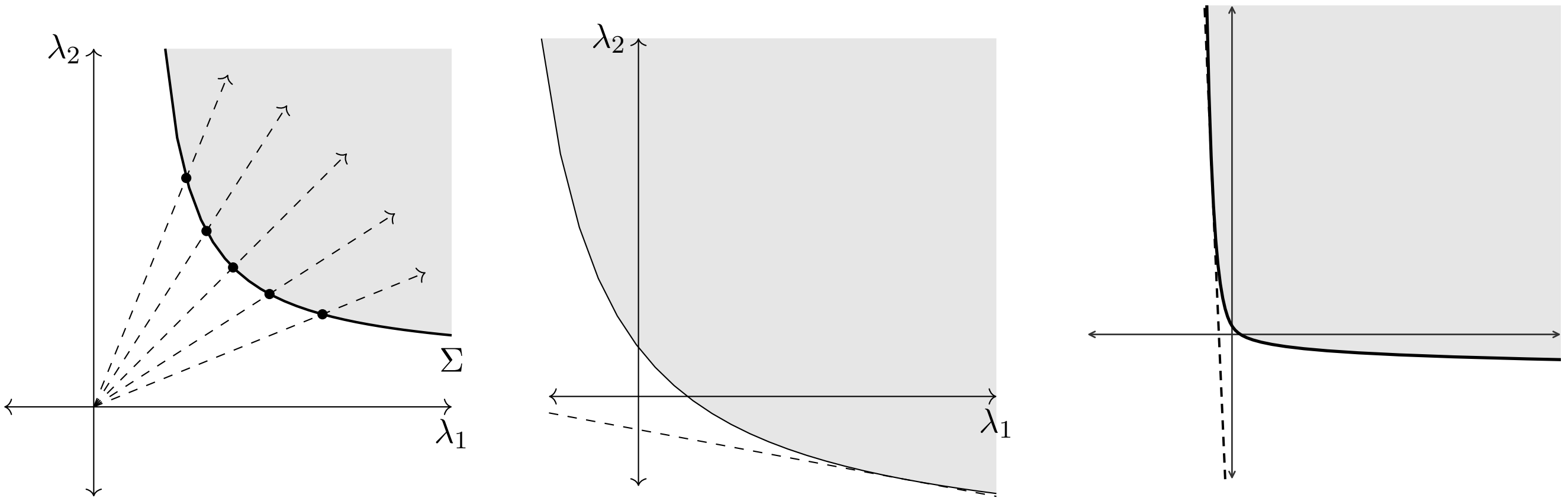
\begin{figure}[t]
    \centering
    \begin{subfigure}[t]{0.32\textwidth}
        \centering
        \begin{tikzpicture}[scale=1]
            \fill[white!90!black] plot[domain=0.7991:4] ({\x}, {tan(deg(2)-atan(\x))}) -- (4,4) -- (0.7991,4);
            \draw[<->, thin] (-1,0) -- (4,0) node[below] {$\lambda_1$};
            \draw[<->, thin] (0,-1) -- (0,4) node[left] {$\lambda_2$};
            \draw[thick] plot[domain=0.7991:4] ({\x}, {tan(deg(2)-atan(\x))}) node[below] {$\Sigma$};
            \draw[dashed, ->] (0,0) -- (3.708, 1.5);
            \draw[dashed, ->] (0,0) -- (3.3666, 2.16);
            \draw[dashed, ->] (0,0) -- (2.8284, 2.8284);
            \draw[dashed, ->] (0,0) -- (2.16, 3.3666);
            \draw[dashed, ->] (0,0) -- (1.5,3.709);
            \filldraw (2.556, 1.034) circle (0.05);
            \filldraw (1.964, 1.260) circle (0.05);
            \filldraw (1.557, 1.557) circle (0.05);
            \filldraw (1.260, 1.964) circle (0.05);
            \filldraw (1.034, 2.556) circle (0.05);
        \end{tikzpicture}
        \caption{}
    \end{subfigure}
    \hfill
    \begin{subfigure}[t]{0.32\textwidth}
        \centering
        \begin{tikzpicture}[scale=1]
            \fill[white!90!black] plot[domain=-1.084:4] ({\x}, {tan(deg(0.5)-atan(\x))}) -- (4,4) -- (-1.084,4);
            \draw[<->, thin] (-1,0) -- (4,0) node[below] {$\lambda_1$};
            \draw[<->, thin] (0,-1) -- (0,4) node[left] {$\lambda_2$};
            \draw plot[domain=-1.084:4] ({\x}, {tan(deg(0.5)-atan(\x))});
            \draw[thin, dashed] (-1, -0.184) -- (4,-1.1164);
        \end{tikzpicture}
        \caption{}
    \end{subfigure}
    \hfill
    \begin{subfigure}[t]{0.32\textwidth}
        \centering
        \includegraphics[width=\linewidth]{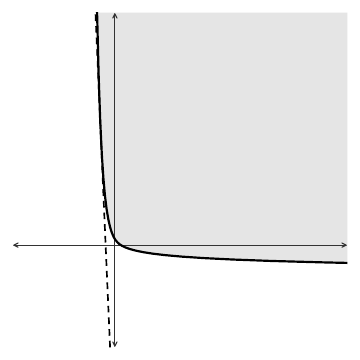}
        \caption{}
    \end{subfigure}
    \caption{(a) By \eqref{eq:Subeq-Reg-intro-CNS-III-11-condition}--\eqref{eq:Subeq-Reg-intro-CNS-III-12-condition}, each ray from the origin intersects $\Sigma$ at most once. There exists a positively 1-homogeneous $f\in C^2(\RR_{\ge0}\Sigma)$ satisfying $f|_{\Sigma} \equiv 1$. (b) $ {\Sigma_{\frac{n-2}{2} + \delta}} := \lambda (F_{\frac{n-2}{2}\pi + \delta}) \subset \RR^n$ \eqref{eq:Subeq-Reg-subequations-for-SL} fails \eqref{eq:Subeq-Reg-intro-CNS-III-12-condition}. (c) Every translate of the $\Sigma_{\exp}$ fails \eqref{eq:Subeq-Reg-intro-CNS-III-12-condition} (Example~\ref{example:Subeq-Reg-fails-CNS-for-all-translations}).}
    \label{fig:Subeq-Reg-CNS-construction}
\end{figure}

Caffarelli--Nirenberg--Spruck additionally studied a BVP somewhat reminiscent of \eqref{eq:Subeq-Reg-const-coeff-subequation-PDE}~\cite[(1)'']{caffarelli-nirenberg-spruck-dirichlet-problem-III-MR0806416},
\[\begin{aligned}
    \lambda(D^2u) \in \Sigma& \qquad \text{on } \Omega,\\
    u = \varphi& \qquad \text{on } \partial \Omega,
\end{aligned}\]
where $\Sigma \subset \RR^n$ is a smooth, permutation invariant hypersurface bounding a convex set, whose inward unit normals have positive components~\cite[(11)]{caffarelli-nirenberg-spruck-dirichlet-problem-III-MR0806416},
\begin{equation}
    \label{eq:Subeq-Reg-intro-CNS-III-11-condition}
    \nu^{\Sigma}(y) \in (0,\infty)^n, \qquad \forall\ y\in \Sigma.
\end{equation}
Their assumptions on $\Sigma$ are equivalent to asserting that
\[ F_{\Sigma} = \left\{A +P: A,P\in \Sym^2(\RR^n),\ \lambda(A) \in \Sigma,\ P\ge 0\right\}\]
is an $O(n)$-invariant, convex, non-degenerate subequation with smooth boundary $\partial F$ 
(we do not prove this equivalence, but it is not hard to show from the arguments of
Theorem~\ref{thm:Subeq-Reg-a-priori-estimates}), where $\lambda(A)$ denotes the vector of eigenvalues of $A$. However, their existence theorem \cite[Theorem~4]{caffarelli-nirenberg-spruck-dirichlet-problem-III-MR0806416} needed one additional assumption \cite[(12)]{caffarelli-nirenberg-spruck-dirichlet-problem-III-MR0806416}:
\begin{equation}
    \label{eq:Subeq-Reg-intro-CNS-III-12-condition}
    \text{each tangent plane to $\Sigma$ separates $\Sigma$ from the origin $0\in \RR^n$.}
\end{equation}
Assumption \eqref{eq:Subeq-Reg-intro-CNS-III-12-condition} was used to construct a smooth, concave, $O(n)$-invariant partial differential operator $\fcal$ on $\Sym^2(\RR^n)$ (more precisely, on the cone $\RR_{\ge0} F_{\Sigma}$), for which $\partial F$ is a regular level set ($\partial F_\Sigma = \{\fcal =1\}$) (Figure~\ref{fig:Subeq-Reg-CNS-construction}~(a)). Then, they apply their existence theorem for partial differential operators \cite[Theorem~2]{caffarelli-nirenberg-spruck-dirichlet-problem-III-MR0806416}.

\noindent It turns out that the condition\eqref{eq:Subeq-Reg-intro-CNS-III-12-condition} fails for several natural subequations. The subequations $F_c$ \eqref{eq:Subeq-Reg-subequations-for-SL} are convex, non-degenerate and $O(n)$-invariant for $c\in [\frac{n-2}{2}\pi, \frac{n}{2}\pi)$~\cite{yuan-global-solutions-MR2199179} (cf. Theorem~\ref{thm:Subeq-Reg-a-priori-estimates}), but $\Sigma_c := \lambda(F_c)$ does not satisfy \eqref{eq:Subeq-Reg-intro-CNS-III-12-condition}, for $c\in (\frac{n-2}{2}\pi, \frac{n-1}{2}\pi)$ (Figure~\ref{fig:Subeq-Reg-CNS-construction}~(b)). It is possible to show that some translation $\Sigma_c+ t_c,\ t_c\in \RR^n$, does satisfy \eqref{eq:Subeq-Reg-intro-CNS-III-12-condition} even for $c\in (\frac{n-2}{2}\pi, \frac{n-1}{2}\pi)$. So in principle \cite[Theorem~4]{caffarelli-nirenberg-spruck-dirichlet-problem-III-MR0806416} can be applied to $F_c$. However, there are many examples of subequations $F$ such that \eqref{eq:Subeq-Reg-intro-CNS-III-12-condition} fails to hold for \textit{any} translate $\Sigma_B:= \lambda(F)+\lambda(B)$, for $B\in \Sym^2(\RR^n)$:
\begin{example}
\label{example:Subeq-Reg-fails-CNS-for-all-translations}
    Consider the function
    \[f:\RR^2 \rightarrow \RR,\qquad f(x,y) := x - e^{-x} + y -e^{-y},\]
    and
    \[\Sigma_{\exp} := \left\{ (x,y)\in \RR^2:  f(x,y) = 0 \right\}.\]
    By a simple calculation, every translate of $\Sigma_{\exp}$ fails to satisfy \eqref{eq:Subeq-Reg-intro-CNS-III-12-condition}. The key point is that the tangent plane to the level set $ \lambda(\partial F) = \{f= 0\}$ at $(-R-1, e^{R+1} +o_R(1)) \in \{f=0\}$ is nearly parallel to the $y$-axis, for $R>>1$. Hence, $(-R,-R)$ will not be separated from the surface by this tangent plane (Figure~\ref{fig:Subeq-Reg-CNS-construction}~(c)).
\end{example}

Our first main observation is that \eqref{eq:Subeq-Reg-intro-CNS-III-12-condition} is not needed for existence, thus considerably (and essentially optimally) enlarging the class of equations \eqref{eq:Subeq-Reg-const-coeff-subequation-PDE} that admit solutions. Without \eqref{eq:Subeq-Reg-intro-CNS-III-12-condition}, we cannot construct a smooth, concave defining partial differential operator $\fcal$ as done in~\cite[Theorem~4]{caffarelli-nirenberg-spruck-dirichlet-problem-III-MR0806416}. Instead, we use the signed-distance $\signdist^{F}$ \eqref{eq:Subeq-Reg-signed-dist-fcal-PDE} to the boundary of the subequation $\partial F$, and the corresponding formulation of \eqref{eq:Subeq-Reg-const-coeff-subequation-PDE}, 
\begin{equation}
    \label{eq:Subeq-Reg-signed-dist-fcal-PDE}
    \begin{aligned}
        \signdist^{F}(D^2u) = 0 \qquad &\text{on } \Omega,\\
        u = \varphi \qquad &\text{on }\partial \Omega.
    \end{aligned}
\end{equation}
By Lemma~\ref{lemma:Subeq-Reg-signed-distance-lemma}~(ii)--(iv), $\signdist^{F}$ satisfies all the assumptions for an existence theorem (e.g., \cite[(1.10)--(1.13)]{Trudinger1995-DHE}) except that $\signdist^{F}$---while defined on all of $\Sym^2(\RR^n)$---is only smooth in a neighborhood of the boundary $\partial F$ (Lemma~\ref{lemma:Subeq-Reg-signed-distance-lemma}~(iii)). We streamline the arguments of Caffarelli--Nirenberg--Spruck~\cite[Theorem~2]{caffarelli-nirenberg-spruck-dirichlet-problem-III-MR0806416}, Trudinger~\cite[Theorem~1.2]{Trudinger1995-DHE} and Guan~\cite[Theorem~1.1]{guan-subsol-MR1265805}, and demonstrate that smoothness of $\signdist^F$ near $\partial F$ is sufficient for many of their arguments to hold.

Our second main observation is that $O(n)$-invariance can be considerably generalized. Most of the general existence theorems for \eqref{eq:Subeq-Reg-const-coeff-subequation-PDE} follow the lead of Caffarelli--Nirenberg--Spruck~\cite{caffarelli-nirenberg-spruck-dirichlet-problem-III-MR0806416} and work with $O(n)$-invariant PDE $\fcal$ \cite{Trudinger1995-DHE, guan-subsol-MR1265805}. The $U(m)$-invariant case has also received some attention: the simplest $U(m)$-invariant equation is the complex Monge--Amp\`ere (CMA). Bedford--Taylor constructed weak solutions to the CMA~\cite[Theorem~D]{Bedford-Taylor-1976-MR445006}, and Caffarelli--Kohn--Nirenberg--Spruck~\cite[Theorem~1]{CKNS-II-MR780073} proved regularity of solutions on uniformly pseudoconvex domains. S.-Y.~Li~\cite[Theorem~1.1]{Li-song-ying-MR2128299} extended this regularity to general $U(m)$-invariant partial differential operators of the complex Hessian. Similarly, the quaternionic Monge--Amp\`ere was introduced by Alesker~\cite[Theorems~1.3--1.4]{Alesker-quaternionic-MA-MR1967025} introduced the quaternionic Monge--Amp\`ere (QMA) and proved existence of weak solutions as well as regularity when the domain is a ball. Zhu~\cite[Theorem~1.1]{Zhu-quaternionic-MA-MR3544695} proved regularity for the QMA on uniformly quaternionic pseudoconvex domains. Though neither of them studied general equations invariant under the symmetries of the quaternions.

In the approach of Caffarelli--Nirenberg--Spruck~\cite{caffarelli-nirenberg-spruck-dirichlet-problem-III-MR0806416}, $O(n)$-invariance is primarily needed to estimate the boundary tangent-normal (mixed) Hessian components. In the $U(m)$-invariant case, Caffarelli--Kohn--Nirenberg--Spruck~\cite[Lemma~1.3]{CKNS-II-MR780073}, and its generalization by S.-Y.~Li~\cite[Lemma~4.3]{Li-song-ying-MR2128299}, replace this mixed Hessian estimate by a lengthy and technical argument using a barrier involving a quadratic in $Du$. Instead, we modify the original argument of Caffarelli--Nirenberg--Spruck~\cite{caffarelli-nirenberg-spruck-dirichlet-problem-III-MR0806416} for the $G$-invariant case. In particular, our argument considerably simplifies the a priori estimates for $U(m)$-invariant equations. Furthermore, it allows Theorem~\ref{thm:Subeq-Reg-a-priori-estimates} and Corollary~\ref{corollary:Subeq-Reg-existence-thm} to cover essentially all known equations of type \eqref{eq:Subeq-Reg-const-coeff-subequation-PDE}, whether they are real, complex, quaternionic, or even more general. 



\bigskip
\noindent\textbf{AI Declaration.} This article is part of the first-named author's
PhD thesis and is mostly a product of the pre-AI boom era, being the culmination of work done over the years 2023-2026.
Having said that, AI had two simplifying inputs.
First, we had originally proved the {\it local} concavity of $\signdist^F$ \eqref{eq:Subeq-Reg-signed-dist-def}, and this was sufficient for our main results; the proof of {\it global} concavity was done by AI (Lemma~\ref{lemma:Subeq-Reg-signed-distance-lemma}~(v)). 
Second, AI was used to digest 
the classical proof of the mixed Hessian estimate from Caffarelli--Kohn--Nirenberg--Spruck~\cite[Lemma~1.3]{CKNS-II-MR780073} and S.-Y.~Li~\cite[Lemma~4.3]{Li-song-ying-MR2128299} for PDE that are $U(m)$-invariant functions of the complex Hessian. This digestion ultimately helped us see how the original argument of Caffarelli--Nirenberg--Spruck~\cite[pp.~281--284]{caffarelli-nirenberg-spruck-dirichlet-problem-III-MR0806416} in the $O(n)$-invariant case also generalizes to $G$-invariant subequations for any $G\subset O(n)$ acting transitively on $S^{n-1}$.

\noindent
Additionally, AI was used to simplify a few technical steps in \S\ref{section:Subeq-Reg-signed-distance} and \S\S\ref{sec:Subeq-Reg-mixed-hessian-estimates}--\ref{sec:Subeq-Reg-normal-hessian-estimates}, to help produce the Example~\ref{example:Subeq-Reg-fails-CNS-for-all-translations}, and to catch typographical errors throughout the article.

\noindent
All the text and figures in this article are written by us, and we take full responsibility for any errors.

\bigskip
\noindent\textbf{Acknowledgments.} 
The first-named author was supported by the Patrick and Marguerite Sung Graduate Research Fellowship in Mathematics
in Fall 2025.
The second-named author was supported by the Bergman Distinguished Visiting Professorship at Stanford
University, and a Simons Fellowship in Mathematics.
Both authors
were supported 
by
NSF grants DMS-1906370, 2204347, 2506872
and by BSF grants 2016173, 2020329, 2024293.
The authors are grateful to the Mathematics Department and the MRC at Stanford University for the hospitality during the academic year 2025--2026.

\section{Non-degenerate subequations}
\label{section:Subeq-Reg-subequations-prelim}

This section characterizes non-degeneracy \eqref{eq:Subeq-Reg-non-degenerate-subequation} of a subequation $F$ in terms of the inward unit normals to the boundary. 
Note that a subequation $F$ admits inward normals to its boundary $\partial F$ if either $F$ is convex (in which case these normals are the supporting hyperplanes) or if $\partial F$ is a smooth manifold. The former are studied in \eqref{subsection:Subeq-Reg-convex-non-deg-subequations-normals} and the latter in \eqref{subsection:Subeq-Reg-smooth-subequations-normals}

\subsection{Convex non-degenerate subequations}
\label{subsection:Subeq-Reg-convex-non-deg-subequations-normals}

Denote the Euclidean inner product on $\Sym^2(\RR^n)$, also called the \textit{Frobenius} or \textit{trace inner product}, by
\begin{equation}
    \label{eq:Subeq-Reg-trace-inner-prod}
    \langle A,B \rangle := \tr(AB), \qquad \forall\ A,B\in \Sym^2(\RR^n).
\end{equation}
A convex set $F\subset \Sym^2(\RR^n)$ admits a supporting hyperplane at each boundary point $A\in \partial F$~\cite[Theorem~11.2]{rockafeller-convex-analysis-MR1451876}, i.e., for each $A_0\in \partial F$, there exists $B_{A_0}\in \Sym^2(\RR^n)$ such that
\begin{equation}
    \label{eq:Subeq-Reg-inward-supporting-hyperplane}
    F \subset \left\{ A\in \Sym^2(\RR^n): \langle A, B_{A_0} \rangle \ge \langle A_0, B_{A_0} \rangle \right\},
\end{equation}
and in particular,
\begin{equation}
    \label{eq:Subeq-Reg-inward-supporting-hyperplane-interior}
    \intr F \subset \left\{ A\in \Sym^2(\RR^n): \langle A, B_{A_0} \rangle > \langle A_0, B_{A_0} \rangle \right\}.
\end{equation}
Note that $B_{A_0}$ is an \textit{inward} normal to a supporting hyperplane to $\partial F$ at $A_0$.

\begin{lemma}
    \label{lemma:Subeq-Reg-convex-subequation-normals}
    For a non-empty, closed, proper convex subset $F\subset \Sym^2(\RR^n)$ with inward normals denoted by \eqref{eq:Subeq-Reg-inward-supporting-hyperplane},
    \[\begin{aligned}
        F \quad \text{is a subequation \eqref{eq:Subeq-Reg-subequation-def-prop-basic}} \qquad \iff& \qquad B_{A_0} \ge 0, \quad \forall\ A_0 \in \partial F,\\
        F \quad \text{is non-degenerate \eqref{eq:Subeq-Reg-non-degenerate-subequation}} \qquad \iff& \qquad B_{A_0} > 0, \quad \forall\ A_0 \in \partial F.
    \end{aligned}\]
\end{lemma}
\begin{proof}
    \noindent\textit{(i) }\textit{\eqref{eq:Subeq-Reg-subequation-def-prop-basic} (resp. \eqref{eq:Subeq-Reg-non-degenerate-subequation}) $\implies\ B_{A_0} \ge 0$ (resp. $>0$) for all $A_0 \in \partial F$:} 
    
    \vspace{0.5em}\noindent Fix $A_0\in \partial F$ and a corresponding inward normal $B_{A_0}$ \eqref{eq:Subeq-Reg-subequation-inward-normal}. By \eqref{eq:Subeq-Reg-subequation-def-prop-basic} (resp. \eqref{eq:Subeq-Reg-non-degenerate-subequation}),
    \[A_0 +P \in F\ (\text{resp. } \intr F), \qquad \forall\ A_0 \in \partial F,\ P\in \pcal\setminus\{0\}.\]
    Hence, by \eqref{eq:Subeq-Reg-inward-supporting-hyperplane}--\eqref{eq:Subeq-Reg-inward-supporting-hyperplane-interior},
    \[\begin{aligned}
        \langle A_0+P, B_{A_0}\rangle - \langle A_0, B_{A_0}\rangle \ge&\ 0 \ (\text{resp. } > 0) \\
        \implies \qquad \langle P, B_{A_0} \rangle \ge&\ 0 \ (\text{resp. } > 0), \qquad \forall\ A_0 \in \partial F,\ P\in \pcal\setminus\{0\}.
    \end{aligned}\]
    Taking $P = vv^T$,
    \[0 \stackrel{(<)}{\le} \langle vv^T, B_{A_0} \rangle = \tr(vv^T B_{A_0}) = \langle B_{A_0} v,v\rangle, \qquad \forall\ A_0\in \partial F,\ v \in \RR^n\setminus \{0\}.\]
    Conclude that $B_{A_0} \ge 0$ (resp. $>0$).

    \vspace{0.5em}\noindent\textit{(ii) }\textit{ $B_{A_0} \ge 0$ (resp. $>0$) for all $A_0 \in \partial F$ $\implies$ \eqref{eq:Subeq-Reg-subequation-def-prop-basic} (resp. \eqref{eq:Subeq-Reg-non-degenerate-subequation}):}  
    
    \vspace{0.5em}\noindent Assume on the contrary that \eqref{eq:Subeq-Reg-subequation-def-prop-basic} (resp. \eqref{eq:Subeq-Reg-non-degenerate-subequation}) fails. Then there exists $A,P \in \Sym^2(\RR^n)$ such that
    \[A\in F,\quad P\in \pcal\setminus\{0\},\quad \text{and}\quad  A+P \notin F \text{ (resp. $\notin \intr F$)}.\]
    Consider
    \[t_0:= \sup \left\{t\in [0,1]: A+tP \in F \text{ (resp. $\in \intr F$)} \right\}.\]
    Since $F$ is closed and $A+P\notin F$,
    \begin{equation}
        \label{eq:Subeq-Reg-convex-pos-supp-hyper-sueq-eq1}
        0\le t_0<1, \qquad A+t_0P \in \partial F, \qquad A+P \notin F  \text{ (resp. $\notin \intr F$)}.
    \end{equation}
    Consider a supporting hyperplane to $\partial F$ at $A + t_0P$, say with an inward normal $B=B_{A+t_0P} \in \Sym^2(\RR^n)$ satisfying \eqref{eq:Subeq-Reg-inward-supporting-hyperplane} at $A_0 = A+t_0P$. Then, by \eqref{eq:Subeq-Reg-convex-pos-supp-hyper-sueq-eq1},
    \[\langle B, (1-t_0)P \rangle = \langle B, (A +P) - (A + t_0P)\rangle \stackrel{(\le)}{<}  0 .\]
    Since $t_0<1$, it follows that $\langle B,P\rangle <0$ (resp. $\le  0$). Hence, $B = B_{A+t_0P}$ has a negative (resp. non-positive) eigenvalue. This contradicts the assumption that every supporting hyperplane has to $F$ has positive semi-definite (resp. positive definite) inward normals.
\end{proof}

\subsection{Smooth subequations}
\label{subsection:Subeq-Reg-smooth-subequations-normals}

The main observation of this section is that a smooth, convex, non-degenerate subequation has positive definite inward unit normals (Lemma~\ref{lemma:Subeq-Reg-subequation-is-elliptic}~(iii)).

For a subequation with smooth boundary $\partial F$, denote \textit{the} inward unit normal to $\partial F$ with respect to the trace inner product \eqref{eq:Subeq-Reg-trace-inner-prod} by 
\begin{equation}
    \label{eq:Subeq-Reg-subequation-inward-normal}
    \nu^{\partial F}:\partial F \rightarrow \Sym^2(\RR^n).
\end{equation}
As it turns out, non-negativity of the unit normals $\nu^{\partial F}\ge 0$ is equivalent to $F$ being a subequation. And when $F$ is convex, by Lemma~\ref{lemma:Subeq-Reg-convex-subequation-normals}, positivity $\nu^{\partial F}> 0$ is equivalent to $F$ being a non-degenerate subequation:

\begin{lemma}
\label{lemma:Subeq-Reg-subequation-is-elliptic}
    For a non-empty, proper, closed subset $F\subsetneq \Sym^2(\RR^n)$ with smooth boundary $\partial F$,
    \begin{enumerate}[label={(\roman*)}]
        
        \item $F$ is a subequation \eqref{eq:Subeq-Reg-subequation-def-prop-basic} $\quad \iff \quad \nu^{\partial F}(A) \ge 0, \quad \forall\ A\in \partial F$,
        
        \item $\nu^{\partial F}(A) > 0, \quad \forall\ A\in \partial F \quad \implies \quad $ $F$ is non-degenerate \eqref{eq:Subeq-Reg-non-degenerate-subequation},
        
        \item $F$ is convex and non-degenerate \eqref{eq:Subeq-Reg-non-degenerate-subequation} $\quad \implies \quad \nu^{\partial F}(A) > 0, \quad \forall\ A\in \partial F$.
    \end{enumerate}
\end{lemma}

A subequation $F$ is called \textit{strictly elliptic} if it has positive definite inward unit normals,
\begin{equation}
    \label{eq:Subeq-Reg-strictly-elliptic-subequation-pos-normals}
    \nu^{\partial F}(A) >0, \qquad \forall\ A\in \partial F.
\end{equation}
When $F$ is the super-level set of a smooth elliptic partial differential operator operator $\fcal$ \eqref{eq:Subeq-Reg-intro-purely-second-order-PDE-elliptic}, $F = \{\fcal \ge 0\}$, \eqref{eq:Subeq-Reg-strictly-elliptic-subequation-pos-normals} follows from strictly ellipticity of $\fcal$, i.e, positivity of its symbol $D\fcal >0$. 

\begin{proof}
    When $F$ is \textit{convex} with smooth boundary $\partial F$, $\nu^{\partial F}(A_0)$ defines the unique supporting hyperplane to $\partial F$ at $A_0$, i.e., $B_{A_0} = \nu^{\partial F}(A_0)$ \eqref{eq:Subeq-Reg-inward-supporting-hyperplane}. Hence, (iii) follows from Lemma~\ref{lemma:Subeq-Reg-convex-subequation-normals}.

    \vspace{0.5em}\noindent\textit{(i) $\implies$:} Fix a boundary point
    \[A\in \partial F.\]
    Since $\partial F$ is a smooth manifold, there exists a neighborhood $ U  \subset \Sym^2(\RR^n)$ of $A\in \partial F$, and a smooth local defining function $f:U \rightarrow \RR$ satisfying
    \[F  \cap U =\{f \ge 0\}, \quad \partial F \cap U = \{f = 0\}, \quad \text{and} \quad D f(A) = \nu^{\partial F}(A) \neq 0.\]
    By \eqref{eq:Subeq-Reg-subequation-def-prop-basic},
    \[A+tP \in F \quad \implies \quad f(A+tP) \ge 0, \qquad \forall\ P\ge 0,\ 0<t\ll 1.\]
    Since $A\in \partial F$, $f(A) = 0$ and so
    \[\langle \nu^{\partial F}(A), P\rangle = \langle Df(A),P\rangle = \lim_{t\rightarrow 0^+} \frac{f(A+tP) - f(A)}{t} \ge 0, \qquad \forall\ P\ge 0. \]
    Setting $P= vv^T \ge 0 $,
    \[0\le \langle \nu^{\partial F}(A), vv^T\rangle = \tr\left( \nu^{\partial F}(A)vv^T \right) = \langle \nu^{\partial F}(A)v,v\rangle \qquad \forall\ v\in \RR^n,\]
   i.e., $\nu^{\partial F}(A) \ge 0$.

    \vspace{0.5em}\noindent\textit{(i) $\impliedby$:} Assume on the contrary that \eqref{eq:Subeq-Reg-subequation-def-prop-basic} fails. Then there exists $A,P \in \Sym^2(\RR^n)$ such that
    \[A\in F,\quad P\ge 0,\quad \text{and}\quad  A+P \notin F.\]
    Since $F^c$ is open, in fact $A+P+\delta I\notin F $ for some $\delta>0$. Setting
    \[Q:= P+\delta I >0,\]
    it follows that
    \begin{equation}
        \label{eq:Subeq-Reg-pos-normals-is-equiv-subeq-eq1}
        A\in F,\quad Q> 0,\quad \text{and}\quad  A+Q \notin F.
    \end{equation}
    Consider
    \[t_0:= \sup \left\{t\in [ 0,1]: A+tQ \in F \right\}.\]
    Since $F$ is closed and $A+Q \notin F$,
    \begin{equation}
        \label{eq:Subeq-Reg-pos-normals-is-equiv-subeq-eq1.25}
        t_0<1, \qquad A+t_0Q \in \partial F, \qquad A+tQ \notin F , \quad \forall\ t\in (t_0,1].
    \end{equation}
    Since $\partial F$ is a smooth manifold, there exists a neighborhood $ U  \subset \Sym^2(\RR^n)$ of $A\in \partial F$, and a smooth function $f:U \rightarrow \RR$ such that
    \[F  \cap U =\{f \ge 0\}, \quad \partial F \cap U \cap \{f = 0\}, \quad \text{and} \quad D f(A+t_0Q) = \nu^{\partial F}(A+t_0 Q) \neq 0.\]
    For some small $\epsilon>0$,
    \[A+tQ \in U, \qquad \forall\ t\in (t_0-\epsilon, t_0+\epsilon),\]
    and by \eqref{eq:Subeq-Reg-pos-normals-is-equiv-subeq-eq1.25},
    \[f(A+t_0 Q) = 0,\qquad f(A+tQ) <0 , \quad \forall\ t\in (t_0,t_0 +\epsilon).\]
    Conclude that
    \[\langle  \nu^{\partial F}(A+t_0Q), Q \rangle = \langle Df(A+t_0Q), Q\rangle  = \lim_{h\rightarrow 0} \frac{f(A+t_0Q+hQ) -f(A+t_0Q)}{h} \le 0.\]
    Since $Q>0$ \eqref{eq:Subeq-Reg-pos-normals-is-equiv-subeq-eq1} and $\langle  \nu^{\partial F}(A+t_0Q), Q \rangle\le 0$, it follows that $\nu^{\partial F}(A+t_0Q)$ has a negative eigenvalue. This contradicts the assumption that $\nu^{\partial F}(A')\ge 0$ for all $A'\in \partial F$. Hence, \eqref{eq:Subeq-Reg-subequation-def-prop-basic}.

    \vspace{0.5em}\noindent\textit{(ii):} Assume on the contrary that \eqref{eq:Subeq-Reg-non-degenerate-subequation} fails, there exists $A,P\in \Sym^2(\RR^n)$ such that
    \begin{equation}
        \label{eq:Subeq-Reg-pos-normals-is-equiv-subeq-eq2.0}
        A\in F,\quad P\ge 0, P\neq 0,\quad \text{and}\quad  A+P \notin \intr F.
    \end{equation}
    Repeating the argument of the reverse direction ($\impliedby$) of part (i), with $P$ replacing $Q$, there exists $t_0\in [0,1)$ such that $A + t_0P \in \partial F$ and
    \[\langle  \nu^{\partial F}(A+t_0P), P \rangle = \langle Df(A+t_0P), P\rangle  = \lim_{h\rightarrow 0} \frac{f(A+t_0P+hP) -f(A+t_0P)}{h} \le 0.\]
    The paragraph of (i) $\impliedby$ now does not apply directly because $P\not>0$. However, $\nu^{\partial F}(A+t_0P)>0$ by assumption and $P\ge 0,\ P\neq 0$ \eqref{eq:Subeq-Reg-pos-normals-is-equiv-subeq-eq2.0}. This contradicts $\langle  \nu^{\partial F}(A+t_0P), P \rangle \le 0$. Hence, \eqref{eq:Subeq-Reg-non-degenerate-subequation} holds.
\end{proof}

\begin{remark}
\label{remark:Subeq-Reg-converse-non-deg-not-imply-strictly-elliptic}
    The converse to Lemma~\ref{lemma:Subeq-Reg-subequation-is-elliptic}~(ii) is false. Consider the subequation
    \[F:= \left\{A\in \Sym^2(\RR^n): \fcal(A):=\tr(A^3) = \sum_{i=1}^n \lambda_i(A)^3 \ge  1 \right\},\]
    where $\lambda_1(A) \ge \dots \ge \lambda_n(A)$ denote the eigenvalues of $A$. $F$ is a non-degenerate subequation \eqref{eq:Subeq-Reg-non-degenerate-subequation} because $P\ge 0 ,\ P\neq 0$ implies $\lambda_i(A+P) > \lambda_i(A)$ for some $1\le i\le n$. and since $x^3$ is an increasing function. Computing the inward normal at $ A\in \partial F$,
    \[\begin{aligned}
        \langle D\fcal(A), B\rangle  =&\ \lim_{h\rightarrow 0}\frac{\tr((A+hB)^3) - \tr (A^3)}{h} \\
        =&\ \tr\left(\lim_{h\rightarrow 0}\frac{(A+hB)^3 - A^3}{h}\right) \\
        =&\ \tr \left( A^2B + ABA + BA^2 \right) = 3\tr(A^2B),\\
        \implies\quad D\fcal(A) =&\ 3A^2.
    \end{aligned} \]
    The inward unit normal to $\partial F$, $\nu^{\partial F}$ \eqref{eq:Subeq-Reg-subequation-inward-normal}, is proportional to $D\fcal$, and so,
    \[\nu^{\partial F}(A) = \frac{3A^2}{|3A^2|} \not>0 \qquad \text{for }A = \diag(1,0,\dots, 0) \in \partial F.\]
    Thus $F$ is smooth and non-degenerate but does \textit{not} have positive definite inward unit normals \eqref{eq:Subeq-Reg-strictly-elliptic-subequation-pos-normals}. In other words, the converse to Lemma~\ref{lemma:Subeq-Reg-subequation-is-elliptic}~(ii) is false.
\end{remark}

\section{Signed-distance}
\label{section:Subeq-Reg-signed-distance}

The strategy to solve the BVP~\eqref{eq:Subeq-Reg-const-coeff-subequation-PDE} is to reduce it to a concave partial differential operator in a neighborhood of the boundary of the subequation $\partial F$. The simplest such operator is the signed-distance $\signdist^F$ \eqref{eq:Subeq-Reg-signed-dist-def} to the boundary $\partial F$. This section explores the properties of the signed-distance. For more details, the reader is referred to Gilbarg--Trudinger~\cite[\S14.6]{gilbarg-trudinger-MR1814364} and Petersen~\cite[p.42]{Petersen-Riem-Geom-MR3469435}.

For a subequation $F\subset \Sym^2(\RR^n)$, the \textit{signed-distance to the boundary of $F$} $\signdist^{F}:\Sym^2(\RR^n) \rightarrow \RR$ is the function \cite[p.42]{Petersen-Riem-Geom-MR3469435},
    \begin{equation}
        \label{eq:Subeq-Reg-signed-dist-def}
        \signdist^{F}(B) := \begin{cases}
        \displaystyle{\inf_{A\in \partial F} |B - A|,} \qquad &\text{if } B\in  F,\\
        \displaystyle{-\inf_{A\in \partial F} |B - A|}, \qquad &\text{if } B\notin  F,
    \end{cases}
    \end{equation}
where $|B|$ is the trace norm with respect to \eqref{eq:Subeq-Reg-trace-inner-prod}. Note that, in our convention, $\signdist^{F}$ is positive on the chosen subequation $F$ and negative on its complement.


The following lemma interprets various assumptions on the subequation $F$ in terms of the signed-distance $\signdist^{F}$. For more details, the reader is referred to Petersen~\cite[p.42]{Petersen-Riem-Geom-MR3469435}. 

\begin{lemma}
\label{lemma:Subeq-Reg-signed-distance-lemma}
    Fix a subequation $F\subset \Sym^2(\RR^n)$ with signed-distance $\signdist^{F}$ \eqref{eq:Subeq-Reg-signed-dist-def}.
    \begin{enumerate}[label={(\roman*)}]
        \item (1-Lipschitz) $\signdist^{F}$ is continuous and 1-Lipschitz on $\Sym^2(\RR^n)$.
        
        \item (Linear growth) There exists $M>0$ such that $\signdist^{F}$ has the linear growth rate
        \begin{equation}
            \label{eq:Subeq-Reg-signed-distance-lemma-linear-growth}
            t - M\le \signdist^{F}(tI) \le \sqrt{n} t +M, \qquad \forall\ t\ge 0.
        \end{equation}
        
        \item (Smooth) If $\partial F$ is a smooth manifold, then $\signdist^{F}$ is smooth on a neighborhood $U$ of the boundary $\partial F$ and has unit gradient
    \begin{equation}
        \label{eq:Subeq-Reg-signed-distance-unit-norm-derivative}
        |D \signdist^{F}| = 1 \quad \text{on } U.
    \end{equation}
    In particular $D \signdist^{F} = \nu^{\partial F}$ \eqref{eq:Subeq-Reg-subequation-inward-normal} on $\partial F$

        \item ($G$-invariant) If $F$ is $G$-invariant for a group $G\subset O(n)$ \eqref{eq:Subeq-Reg-group-invariant-subequation}, then so is $\signdist^{F}$, i.e.,
        \begin{equation}
            \label{eq:Subeq-Reg-sign-dist-is-orthogonal-invariant}
            \signdist^{F}(O^TAO) = \signdist^{F}(A), \qquad \forall\ A\in F,\ O\in G.
        \end{equation}
    
        \item (Concave) If $F$ is a convex set, then $\signdist^{F}$ is concave on $\Sym^2(\RR^n)$.

        \item (Strictly elliptic) If $\partial F$ is smooth and strictly elliptic \eqref{eq:Subeq-Reg-strictly-elliptic-subequation-pos-normals}, then
        \begin{equation}
            \label{eq:Subeq-Reg-sign-dist-has-pos-grad}
            D\signdist^{F}(A)>0 \qquad \forall\ A\in \partial F.
        \end{equation}

        \item (Monotonicity) By \eqref{eq:Subeq-Reg-subequation-def-prop-basic},
        \begin{equation}
            \label{eq:Subeq-Reg-basic-monotonicity-for-sign-dist}
            \signdist^F(B+P) \ge \signdist^F(B), \qquad \forall\ B,P\in \Sym^2(\RR^n),\ P \ge 0,
        \end{equation}

        \item (Non-degenerate) If $F$ is non-degenerate \eqref{eq:Subeq-Reg-non-degenerate-subequation}, then
        \begin{equation}
            \label{eq:Subeq-Reg-sign-dist-is-strictly-elliptic}
            \signdist^{F}(B+P) >\signdist^{F}(B), \qquad \forall\ B,P\in \Sym^2(\RR^n),\ P \ge 0,\ P\neq 0. 
        \end{equation}
    \end{enumerate}
\end{lemma}


\begin{proof}[Proof of Lemma~\ref{lemma:Subeq-Reg-signed-distance-lemma}]
\noindent\textit{(i):} Being a distance to subset of $\Sym^2(\RR^n)$, $\signdist^{F}$ is continuous, and is 1-Lipschitz by triangle inequality on $\Sym^2(\RR^n)$.

\vspace{0.5em}\noindent\textit{(ii):} The upper bound in \eqref{eq:Subeq-Reg-signed-distance-lemma-linear-growth} follows easily since $\signdist^{F}$ is 1-Lipschitz by (i),
\begin{equation}
    \label{eq:Subeq-Reg-sign-dist-linear-growth-eq2}
    \signdist^{F}(tI) \le \signdist^{F}(0) + |tI| = \signdist^{F}(0) + \sqrt{n}t.
\end{equation}
It remains to prove the lower bound in \eqref{eq:Subeq-Reg-signed-distance-lemma-linear-growth}. Since $F$ is a subequation, there exists $t_0>0$ such that \cite[p.~11, (4)]{HarLaw-dirichletduality}
\[t_0 I \in F.\]
Then, by \eqref{eq:Subeq-Reg-subequation-def-prop-basic},
\[tI \in t_0I + \pcal  \subset F, \qquad \forall\ t\ge t_0.\]
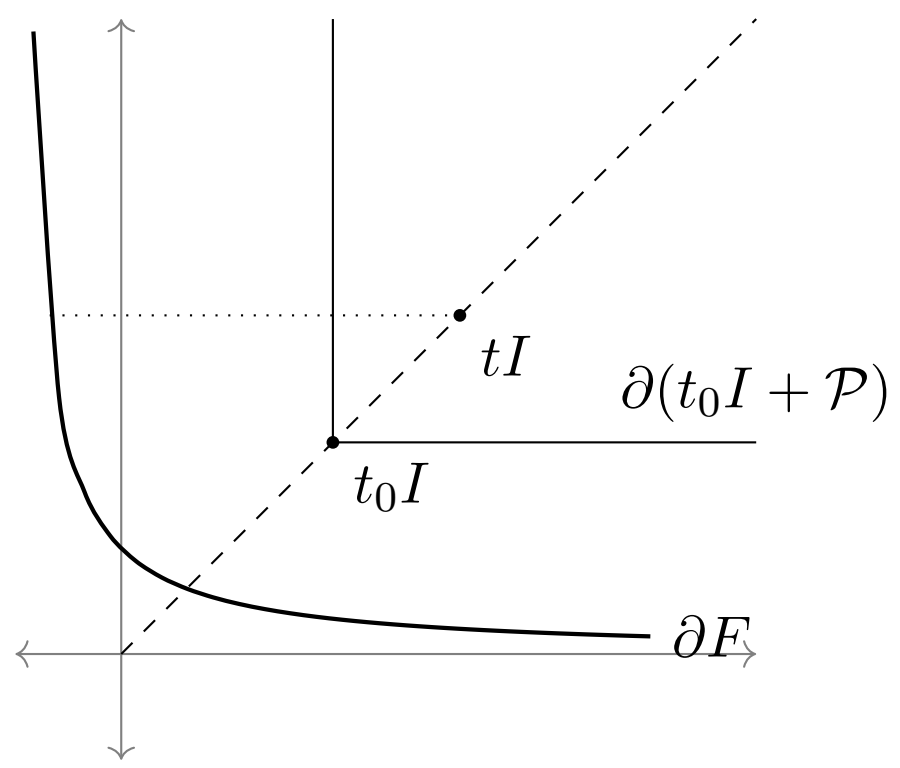
\begin{figure}[t]
    \centering
    \begin{tikzpicture}[scale=0.7]
        \draw[<->, white!50!black] (-1,0) -- (6,0);
        \draw[<->, white!50!black] (0,-1) -- (0,6); 
        \draw[domain=-0.83:5, thick, smooth, variable=\x] plot ({\x}, {1/(\x+1)}) node[right] {$\partial F$};
        \draw (2,6) -- (2,2) -- (6,2) node[above] {$\partial(t_0I + \pcal )$};
        \draw[dashed] (0,0) -- (6,6);
        \filldraw (2,2) circle (0.05) node[anchor=north west] {$t_0I$};
        \filldraw (3.2,3.2) circle (0.05) node[anchor=north west] {$tI$};
        \draw[dotted] (3.2,3.2) -- (-0.6875,3.2);
    \end{tikzpicture}
    \caption{Since $t_0I + \pcal  \subset F$, the point $tI$ is closer to $\partial(t_0 I + \pcal )$ than $\partial F$.}
    \label{fig:Subeq-Reg-sign-dist-linear-growth}
\end{figure}
Hence $tI$ is closer to $\partial(t_0+\pcal )$ than to $\partial F$ (see Figure~\ref{fig:Subeq-Reg-sign-dist-linear-growth}),
\begin{align}
    \signdist^{F}(tI) \ge&\ \signdist^{t_0I + \partial \pcal }(tI) \nonumber \\
    =&\ \inf_{A\in t_0I + \partial\pcal }|tI-A| \nonumber  \\
    =&\ \inf_{A\in \partial \pcal }|(t-t_0)I- A| \nonumber \\
    =&\ \signdist^{ \pcal }((t-t_0)I) \nonumber \\
    =&\ t-t_0, \qquad \forall\ t\ge t_0, \label{eq:Subeq-Reg-sign-dist-linear-growth-eq1}
\end{align}
where $\signdist^{\pcal}$ is the signed-distance to $\partial \pcal$ \eqref{eq:Subeq-Reg-positive-semi-def-subeq}, and the final line follows from Claim~\ref{claim:Subeq-Reg-linear-growth-claim}:
\begin{claim}
\label{claim:Subeq-Reg-linear-growth-claim}
     $\signdist^{\pcal}(P) = \lambda_{\min}(P)$, the smallest eigenvalue of $P$, for any $P\ge 0$.
\end{claim}
By continuity of $\signdist^F$ (part (i)),
\begin{equation}
    \label{eq:Subeq-Reg-sign-dist-linear-growth-eq1.00}
    \signdist^F(tI) -t \ge - \max_{s\in [0,t_0] }\left[ t- \signdist^F(sI) \right] , \qquad \forall\ t\in [0,t_0].
\end{equation}
Combining \eqref{eq:Subeq-Reg-sign-dist-linear-growth-eq2}--\eqref{eq:Subeq-Reg-sign-dist-linear-growth-eq1.00}, \eqref{eq:Subeq-Reg-signed-distance-lemma-linear-growth} holds for
\[M := \max\left\{\signdist^{F}(0),\ t_0,\ \max_{s\in [0,t_0] }\left[ t- \signdist^F(sI) \right] \right\}\]
It remains to prove Claim~\ref{claim:Subeq-Reg-linear-growth-claim}:
\begin{proof}[Proof of Claim~\ref{claim:Subeq-Reg-linear-growth-claim}]
    Suppose $\lambda\ge0$ is the smallest eigenvalue of $P$ with a corresponding unit eigenvector $v\in \RR^n$, $|v|=1$. In other words,
    \[Pv = \lambda v, \qquad \text{and} \qquad \langle P w, w \rangle \ge \lambda |w|^2, \quad \forall\ w\in \RR^n.\]
    
    \vspace{0.5em}\noindent\textit{Upper bound:} Consider the matrix
    \[A_0:= P - \lambda vv^T.\]
    Since $Pv = \lambda v$, it follows that $A_0v = 0$ making $A_0$ degenerate. Furthermore
    \[\langle A_0w, w\rangle = \langle Pw,w\rangle - \lambda |\langle v,w\rangle|^2 \ge \langle Pw,w\rangle - \lambda |w|^2 \ge 0, \qquad \forall\ w\in \RR^n.\]
    Hence, $A_0$ is also positive semi-definite. Conclude that
    \[A_0 \in \partial \pcal .\]
    Finally,
    \[\signdist_{\pcal }(P) = \inf_{A\in \partial \pcal }|P-A| \le |P - A_0| = |\lambda vv^T| = \lambda.\]

    \vspace{0.5em}\noindent\textit{Lower bound:} Note that the trace inner norm \eqref{eq:Subeq-Reg-trace-inner-prod} bounds the quadratic forms of any matrix,
    \[\langle Bw_1,w_1\rangle \le |Bw_1|\le   |B| = \sqrt{\tr (BB^T)}, \qquad \forall \ w_1\in \RR^n, |w_1|=1\]
    This can be seen for example by extending $w_1$ to an orthonormal basis $w_1,\dots, w_n$ of $\RR^n$ and computing the trace in this basis to get
    \[|B| = \sqrt{\sum_{i=1}^n|Bw_i|^2}.\]
    Pick $A\in \partial \pcal $ and apply this bound to $B = P-A$,
    \[\langle (P-A)w_1,w_1\rangle \le |P-A|, \qquad \forall\ w_1\in \RR^n,\ |w_1|=1,\ A\in \partial \pcal .\]
    Since $A\in \partial \pcal $ it has a null eigenvector, say $Av_1 = 0$, $|v_1|=1$. Then,
    \[|P-A| \ge \langle Pv_1, v_1\rangle \ge \lambda, \qquad \forall\ A\in \partial \pcal ,\ Av_1 = 0. \]
    Taking an infimum over $A\in \partial \pcal $,
    \[\signdist_{\pcal }(P) \ge \lambda.\tag*{\qedhere}\]
\end{proof}

\vspace{0.5em}\noindent\textit{(iii):} The signed-distance $\signdist^{F}$ \eqref{eq:Subeq-Reg-signed-dist-def} is smooth on some neighborhood $U$ of the boundary $\partial F$ with unit length gradients \eqref{eq:Subeq-Reg-signed-distance-unit-norm-derivative} by a computation in local charts near $\partial F$. The reader is referred to Petersen \cite[p. 42]{Petersen-Riem-Geom-MR3469435} for a complete argument. Finally, since $\signdist^F \equiv 0$ on $\partial F$, $D\signdist^{F}(A)$ is a unit vector normal to the tangent plane $T_{A}\partial F$ and points inwards. In other words, 
\[D\signdist^{F}(A) = \nu^{\partial F}(A),\qquad \forall\ A\in \partial F.\]

    \vspace{0.5em}\noindent\textit{(iv):} Note that the norm $|\cdot|$ in \eqref{eq:Subeq-Reg-signed-dist-def} is with respect to the trace inner product \eqref{eq:Subeq-Reg-trace-inner-prod}, which is $O(n)$-invariant, and in particular $G$-invariant for any $G\subset O(n)$. Thus, when $F$ is $G$-invariant, so is $\signdist^{F}$ by \eqref{eq:Subeq-Reg-signed-dist-def}.

    \vspace{0.5em}\noindent\textit{(v):} First, by convexity, $F$ admits a supporting half-space through each boundary point \eqref{eq:Subeq-Reg-inward-supporting-hyperplane-interior}. Denote the set of supporting half-spaces for $F$,
    \begin{equation}
        \label{eq:Subeq-Reg-set-of-supporting-half-spaces}
        \SHS_F:= \left\{H\subset \Sym^2(\RR^n): H \text{ is a half-space},\ F\subset H,\  \partial  F \cap \partial H \neq \varnothing   \right\}.
    \end{equation}
    Half-spaces are parametrized by their inward normal, and a constant,
    \begin{equation}
        \label{eq:Subeq-Reg-half-space-notation}
        H_{B,c} := \left\{ A\in \Sym^2(\RR^n): \langle B,A\rangle \ge c \right\}, \qquad B\in \Sym^2(\RR^n)\setminus \{0\},\ c\in \RR,
    \end{equation}
    Note that $H_{B,c}$ is a subequation \eqref{eq:Subeq-Reg-subequation-def-prop-basic} if $B\ge 0$, and by Lemma~\ref{lemma:Subeq-Reg-convex-subequation-normals}~(i), $B\ge 0$ whenever $H_{B,c} \in \SHS_F$. And by the supporting hyperplane theorem~\cite[Theorem~11.2]{rockafeller-convex-analysis-MR1451876},
    \begin{equation}
        \label{eq:Subeq-Reg-signed-distance-convex-as-intersection}
        F = \bigcap_{H \in \SHS_F} H .
    \end{equation}
    It turns out that the signed-distance to $\partial F$ can be expressed as the infimum of signed-distances to the boundaries of its supporting half-spaces $\partial H$, $H\in \SHS_F$, each of which is an affine (and hence convex) function \eqref{eq:Subeq-Reg-signed-ditance-to-half-space-affine}:
    \begin{claim}
    \begin{equation}
        \label{eq:Subeq-Reg-signed-distance-to-convex-is-concave}
        \signdist^{F} = \inf_{H \in \SHS_F} \signdist^{ H } \qquad \text{on } \Sym^2(\RR^n).
    \end{equation}
    \end{claim}
    \begin{proof}
    For the whole argument, fix a matrix
    \[B\in \Sym^2(\RR^n).\]
    Since $\partial F$ is closed, the distance from a given matrix $B\in \Sym^2(\RR^n)$ to the boundary $\partial F$ is achieved at some matrix $A_0\in \partial F$, i.e.,
    \begin{equation}
        \label{eq:Subeq-Reg-signed-dist-achieved-at-perp}
        \signdist^{F}(B) = \pm |B-A_0|.
    \end{equation}
    Also, for a half-space of the form $H_{B_0,c}$ \eqref{eq:Subeq-Reg-half-space-notation}, with $B_0\neq 0$, the signed-distance $\signdist^{H_{B_0,c}}(B)$ is the length of the component of $B$ along the inward normal $B_0$, and in particular is an affine map,
    \begin{equation}
        \label{eq:Subeq-Reg-signed-ditance-to-half-space-affine}
        \signdist^{H_{B_0,c}}(B)
        =  \frac{\langle B_0, B \rangle -c}{|B_0|}, \qquad \forall\ B\in \Sym^2(\RR^n),\ B_0 \neq 0,\ c\in \RR.
    \end{equation}
    Since $B\in \Sym^2(\RR^n)$ is arbitrary, \eqref{eq:Subeq-Reg-signed-distance-to-convex-is-concave} is broken down into two statements,
    \begin{align}
        \signdist^F (B) =&\ \signdist^H (B), \qquad \text{for some } H\in \SHS_F, \label{eq:Subeq-Reg-signed-distance-to-convex-is-concave-part-1} \\ 
        \signdist^F (B) \le&\ \signdist^H (B), \qquad \forall\ H\in \SHS_F. \label{eq:Subeq-Reg-signed-distance-to-convex-is-concave-part-2} 
    \end{align}
    Equations \eqref{eq:Subeq-Reg-signed-distance-to-convex-is-concave-part-1}--\eqref{eq:Subeq-Reg-signed-distance-to-convex-is-concave-part-2} are proved separately.

    \vspace{0.5em}\noindent\textit{Part 1. \eqref{eq:Subeq-Reg-signed-distance-to-convex-is-concave-part-1}:} Recall the arbitrarily fixed matrix $B\in \Sym^2(\RR^n)$. There are three cases to consider:
    \[\text{(I)}\ B \notin F, \quad \text{(II)}\ B \in \partial F,\quad \text{(III)}\ B \in \intr F.\]
    \vspace{0.5em}\noindent\textit{Case (I):} Suppose $B\notin F$. Recall that the distance $\signdist^F(B)$ is achieved at $A_0\in \partial F$ \eqref{eq:Subeq-Reg-signed-dist-achieved-at-perp}. In other words, $A_0$ minimizes the distance to $B$ among all $A\in F$. In particular, for any $A\in F$, the quadratic
    \[[0,1]\ni t\mapsto |B-A_0 - t(A-A_0)|^2 = |B-A_0|^2 - 2t \langle B-A_0, A-A_0\rangle + t^2|A-A_0|^2\]
    achieves its minimum at $t=0$, since $A_0 + t(A-A_0) = (1-t)A_0 + tA \in F$ by convexity of $F$. Hence, its derivative at $t=0$ is non-negative, i.e.,
    \[ \langle A_0-B, A-A_0\rangle  = -\langle B- A_0, A-A_0\rangle \ge 0, \qquad \forall\ A\in F, \]
    with equality at $A=A_0$. In other words, the half-space $H_{A_0-B, \langle A_0-B,A_0\rangle }$ \eqref{eq:Subeq-Reg-half-space-notation}, is a supporting half-space,
    \[H_{A_0-B, \langle A_0-B,A_0\rangle } \in \SHS_F.\]
    But observe that, by \eqref{eq:Subeq-Reg-signed-ditance-to-half-space-affine},
    \[\signdist^{H_{A_0-B, \langle A_0-B,A_0\rangle }}(B) = \frac{\langle A_0-B, B \rangle - \langle A_0-B,A_0\rangle}{|A_0-B|} = \frac{\langle A_0 -B, B - A_0 \rangle}{|A_0 - B| } = -|A_0 - B|.\]
    So, by \eqref{eq:Subeq-Reg-signed-dist-achieved-at-perp},
    \[\signdist^F(B) = \signdist^{H_{A_0-B, \langle A_0-B,A_0\rangle }}(B), \qquad H_{A_0-B, \langle A_0-B,A_0\rangle } \in \SHS_F.\]
    Thus, \eqref{eq:Subeq-Reg-signed-distance-to-convex-is-concave-part-1} holds in Case (I).
    
    \vspace{0.5em}\noindent\textit{Case (II):} If $B\in \partial F$, by the supporting hyperplane theorem~\cite[Theorem~11.2]{rockafeller-convex-analysis-MR1451876}, there is a supporting half-space to $F$ at $B\in \partial F$, i.e., $H\in \SHS_F$ \eqref{eq:Subeq-Reg-set-of-supporting-half-spaces} such that $B\in \partial F \cap \partial H$. So by \eqref{eq:Subeq-Reg-signed-dist-def},
    \[\signdist^F(B) = \signdist^H(B) = 0.\]
    Thus, \eqref{eq:Subeq-Reg-signed-distance-to-convex-is-concave-part-1} holds in Case (II).
    
    \vspace{0.5em}\noindent\textit{Case (III):} Suppose $B\in \intr F$ and $\signdist^F(B)$ is achieved at $A_0 \in \partial F$ \eqref{eq:Subeq-Reg-signed-dist-achieved-at-perp}. By the supporting hyperplane  theorem~\cite[Theorem~11.2]{rockafeller-convex-analysis-MR1451876}, there is a supporting half-space $H\in \SHS_F$ \eqref{eq:Subeq-Reg-set-of-supporting-half-spaces} to $F$ at $A_0\in \partial F$, i.e.,
    \[F\subset H, \qquad A_0 \in \partial F\cap \partial H.\]
    Since $B\in F\subset H$, the distance from $B$ to $\partial F$ is smaller than or equal to the distance from $B$ to $\partial H$, i.e.,
    \[|B-A_0| = \signdist^F(B) = \inf_{A\in \partial F}|B - A|  \le \inf_{A\in \partial H} |B-A| = \signdist^H(B) \le |B-A_0|, \]
    by \eqref{eq:Subeq-Reg-signed-dist-def} and \eqref{eq:Subeq-Reg-signed-dist-achieved-at-perp}, since $A_0 \in \partial F\cap \partial H$. Thus, \eqref{eq:Subeq-Reg-signed-distance-to-convex-is-concave-part-1} holds in Case (III).

    \vspace{0.5em}\noindent\textit{Part 2. \eqref{eq:Subeq-Reg-signed-distance-to-convex-is-concave-part-2}:} Fix some $H\in \SHS_F$. By \eqref{eq:Subeq-Reg-set-of-supporting-half-spaces}, $F\subset H$. So, the argument is again divided into three cases,
    \[\text{(I)}\ B \in F, \quad \text{(II)}\ B \in H\setminus F,\quad \text{(III)}\ B \notin H.\]
    If $B\in F\subset H$, then the distance from $B\in F$ to $\partial F$ is smaller than the distance from $F$ to $\partial H$, i.e.,
    \begin{equation}
        \label{eq:Subeq-Reg-signed-dist-convex-eq1.1}
        \signdist^{F}(B)\le \signdist^{H}(B), \qquad \text{if } B\in F,\ H\in \SHS_F.
    \end{equation}
    Next, by the sign convention in \eqref{eq:Subeq-Reg-signed-dist-def},
    \begin{equation}
        \label{eq:Subeq-Reg-signed-dist-convex-eq1.2}
        \signdist^{F}(B)<0 \le \signdist^{H}(B),
        \qquad \text{if } B\in H\setminus F,\  H\in\SHS_F.
    \end{equation}
    And finally if $B\notin H$, since $F\subset H$,
    \[ \begin{aligned}
            \inf_{A\in\partial H}|B-A|
            =&\ \inf_{A\in H}|B-A| \\
            \le&\ \inf_{A\in F}|B-A|
            =\inf_{A\in\partial F}|B-A|,
            \qquad  \text{if } B\notin H,\ H\in \SHS_F.
        \end{aligned} \]
    So, by the sign convention in \eqref{eq:Subeq-Reg-signed-dist-def},
    \begin{equation}
        \label{eq:Subeq-Reg-signed-dist-convex-eq1.3}
        \signdist^{F}(B)\le \signdist^{H}(B),
        \qquad \text{if } B\notin H,\ H\in \SHS_F.
    \end{equation}
    Then \eqref{eq:Subeq-Reg-signed-distance-to-convex-is-concave-part-2} follows from \eqref{eq:Subeq-Reg-signed-dist-convex-eq1.1}--\eqref{eq:Subeq-Reg-signed-dist-convex-eq1.3}, and \eqref{eq:Subeq-Reg-signed-distance-to-convex-is-concave} follows from \eqref{eq:Subeq-Reg-signed-distance-to-convex-is-concave-part-1}--\eqref{eq:Subeq-Reg-signed-distance-to-convex-is-concave-part-2}.
    \end{proof}
     
    By \eqref{eq:Subeq-Reg-signed-distance-to-convex-is-concave} and \eqref{eq:Subeq-Reg-signed-ditance-to-half-space-affine}, $\signdist^{F}$ is an infimum of affine functions on $\Sym^2(\RR^n)$, and hence is concave.

    \vspace{0.5em}\noindent\textit{(vi):} Since $\signdist^{F}$ is smooth near $\partial F$, constant on $\partial F$, positive on $F$, and has unit length gradient \eqref{eq:Subeq-Reg-signed-distance-unit-norm-derivative}, its derivative equals the inward unit normal to $F$ on $\partial F$,
    \[D\signdist^{F}(A) = \nu^{\partial F}(A), \qquad \forall\ A\in \partial F.\]
    Thus \eqref{eq:Subeq-Reg-sign-dist-has-pos-grad} follows from strict positivity of the inward unit normal of $F$ \eqref{eq:Subeq-Reg-strictly-elliptic-subequation-pos-normals}.
    
    \vspace{0.5em}\noindent\textit{(vii)--(viii):} For the whole argument fix arbitrary matrices
    \[B,P\in \Sym^2(\RR^n), \qquad P\ge 0,\ P\neq 0.\]
    By \eqref{eq:Subeq-Reg-subequation-def-prop-basic} there are three cases: (a) $B, B+P \in F$, (b) $B,B+P \notin F$, and (c) $B \notin F,\ B+P \in F$.
    
    \vspace{0.5em}\noindent\textit{Case (a)} $B, B+P \in F$: Since $\partial F$ is closed, there exists $A_0\in \partial F$ attaining the infimum defining $\signdist^F(B+P)$ \eqref{eq:Subeq-Reg-signed-dist-def}, i.e.,
    \[\signdist^F(B+P) = \inf_{A\in \partial F} |B+P - A| = |B+P-A_0|. \]
    By \eqref{eq:Subeq-Reg-subequation-def-prop-basic} observe that $\intr F + \pcal \subset \intr F$ \cite[p.~11, (3)]{HarLaw-dirichletduality}. So since $A_0 \in \partial F$,
    \[A_0-P' \notin \intr F , \qquad \forall\ P'\ge 0.\]
    In particular, $A_0-P \in \overline{F^c}$. So by \eqref{eq:Subeq-Reg-signed-dist-def},
    \[\begin{aligned}
        \signdist^F(B+P) =&\ |B+P-A_0|\\
        =& |B-(A_0-P)|\\
        \ge&\ \inf_{A\in \overline{F^c}} |B-A|\\
        =&\ \inf_{A\in \partial F}|B-A| = \signdist^F(B).
    \end{aligned}\]
    Hence \eqref{eq:Subeq-Reg-basic-monotonicity-for-sign-dist} holds in case (a).

    If in addition \eqref{eq:Subeq-Reg-non-degenerate-subequation} holds, then
    \[A_0-P \notin F.\]
    The distance from $B\in F$ to $A_0-P\in F^c$ is strictly larger than the distance from $B$ to $\partial F$ since $F^c$ is open. To see this, note that the line segment joining $B$ to $A_0-P$ intersects $\partial F$. Hence,
    \[\begin{aligned}
        \signdist^F(B+P) =&\ |B+P-A_0|\\
        =& |B-(A_0-P)|\\
        >&\ \inf_{A\in \partial F}|B-A| = \signdist^F(B).
    \end{aligned}\]
    Hence, in case (a), \eqref{eq:Subeq-Reg-sign-dist-is-strictly-elliptic} holds when \eqref{eq:Subeq-Reg-non-degenerate-subequation} does.
    
    \vspace{0.5em}\noindent\textit{Case (b)} $B,B+P\notin F$: Once again, since $\partial F$ is closed, there exists $A_0\in \partial F$ achieving the distance defining $\signdist^F(B)$ \eqref{eq:Subeq-Reg-signed-dist-def},
    \[ \signdist^F(B) = -\inf_{A\in \partial F} |B -A| = -|B-A_0|.\]
    Since $A_0\in F$, $A_0+P\in F$ by \eqref{eq:Subeq-Reg-subequation-def-prop-basic}. So,
    \[\begin{aligned}
        \signdist^F(B) =&\ -|B-A_0|\\
        =&\ -|(B+P) - (A_0+P)|\\
        \le&\ -\inf_{A\in F}|B+P - A|\\
        =&\ -\inf_{A\in \partial F}|B+P - A| = \signdist^F(B+P),
    \end{aligned}\]
   Thus, \eqref{eq:Subeq-Reg-basic-monotonicity-for-sign-dist} follows in case (b).

   If in addition \eqref{eq:Subeq-Reg-non-degenerate-subequation} holds, then
   \[A_0 + P \in \intr F.\]
   Since $B+P \in F^c$ and $A_0+P\in \intr F$, $B+P$ is strictly closer to $\partial F$ than to $A_0+P$. So,
   \[\begin{aligned}
        \signdist^F(B) =&\ -|B-A_0|\\
        =&\ -|(B+P) - (A_0+P)|\\
        <&\  - \inf_{A\in \partial F} |B+P - A| = \signdist^F(B+P).
   \end{aligned}\]
    Hence, in case (b), \eqref{eq:Subeq-Reg-sign-dist-is-strictly-elliptic} holds when \eqref{eq:Subeq-Reg-non-degenerate-subequation} does.
    
    \vspace{0.5em}\noindent\textit{Case (c)} $B\notin F,\ B+P \in F$: In this case, by \eqref{eq:Subeq-Reg-signed-dist-def}
    \begin{equation}
        \label{eq:Subeq-Reg-monotonicity-of-sign-dist-vii-case-c-eq1}
        \signdist^F(B)<0\le \signdist^F(B+P).
    \end{equation}
    So, both \eqref{eq:Subeq-Reg-basic-monotonicity-for-sign-dist} and \eqref{eq:Subeq-Reg-sign-dist-is-strictly-elliptic} holds in case (c).
\end{proof}

\subsection{Linearization of signed-distance}
\label{subsection:Subeq-Reg-lineqrization-of-signed-distance}

Recall that the fully non-linear PDE \eqref{eq:Subeq-Reg-signed-dist-fcal-PDE} is said to be strictly elliptic, if the associated linear PDE at any solution $u\in C^2(\overline{\Omega})$,
\[v\mapsto \left. \frac{d}{dt}\right|_{t=0} \signdist^{F} (D^2u + tD^2v) = \langle D\signdist^F(D^2u), D^2v \rangle,\]
is a strictly elliptic linear PDE. The matrix valued function $D\signdist^F(D^2u)$ is called the \textit{symbol} of the non-linear PDE \eqref{eq:Subeq-Reg-signed-dist-fcal-PDE}. The symbol of \eqref{eq:Subeq-Reg-signed-dist-fcal-PDE} is positive definite when convex and non-degenerate, by Lemma~\ref{lemma:Subeq-Reg-subequation-is-elliptic}~(iii):

\begin{claim}
        \label{claim:Subeq-Reg-dij-is-elliptic-PDE}
        For a strictly elliptic \eqref{eq:Subeq-Reg-strictly-elliptic-subequation-pos-normals} subequation $F \subset \Sym^2(\RR^n)$ with smooth boundary, a solution $u \in {C^2}({\Omega})$ to \eqref{eq:Subeq-Reg-const-coeff-subequation-PDE},
        and signed-distance $\signdist^{F}$ \eqref{eq:Subeq-Reg-signed-dist-def}, the linear partial differential operator $\Delta_{D\signdist^{F}(D^2u)}$
        \begin{equation}
            \label{eq:Subeq-Reg-linearization-of-dij-elliptic-pde}
            \Delta_{D\signdist^{F}(D^2u)} v:= \left\langle   D\signdist^{F}(D^2u), D^2v \right\rangle =  \sum_{i,j=1}^n \signdist^{F}_{ij}(D^2u) \partial_i \partial_j v, \qquad v\in C^2({\Omega})
        \end{equation}
        is a linear elliptic PDE on $\overline{\Omega}$, where $D\signdist^{F} = (\signdist^{F}_{ij})$ denotes the coordinates of the symmetric matrix $D\signdist^{F}$.
    \end{claim}
    \begin{proof}
        By strict ellipticity \eqref{eq:Subeq-Reg-strictly-elliptic-subequation-pos-normals} and Lemma~\ref{lemma:Subeq-Reg-signed-distance-lemma}~(vi),
    \begin{equation}
        \label{eq:Subeq-Reg-interior-grad-est-eq1}
        D\signdist^{F}(A) = \nu^{\partial F}(A) >0, \qquad \forall\ A\in \partial F,
    \end{equation}
    where $\nu^{\partial F}(A)$ \eqref{eq:Subeq-Reg-subequation-inward-normal}---the inward unit normal to $\partial F \subset \Sym^2(\RR^n)$---is a symmetric matrix. Hence, the linear PDE $\Delta_{D\signdist(D^2u)}$ \eqref{eq:Subeq-Reg-linearization-of-dij-elliptic-pde} has a positive definite symbol $D\signdist(D^2u) = (\signdist_{ij}(D^2u))$, making it elliptic.
    \end{proof}

    

\subsection{Linearization and symmetry}

The BVP \eqref{eq:Subeq-Reg-const-coeff-subequation-PDE} is purely second order, and in particular is invariant translations of the domain $\Omega$. Similarly, when $F$ is $G$-invariant for a group $G\subset O(n)$, the BVP \eqref{eq:Subeq-Reg-const-coeff-subequation-PDE} is invariant under rotations of $\Omega$ under elements of $G$. If the vector field $\tau:\RR^n \rightarrow \RR^n$ is an infinitesimal generator of such a symmetry of $F$, then differentiating \eqref{eq:Subeq-Reg-signed-dist-fcal-PDE} along these symmetries, the derivative
\begin{equation}
    \label{eq:Subeq-Reg-C2-TN-utau-def}
    \partial_{\tau} u(x) := \langle Du(x), \tau(x)\rangle,
\end{equation}
is a solution to the linearized equation:

\begin{lemma}
    \label{lemma:Subeq-Reg-symmetries-of-linearization}
    For a group $G \subset O(n)$, a $G$-invariant subequation $F\subset \Sym^2(\RR^n)$ with smooth boundary $\partial F$, and a solution $u\in C^3(\overline{\Omega})$ to the BVP \eqref{eq:Subeq-Reg-const-coeff-subequation-PDE},
    \begin{equation}
        \label{eq:Subeq-Reg-symmetries-of-linearization}
        \Delta_{D\signdist^F(D^2 u)} \partial_{\tau_{e,Q}} u = 0,
    \end{equation}
    where $\Delta_{D\signdist^F(D^2u)}$ \eqref{eq:Subeq-Reg-linearization-of-dij-elliptic-pde} is the linearized operator, and $\tau_{e,Q}:\RR^n \rightarrow \RR^n$ is the vector field
    \begin{equation}
        \label{eq:Subeq-Reg-tau-VQ-vector-field-def}
        \tau_{e,Q}(x) := e +Qx, \qquad \forall\  e\in \RR^n,\ Q\in \mathfrak{g}.
    \end{equation}
\end{lemma}
A slightly less general statement is due to Caffarelli--Nirenberg--Spruck~\cite[(5.4)]{caffarelli-nirenberg-spruck-dirichlet-problem-III-MR0806416}.
\begin{proof}
    First, differentiating \eqref{eq:Subeq-Reg-signed-dist-fcal-PDE} and applying \eqref{eq:Subeq-Reg-linearization-of-dij-elliptic-pde}, observe that
    \begin{equation}
    \label{eq:Subeq-Reg-constant-vector-field-and-linearization}
        \begin{aligned}
            0=\partial_e(\signdist^F(D^2u)) =&\ \langle D\signdist^F(D^2u), \partial_e  D^2u\rangle \\
            =&\ \langle D\signdist^F(D^2u), D^2\partial_e u\rangle = \Delta_{D\signdist^F(D^2u)} \partial_e u, \qquad \forall\ e\in \RR^n.
        \end{aligned}
    \end{equation}
    
    Next, fix $e\in \RR^n$ and $Q\in \mathfrak{g} = T_IG$. By the definition of the Lie algebra $\mathfrak{g}$, the exponential flow of $Q$ lies in $G$, i.e.,
    \[e^{sQ} \in G, \qquad \forall\ s\in \RR.\]
    The corresponding family of rotations on $\RR^n$ is the integral flow of the vector field $x\mapsto Qx$,
    \[\left. \frac{d}{ds} \right|_{s=0} e^{sQ}x = Qx = \tau_{e,Q}(x) - e,\]
    since $\tau_{e,Q} = e +Qx$ \eqref{eq:Subeq-Reg-tau-VQ-vector-field-def}. Let
    \[u_s(x) := u(e^{sQ}x), \qquad x\in \Omega,\ |s|\ll 1,\]
    and observe that
    \[\left. \frac{d}{ds} \right|_{s=0} u_s(x) = \left. \frac{d}{ds} \right|_{s=0} u(e^{sQ}x) = \left\langle Du(x),  \left. \frac{d}{ds} \right|_{s=0} e^{sQ}x \right\rangle = \langle Du(x), Qx\rangle = \partial_{\tau_{e,Q}} u - \partial_{e} u. \]
    Since $F$ is $G$-invariant, so is $\signdist^F$ \eqref{eq:Subeq-Reg-signed-dist-def}. And because $e^{sQ} \in G$,
    \[\signdist^F(D^2u_s(x)) = \signdist^F(e^{-sQ} D^2u(e^{sQ} x) e^{sQ}) = \signdist^F(D^2u(e^{sQ} x)) = 0,\qquad x\in \Omega,\ |s|\ll 1.\]
    Differentiating along $s$ at $s=0$,
    \[\begin{aligned}
        0= \left. \frac{d}{ds}  \right |_{s=0}\signdist^F(D^2u_s) =&\ \left\langle D\signdist^F(D^2u),  \left. \frac{d}{ds} \right |_{s=0} D^2u_s  \right\rangle\\
        =&\ \left\langle D\signdist^F(D^2u), D^2  \left. \frac{d}{ds} \right |_{s=0} u_s  \right\rangle \\
        =&\  \left\langle D\signdist^F(D^2u), D^2  (\partial_{\tau_{e,Q}} u  - \partial_{e} u) \right\rangle\\
        =&\ \Delta_{ D\signdist^F(D^2u)} \partial_{\tau_{e,Q}} u - \Delta_{ D\signdist^F(D^2u)} \partial_{e} u,\\
        =&\ \Delta_{ D\signdist^F(D^2u)} \partial_{_{e,Q}} u + 0 \qquad \text{on } \Omega,
    \end{aligned}\]
    where the last line follows by \eqref{eq:Subeq-Reg-constant-vector-field-and-linearization}.
\end{proof}

\section{Conormal subequation and signed-distance}
\label{sec:Subeq-Reg-projected-subeq-and-signed-distance}

This section is an aside on the orthogonal projections of a subequation. The results of this section underpin the a priori boundary normal Hessian estimate on a solution $u$ to the \eqref{eq:Subeq-Reg-const-coeff-subequation-PDE}.

\begin{figure}[t]
    \centering
    \begin{subfigure}[t]{0.32\textwidth}
        \centering
        \begin{tikzpicture}[scale=1]
            \fill[white!90!black] plot[domain=0.25:4] ({\x}, {1/\x}) -- (4,4) -- (0.25,4);
            \fill[white!60!black] (0,0.1) -- (0,-0.1) -- (4,-0.1) -- (4,0.1) -- (0,0.1);
            \filldraw[fill=white] (0,0) circle (0.1);
            \draw[<->] (-1,0) -- (4,0) node[below] {$\lambda_1$};
            \draw[<->] (0,-1) -- (0,4) node[left] {$\lambda_2$};
            \draw[thick] plot[domain=0.25:4] ({\x}, {1/\x});
            \draw (4,0.5) node[anchor=south east] {$F$};
            \draw (2,0) node[below] {$\pi_+(F)$};
        \end{tikzpicture}
        \caption{$F=\{A\ge,\, \det A\ge 1\}$}
    \end{subfigure}
    \hfill
    \begin{subfigure}[t]{0.32\textwidth}
        \centering
        \begin{tikzpicture}[scale=1]
            \fill[white!90!black] (0,4) -- (0,0) -- (4,0) --  (4,4) -- (0,4);
            \fill[white!60!black] (0,0.1) -- (0,-0.1) -- (4,-0.1) -- (4,0.1) -- (0,0.1);
            \fill[fill=white!60!black] (0,0) circle (0.1);
            \draw[<->] (-1,0) -- (4,0) node[below] {$\lambda_1$};
            \draw[<->] (0,-1) -- (0,4) node[left] {$\lambda_2$};
            \draw[thick] (0,4) -- (0,0) -- (4,0);
            \draw (4,0.5) node[anchor=south east] {$F$};
            \draw (2,0) node[below] {$\pi_+(F)$};
        \end{tikzpicture}
        \caption{$F=\pcal$ \eqref{eq:Subeq-Reg-positive-semi-def-subeq}}
    \end{subfigure}
    \hfill
    \begin{subfigure}[t]{0.32\textwidth}
        \centering
        \begin{tikzpicture}[scale=1]
            \fill[white!90!black] plot[domain=-1:4] ({\x}, {3-\x}) -- (4,4) -- (-1,4);
            \fill[white!60!black] (-1,0.1) -- (-1,-0.1) -- (4,-0.1) -- (4,0.1) -- (-1,0.1);
            \draw[<->] (-1,0) -- (4,0) node[below] {$\lambda_1$};
            \draw[<->] (0,-1) -- (0,4) node[left] {$\lambda_2$};
            \draw[thick] plot[domain=-1:4] ({\x}, {3-\x});
            \draw (4,0.5) node[anchor=south east] {$F$};
            \draw (2,0) node[below] {$\pi_+(F)$};
        \end{tikzpicture}
        \caption{$F = \{\tr A\ge 1\}$}
    \end{subfigure}
    \caption{Subequation $F$ (light gray) and its conormal projection $\pi_+(F)$ \eqref{eq:Subeq-Reg-projection-onto-the-principal-minor} (dark gray). (a) $F$ is convex and non-degenerate, and $\pi_+(F) \subset \Sym^2(\RR^{n-1})$ is a proper open subset. (b) $F$ is convex but not non-degenerate. $\pi_+(F)$ is a closed. (c) $F$ is convex and non-degenerate, and $\pi_+(F) = \Sym^2(\RR^{n-1})$.} 
    \label{fig:Subeq-Reg-Trudinger-projection-operator}
\end{figure}
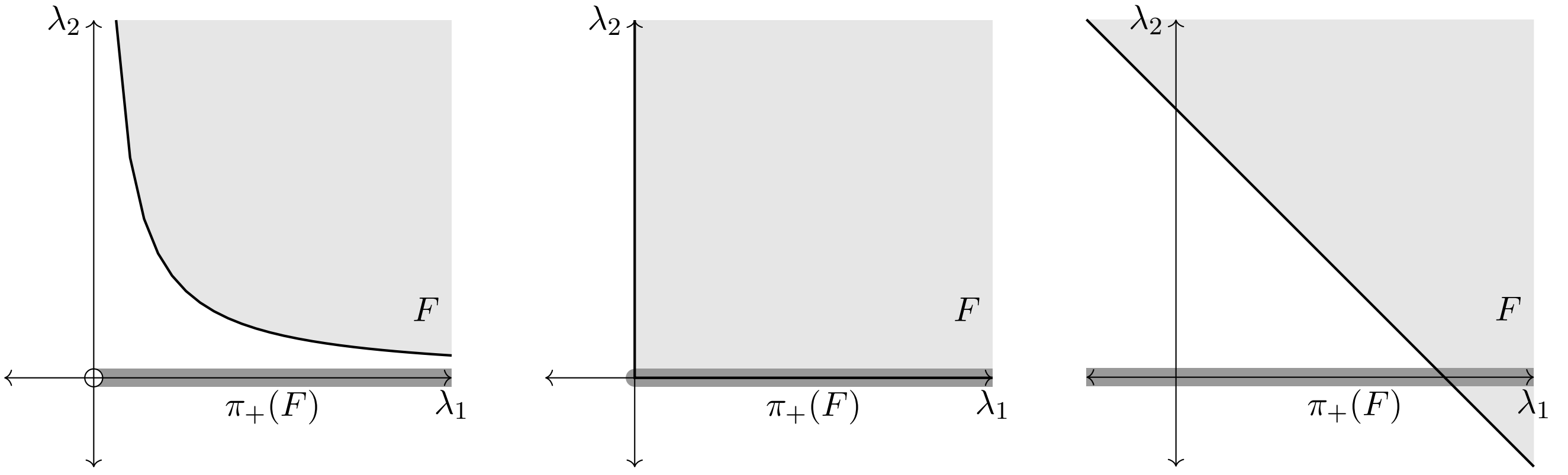

Consider the \textit{orthogonal projection} $\pi_+:\Sym^2(\RR^n) \rightarrow \Sym^2(\RR^{n-1})$,
\begin{equation}
    \label{eq:Subeq-Reg-projection-onto-the-principal-minor}
    \pi_+(A) :=  i_{e_n^\perp}^T \circ A \circ i_{e_n^\perp}, \qquad \forall\ A \in \Sym^2(\RR^{n}),
\end{equation}
where
\[i_{e_n^\perp}(x_1,\dots, x_{n-1}) = (x_1,\dots, x_{n-1},0),\]
and $i_{e_n^\perp}^T: \RR^n \rightarrow \RR^{n-1}$ is the orthogonal projection. In terms of a block decomposition of $A$,
\begin{equation}
    \label{eq:Subeq-Reg-projection-onto-the-principal-minor-block-repn}
    \pi_+\begin{pmatrix}
        A_+ & \vec{a} \\
        \vec{a}^T & a_{nn}
    \end{pmatrix} =  A_+, \qquad \forall\ A_+ \in \Sym^2(\RR^{n-1}),\ \vec{a} \in \RR^n,\ a_{nn} \in \RR.
\end{equation}

The \textit{conormal subequation} associated to a subequation $F\subset \Sym^2(\RR^n)$ is
\begin{equation}
    \label{eq:Subeq-Reg-F+-projected-subequation}
    F_+ := \overline{\pi_+(F)} \subset \Sym^2(\RR^{n-1}).
\end{equation}
The conormal subequation might be the whole space $F_+ = \Sym^2(\RR^{n-1})$. But when it is a proper subset, $F_+$ is a subequation, and is convex when $F$ is so (Lemma~\ref{lemma:Subeq-Reg-projection-F+-is-subeq}). Interestingly, when $F$ is convex and non-degenerate, $\pi_+(F) = \intr F_+$ (Lemma~\ref{lemma:Subeq-Reg-projection-of-conv-non-deg-is-open}). In particular,
\[\signdist^{F_+}(\pi_+(A)) > 0, \qquad \forall\ A\in \partial F,\]
where $\signdist^{F_+}:\Sym^2(\RR^{n-1}) \rightarrow \RR$ is the signed-distance to $\partial F^+$ \eqref{eq:Subeq-Reg-signed-dist-def}, and we note that $\signdist^{F_+} \equiv \infty$ if $F_+ = \Sym^2(\RR^{n-1})$. The main result of this section is a quantitative equivalence between bounds on $A\in \partial F$ and lower bounds on the signed-distance to the boundary of the conormal subequation, $\signdist^{F_+}\circ \pi_+(A)$ (Proposition~\ref{prop:Subeq-Reg-C2-NN-trudinger-prop}). These estimates are applied to estimate the boundary normal component of the Hessian $D^2u$ to a solution $u$ to \eqref{eq:Subeq-Reg-const-coeff-subequation-PDE} (see \S\ref{sec:Subeq-Reg-normal-hessian-estimates}).

\subsection{Conormal subequation}

This section shows that the conormal subequation $F_+$ \eqref{eq:Subeq-Reg-F+-projected-subequation} is a subequation, and $\pi_+(F) = \intr F_+$ when $F$ is convex and non-degenerate (Lemma~\ref{lemma:Subeq-Reg-projection-of-conv-non-deg-is-open}).

\begin{lemma}
    \label{lemma:Subeq-Reg-projection-F+-is-subeq}
    For a subequation $F\subset \Sym^2(\RR^n)$, and $F_+ = \overline{\pi_+(F)}$ \eqref{eq:Subeq-Reg-F+-projected-subequation} (see Definition~\ref{def:Subeq-Reg-subequations}),
    \begin{equation}
        \label{eq:Subeq-Reg-F+everything-or-subeq}
        \text{either} \qquad F_+ = \Sym^2(\RR^{n-1}) \qquad \text{or} \qquad F_+\text{ is a subequation.}
    \end{equation}
    Furthermore, if $F$ is convex, then so is $F$.
\end{lemma}
If $F_+ = \Sym^2(\RR^{n-1})$, then $\partial F_+ = \varnothing$ and by definition \eqref{eq:Subeq-Reg-signed-dist-def}, $\signdist^{F_+} \equiv \infty$. So if $F\subset \Sym^2(\RR^n)$ is a convex subequation, then by Lemma~\ref{lemma:Subeq-Reg-signed-distance-lemma}~(v) and Lemma~\ref{lemma:Subeq-Reg-projection-F+-is-subeq},
\begin{equation}
    \label{eq:Subeq-Reg-signdist-F+-inf-or-convex}
    \text{either} \qquad \signdist^{F_+}\equiv \infty \qquad \text{or}\qquad   \signdist^{F_+}:\Sym^2(\RR^{n-1})\rightarrow \RR \quad \text{is concave,}
\end{equation}
in the two respective cases of \eqref{eq:Subeq-Reg-F+everything-or-subeq}.
\begin{proof}
    Fix any point $A_+ \in \pi(F)$, say $A_+ = \pi_+(A)$ for some $A\in  F$. Since $F$ is a subequation \eqref{eq:Subeq-Reg-subequation-def-prop-basic},
    \[A + \begin{pmatrix}
        P_+ & 0 \\
        0 & 0
    \end{pmatrix} \in F \qquad \forall\ P_+ \in \Sym^2(\RR^{n-1}),\ P_+\ge 0. \]
    Hence,
    \[A_+ +P' = \pi_+\left(A + \begin{pmatrix}
        P' & 0 \\
        0 & 0
    \end{pmatrix} \right) \in \pi_+(F), \qquad \forall\ P'\ge 0.\]
    In other words, $\pi_+(F) + \pcal \subset \pi_+(F)$. Taking closures, conclude that $F_+ = \overline{\pi_+(F)}$ satisfies \eqref{eq:Subeq-Reg-subequation-def-prop-basic}. Being a closed, non-empty subset of $\Sym^2(\RR^{n-1})$ satisfying \eqref{eq:Subeq-Reg-subequation-def-prop-basic}, $F_+$ is a subequation when it is proper, i.e., when $F\neq \Sym^2(\RR^{n-1})$.
\end{proof}

\begin{remark}
    If we drop the properness assumption that subequations are proper subsets in Definition~\ref{def:Subeq-Reg-subequations}, not change much changes in the theory. The reason to exclude them is that the boundary $\partial F$, and hence the PDE \eqref{eq:Subeq-Reg-const-coeff-subequation-PDE}, are not well-defined when $F=\Sym^2(\RR^n)$.
\end{remark}

\begin{lemma}
    \label{lemma:Subeq-Reg-projection-of-conv-non-deg-is-open}
    For a convex, non-degenerate subequation $F\subset \Sym^2(\RR^n)$  \eqref{eq:Subeq-Reg-non-degenerate-subequation},
    \begin{equation}
        \label{eq:Subeq-Reg-projection-of-conv-non-deg-is-open-main}
        \pi_+(\partial F) = \pi_+(\intr F) = \pi_+ (F) = \intr F_+.
    \end{equation}
\end{lemma}
\begin{proof}
    Observe that the rank one orthogonal projection operator $P_{e_n} \in \Sym^2(\RR^n)$,
    \[P_{e_n} := e_n\circ e_n^T,\]
    lies in the kernel of $\pi_+$ \eqref{eq:Subeq-Reg-projection-onto-the-principal-minor},
    \begin{equation}
        \label{eq:Subeq-Reg-projection-of-conv-non-deg-is-open-eq0}
        \pi_+(P_{e_n}) = 0.
    \end{equation}
    But by non-degeneracy \eqref{eq:Subeq-Reg-non-degenerate-subequation} and \eqref{eq:Subeq-Reg-projection-of-conv-non-deg-is-open-eq0},
    \begin{equation}
        \label{eq:Subeq-Reg-projection-of-conv-non-deg-is-open-eq0.5}
        \pi_+(\partial F) = \pi_+(\partial F + P_{e_n}) \subset \pi_+( \intr F).
    \end{equation}
    It remains to prove the reverse inclusion. Fix an arbitrary element
    \begin{equation}
        \label{eq:Subeq-Reg-projection-of-conv-non-deg-is-open-eq0.75}
        A\in \intr F.
    \end{equation}
    Also pick any $A_0 \in \partial F$. By convexity of $F$, it admits a supporting half-space at $A_0$, i.e., $B_{A_0} \in \Sym^2(\RR^n)$ such that \eqref{eq:Subeq-Reg-inward-supporting-hyperplane},
    \begin{equation}
        \label{eq:Subeq-Reg-projection-of-conv-non-deg-is-open-eq1}
        F \subset \left\{B\in \Sym^2(\RR^n): \langle B, B_{A_0} \rangle \ge \langle A_0, B_{A_0} \rangle \right\}.
    \end{equation}
    Since $F$ is non-degenerate, by Lemma~\ref{lemma:Subeq-Reg-convex-subequation-normals} $B_{A_0} >0$. In particular, $\langle P_{e_n}, B_{A_0}\rangle >0$,
    \[\lim_{t\rightarrow \infty} \langle A-tP_{e_n}, B_{A_0}\rangle  = \langle A, B_{A_0} \rangle - \lim_{t\rightarrow \infty} t\langle P_{e_n}, B_{A_0}\rangle = -\infty < \langle A_0, B_{A_0}\rangle,\]
    where $A$ is the chosen matrix \eqref{eq:Subeq-Reg-projection-of-conv-non-deg-is-open-eq0.75}. So, by \eqref{eq:Subeq-Reg-projection-of-conv-non-deg-is-open-eq1}, there exists $t_0>0$ such that
    \[A-t_0P_{e_n} \notin F.\]
    But since $A\in \intr F$ \eqref{eq:Subeq-Reg-projection-of-conv-non-deg-is-open-eq0.75}, there exists $t_1\in [0,t_0)$ such that
    \[A-t_1P_{e_n} \in \partial F.\]
    So by \eqref{eq:Subeq-Reg-projection-of-conv-non-deg-is-open-eq0.5},
    \[\pi_+(A) = \pi_+(A-t_1P_{e_n}) \in \pi_+(\partial F).\]
    Since $A\in \intr F$ \eqref{eq:Subeq-Reg-projection-of-conv-non-deg-is-open-eq0.75} is arbitrary,
    \begin{equation}
        \label{eq:Subeq-Reg-projection-of-conv-non-deg-is-open-eq1.5}
        \pi_+(\partial F) \subset \pi_+(\intr F).
    \end{equation}
    Then \eqref{eq:Subeq-Reg-projection-of-conv-non-deg-is-open-main} follows from \eqref{eq:Subeq-Reg-projection-of-conv-non-deg-is-open-eq0.5}--\eqref{eq:Subeq-Reg-projection-of-conv-non-deg-is-open-eq1.5}, and because $F_+$ \eqref{eq:Subeq-Reg-F+-projected-subequation} is the closure of the open set $\pi_+(F)$.
\end{proof}

\subsection{Formula for conormal signed-distance}

This section proves a formula for the signed-distance $\signdist^{F_+}$ to $\partial F_+$ \eqref{eq:Subeq-Reg-F+-projected-subequation} in terms of $\signdist^F$ \eqref{eq:Subeq-Reg-formula-for-sign-dist-to-F+}. We begin with a simple characterization of $\intr F_+$ \eqref{eq:Subeq-Reg-F+-projected-subequation} in terms of $F$:

\begin{lemma}
\label{lemma:Subeq-Reg-proj-of-A-in-pi+-of-F}
    For a subequation $F \Sym^2(\RR^n)$, if $A\in \Sym^2(\RR^n)$ is such that $\pi_+(A) \in \intr F_+$, then there exists $R>0$ such that $A+RP_{e_n} \in F$.
\end{lemma}
\begin{proof}
     Suppose $\pi_+(A) \in \intr F_+$ \eqref{eq:Subeq-Reg-projection-of-conv-non-deg-is-open-main}. Then, there exists $\epsilon>0$ such that
     \[\pi_+(A-\epsilon I_n) = \pi_+(A) - \epsilon I_{n-1}  \in \intr F_+, \]
     where $I_n$ and $I_{n-1}$ are respective the $n$ and $n-1$ dimensional identity matrices. Then,
    \[\pi_+(A-\epsilon I_n) = \pi_+(B) \qquad \text{for some } B\in F.\]
    Use the block notation
    \[A = \begin{pmatrix}
        A_+ & \vec{a} \\
        \vec{a}^T & a_{nn}
    \end{pmatrix} \quad \text{and} \quad B := \begin{pmatrix}
        A_+ - \epsilon I_{n-1} & \vec{b} \\
        \vec{b}^T & b_{nn}
    \end{pmatrix}, \]
    for some $A_+\in \Sym^2(\RR^{n-1}),\ \vec{a}, \vec{b} \in \RR^n,\ a_{nn}, b_{nn} \in \RR$. Computing the quadratic form of the matrix $A +RP_{e_n} - B$, observe that
    \[A +RP_{e_n} - B = \begin{pmatrix}
        \epsilon I_{n-1} & \vec{a}- \vec{b} \\
        \vec{a}^T - \vec{b}^T & a_{nn} - b_{nn} + R
    \end{pmatrix} \ge 0,\]
    whenever
    \[R\ge R_0: = \epsilon^{-1}|\vec{a}- \vec{b}|^2  + |a_{nn}|+|b_{nn}|.\]
    By \eqref{eq:Subeq-Reg-subequation-def-prop-basic},
    \[A+ RP_{e_n} = B + (A+ RP_{e_n}  - B) \in F+ \pcal \subset F, \qquad \forall\ R\ge R_0.\tag*{\qedhere}\]
\end{proof}

Lemma~\ref{lemma:Subeq-Reg-proj-of-A-in-pi+-of-F} indicates that the asymptotics of $ R\mapsto \signdist^F(A+RP_{e_n})$ captures a lot of information about the projection $\pi_+(A)$. When $F$ is convex and non-degenerate, this can be made more precise: 

\begin{lemma}
\label{lemma:Subeq-Reg-C2-NN-projected-subeq-formula}
    For a convex and non-degenerate subequation $F\subset \Sym^2(\RR^{n-1})$ \eqref{eq:Subeq-Reg-non-degenerate-subequation}, the signed-distance $\signdist^{F_+}$ \eqref{eq:Subeq-Reg-signed-dist-def} to $\partial F_+$ satisfies
    \begin{equation}
        \label{eq:Subeq-Reg-formula-for-sign-dist-to-F+}
        \signdist^{F_+}\circ \pi_+ = \lim_{R\rightarrow \infty} \signdist^F(\, \cdot\, + RP_{e}).
    \end{equation}
\end{lemma}
In the block notation \eqref{eq:Subeq-Reg-projection-onto-the-principal-minor-block-repn}, \eqref{eq:Subeq-Reg-formula-for-sign-dist-to-F+} asserts that
\[\signdist^{F_+}(A_+) = \lim_{R\rightarrow \infty} \signdist^F \begin{pmatrix}
    A_+ & \vec{a} \\
    \vec{a}^T & R
\end{pmatrix}, \qquad \forall\ A_+\in \Sym^2(\RR^{n-1}),\ \vec{a} \in \RR^n.\]

\begin{proof}
    By \eqref{eq:Subeq-Reg-signdist-F+-inf-or-convex}, there are two cases to consider,
    \[F_+ = \Sym^2(\RR^{n-1}) \qquad \text{and} \qquad F_+ \subsetneq \Sym^2(\RR^{n-1}).\]
    When $F_+ = \Sym^2(\RR^{n-1})$, then $\partial F_+ = \varnothing$, and $\signdist^{F_+}\equiv \infty$ \eqref{eq:Subeq-Reg-signdist-F+-inf-or-convex}. We first prove \eqref{eq:Subeq-Reg-formula-for-sign-dist-to-F+} in this case.
    
    \vspace{0.5em}\noindent \textit{Step 1.} \textit{If $F_+ = \Sym^2(\RR^{n-1})$, then $ \lim_{R\rightarrow \infty} \signdist^F(\, \cdot\, + RP_{e_n}) \equiv \infty$:}

    \vspace{0.5em}\noindent Fix an arbitrary matrix $A\in \Sym^2(\RR^n)$. By the assumption $F_+ = \Sym^2(\RR^{n-1})$,
    \[\pi_+(A-tI_n) = \pi_+(A) - tI_{n-1} \in \Sym^2(\RR^{n-1}) = \intr F_+ , \qquad \forall\ t\in \RR.\]
    By Lemma~\ref{lemma:Subeq-Reg-proj-of-A-in-pi+-of-F} applied to $A-tI_n$, there exists $R_t$ such that
    \begin{equation}
        \label{eq:Subeq-Reg-formula-for-sign-dist-to-F+-eq0.4}
        A- tI_n + RP_{e_n} \in F, \qquad \forall\ R\ge R_t,\ t\in \RR. 
    \end{equation}
    Observe that $A' + tI_n \ge 0$ for any $|A'|\le t$. So by \eqref{eq:Subeq-Reg-formula-for-sign-dist-to-F+-eq0.4} and \eqref{eq:Subeq-Reg-subequation-def-prop-basic},
    \begin{equation}
        \label{eq:Subeq-Reg-formula-for-sign-dist-to-F+-eq0.5}
        A+RP_{e_n} -A' = A+RP_{e_n} -tI_n + (tI_n -A') \in F + \pcal \subset F, \qquad\ \forall\ |A'|\le t,\ R\ge R_t,\ t\ge 0.
    \end{equation}
    In other words, the ball of radius $t$ about $A+R_tP_{e_n}$ lies in $F$ for each $t\ge 0$. So the distance from $A+R_tP_{e_n}$ to $\partial F$ is at least $t$,
    \[\signdist^F(A+R_tP_{e_n}) \ge t, \qquad \forall\ t \ge 0.\]
    Since $A\in \Sym^2(\RR^n)$ is arbitrary,
    \[\lim_{R\rightarrow \infty} \signdist^F(\, \cdot\, + RP_{e_n}) \equiv \infty \qquad \text{on } \Sym^2(\RR^n).\]
    This completes Step 1. For the rest of the argument, assume that $F_+ \subsetneq \Sym^2(\RR^{n-1})$.
    
    \vspace{0.5em}\noindent\textit{Step 2.} \textit{If $F_+ \subsetneq \Sym^2(\RR^{n-1})$, then $\signdist^{F_+}\circ \pi_+ \ge \signdist^F$:}

    \vspace{0.5em}\noindent Fix a matrix $A\in F$, and use the block decomposition
    \[A := \begin{pmatrix}
        A_+ & \vec{a} \\
        \vec{a}^T & a_{nn}
    \end{pmatrix} \in F \subset \Sym^2(\RR^n), \qquad A_+\in \Sym^2(\RR^{n-1}),\ \vec{a}\in \RR^n,\ a_{nn} \in \RR.\]
    Since $F_+\subsetneq \Sym^2(\RR^{n-1})$, by \eqref{eq:Subeq-Reg-signed-dist-def},
    \[\signdist^{F_+}(\pi_+(A)) = |B_+-\pi_+(A)|, \qquad \text{for some } B_+\in \partial F_+.\]
    Since $\pi_+(F) = \intr F_+$ \eqref{eq:Subeq-Reg-projection-of-conv-non-deg-is-open-main},
    \[  B:= \begin{pmatrix}
        B_+ & \vec{a} \\
        \vec{a}^T & a_{nn}
    \end{pmatrix} \notin F,\]
    since $\pi_+(B) = B_+ \in \partial F_+$. Hence,
    \[\signdist^{F_+}(\pi_+(A)) = |B_+-\pi_+(A)| = |B-A| \ge \signdist^F(A).\]
    Conclude that
    \begin{equation}
        \label{eq:Subeq-Reg-F+-signdist-formula-step-1-eq-1}
        \signdist^{F_+}\circ \pi_+ \ge \signdist^F \qquad \text{on } F.
    \end{equation}

    Next, consider a matrix $A\notin F$. If $\pi_+(A) \in F_+$, then by the sign convention in \eqref{eq:Subeq-Reg-signed-dist-def}, $\signdist^{F_+}(\pi_+(A)) \ge 0 > \signdist^F(A)$. And if $\pi_+(A) \notin F_+$, then in particular, $A\notin F$ by \eqref{eq:Subeq-Reg-projection-of-conv-non-deg-is-open-main}, so by \eqref{eq:Subeq-Reg-signed-dist-def} and Lemma~\ref{lemma:Subeq-Reg-projection-of-conv-non-deg-is-open},
    \[\begin{aligned}
        \signdist^{F}(A) = -\inf_{B\in F} |A-B| \le&\ -\inf_{B\in F}\left| \pi_+(A) - \pi_+(B) \right| \\
        =&\ -\inf_{B_+\in \intr F_+}   \left| \pi_+(A) - B_+ \right| \\
        =&\ \signdist^{F_+}(\pi_+(A)).
    \end{aligned} \]
    Conclude that
    \begin{equation}
        \label{eq:Subeq-Reg-F+-signdist-formula-step-1-eq-2}
        \signdist^{F_+} \circ \pi_+ \ge \signdist^F \qquad \text{on } \Sym^2(\RR^n) \setminus F.
    \end{equation}
    Step 2 follows from \eqref{eq:Subeq-Reg-F+-signdist-formula-step-1-eq-1}--\eqref{eq:Subeq-Reg-F+-signdist-formula-step-1-eq-2}.

    \vspace{0.5em}\noindent \textit{Step 3.} \textit{If $F_+ \subsetneq \Sym^2(\RR^{n-1})$, then $ \signdist^{F_+}\circ \pi_+ \le  \lim_{R\rightarrow \infty} \signdist^F(\, \cdot\, + RP_{e_n})$:}

    \vspace{0.5em}\noindent Fix $A\in \Sym^2(\RR^n)$. We work in two cases: (i) $\pi_+(A) \in \intr F_+$ and (ii) $\pi_+(A) \notin \intr F_+$.
    
    \vspace{0.5em}\noindent\textit{(i):} Suppose $\pi_+(A) \in \intr F_+ = \pi_+(F)$ \eqref{eq:Subeq-Reg-projection-of-conv-non-deg-is-open-main}. By Lemma~\ref{lemma:Subeq-Reg-proj-of-A-in-pi+-of-F},
    \[\signdist^F(A+RP_{e_n}) \ge 0, \qquad \forall\ R\ge R_0,\]
    for some $R_0 \in \RR$, and by \eqref{eq:Subeq-Reg-signed-dist-def},
    \[\signdist^{F_+}(\pi_+(A)) \ge 0.\]
    For each $R\ge 0$ and $A\in \Sym^2(\RR^n)$, $\signdist^F(A+RP_{e_n})$ \eqref{eq:Subeq-Reg-signed-dist-def} is achieved at some $C_{A,R} \in \partial F$,
    \[\signdist^F(A+RP_{e_n}) =  |C_{A,R} - A- RP_{e_n}| \qquad \text{for some } C_{A,R} \in \partial F.\]
    By Step 2 and \eqref{eq:Subeq-Reg-projection-of-conv-non-deg-is-open-eq0},
    \[\begin{aligned}
        |C_{A,R} - A- RP_{e_n}| =&\ \signdist^F(A+RP_{e_n})  \\
        \le&\ \signdist^{F_+}\circ \pi_+(A+RP_{e_n}) = \signdist^{F_+} \circ \pi_+(A)<\infty, \quad R\ge R_0.
    \end{aligned}\]
    In particular, $C_{A,R} -  RP_{e_n}$, $R\ge R_0$, is a bounded family of matrices, and admits a subsequential limit,
    \[ C_{A} := \lim_{j\rightarrow \infty } C_{A,R_j} - R_jP_{e_n},\]
    for some sequence $R_j \rightarrow \infty$. In particular,
    \begin{align}
        \lim_{j\rightarrow \infty} \signdist^F(A+R_jP_{e_n}) =&\ \lim_{j\rightarrow \infty} |C_{A,R_j} - A- R_jP_{e_n}| \nonumber \\
        =&\ |C_A - A| \ge |\pi_+(C_A) - \pi_+(A)|. \label{eq:Subeq-Reg-formula-for-sign-dist-to-F+-eq1}
    \end{align}
    \begin{claim}
    \label{claim:Subeq-Reg-fomrula-for-F+-claim-in-bdry}
        $\pi_+(C_A) \in \partial \intr F_+$.
    \end{claim}
    \begin{proof}
        First, by \eqref{eq:Subeq-Reg-projection-of-conv-non-deg-is-open-main} and construction of $C_{A,R}$,
        \[ \pi_+(C_{A,R_j}- R_jP_{e_n}) = \pi_+(C_{A,R_j}) \in  \pi_+(\partial F) \subset F_+,\]
        and
        \[\lim_{j\rightarrow \infty} \pi_+(C_{A,R_j}- R_jP_{e_n}) = \pi_+(C_A).\]
        Conclude that $\pi_+(C_A) \in F_+$. Next, suppose on the contrary that
        \[\pi_+(C_A) \in \intr F_+ = \pi_+(\partial F).\]
        Then by Lemma~\ref{lemma:Subeq-Reg-proj-of-A-in-pi+-of-F} applied to $C_A-\epsilon I$, $C_A+RP_{e_n} \in \intr F$ for $R \gg 1$, and so,
        \[ \lim_{j\rightarrow \infty } \signdist^F(C_A + R_jP_{e_n})  >0. \]
        However, since $\signdist^F$ is 1-Lipschitz (Lemma~\ref{lemma:Subeq-Reg-signed-distance-lemma}~(i)) and $C_{A,R_j} \in \partial F$,
        \[\begin{aligned}
            \lim_{j\rightarrow \infty } \signdist^F(C_A + R_jP_{e_n}) \le&\ \lim_{j\rightarrow \infty } \signdist^F(C_{A,R_j}) + \lim_{j\rightarrow \infty }|C_{A} + R_jP_{e_n}- C_{A,R_j}|\\
            =&\ \lim_{j\rightarrow \infty } 0 + \lim_{j\rightarrow \infty }|C_{A} + R_jP_{e_n}- C_{A,R_j}| = 0. 
        \end{aligned}\]
        This contradicts the previous equation. Conclude that $\pi_+(C_A) \in \partial \intr F_+$.
    \end{proof}
    By \eqref{eq:Subeq-Reg-signed-dist-def}, \eqref{eq:Subeq-Reg-formula-for-sign-dist-to-F+-eq1}, and since $\pi_+(C_A)\in \partial F_+$ (Claim~\ref{claim:Subeq-Reg-fomrula-for-F+-claim-in-bdry}),
    \[\begin{aligned}
        \lim_{j\rightarrow \infty} \signdist^F(A+R_jP_{e_n}) \ge&\ |\pi_+(C_A) - \pi_+(A)| \\
        \ge&\ \inf_{B_+\in \partial F_+} |B_+ - \pi_+(A)| =  \signdist^{F_+}(\pi_+(A)).
    \end{aligned}\]
    But the subsequence $R_j \rightarrow \infty$ is arbitrary. Thus
    \begin{equation}
        \label{eq:Subeq-Reg-F+-formula-final-step-case-i}
        \signdist^{F_+}\circ \pi_+(A) \le  \lim_{R\rightarrow \infty} \signdist^F(A + RP_{e_n}), \qquad \forall\ A\in \Sym^2(\RR^n),\ \pi_+(A) \in \intr F_+.
    \end{equation}

    \vspace{0.5em}\noindent\textit{(ii):} Suppose
    $\pi_+(A)\notin \intr F_+$. Then signed-distance $\signdist^{F_+}(\pi_+(A))$ is achieved \eqref{eq:Subeq-Reg-signed-dist-def},
    \[\signdist^{F_+}(\pi_+(A))
    =-|B_+-\pi_+(A)| \qquad \text{for some } B_+\in \partial F_+.\]
    Note that, $\pi_+(A) \in \partial F_+$, then $B_+ = \pi_+(A)$. Use the block notation
    \[A=\begin{pmatrix}
        A_+ & \vec{a}\\
        \vec{a}^T & a_{nn}
    \end{pmatrix}, \qquad A_+=\pi_+(A) \notin F_+.\]
    For each $\epsilon>0$, let
    \[B_\epsilon:=\begin{pmatrix}
        B_+ +\epsilon I_{n-1} & \vec{a}\\
        \vec{a}^T & a_{nn}
    \end{pmatrix}.\]
    Since
    \[\pi_+(B_\epsilon)=B_+ +\epsilon I_{n-1}
    \in\intr F_+,\]
    by Lemma~\ref{lemma:Subeq-Reg-proj-of-A-in-pi+-of-F}, there exists $R_\epsilon \ge 0$ such that,
    \begin{equation}
        \label{eq:Subeq-Reg-formula-for-F+-final-step-eq1}
        B_\epsilon+RP_{e_n}\in F,
    \qquad \forall\ R\ge R_\epsilon.
    \end{equation}
    And since $\pi_+(A) \notin \intr F_+$, by \eqref{eq:Subeq-Reg-projection-of-conv-non-deg-is-open-main} and \eqref{eq:Subeq-Reg-projection-of-conv-non-deg-is-open-eq0},
    \[\pi_+(A) = \pi_+(A+RP_{e_n}) \notin \pi_+(F) = \intr F_+, \qquad \forall\ R\ge 0.\]
    So by Lemma~\ref{lemma:Subeq-Reg-proj-of-A-in-pi+-of-F}, $A+ RP_{e_n}\notin F$ for every $R\ge 0$. So by \eqref{eq:Subeq-Reg-signed-dist-def} and \eqref{eq:Subeq-Reg-formula-for-F+-final-step-eq1},
    \[\begin{aligned}
        \signdist^F(A+RP_{e_n})
        =&\ -\inf_{B\in F}|A+RP_{e_n}-B|\\
        \ge&\ -|A+RP_{e_n}-(B_\epsilon+RP_{e_n})|\\
        =&\ -|A_+-B_+-\epsilon I_{n-1}|,
        \qquad \forall\ R\ge R_\epsilon.
    \end{aligned}\]
    Letting $R\rightarrow\infty$ and then $\epsilon\rightarrow 0$,
    \[\lim_{R\rightarrow\infty}\signdist^F(A+RP_{e_n})
    \ge -|A_+-B_+| = \signdist^{F_+}(\pi_+(A)).\]
    Combined with Steps 1--2 and \eqref{eq:Subeq-Reg-F+-formula-final-step-case-i}, this proves
    \eqref{eq:Subeq-Reg-formula-for-sign-dist-to-F+}.
\end{proof}

\subsection{Lower bound}


This section is an aside on lower bounds to components of $A\in F$ for a convex subequation $F$. If $F$ is a convex, non-degenerate subequation and $A\in F$, then each diagonal element of $A$ can be bounded from below in terms of the remaining components of $A$. This plays a role in the proof of Proposition~\ref{prop:Subeq-Reg-C2-NN-trudinger-prop}, and is also used to obtain a lower bound for the normal component of the Hessian $D^2u$ in \S\ref{sec:Subeq-Reg-normal-hessian-estimates}.

\begin{lemma}
\label{lemma:Subeq-Reg-C2-NN-lower-bound-lin-alg}
    For a convex, non-degenerate subequation $F\subset \Sym^2(\RR^n)$ \eqref{eq:Subeq-Reg-non-degenerate-subequation} and $C_1>0$, there exists $C_2 = C_2(n, F,C_1)>0$ such that for any $e\in S^{n-1}$ and $A\in F$,
    \[ | \langle A v, w\rangle| \le C_1,\quad \forall\ v,w\in S^{n-1},\ w\perp e, \qquad \implies \qquad \langle Ae, e\rangle \ge -C_2.\]
\end{lemma}


\begin{proof}
    Fix $e\in S^{n-1}$. Extend $e$ to an orthonormal basis
    \[e_1,\dots, e_{n-1}, e_n := e.\]
    Fix $A_0 \in \partial F$. Since $F$ is convex, it admits a supporting half-space at $A_0$ \eqref{eq:Subeq-Reg-inward-supporting-hyperplane},
    \begin{equation}
        \label{eq:Subeq-Reg-C2-NN-lower-bound-matrix-eq0.25}
        F \subset \left\{B\in \Sym^2(\RR^n): \langle B,B_{A_0} \rangle \ge \langle A_0 ,B_{A_0} \rangle \right\}.
    \end{equation}
    By Lemma~\ref{lemma:Subeq-Reg-convex-subequation-normals},
    \begin{equation}
        \label{eq:Subeq-Reg-C2-NN-lower-bound-matrix-eq0.5}
         \left( b_{ij} \right)_{i,j=1}^n = B_{A_0} >0,
    \end{equation}
    where
    \[b_{ij} := \langle B_{A_0} e_i, e_j \rangle, \qquad i,j=1,\dots, n.\]
    By \eqref{eq:Subeq-Reg-C2-NN-lower-bound-matrix-eq0.25},
    \[\langle A , B_{A_0} \rangle  = \sum_{i,j=1}^n a_{ij} b_{ij} \ge \langle A_0, B_{A_0} \rangle, \qquad \forall\ A = (a_{ij})_{i,j=1}^n \in  F.\]
    Since $B_{A_0} >0$ and $e_n = e$ is a unit vector,
    \[b_{nn} = \langle B_{A_0} e_n, e_n \rangle = \langle B_{A_0} e, e\rangle \ge \lambda_{\min}(B_{A_0}) >0.\]
    Hence,
    \begin{equation*}
        a_{nn} \ge - \frac{1}{b_{nn}} \sum_{\substack{i,j=1\\(i,j)\neq (n,n)}}^n a_{ij}b_{ij} + \frac{\langle A_0,B_{A_0}\rangle}{b_{nn}}, \qquad \forall\ A = (a_{ij})_{i,j=1}^n \in F. 
    \end{equation*}
    In particular, if
    \[|\langle Av,w\rangle | \le C_1 \qquad \forall\ v,w \in S^{n-1}, w\perp e= e_n,\]
    when $|a_{ij}|\le C_1$ for $i,j=1,\dots, n$, such that $(i,j)\neq (n,n)$, and
    \[\begin{aligned}
        \implies \quad \langle Ae,e \rangle = a_{nn} \ge&\ -\frac{1}{b_{nn}} \sum_{\substack{i,j=1\\(i,j)\neq (n,n)}}^n |b_{ij}| C_1 - \frac{|\langle A_0, B_{A_0}\rangle|}{b_{nn}} \\
        \ge&\ -\frac{ n^2C_1|B_{A_0}| +|A_0||B_{A_0}|}{\lambda_{\min}(B_{A_0})}.
    \end{aligned}\]
    The constant in the final line depends only on $n$, $C_1$, $|A_0|$, $|B_{A_0}|$ and $\lambda_{\min}(B_{A_0})$, and in particular is independent of the vector $e\in S^{n-1}$ and chosen basis $e_1,\dots, e_n$. This completes the proof of Lemma~\ref{lemma:Subeq-Reg-C2-NN-lower-bound-lin-alg}.
\end{proof}

\subsection{Conormal signed-distance and estimates}

By Lemma~\ref{lemma:Subeq-Reg-projection-of-conv-non-deg-is-open}, $\pi_+(\partial F) = \intr F_+$. In particular, $\signdist^{F_+}\circ \pi_+(A)>0$ for any $A\in \partial F$. The goal of this section is to provide a quantitative version of this estimate: $A\in \partial F$ is uniformly bounded if and only if $\signdist^{F_+}\circ \pi_+(A)$ is uniformly positive:

\begin{proposition}
\label{prop:Subeq-Reg-C2-NN-trudinger-prop}
    For a convex, non-degenerate subequation $F\subset \Sym^2(\RR^n)$  \eqref{eq:Subeq-Reg-non-degenerate-subequation}, and matrix
    \[A = \begin{pmatrix}
        A_+ & \vec{a} \\
        \vec{a}^T & a_{nn}
    \end{pmatrix}\in \partial F, \]
    satisfying
    \begin{equation}
        \label{eq:Subeq-Reg-C2-NN-trudinger-prop-tan-mix-bound}
        |A_+|, |\vec{a}|\le C_1,
    \end{equation}
    for some constant $C_1>0$, if
    \begin{equation}
        \label{eq:Subeq-Reg-C2-NN-trudinger-prop-1}
        \signdist^{F_+} (A') \ge \delta
    \end{equation}
    for some $\delta>0$, then there exists a constant $C=C(\delta, n, F, C_1)>0$ such that
    \begin{equation}
        \label{eq:Subeq-Reg-C2-NN-trudinger-prop-2}
        |A|\le C.
    \end{equation}
    Conversely, if \eqref{eq:Subeq-Reg-C2-NN-trudinger-prop-2} holds for some $C>0$, then there exists $\delta = \delta(F,C)>0$ such that \eqref{eq:Subeq-Reg-C2-NN-trudinger-prop-1} holds.
\end{proposition}

Proposition~\ref{prop:Subeq-Reg-C2-NN-trudinger-prop} is adapted from a brief discussion of Trudinger~\cite[p. 158]{Trudinger1995-DHE}. It is used in \S\ref{sec:Subeq-Reg-normal-hessian-estimates} to reduce the boundary normal Hessian estimate to a special boundary point $x_{\min} \in \partial \Omega$ \eqref{eq:Subeq-Reg-Trudingers-point} (see Corollary~\ref{corollary:Subeq-Reg-C2-NN-trudinger-prop}).

\begin{proof}
    \noindent\underline{\eqref{eq:Subeq-Reg-C2-NN-trudinger-prop-1}$\implies$\eqref{eq:Subeq-Reg-C2-NN-trudinger-prop-2}:} First, note that by Lemma~\ref{lemma:Subeq-Reg-C2-NN-lower-bound-lin-alg}, there exists $C_2 = C_2(n,F,C_1)$ \eqref{eq:Subeq-Reg-C2-NN-trudinger-prop-tan-mix-bound} such that
    \begin{equation}
        \label{eq:Subeq-Reg-C2-NN-trud-prop-ann-lower-main}
        a_{nn} \ge -C_2.
    \end{equation}
    Combined with \eqref{eq:Subeq-Reg-C2-NN-trudinger-prop-tan-mix-bound}, \eqref{eq:Subeq-Reg-C2-NN-trudinger-prop-2} would follow from an upper bound on $a_{nn}$.
    
    By \eqref{eq:Subeq-Reg-formula-for-sign-dist-to-F+} Since $\signdist^{F_+}\circ \pi_+(A) = \signdist^{F_+}(A_+) \ge \delta$, by \eqref{eq:Subeq-Reg-C2-NN-trudinger-prop-1} there exists $R_1 = R_1(C_1)>0$ \eqref{eq:Subeq-Reg-C2-NN-trudinger-prop-tan-mix-bound} such that
    \begin{equation}
        \label{eq:Subeq-Reg-C2-NN-formula-for-fcal-tilde-eq1}
         \signdist^{F}\begin{pmatrix}
            A_+ & \vec{a} \\
            \vec{a}^T &  R
        \end{pmatrix} \ge \frac{\delta}{2}, \qquad \forall\ R\ge R_1.
    \end{equation}
    Note that $R_1$ can be chosen independent of $A_+$ depending only on the bound $|A_+|\le C_1$ \eqref{eq:Subeq-Reg-C2-NN-trudinger-prop-tan-mix-bound}, by continuity of $\signdist^{F_+}$ (Lemma~\ref{lemma:Subeq-Reg-signed-distance-lemma}~(i)).
    
    Since $A\in \partial F$, $\signdist^F(A) = 0$ by \eqref{eq:Subeq-Reg-signed-dist-def}. So by \eqref{eq:Subeq-Reg-C2-NN-formula-for-fcal-tilde-eq1},
    \[a_{nn} \le R_1 .\]
    Combined with \eqref{eq:Subeq-Reg-C2-NN-trudinger-prop-tan-mix-bound} and \eqref{eq:Subeq-Reg-C2-NN-trud-prop-ann-lower-main}, \eqref{eq:Subeq-Reg-C2-NN-trudinger-prop-2} follows.

    \vspace{0.5em}\noindent\underline{\eqref{eq:Subeq-Reg-C2-NN-trudinger-prop-2}$\implies$\eqref{eq:Subeq-Reg-C2-NN-trudinger-prop-1}:} By strict ellipticity of $\signdist^{F}$ \eqref{eq:Subeq-Reg-sign-dist-is-strictly-elliptic},
    \[\signdist^{F}(B+P_{e_n}) > \signdist^{F}(B), \qquad \forall\ B\in \Sym^2(\RR^n),\]
    where $P_{e_n} = e_n\circ e_n^T \in \Sym^2(\RR^n)$ is the rank one projection along $e_n$. For $C$ as in \eqref{eq:Subeq-Reg-C2-NN-trudinger-prop-2}, let
    \begin{equation}
        \label{eq:Subeq-Reg-C2-NN-trudinger-prop-eq3}
        \delta := \max \left\{ \signdist^{F}(B+P_{e_n}) - \signdist^{F}(B): B\in \Sym^2(\RR^n),\quad |B|\le C\right\}.
    \end{equation}
    Being the maximum of a continuous and positive function on a compact set,
    \[\delta>0.\]
    Furthermore, by \eqref{eq:Subeq-Reg-formula-for-sign-dist-to-F+} and \eqref{eq:Subeq-Reg-sign-dist-is-strictly-elliptic},
    \[\signdist^{F_+}(A_+) = \lim_{R\rightarrow \infty} \signdist^{F}(A+RP_{e_n})  \ge \signdist^{F}(A+P_{e_n}).\]
    And since $|A|<C$ \eqref{eq:Subeq-Reg-C2-NN-trudinger-prop-2}, and since $\signdist^F(A) = 0$,
    \[\signdist^{F_+}(A_+)  \ge \signdist^{F}(A+P_{e_n}) \ge \signdist^{F}(A) + \delta =0+\delta>0,\]
Thus \eqref{eq:Subeq-Reg-C2-NN-trudinger-prop-1} holds for $\delta = \delta(F,C)>0$ \eqref{eq:Subeq-Reg-C2-NN-trudinger-prop-eq3}.
\end{proof}

\section{Barriers}
\label{section:Subeq-Reg-barriers}

The primary technique for producing a priori boundary estimates to non-linear PDE is to construct barrier functions that glance the solution at the boundary of the domain and stay either above or below the solution on the rest of the domain. Such barriers immediately imply the existence of a $C^0$ and $C^1$ estimate.

\subsection{Harvey--Lawson boundary convexity}

This section introduces strict $\vec{F}$-convexity of the boundary $\partial \Omega$ of a domain $\Omega$. Proposition~\ref{prop:Subeq-Reg-HL-str-F-convex-implies-str-subsolution} is the only consequence of strict $\vec{F}$-convexity needed in this article.

\begin{definition}
    \label{def:Subeq-Reg-HL-ray-set}
    {\textrm{(Harvey--Lawson~\cite[Definition~5.8]{HarLaw-dirichletduality})}} The \textit{asymptotic  interior} of a subequation $F\subset \Sym^2(\RR^n)$ (Definition~\ref{def:Subeq-Reg-subequations}) is
    \begin{equation}
        \label{eq:Subeq-Reg-Fvec-set-def}
        \intr \vec{F} := \left\{ A\in \Sym^2(\RR^n): \exists \ t_0, \epsilon >0 \quad \text{such that} \quad t(A-\epsilon I) \in F\quad \forall\ t\ge t_0 \right\},
    \end{equation}
    whose closure is called the \textit{associated ray subequation} of $F$,
    \begin{equation}
        \label{eq:Subeq-Reg-associated-ray-subequation}
        \vec{F} := \overline{\intr \vec{F} }.
    \end{equation}
    $\vec{F}$ is subequation by \cite[Proposition~5.11]{HarLaw-dirichletduality}.
\end{definition}

\begin{definition}[Strict $\vec{F}$-convexity]
\label{def:Subeq-Reg-HL-boundary-convexity}
Fix a domain $\Omega\subset \RR^n$ with $C^2$ boundary, and a subequation $F\subset \Sym^2(\RR^n)$ with associated ray subequation $\vec{F}$ \eqref{eq:Subeq-Reg-associated-ray-subequation}.

\begin{enumerate}[label={(\roman*)}]
    \item $\rho \in C^2(\overline{\Omega})$ is called a \textit{global boundary defining function (bdf)} if
    \begin{equation}
        \label{eq:Subeq-Reg-bdf-def}
        \begin{aligned}
        \rho <0& \qquad \text{on } \Omega,\\
        \rho = 0& \qquad \text{on } \partial \Omega,\\
        D\rho(x)\neq0& \qquad \forall\ x\in \partial \Omega.
    \end{aligned} 
    \end{equation}
    
    \item {\textrm{(Harvey--Lawson~\cite[Definition~5.1]{HarLaw-dirichletduality})}} $\partial \Omega$ is \textit{strictly $\vec{F}$-convex} if there exists a global bdf $\rho \in C^2(\overline{\Omega})$ \eqref{eq:Subeq-Reg-bdf-def} such that
    \begin{equation}
        \label{eq:Subeq-Reg-str-F-convex-bdry-def}
        \forall\ x\in \partial \Omega \quad \exists\ B_x\in \intr\vec{F}\qquad \text{such that}\qquad  D^2\rho(x)|_{T_x\partial \Omega} = B_x|_{T_x\partial \Omega},
    \end{equation}
    where $B|_{T_x\partial \Omega} : {T_x\partial \Omega} \rightarrow {T_x\partial \Omega} $ denotes the composition of $B$ with the orthogonal projection $\pi_x : \RR^n \rightarrow T_x\partial \Omega$.
\end{enumerate}
\end{definition}

It turns out that if \eqref{eq:Subeq-Reg-str-F-convex-bdry-def} is true for one global bdf \eqref{eq:Subeq-Reg-bdf-def}, it is true for every global bdf, and strict $\vec{F}$-convexity of $\partial \Omega$ depends only on its second fundamental form $\II_{\partial \Omega}$ \cite[Lemma~5.3, Corollary~5.4]{HarLaw-dirichletduality}. However, we do not need these facts. The only consequence of strict $\vec{F}$-convexity of $\partial \Omega$ used in this article is:

\begin{proposition}
    \label{prop:Subeq-Reg-HL-str-F-convex-implies-str-subsolution}
    {\rm{(Harvey--Lawson~\cite[Theorem~5.12]{HarLaw-dirichletduality})}} For a subequation $F\subset \Sym^2(\RR^n)$ with associated ray subequation $\vec{F}$ \eqref{eq:Subeq-Reg-associated-ray-subequation}, and a bounded domain $\Omega\subset \RR^n$ with $C^2$, strictly $\vec{F}$-convex boundary $\partial \Omega$, there exists a global bdf $\rho \in C^2(\overline{\Omega})$ \eqref{eq:Subeq-Reg-bdf-def}, and constants $\epsilon, R_0>0$ such that
    \begin{equation}
        \label{eq:Subeq-Reg-str-F-convexity-implies-strong-F-bdf}
        R \left( D^2\rho(x) - \epsilon I \right) \in F,\qquad \forall\ x\in \overline{\Omega},\ R\ge R_0. 
    \end{equation}
\end{proposition}

\subsection{Strict subsolutions and supersolution}

This section constructs admissible strict subsolutions and supersolutions to \eqref{eq:Subeq-Reg-const-coeff-subequation-PDE}. Strict subsolutions exist when $\partial \Omega$ is strict $\vec{F}$-convex, by Proposition~\ref{prop:Subeq-Reg-HL-str-F-convex-implies-str-subsolution} (Proposition~\ref{prop:Subeq-Reg-existence-of-delta-subsol-and-supersol}~(i)), and strict supersolutions exist whenever $F$ is convex and non-degenerate \eqref{eq:Subeq-Reg-non-degenerate-subequation} (Proposition~\ref{prop:Subeq-Reg-existence-of-delta-subsol-and-supersol}~(ii)).

$u\in C^2(\Omega)$ is called an \textit{$F$-subsolution} if
\begin{equation}
    \label{eq:Subeq-Reg-F-subsolution-def}
    D^2u(x) \in F, \qquad \forall\ x\in \Omega,
\end{equation}
and an $F$-supersolution if
\begin{equation}
    \label{eq:Subeq-Reg-F-supersolution-def}
    -D^2u(x) \in \widetilde{F}, \qquad \forall\ x\in \Omega,
\end{equation}
i.e., if $-u$ is an $\widetilde{F}$-subsolution, where $\widetilde{F}$ \eqref{eq:Subeq-Reg-dual-to-a-subequation} is the dual subequation.

When $\partial \Omega$ is strictly $\vec{F}$-convex, by \eqref{eq:Subeq-Reg-str-F-convexity-implies-strong-F-bdf} it admits a global bdf $\rho \in C^2(\overline{\Omega})$ that is not just an $F$-subsolution, but is a \textit{strict} subsolution. In Theorem~\ref{thm:Subeq-Reg-Har-Law-exist-thm}, both strict $\vec{F}$- and $\vec{\widetilde{F}}$-convexity are assumed to ensure existence of a strict subsolution and a strict supersolution.

\begin{definition}
    \label{def:Subeq-Reg-str-subsol-supersol}
    Fix $\delta>0$, a domain $\Omega\subset \RR^n$, a boundary condition $\varphi \in C^2(\partial \Omega)$, and a subequation ${F} \subset \Sym^2 (\RR^n)$ with dual $\widetilde{F}$ \eqref{eq:Subeq-Reg-dual-to-a-subequation}.
    \begin{enumerate}[label={(\roman*)}]
        \item $\underline{u}\in C^2(\overline{\Omega})$ is called an \textit{admissible $\delta$-subsolution} to \eqref{eq:Subeq-Reg-const-coeff-subequation-PDE}
        \begin{equation}
            \label{eq:Subeq-Reg-strongly-admissible-subsolution}
            \begin{aligned}
                D^2\underline{u} -\delta I\in F\qquad &\text{on } \Omega,\\
                \underline{u} = \varphi\qquad& \text{on } \partial \Omega,
            \end{aligned}
        \end{equation}
        where $I$ is the identity matrix.

        \item $\overline{u}\in C^2(\overline{\Omega})$ is called an \textit{admissible $\delta$-supersolution} to \eqref{eq:Subeq-Reg-const-coeff-subequation-PDE} if
        \begin{equation}
            \label{eq:Subeq-Reg-strongly-admissible-supersolution}
            \begin{aligned}
                -D^2\overline{u} - \delta I\in \widetilde{F}\qquad &\text{on } \Omega,\\
                \overline{u} = \varphi\qquad& \text{on } \partial \Omega.
            \end{aligned}
        \end{equation}
    \end{enumerate}
\end{definition}

\begin{proposition}
\label{prop:Subeq-Reg-existence-of-delta-subsol-and-supersol}
    For $\delta>0$, a bounded domain $\Omega\subset \RR^n$ with $C^3$ boundary, a boundary condition $\varphi \in C^3(\partial \Omega)$, and a subequation ${F} \subset \Sym^2 (\RR^n)$, there exists a constant $C = C(F,\Omega, \|\varphi\|_{C^3})$ such that
    \begin{enumerate}[label={(\roman*)}]
        \item if $\partial \Omega$ is strictly $\vec{F}$-convex \eqref{eq:Subeq-Reg-str-F-convex-bdry-def}, then \eqref{eq:Subeq-Reg-const-coeff-subequation-PDE} admits an admissible $\delta$-subsolution $\underline{u}\in C^2(\overline{\Omega})$ \eqref{eq:Subeq-Reg-strongly-admissible-subsolution}, satisfying
        \[\|\underline{u}\|_{C^2} \le C.\]

        \item if $F$ is convex and non-degenerate \eqref{eq:Subeq-Reg-non-degenerate-subequation}, then \eqref{eq:Subeq-Reg-const-coeff-subequation-PDE} admits an admissible $\delta$-supersolution $\overline{u}\in C^2(\overline{\Omega})$ \eqref{eq:Subeq-Reg-strongly-admissible-supersolution}, satisfying
        \[\|\overline{u}\|_{C^2} \le C.\]
    \end{enumerate}
\end{proposition}

\begin{proof}
    \textit{(i):} By Proposition~\ref{prop:Subeq-Reg-HL-str-F-convex-implies-str-subsolution}, there exists a global bdf $\rho \in C^2(\overline{\Omega})$ \eqref{eq:Subeq-Reg-bdf-def} satisfying \eqref{eq:Subeq-Reg-str-F-convexity-implies-strong-F-bdf}, i.e., there exist $\epsilon, R_0>0$ such that
    \[\begin{aligned}
        \rho <0& \qquad \text{on } \Omega,\\
        \rho = 0& \qquad \text{on } \partial \Omega,\\
        D\rho(x)\neq0& \qquad \forall\ x\in \partial \Omega,\\
        R\left(D^2\rho(x) - \epsilon I \right)\in F& \qquad \forall\ x\in \partial \Omega, \ R\ge R_0.
    \end{aligned} \]
    Fix an arbitrary extension $\varphi \in C^2(\overline{\Omega})$, and a constant
    \[C_1>0\]
    to be determined, and consider the function
    \[\underline{u} := C_1\rho + \varphi \in C^2(\overline{\Omega}).\]
    First, since $\rho \equiv 0$ on $\partial \Omega$,
    \begin{equation}
        \label{eq:Subeq-Reg-str-F-convex-subol-eq1}
        \underline{u} = \varphi \qquad \text{on } \partial \Omega.
    \end{equation}
    Next, observe that
    \[D^2\underline{u} -\delta I = C_1 D^2\rho +D^2\varphi - \delta I.\]
    Choosing
    \begin{equation}
        \label{eq:Subeq-Reg-str-F-convex-subol-eq2}
        C_1 = \max \left\{ \frac{ \|\varphi\|_{C^2(\overline{\Omega})} + \delta}{\epsilon } ,\ R_0 \right\}
    \end{equation}
    observe that
    \begin{align}
        D^2\underline{u} -\delta I =&\ C_1 D^2\rho +D^2\varphi - \delta I \nonumber \\
        \ge&\ C_1D^2\rho - \left(\|\varphi\|_{C^2} + \delta \right)I \nonumber \\
        =&\ C_1 \left( D^2\rho - \frac{\|\varphi\|_{C^2} + \delta }{C} I \right)  \nonumber \\
        \ge&\ C_1 \left( D^2\rho - \epsilon I \right) \in F \qquad \text{on } \overline{\Omega}, \label{eq:Subeq-Reg-str-F-convex-subol-eq3}
    \end{align}
    where the last line follows since $C_1\ge R_0$ and $R(D^2\rho - \epsilon I) \in F$ for any $R\ge R_0$ \eqref{eq:Subeq-Reg-str-F-convexity-implies-strong-F-bdf}.
    
    By \eqref{eq:Subeq-Reg-subequation-def-prop-basic}, $A\ge B$ and $B\in F$ implies $A\in F$. Conclude by \eqref{eq:Subeq-Reg-str-F-convex-subol-eq3}, that
    \begin{equation}
    \label{eq:Subeq-Reg-str-F-convex-subol-eq4}
        D^2\underline{u} -\delta I \in F \qquad \text{on } \Omega.
    \end{equation}
    Hence, \eqref{eq:Subeq-Reg-strongly-admissible-subsolution} follows by \eqref{eq:Subeq-Reg-str-F-convex-subol-eq1} and \eqref{eq:Subeq-Reg-str-F-convex-subol-eq4}. Furthermore, since $\underline{u} = C_1\rho +\varphi$,
    \[\|\underline{u}\|_{C^2} \le C_1 \|\rho\|_{C^2} + \|\varphi\|_{C^2}.\]

    \vspace{0.5em}\noindent\textit{(ii):} Fix $A_0 \in \partial F$. By convexity, $F$ admits a supporting hyperplane at $A_0$, i.e., there exists $B_{A_0} \in \Sym^2(\RR^n)$ satisfying \eqref{eq:Subeq-Reg-inward-supporting-hyperplane} 
    \begin{equation}
        \label{eq:Subeq-Reg-str-F-convex-subol-eq5}
        \left\langle A, B_{A_0} \right\rangle \ge \left\langle A_0, B_{A_0} \right\rangle  \qquad \forall\ A\in F.
    \end{equation}
    Since $F$ is convex and non-degenerate \eqref{eq:Subeq-Reg-non-degenerate-subequation}, by Lemma~\ref{lemma:Subeq-Reg-convex-subequation-normals},
    \[B_{A_0}>0.\]
    In particular, the linear PDE
    \[\Delta_{B_{A_0}} v := \langle B_{A_0}, D^2v \rangle = \tr(B_{A_0} D^2v)\]
    is uniformly elliptic. By existence of solutions to uniformly linear elliptic PDE~\cite[Theorem~6.14]{gilbarg-trudinger-MR1814364}, there exists $\overline{u}\in C^2(\overline{\Omega})$ satisfying
    \begin{equation}
        \label{eq:Subeq-Reg-str-F-convex-subol-eq6}
        \begin{aligned}
        \Delta_{B_{A_0}} \overline{u} = \langle  A_0, B_{A_0} \rangle - \langle \delta I, B_{A_0}\rangle - 1& \qquad \text{on } \Omega,\\
        \overline{u} = \varphi,& \qquad \text{on } \partial \Omega,
    \end{aligned}
    \end{equation}
    along with a $C^2$ estimate
    \[\|\overline{u}\|_{C^2} \le C_2 = C_2(\Omega, \|\varphi\|_{C^3}, |A_0|,|B_{A_0}|,\delta).\]
    But by \eqref{eq:Subeq-Reg-str-F-convex-subol-eq6},
    \[ \langle D^2\overline{u} +\delta I , B_{A_0}\rangle = \Delta_{B_{A_0}}\overline{u} + \langle \delta I, B_{A_0}\rangle = \langle  A_0, B_{A_0} \rangle  - 1 < \langle  A_0, B_{A_0} \rangle.\]
    Conclude by \eqref{eq:Subeq-Reg-dual-to-a-subequation} and \eqref{eq:Subeq-Reg-str-F-convex-subol-eq5} that
    \begin{equation}
        \label{eq:Subeq-Reg-str-F-convex-subol-eq7}
        D^2\overline{u} +\delta I \in \Sym^2(\RR^n)\setminus F \subset -\widetilde{F}.
    \end{equation}
    Then, \eqref{eq:Subeq-Reg-strongly-admissible-supersolution} follows from \eqref{eq:Subeq-Reg-str-F-convex-subol-eq6} and \eqref{eq:Subeq-Reg-str-F-convex-subol-eq7}.
\end{proof}

\subsection{Barriers to the linearization}

    Proposition~\ref{prop:Subeq-Reg-existence-of-delta-subsol-and-supersol} constructs subsolutions and supersolutions to the PDE \eqref{eq:Subeq-Reg-signed-dist-fcal-PDE}. These subsolutions and supersolutions easily yield $C^1$ estimates (see Proposition~\ref{prop:Subeq-Reg-unif-grad-est}). However, to obtain estimates on higher derivatives of $u$, we need to differentiate \eqref{eq:Subeq-Reg-signed-dist-fcal-PDE}, and construct subsolution and supersolutions to the linearization $\Delta_{D\signdist^F(D^2u)}$ \eqref{eq:Subeq-Reg-linearization-of-dij-elliptic-pde}. The following construction of Guan~\cite[Lemma~3.1]{guan-subsol-MR1265805} shows that an admissible $\delta$-subsolution \eqref{eq:Subeq-Reg-strongly-admissible-subsolution} gives us subsolutions to the linearization \eqref{eq:Subeq-Reg-linearization-of-dij-elliptic-pde}:
    \begin{proposition}
    \label{prop:Subeq-Reg-domain-with-adm-subsol-guan-prop}
        For $\delta>0$, a bounded domain $\Omega \subset \RR^n$, a boundary condition $\varphi$, and a convex subequation $F\subset \Sym^2(\RR^n)$ with smooth boundary $\partial F$, if $\underline{u}\in C^2(\overline{\Omega})$ is an admissible $\delta$-subsolution \eqref{eq:Subeq-Reg-strongly-admissible-subsolution}, there exists $\epsilon = \epsilon(n,\Omega, F, \|\underline{u}\|_{C^2( \overline{\Omega})}, \delta)>0$ such that if $u\in C^2(\overline{\Omega})$ is a solution to \eqref{eq:Subeq-Reg-const-coeff-subequation-PDE}, then
        \begin{equation}
            \label{eq:Subeq-Reg-doman-with-adm-subsol-guan-prop-conclusion}
            \Delta_{D \signdist^{F}(D^2u)}(u - \underline{u}) \le -\epsilon <0, \qquad \text{on } \Omega,
        \end{equation}
        where $\Delta_{D \signdist^{F}(D^2u)}$ is the linearization \eqref{eq:Subeq-Reg-linearization-of-dij-elliptic-pde}.
    \end{proposition}
    \begin{proof}
        By \eqref{eq:Subeq-Reg-signed-dist-def} and \eqref{eq:Subeq-Reg-strongly-admissible-subsolution},
        \begin{equation}
        \label{eq:Subeq-Reg-adm-subsol-sign-dist-of-subsol-non-neg}
            \signdist^{F}(D^2\underline{u} - \delta I) \ge 0 \qquad \text{on } \overline{\Omega}.
        \end{equation}
        Hence, by \eqref{eq:Subeq-Reg-sign-dist-is-strictly-elliptic},
        \begin{equation}
            \label{eq:Subeq-Reg-domain-with-adm-subsol-guan-prop-eq1}
            \signdist^{F}(D^2\underline{u}) \ge   \signdist^{F}(D^2\underline{u} - \delta I) + \epsilon \ge \epsilon \qquad \text{on } \overline{\Omega},
        \end{equation}
        where
        \[\epsilon := \min_{x\in \overline{\Omega}} \left[ \signdist^{F}(D^2\underline{u}) - \signdist^{F}(D^2\underline{u} - \delta I) \right] >0.\]
        By concavity of $\signdist^{F}$ (Lemma~\ref{lemma:Subeq-Reg-signed-distance-lemma}~(v)) and smoothness near $\partial F$, and by \eqref{eq:Subeq-Reg-domain-with-adm-subsol-guan-prop-eq1},
        \[ 0<\epsilon \le \signdist^{F}(D^2\underline{u}) \le \signdist^{F}(D^2{u}) + \left\langle D\signdist^{F}(D^2u), D^2\underline{u} - D^2 u \right\rangle.\]
        Conclude by \eqref{eq:Subeq-Reg-signed-dist-fcal-PDE} and \eqref{eq:Subeq-Reg-linearization-of-dij-elliptic-pde}, that
        \[\Delta_{D\signdist^{F}(D^2u)}(u - \underline{u}) \le -\epsilon. \tag*{\qedhere}\]
    \end{proof}

    Proposition~\ref{prop:Subeq-Reg-domain-with-adm-subsol-guan-prop} is used to construct barriers to $\Delta_{D\signdist^F(D^2u)}$ \eqref{eq:Subeq-Reg-linearization-of-dij-elliptic-pde} \cite[pp. 410--411]{guan-subsol-MR1265805}:
    \begin{corollary}
    \label{corollary:Subeq-Reg-guan-corollary}
        Under the assumptions of Proposition~\ref{prop:Subeq-Reg-domain-with-adm-subsol-guan-prop}, if $\underline{u}\in C^2(\overline{\Omega})$ is an admissible $\delta$-subsolution \eqref{eq:Subeq-Reg-strongly-admissible-subsolution}, $u\in C^2(\overline{\Omega})$ a solution to \eqref{eq:Subeq-Reg-const-coeff-subequation-PDE}, $x_0\in \partial \Omega$, and $h\in C^2(\partial{\Omega})$, then there exists $C = C(n, \Omega, F, \|\underline{u}\|_{C^2}, \delta, \|h\|_{C^2})  >0$, and a function $w\in C^2(\overline{\Omega})$ such that
        \begin{equation}
            \label{eq:Subeq-Reg-guan-corollary-assump-conclusion}
            \begin{aligned}
            \Delta_{D\signdist^{F}(D^2u)}w \le&\ 0\qquad \text{on } \Omega\\
            w \ge&\ h \qquad \text{on } \partial \Omega,\\
            w(x_0) =&\ h(x_0),\\
            \|w\|_{C^1(\overline{\Omega})} \le&\ C.
        \end{aligned}
        \end{equation}
    \end{corollary}
    \begin{proof}
        Consider an arbitrary extension of $h\in C^2(\overline{\Omega})$ to $C^2(\overline{\Omega})$, and denote it again by $h\in C^2(\overline{\Omega})$. By Taylor's theorem, there exists $\delta_1>0$ such that
        \[h(x) \le h(x_0) + \langle Dh(x_0), x-x_0\rangle + \|h\|_{C_2}|x-x_0|^2, \qquad \forall\ |x-x_0|<  \delta_1. \]
        And if $|x-x_0|\ge\delta_1$, then
        \[\begin{aligned}
            h(x) - h(x_0) - \langle Dh(x_0), x-x_0\rangle \le&\  2\|h\|_{C^0} + 2\|h\|_{C^1}\sup_{\overline{\Omega}}|x| \\
            \le&\ \frac{2\|h\|_{C^0} + 2\|h\|_{C^1}\sup_{\overline{\Omega}}|x|}{\delta_1^2} |x-x_0|^2, \qquad \forall\ |x-x_0|\ge \delta_1.
        \end{aligned}\]
        So letting
        \[C_1 := \max\left\{ \|h\|_{C^2},  \frac{2\|h\|_{C^0} + 2\|h\|_{C^1}\sup_{\overline{\Omega}}|x|}{\delta_1^2}\right\},\]
        conclude that
        \begin{equation}
            \label{eq:Subeq-Reg-guan-cor-upper-quad-bound-on-h}
            h(x) \le h(x_0) + \langle Dh(x_0), x-x_0\rangle + C_1|x-x_0|^2.
        \end{equation}
        Consider the function
        \[w_A(x) := h(x_0)  + \langle Dh(x_0),x-x_0\rangle + C_1|x-x_0|^2 + A(u(x)-\underline{u}(x))  \qquad\forall\ x\in \overline{\Omega},\ A\ge 0.\]
        By \eqref{eq:Subeq-Reg-const-coeff-subequation-PDE} and \eqref{eq:Subeq-Reg-strongly-admissible-subsolution}, $u=\underline{u}$ on $\partial \Omega$. So by \eqref{eq:Subeq-Reg-guan-cor-upper-quad-bound-on-h},
        \[w_A(x) = h(x_0)  + \langle Dh(x_0),x-x_0\rangle + C_1|x-x_0|^2 \ge h(x), \qquad \forall\ x\in \partial \Omega,\ A\ge 0.\]
        Hence, the last three lines of \eqref{eq:Subeq-Reg-guan-corollary-assump-conclusion} immediately follow. Furthermore, by \eqref{eq:Subeq-Reg-trace-inner-prod}, \eqref{eq:Subeq-Reg-linearization-of-dij-elliptic-pde} and \eqref{eq:Subeq-Reg-doman-with-adm-subsol-guan-prop-conclusion},
        \[\begin{aligned}
            \Delta_{D\signdist^{F}(D^2u)}w_A =&\ A\Delta_{D\signdist^{F}(D^2u)}(u-\underline{u}) + C_1 \Delta_{D\signdist^{F}(D^2u)}|x-x_0|^2 \\
            =&\ A\Delta_{D\signdist^{F}(D^2u)}(u-\underline{u}) + 2C_1\langle D\signdist^{F}(D^2u), I \rangle \\
            =&\ A\Delta_{D\signdist^{F}(D^2u)}(u-\underline{u}) + 2C_1\tr(D\signdist^{F}(D^2u)) \\
            \le&\ -A\epsilon + 2C_1\tr(D\signdist^{F}(D^2u)).
        \end{aligned}\]
        By \eqref{eq:Subeq-Reg-signed-distance-unit-norm-derivative} and \eqref{eq:Subeq-Reg-sign-dist-is-strictly-elliptic},
        \[0< \tr(D\signdist^{F}(D^2u)) = \langle D\signdist^{F}(D^2u), I\rangle \le |D\signdist^{F}(D^2u)||I| =\sqrt{n}.\]
        Conclude that
        \[\Delta_{D\signdist^{F}(D^2u)}w_A \le -A\epsilon + 2C_1\sqrt{n}.\]
        Hence, for $A\ge 2C_1 \sqrt{n} \epsilon^{-1}$, $w_A$ satisfies \eqref{eq:Subeq-Reg-guan-corollary-assump-conclusion}.
    \end{proof}

\subsection{Uniformly convex domains}
\label{subsection:Subeq-Reg-unif-convex-domains}

This subsection constructs {admissible $\delta$-subsolutions} and supersolutions \eqref{eq:Subeq-Reg-strongly-admissible-subsolution}--\eqref{eq:Subeq-Reg-strongly-admissible-supersolution} on a uniformly convex domain.

\begin{definition}[Uniformly convex domains]
\label{def:Subeq-Reg-uniformly-convex-domain}
    A domain $\Omega \subset \RR^n$ with $C^2$ boundary $\partial \Omega$ is said to be \textit{uniformly convex} if there exists $\rho \in C^k(\overline{\Omega})$ and $\kappa>0$ such that
    \begin{equation}
        \label{eq:Subeq-Reg-unif-conv-dom-rho-prop}
        \begin{aligned}
        \rho<0 \quad &\text{on } \Omega,\\
        \rho=0 \quad &\text{on } \partial \Omega,\\
        D\rho \neq 0 \quad &\text{on } \partial \Omega, \\
        D^2 \rho \ge \kappa I \quad &\text{on } \Omega.
    \end{aligned} 
    \end{equation}
    In other words, $\rho$ is a uniformly convex bdf \eqref{eq:Subeq-Reg-bdf-def}.
\end{definition}

\begin{proposition}
\label{prop:Subeq-Reg-unif-convex-is-F-vec-for-all-F}
    A bounded, $C^2$ uniformly convex domain $\Omega \subset \RR^n$ is strictly $\vec{F}$-convex for any subequation $F\subset \Sym^2(\RR^n)$.
\end{proposition}
\begin{proof}
     Since $F$ is a subequation, there exists $t_0>0$ such that \cite[p.~11, (4)]{HarLaw-dirichletduality}
    \begin{equation}
        \label{eq:Subeq-Reg-subsol-supersol-eq1.0}
        \begin{aligned}
        tI \in F, \qquad \forall\ t\ge t_0.
    \end{aligned} 
    \end{equation}
    By uniform convexity of $\Omega$, there exists $\rho \in C^2(\overline{\Omega})$ and $\kappa>0$ such that \eqref{eq:Subeq-Reg-unif-conv-dom-rho-prop} holds. Then, for any $\epsilon \in (0,\kappa)$,
    \[C(D^2\rho(x)- \epsilon I) \ge C(\kappa I - \epsilon I) \ge t_0I \in F, \qquad \forall\ C\ge t_0/(\kappa-\epsilon), \ x\in \overline{\Omega}.\]
    By \eqref{eq:Subeq-Reg-Fvec-set-def},
    \[D^2\rho(x) \in \intr \vec{F} ,\qquad \forall\ x\in \overline{\Omega}.\]
    By \eqref{eq:Subeq-Reg-str-F-convex-bdry-def}, $\partial \Omega$ is strictly $\vec{F}$-convex.
\end{proof}

\begin{allproofs}
\subsection{Principal coordinate system}

The boundary of a $C^2$ domain can locally be parametrized as the graph of a smooth function:
\begin{lemma}
    \label{lemma:Subeq-Reg-unif-convex-boundary-coordinates}
    For a domain $ \Omega \subset \RR^n$ with $C^3$ boundary $\partial \Omega $ and origin on the boundary $0 \in \partial \Omega$, there exists an orthonormal basis $e_1,\dots, e_n$, open neighborhoods of the origin
    \[0\in U' \subset \RR^{n-1} \quad \text{and}\quad 0\in U \subset \RR^n,\]
    and a function $\chi\in C^3(U') $ such that $\partial \Omega$ is parametrized as the graph of $\chi$ near the origin,
    \[(\id,\chi)(U') = U\cap \partial \Omega,\]
    and such that $\chi$ has the form
    \begin{equation}
        \label{eq:Subeq-Reg-principal-coordinate-graph-sigma} 
        \chi(x') = \frac{1}{2}\sum_{i=1}^{n-1} \kappa_i x_i^2 + O(|x'|^3),
    \end{equation}
    for some real numbers $\kappa_1,\dots, \kappa_{n-1} \in \RR$. If $\Omega$ is in addition uniformly convex (Definition~\ref{def:Subeq-Reg-uniformly-convex-domain}), then
    \begin{equation}
        \label{eq:Subeq-Reg-principal-coordinate-unif-convex-kappa-positive}
        \kappa_1,\dots, \kappa_{n-1} >0.
    \end{equation}
\end{lemma}

\begin{figure}[t]
    \centering
    \begin{tikzpicture}[scale=3.5, domain=-1:1]
        \foreach \x in {-0.95,-0.90,...,1.00}
        {
            \draw[very thin, color=black!40!white] (\x,+0.5*\x*\x + 0.125*\x*\x*\x*\x ) -- (\x, 0.75);
        }
        \draw (0.2,-0.07) node {$x_0=0$};
        \draw (0.8,0.81) node {${\Omega}$};
        \draw[thick, <->] (-1,0) -- (1,0) node[anchor=north west] {$(x_1,\dots, x_{n-1})=x'$};
        \draw[thick, <->] (0, -0.1) -- (0, 0.85) node[above] {$x_n$}; 
        \draw[thin] (-0.7,0.05) -- (-0.7,-0.05) node[below] {$-\delta$};
        \draw[thin] (0.7,0.05) -- (0.7,-0.05) node[below] {$\delta$};
        \draw[domain=-1:1, thick, smooth, variable=\x] plot ({\x}, {0.5*\x*\x + 0.125*\x*\x*\x*\x}) node[right] {$\operatorname{Gr}(\chi)$};
    \end{tikzpicture}
    
    \caption{The domain ${\Omega}$ pictured in principal coordinates centered at $x_0 \in \partial \Omega$.}
    \label{fig:Subeq-Reg-uniform-conv-lin-barrier}
\end{figure}

\begin{proof}[Sketch of proof of Lemma~\ref{lemma:Subeq-Reg-unif-convex-boundary-coordinates}]
Since $\partial \Omega$ is a smooth manifold, it can locally be expressed as the graph of a smooth function $\chi$ over the last coordinate near $0\in \partial \Omega$. After a further rotation in the tangent plane $T_0\partial \Omega$ to diagonalize $D^2\chi(0)$, \eqref{eq:Subeq-Reg-principal-coordinate-graph-sigma} holds. Finally \eqref{eq:Subeq-Reg-principal-coordinate-unif-convex-kappa-positive} follows by differentiating the equation $\rho\circ (\id, \chi) = 0$ where $\rho$ is a uniformly convex boundary defining function for $\Omega$ \eqref{eq:Subeq-Reg-unif-conv-dom-rho-prop}. For more details, the reader is referred to Gilbarg--Trudinger~\cite[\S14.6]{gilbarg-trudinger-MR1814364}.
\end{proof}

\begin{definition}[Principal coordinate system]
    \label{def:Subeq-Reg-principal-coordinate-system}
    For a domain $\Omega\subset \RR^n$ with $C^3$ boundary $\partial \Omega$, \textit{a principal coordinate system at $x_0\in \partial \Omega$} is a choice of orthonormal coordinates $(x_1,\dots, x_n)$ on $\RR^n$ which places the origin at $x_0$, in which the inward unit normal to $\partial \Omega$ at $x_0$ is the coordinate vector $e_n=(0,\dots , 0,1)$, and where $\partial \Omega$ is expressed as the graph $x_n = \chi(x')$ for a $C^3$ function $\chi$ satisfying \eqref{eq:Subeq-Reg-principal-coordinate-graph-sigma}.
\end{definition}

\begin{proof}[Proof of Lemma~\ref{lemma:Subeq-Reg-unif-convex-boundary-coordinates}]
    Fix an extension $\varphi \in C^2(\overline{\Omega})$. Since $\partial \Omega$ is a $C^3$ manifold, it can be expressed as the graph of a $C^3$ function $\chi$ over some sets $U' \subset \RR^{n-1}$ and $U\subset \RR^n$. Using the coordinates described in the paragraphs right before Lemma~\ref{lemma:Subeq-Reg-unif-convex-boundary-coordinates}, it is clear that
    \[\chi(0) = 0. \]
    Furthermore, since the tangent space at the origin to the graph of $\chi$ are spanned by
    \[e_i + \partial_i\chi(0) e_n, \qquad i=1,\dots, n-1,\]
    and since $e_n$ is orthogonal to this tangent plane,
    \[D\chi(0) = 0.\]
    Finally, after an orthogonal rotation in $e_1, \dots, e_{n-1}$, $D^2\chi(0)$ can be assumed to be diagonal yielding \eqref{eq:Subeq-Reg-principal-coordinate-graph-sigma} . Note that $\kappa_1,\dots, \kappa_{n-1}$ are the eigenvalue of $D^2\chi(0)$,
    \begin{equation}
        \label{eq:Subeq-Reg-principal-coordinates-sigma-hess-coord-eq0.1}
        D^2\chi(0) = \diag(\kappa_1,\dots, \kappa_{n-1}).
    \end{equation}
    
    Finally, if $\Omega$ is uniformly convex and $\rho \in C^3(\overline{\Omega})$ is a uniformly convex boundary defining function satisfying \eqref{eq:Subeq-Reg-unif-conv-dom-rho-prop}, then
    \[\rho(x',\chi(x')) = 0, \qquad \forall\ x'\in U'.\]
    Compute the derivative and Hessian of this equation. For simplicity of notation, represent partial derivatives by subscripts, e.g.,
    \[\chi_i:= \partial_i\chi \quad \text{and}\quad \chi_{ij} := \partial_i\partial_j \chi.\]
    First, the derivative is given by
    \[ \rho_i(x',\chi(x')) + \rho_n(x',\chi(x')) \chi_i(x') = 0, \qquad i=1,\dots, n-1.\]
    Differentiating again,
    \[ \rho_{ij}(x',\chi(x')) +   \rho_{in} (x',\chi(x')) \chi_j(x') + \rho_{jn}(x',\chi(x')) \chi_i(x') + \rho_n(x',\chi(x'))  \chi_{ij} (x') =0, \qquad i,j=1,\dots, n-1. \]
    Evaluating at $x'=0$, since $\chi(0)= 0$ and $D\chi(0)=0$,
    \begin{equation}
        \label{eq:Subeq-Reg-principal-coordinates-sigma-hess-coord-eq1}
        \rho_{ij}(0) + \rho_n(0)  \chi_{ij}(0) = 0.
    \end{equation}
    But by \eqref{eq:Subeq-Reg-unif-conv-dom-rho-prop}, since $e_1,\dots, e_{n-1} \in T_0\partial \Omega$,
    \[\rho_i(0) =  0,\quad i=1,\dots, n-1, \quad \implies \quad \rho_n(0) \neq 0.\]
    Furthermore by \eqref{eq:Subeq-Reg-unif-conv-dom-rho-prop}, $\rho(0) = 0$, $\rho<0$ on $ \Omega$, and $e_n$ is the inward normal. Hence,
    \[\rho_n(0) < 0. \]
    Combined with \eqref{eq:Subeq-Reg-unif-conv-dom-rho-prop}, \eqref{eq:Subeq-Reg-principal-coordinates-sigma-hess-coord-eq0.1} and \eqref{eq:Subeq-Reg-principal-coordinates-sigma-hess-coord-eq1},
    \[\kappa_i = \chi_{ii}(0) = \frac{\rho_{ii}(0)}{ |\rho_n(0)|} {\ge} \frac{\kappa}{ |\rho_n(0)|}  >0, \qquad i=1,\dots, n-1. \]
    This completes the proof of Lemma~\ref{lemma:Subeq-Reg-unif-convex-boundary-coordinates}.
\end{proof}    
\end{allproofs}

\begin{allproofs}
\subsection{Second fundamental form}

The second fundamental form of $\partial \Omega$ is a symmetric linear map on the tangent bundle $T\partial\Omega$ that captures information about the curvature of $\partial \Omega$~\cite[p.~43]{Petersen-Riem-Geom-MR3469435}.

Denote the inward unit normal to $\Omega$ by
\begin{equation}
\label{eq:Subeq-Reg-inward-normal-to-domain-bdry-notation}
    \nu := \nu^\Omega:\partial \Omega \rightarrow S^{n-1}
\end{equation}
and denote the directional derivative along $\nu$ by
\begin{equation}
    \label{eq:Subeq-Reg-inward-normal-to-domain-bdry-notation-2}
    \partial_{\nu}u(x) := \langle Du(x), \nu(x)\rangle, \qquad \forall\ x\in \partial \Omega.
\end{equation}

\begin{definition}
    For a domain $\Omega \subset \RR^n$ with $C^2$ boundary, the \textit{second fundamental form} of $\partial \Omega$, denoted $\II_{\partial \Omega}$, is a linear map on the tangent bundle $\II_{\partial\Omega} :T\partial \Omega \rightarrow T\partial \Omega$,
    \begin{equation}
        \label{eq:Subeq-Reg-second-fundamental-form}
        \langle \II_{\partial \Omega}(x) V, W \rangle := -\langle \partial_V \nu(x), W\rangle, \qquad \forall\ V,W\in T_x\partial \Omega,
    \end{equation}
    where $\nu \in C^2(\overline{\Omega})$ denotes any $C^2$ extension of the inward unit normal \eqref{eq:Subeq-Reg-inward-normal-to-domain-bdry-notation}.
    
    The linear map $ \II_{\partial \Omega}(x): T_x\partial \Omega \rightarrow T_x\partial \Omega$ is well-defined, symmetric (with respect to the induced inner product on $T_x\partial \Omega \subset \RR^n$), and independent of the chosen extension of $\nu$~\cite[Chapter~5, Exercise (8)]{Petersen-Riem-Geom-MR3469435}.
\end{definition}

\begin{remark}[Second fundamental form]
    The eigenvalues $\kappa_1,\dots, \kappa_{n-1}$ of the Hessian $D^2\chi$ from \eqref{eq:Subeq-Reg-principal-coordinate-graph-sigma} are in fact the principal curvatures of $\partial \Omega$ at the boundary point $0\in \partial \Omega$. Recall that the \textit{second fundamental form} of $\partial \Omega$ is a quadratic form on the tangent plane,
    \begin{equation}
        \label{eq:Subeq-Reg-second-fundamental-form-remark}
        \II_{\partial \Omega}(x) \in \Sym^2(T_x\partial \Omega),\qquad  \forall\ x\in \partial \Omega,\\
    \end{equation}
    whose eigenvalues are the principal curvatures $\kappa_1,\dots, \kappa_{n-1}$ \cite[p.~43]{Petersen-Riem-Geom-MR3469435} (when $\partial \Omega$ is oriented by the inward unit normal). Hence, by \eqref{eq:Subeq-Reg-principal-coordinate-graph-sigma},
    \begin{equation}
        \label{eq:Subeq-Reg-second-fundamental-form-as-hessian-of-sigma}
        \II_{\partial \Omega}(0) = D^2\chi(0).
    \end{equation}
    By Lemma~\ref{lemma:Subeq-Reg-unif-convex-boundary-coordinates}, any boundary point $x_0\in \partial \Omega$ admits a principal coordinate system in which \eqref{eq:Subeq-Reg-principal-coordinate-graph-sigma} holds. So, for the purposes of this article, \eqref{eq:Subeq-Reg-second-fundamental-form-as-hessian-of-sigma} can be treated as the definition of $\II_{\partial \Omega}$.
\end{remark}
\end{allproofs}

\section{Uniform and gradient estimates}
\label{sec:Subeq-Reg-unif-grad-estimates}

Sections~\ref{sec:Subeq-Reg-unif-grad-estimates}--\ref{sec:Subeq-Reg-normal-hessian-estimates} together prove Theorem~\ref{thm:Subeq-Reg-a-priori-estimates}. This section establishes a priori $C^0$ and $C^1$ estimates on a solution to \eqref{eq:Subeq-Reg-const-coeff-subequation-PDE}. These estimates do not need $G$-invariance of $F$.

\begin{proposition}
\label{prop:Subeq-Reg-unif-grad-est}
    For a convex, non-degenerate subequation $F\subset \Sym^2(\RR^n)$ \eqref{eq:Subeq-Reg-non-degenerate-subequation} with smooth boundary $\partial F$, a bounded domain $\Omega \subset \RR^n$ with $C^3$ boundary $\partial \Omega$ that is strictly $\vec{F}$-convex, and a boundary condition $\varphi \in C^3(\overline{\Omega})$ there exists a constant $C = C(F, \Omega, \|\varphi\|_{C^3})$ such that any solution $u \in C^3(\overline{\Omega})$ to \eqref{eq:Subeq-Reg-const-coeff-subequation-PDE} satisfies
    \begin{equation}
        \label{eq:Subeq-Reg-main-grad-est}
        \|u\|_{C^1(\overline{\Omega})} \le C.
    \end{equation}
\end{proposition}

By Lemma~\ref{lemma:Subeq-Reg-subequation-is-elliptic}~(iii), a convex and non-degenerate \eqref{eq:Subeq-Reg-non-degenerate-subequation} subequation with smooth boundary is strictly elliptic \eqref{eq:Subeq-Reg-strictly-elliptic-subequation-pos-normals}. So, Proposition~\ref{prop:Subeq-Reg-unif-grad-est} follows immediately from Lemmas~\ref{lemma:Subeq-Reg-unif-est}, \ref{lemma:Subeq-Reg-interior-grad-est}, and \ref{lemma:Subeq-Reg-boundary-grad-est}.

\begin{remark}[A priori $C^1$ estimate holds for strictly elliptic \eqref{eq:Subeq-Reg-strictly-elliptic-subequation-pos-normals} subequations]
    Convexity and non-degeneracy \eqref{eq:Subeq-Reg-non-degenerate-subequation} of $F$ are used to establish strict ellipticity of $F$ \eqref{eq:Subeq-Reg-strictly-elliptic-subequation-pos-normals}, which is needed for the gradient maximum principle (Lemma~\ref{lemma:Subeq-Reg-interior-grad-est}). Convexity of $F$ is also needed to invoke Proposition~\ref{prop:Subeq-Reg-existence-of-delta-subsol-and-supersol}~(ii) to obtain an admissible $\delta$-supersolution \eqref{eq:Subeq-Reg-strongly-admissible-supersolution}.
    
    If $\partial \Omega$ is both strictly $\vec{F}$- and $\vec{\widetilde{F}}$-convex \eqref{eq:Subeq-Reg-dual-to-a-subequation}, then by the construction of Proposition~\ref{prop:Subeq-Reg-existence-of-delta-subsol-and-supersol}~(i) applied to $\widetilde{F}$, an admissible $\delta$-supersolution \eqref{eq:Subeq-Reg-strongly-admissible-supersolution} exists. Hence, the $C^1$ estimate \eqref{eq:Subeq-Reg-main-grad-est} holds for any subequation $F\subset \Sym^2(\RR^n)$ with smooth boundary $\partial F$ that is strictly elliptic \eqref{eq:Subeq-Reg-strictly-elliptic-subequation-pos-normals}.
    
    In particular, by Proposition~\ref{prop:Subeq-Reg-unif-convex-is-F-vec-for-all-F} the $C^1$ estimate \eqref{eq:Subeq-Reg-main-grad-est} holds for strictly elliptic \eqref{eq:Subeq-Reg-strictly-elliptic-subequation-pos-normals} subequations $F$ with smooth boundary $\partial F$, on uniformly convex domain $\Omega \subset \RR^n$ (Definition~\ref{def:Subeq-Reg-uniformly-convex-domain}).
\end{remark}


\subsection{Uniform estimate}

By Theorem~\ref{thm:Subeq-Reg-Har-Law-exist-thm}, the BVP \eqref{eq:Subeq-Reg-const-coeff-subequation-PDE} admits a unique continuous weak solution $u\in C^0(\overline{\Omega})$ for any subequation $F\subset \Sym^2(\RR^n)$, continuous boundary condition $\varphi \in C^0(\partial \Omega)$, and uniformly convex domain. In particular, that solution has a bounded $C^0$ norm. However, the barriers from Proposition~\ref{prop:Subeq-Reg-existence-of-delta-subsol-and-supersol} provide a simple $C^0$ estimate:

\begin{lemma}
\label{lemma:Subeq-Reg-unif-est}
    For a bounded uniformly convex domain $\Omega \subset \RR^n$, a boundary condition $\varphi \in C^2( \partial \Omega)$ and a subequation $F\subset \Sym^2(\RR^n)$, there exists a constant $C>0$ depending only on $F$, $\Omega$ and $\|\varphi\|_{C^2(\partial \Omega)}$ such that if $u\in C^2(\overline{\Omega})$ solves \eqref{eq:Subeq-Reg-const-coeff-subequation-PDE}, then
    \[ \max_{\overline{\Omega}}|u| \le  C.\]
\end{lemma}

\begin{proof}
    By Proposition~\ref{prop:Subeq-Reg-existence-of-delta-subsol-and-supersol}, there exists $\underline{u}, \overline{u} \in C^2(\overline{\Omega})$, whose $C^2$ norms are controlled by $C_1 = C_1(F,\Omega, \|\varphi\|_{C^3})$, satisfying \eqref{eq:Subeq-Reg-strongly-admissible-subsolution} and \eqref{eq:Subeq-Reg-strongly-admissible-supersolution} respectively. In particular,
    \[\begin{aligned}
        D^2\underline{u} \in \intr F, \quad D^2 u \in \partial F, \quad D^2\overline{u} \in F^c = -\intr \widetilde{F} &\qquad \text{on }  \Omega,\\
        \underline{u} = u = \overline{u} = \varphi  &\qquad \text{on } \partial \Omega.
    \end{aligned}\]
    By the comparison principle~\cite[Theorem~6.4]{HarLaw-dirichletduality},
    \[\underline{u} \le u \le \overline{u} \qquad \text{on } \overline{\Omega}. \tag*{\qedhere}\]
\end{proof}

\subsection{Gradient maximum principle}
\label{subsec:Subeq-Reg-gradient-maximum-principle}

By Lemma~\ref{lemma:Subeq-Reg-signed-distance-lemma}, a solution $u$ to \eqref{eq:Subeq-Reg-const-coeff-subequation-PDE} also satisfies \eqref{eq:Subeq-Reg-signed-dist-fcal-PDE}. Hence, its gradient $Du$ satisfies the linearization of \eqref{eq:Subeq-Reg-signed-dist-fcal-PDE}, which is a linear elliptic PDE. Being a solution to a linear elliptic PDE, the components of $Du$ attain their maximum and minimum on the boundary of the domain $\Omega$ \cite[p. 279]{caffarelli-nirenberg-spruck-dirichlet-problem-III-MR0806416}. Note that $\signdist^{F}$ only needs to be smooth on a neighborhood of $D^2u(\overline{\Omega})$ in this argument.

\begin{lemma}
\label{lemma:Subeq-Reg-interior-grad-est}
    For a bounded domain $\Omega\subset \RR^n$, a strictly elliptic \eqref{eq:Subeq-Reg-strictly-elliptic-subequation-pos-normals} subequation $F\subset \Sym^2(\RR^n)$ with smooth boundary $\partial F$, and any boundary condition $\varphi$, any solution $u\in C^3(\overline{\Omega})$ to \eqref{eq:Subeq-Reg-const-coeff-subequation-PDE} satisfies
    \begin{equation}
        \label{eq:Subeq-Reg-interior-grad-est}
        \max_{\overline{\Omega}}|Du| \le C(n) \max_{\partial \Omega} |Du|.
    \end{equation}
\end{lemma}
\begin{proof}
    Since $u$ satisfies \eqref{eq:Subeq-Reg-signed-dist-fcal-PDE}. Differentiating \eqref{eq:Subeq-Reg-signed-dist-fcal-PDE} along $x\in {\Omega}$,
    \[\frac{\partial}{\partial x_k} \signdist^F(D^2u) = \langle D\signdist^F(D^2u), D^2\, \partial_k u\rangle = \Delta_{D\signdist^F(D^2u)} \partial_k u = 0 \qquad \text{on } \Omega, \quad k=1,\dots, n,
    \]
    where $\Delta_{D\signdist(D^2u)}$ \eqref{eq:Subeq-Reg-linearization-of-dij-elliptic-pde} is a linear elliptic PDE by Claim~\ref{claim:Subeq-Reg-dij-is-elliptic-PDE}. Since the domain $\Omega$ is bounded, by the maximum principle for linear elliptic PDE \cite[Theorem~3.1]{gilbarg-trudinger-MR1814364},
    \[\max_{\overline{\Omega}} |\partial_k u| = \max_{\partial \Omega} |\partial_k u|, \qquad k=1,\dots, n. \tag*{\qedhere}\]
\end{proof}

\subsection{Boundary gradient estimate}
\label{subsection:Subeq-Reg-boundary-gradient-estimate}

The equation \eqref{eq:Subeq-Reg-const-coeff-subequation-PDE} admits a comparison principle~\cite[Theorem~6.4]{HarLaw-dirichletduality}, and also admits both admissible subsolutions and supersolutions (Proposition~\ref{prop:Subeq-Reg-existence-of-delta-subsol-and-supersol}), a boundary gradient estimate follows by adapting an argument of Caffarelli--Nirenberg--Spruck~\cite[p.~278]{caffarelli-nirenberg-spruck-dirichlet-problem-III-MR0806416}.

\begin{lemma}
\label{lemma:Subeq-Reg-boundary-grad-est}
    For convex subequation $F\subset \RR^n$, a bounded domain $\Omega \subset \RR^n$ with $C^3$, strictly $\vec{F}$-convex boundary $\partial \Omega$ (Definition~\ref{def:Subeq-Reg-HL-boundary-convexity}), $\Omega\subset \RR^n$ a subequation $F\subset \Sym^2(\RR^n)$, and a boundary condition $\varphi \in C^3(\partial{\Omega})$, there exists $C = C(F,\Omega, \|\varphi\|_{C^3})$ such that any solution $u\in  C^3(\overline{\Omega})$ to \eqref{eq:Subeq-Reg-const-coeff-subequation-PDE} satisfies 
    \begin{equation}
        \label{eq:Subeq-Reg-boundary-grad-est}
        \max_{\partial \Omega} |Du| \le C.
    \end{equation}
\end{lemma}

\begin{proof} 
    \noindent\textit{Step 1.} \textit{Tangential gradient estimate:}

    \vspace{0.5em}\noindent Since $u=\varphi$ on $\partial \Omega$, differentiating along directions tangent to the smooth boundary $\partial \Omega$,
    \begin{equation}
        \label{eq:Subeq-Reg-boundary-grad-est-tangent}
        \left| \langle Du(x), V\rangle \right| = \left| \langle D\varphi(x), V\rangle \right|  \le \|\varphi\|_{C^1}|V|, \qquad \forall\ x\in \partial \Omega,\ V\in T_x\partial \Omega.
    \end{equation}
    
    \vspace{0.5em}\noindent\textit{Step 2.} \textit{Normal gradient estimate:}

    \vspace{0.5em}\noindent  Denote the inward unit normal to $\Omega$ by
    \begin{equation}
    \label{eq:Subeq-Reg-inward-normal-to-domain-bdry-notation}
        \nu := \nu^\Omega:\partial \Omega \rightarrow S^{n-1}
    \end{equation}
    and denote the directional derivative along $\nu$ by
    \begin{equation}
        \label{eq:Subeq-Reg-inward-normal-to-domain-bdry-notation-2}
        \partial_{\nu}u(x) := \langle Du(x), \nu(x)\rangle, \qquad \forall\ x\in \partial \Omega.
    \end{equation}
    By Proposition~\ref{prop:Subeq-Reg-existence-of-delta-subsol-and-supersol}, there exists a subsolution $\underline{u}, \overline{u} \in {C^2}(\overline{\Omega})$, satisfying \eqref{eq:Subeq-Reg-strongly-admissible-subsolution} and \eqref{eq:Subeq-Reg-strongly-admissible-supersolution} respectively, whose $C^2$ norms are bounded in terms of $F$, $\Omega$, and $\varphi$. By \eqref{eq:Subeq-Reg-strongly-admissible-subsolution}--\eqref{eq:Subeq-Reg-strongly-admissible-supersolution} and the comparison principle~\cite[Theorem~6.4]{HarLaw-dirichletduality},
    \begin{equation}
        \label{eq:Subeq-Reg-boundary-grad-est-normal-eq1}
        \begin{aligned}
        \underline{u} \le u \le \overline{u} \qquad \text{on } \overline{\Omega},\\
        \underline{u} = u = \overline{u} \qquad \text{on } \partial {\Omega}.
    \end{aligned}
    \end{equation}
    The normal gradient estimate follows from \eqref{eq:Subeq-Reg-boundary-grad-est-normal-eq1} by the following elementary observation, stated without proof:

    \begin{claim}
    \label{claim:Subeq-Reg-normal-gradient-claim}
        For a domain $\Omega \subset \RR^n$ with $C^2$ boundary $\partial \Omega$, and $x_0\in \partial \Omega$, if $u,v\in C^1(\overline{\Omega})$ satisfy
        \begin{equation}
        \label{eq:Subeq-Reg-normal-gradient-lemma-assump}
            \begin{aligned}
                 u\le&\ v \qquad \text{on } \Omega,\\
            u(x_0)  =&\ v(x_0),
            \end{aligned}
        \end{equation}
        then using the notation \eqref{eq:Subeq-Reg-inward-normal-to-domain-bdry-notation}--\eqref{eq:Subeq-Reg-inward-normal-to-domain-bdry-notation-2},
        \begin{equation}
            \label{eq:Subeq-Reg-normal-gradient-lemma-conclusion}
            \partial_\nu u(x_0)  \le \partial_\nu v (x_0).
        \end{equation}
    \end{claim}
\begin{allproofs}
    \begin{proof}
        Since $\Omega$ is an open set and $\nu(x_0)$ is the inward normal at $x_0 \in \Omega$, there exists $t_0>0$ such that
        \[x_0 +  t\nu(x_0) \in \Omega, \qquad \forall\ t\in (0,t_0).\]
        So, computing the directional derivative along $\nu(x_0)$ and applying \eqref{eq:Subeq-Reg-normal-gradient-lemma-assump},
        \[\begin{aligned}
            \frac{u(x_0+t\nu(x_0)) - u(x_0)}{t} \le&\  \frac{v(x_0+t\nu(x_0)) - v(x_0)}{t} , \quad \forall\ t\in (0,t_0). 
        \end{aligned}\]
        Taking a limit at $t\rightarrow 0^+$, \eqref{eq:Subeq-Reg-normal-gradient-lemma-conclusion} follows.
    \end{proof}
\end{allproofs}
    Combining \eqref{eq:Subeq-Reg-boundary-grad-est-normal-eq1} and Claim~\ref{claim:Subeq-Reg-normal-gradient-claim},
    \[\partial_\nu \underline{u}(x) \le \partial_\nu {u}(x)  \le \partial_\nu \overline{u}(x) , \qquad \forall\ x\in \partial \Omega.\]
    \[ \left\langle D\underline{u}(x), \nu(x) \right\rangle \le  \left\langle Du(x), \nu(x) \right\rangle \le  \left\langle D\overline{u}(x), \nu(x) \right\rangle, \qquad \forall\ x\in \partial \Omega.\]
    In other words, by \eqref{eq:Subeq-Reg-inward-normal-to-domain-bdry-notation-2},
    \begin{equation}
        \label{eq:Subeq-Reg-boundary-grad-est-normal}
        \implies \qquad \left| \langle Du(x), V\rangle \right| \le |V|\max\left\{ \|\underline{u}\|_{C^1} ,  \|\overline{u}\|_{C^1} \right\} , \qquad \forall\ x\in \partial \Omega,\ V\in (T_x\partial \Omega)^\perp.
    \end{equation}
    Then, \eqref{eq:Subeq-Reg-boundary-grad-est} follows from \eqref{eq:Subeq-Reg-boundary-grad-est-tangent} and \eqref{eq:Subeq-Reg-boundary-grad-est-normal}.
\end{proof}

\section{Hessian maximum principle and tangential Hessian estimate}
\label{sec:Subeq-Reg-interior-and-tangential-hessian-estimates}

Proposition~\ref{prop:Subeq-Reg-unif-grad-est} established a priori $C^1$ estimates on a solution $u$ to \eqref{eq:Subeq-Reg-const-coeff-subequation-PDE}. Sections~\ref{sec:Subeq-Reg-interior-and-tangential-hessian-estimates}--\ref{sec:Subeq-Reg-normal-hessian-estimates} use this estimate and build up to a proof of Theorem~\ref{thm:Subeq-Reg-a-priori-estimates}, eventually presented in Section~\ref{sec:Subeq-Reg-normal-hessian-estimates}. The first two steps are a maximum principle for the Hessian (Lemma~\ref{lemma:Subeq-Reg-interior-hess-est}) and an estimate on the tangential components of the Hessian at the boundary (Lemma~\ref{lemma:Subeq-Reg-C2-TT-estimate}).

\subsection{Hessian maximum principle}

The Hessian maximum principle is adapted from an argument of Caffarelli--Nirenberg--Spruck\cite[pp. 279--280]{caffarelli-nirenberg-spruck-dirichlet-problem-III-MR0806416}. It does not require $G$-invariance of $F$. Note that $\signdist^{F}$ only needs to be smooth in a neighborhood of the solution $D^2u(\overline{\Omega})$:

\begin{lemma}
\label{lemma:Subeq-Reg-interior-hess-est}
    For a bounded domain $\Omega \subset \RR^n$, a convex, non-degenerate subequation $F\subset \Sym^2(\RR^n)$ \eqref{eq:Subeq-Reg-non-degenerate-subequation} with smooth boundary $\partial F$, and any boundary condition $\varphi$, there exists $C = C(n,F)$ such that any solution $u\in C^2(\overline{\Omega}) \cap C^4({\Omega})$ to \eqref{eq:Subeq-Reg-const-coeff-subequation-PDE} satisfies
    \begin{equation}
        \label{eq:Subeq-Reg-interior-hess-est}
        \max_{\overline{\Omega}}|D^2u| \le C(1+ \max_{\partial \Omega}|D^2u|).
    \end{equation}
\end{lemma}

\begin{proof}
     Since all norms on $\Sym^2(\RR^n)$ are equivalent, the norm $|D^2u|$ in \eqref{eq:Subeq-Reg-interior-hess-est} can be replaced by the \textit{spectral norm}
    \begin{equation}
        \label{eq:Subeq-Reg-interior-hess-est-spectral-norm}
        |D^2u| = \max_{V\in S^{n-1}} |\langle D^2u\, V,V\rangle |,
    \end{equation}
    up to changing the constant $C$ in \eqref{eq:Subeq-Reg-interior-hess-est} by a factor depending on $n$.

    \vspace{0.5em}\noindent\textit{Step 1.} \textit{Upper bound on the quadratic form of $D^2u$:}

    \vspace{0.5em}\noindent Use subscripts to denote derivatives of $u$ along $e \in S^{n-1}$
    \[u_{e} := \langle Du, e\rangle\quad \text{and} \quad u_{ee}:= \langle D^2u\, e, e\rangle, \qquad \forall\ e\in S^{n-1} .\]
    By Lemma~\ref{lemma:Subeq-Reg-signed-distance-lemma}~(ii), $\signdist^F$ \eqref{eq:Subeq-Reg-signed-dist-def} is smooth on a neighborhood $U$ of the boundary $\partial F$, which in particular contains $D^2u(\overline{\Omega})$ by \eqref{eq:Subeq-Reg-const-coeff-subequation-PDE}. Since $u$ satisfies \eqref{eq:Subeq-Reg-const-coeff-subequation-PDE} it also satisfies \eqref{eq:Subeq-Reg-signed-dist-fcal-PDE}. Differentiating \eqref{eq:Subeq-Reg-signed-dist-fcal-PDE} along $e$, by \eqref{eq:Subeq-Reg-linearization-of-dij-elliptic-pde},
    \[\Delta_{D\signdist^F(D^2u)}u_e = \sum_{i,j=1}^n \signdist^F_{ij}(D^2u) \partial_i \partial_j u_e = 0 \qquad \text{on } {\Omega}, \quad \forall\ e\in S^{n-1}, \]
    where $\signdist^F_{ij}$ denote the components of the derivative $D\signdist^F$. Differentiating again along $e$, 
    \begin{equation}
        \label{eq:Subeq-Reg-hess-max-prin-new-eq1}
        \sum_{i,j,k,l=1}^n \signdist^F_{ij;kl}(D^2u)\, \partial_i \partial_j u_e\, \partial_k \partial_l u_e + \sum_{i,j=1}^n \signdist^F_{ij}(D^2u) \partial_i \partial_j u_{ee} =0 \qquad \text{on } {\Omega}, \quad \forall\ e\in S^{n-1},
    \end{equation}
    where $D^2\signdist^F = (\signdist^F_{ij;kl})_{i,j,k,l=1}^n$ is a quadratic form on $\Sym^2(\RR^n)$.
    The first term in \eqref{eq:Subeq-Reg-hess-max-prin-new-eq1} can be expressed in terms of this quadratic form as $D^2 \signdist^F(D^2u)(D^2u_e, D^2u_e)$, which is non-positive by concavity of $\signdist^F$ (Lemma~\ref{lemma:Subeq-Reg-signed-distance-lemma}~(v)). Hence, $u_{ee}$ is a subsolution to the linear PDE $\Delta_{D\signdist^F (D^2u)}$ \eqref{eq:Subeq-Reg-linearization-of-dij-elliptic-pde} on $\Omega$,
    \[\Delta_{D\signdist^F(D^2u)} u_{ee} = \sum_{i,j=1}^n \signdist^F_{ij}(D^2u) \partial_i \partial_j u_{ee} \ge 0 \qquad \text{on } {\Omega}.\]
    By Claim~\ref{claim:Subeq-Reg-dij-is-elliptic-PDE}, $\Delta_{D\signdist^F (D^2u)}$ is strictly elliptic. Since the domain $\Omega$ is bounded, by the maximum principle for linear elliptic PDE \cite[Theorem~3.1]{gilbarg-trudinger-MR1814364} applied to $u_{ee} \in C^0(\overline{\Omega}) \cap C^2(\Omega)$,
    \[\max_{\overline{\Omega}} u_{ee} \le \max_{\partial {\Omega}} u_{ee}.\]
    In other words,
    \begin{equation}
        \label{eq:Subeq-Reg-C2-max-prin-Step1-goal}
        \max_{\overline{\Omega}} \langle D^2u\, e, e\rangle \le \max_{\partial {\Omega}} \langle D^2u\, e, e\rangle \le \max_{\partial \Omega} |D^2u|, \qquad \forall\ e\in S^{n-1}.
    \end{equation}

    \vspace{0.5em}\noindent\textit{Step 2.} {Lower bound on the quadratic form of $D^2u$. Follows from Step 1 and Claim~\ref{claim:Subeq-Reg-upper-bound-on-Hess-implies-lower-bound}:}

    \begin{claim}
        \label{claim:Subeq-Reg-upper-bound-on-Hess-implies-lower-bound}
        For a convex, non-degenerate subequation $F\subset \Sym^2(\RR^n)$  \eqref{eq:Subeq-Reg-non-degenerate-subequation}, given $C_1\ge 0$ there exists $C_2\ge 1$ depending only on $n$ and $F$ such that
        \[A\in F \quad \text{and} \quad A\le C_1I \qquad \implies \qquad A \ge -C_2(1+C_1)I.\]
    \end{claim}
    \begin{proof}
        Fix a matrix $A_0 \in \partial F$. By convexity, $F$ admits a supporting half-space at $A_0$, i.e.,
    \begin{equation}
        \label{eq:Subeq-Reg-interior-hess-est-lower-eq1-alt}
        \langle B_{A_0}, A \rangle \ge \langle B_{A_0}, A_0 \rangle, \qquad \forall\ A\in F, 
    \end{equation}
    where $B_{A_0}$ \eqref{eq:Subeq-Reg-inward-supporting-hyperplane} is the inward unit normal to some supporting hyperplane to $\partial F$ at $A_0$. Fix a symmetric matrix $A\in F$ satisfying
    \begin{equation}
        \label{eq:Subeq-Reg-interior-hess-est-lower-eq5-alt}
         A\le C_1I,
    \end{equation}
    and consider its decomposition into its positive and negative semi-definite parts
    \[A = A_+ + A_-, \qquad A_+\ge 0,\qquad A_-\le 0,\qquad A_-A_+ = 0.\]
    Since $A_+\ge0$, and its eigenvalues are either $0$ or equal to those of $A$, all of its eigenvalues are bounded $0\le A_+\le C_1 I$ \eqref{eq:Subeq-Reg-interior-hess-est-lower-eq5-alt}. So, by \eqref{eq:Subeq-Reg-trace-inner-prod}, $|A_+| \le \sqrt{n}C_1$, and by \eqref{eq:Subeq-Reg-interior-hess-est-lower-eq1-alt},
    \begin{equation}
        \label{eq:Subeq-Reg-interior-hess-est-lower-eq5-alt-eq-eq1}
        \langle B_{A_0}, A_- \rangle \ge \langle  B_{A_0} , A \rangle - \langle B_{A_0} , A_+\rangle \ge \langle B_{A_0}, A_0 \rangle - \sqrt{n}C_1|B_{A_0}|.
    \end{equation}
    But $A_-\le 0$ by construction, and since $B_{A_0} \ge \lambda_{\min}\left(B_{A_0}\right) I $,
    \[\langle B_{A_0}, A_- \rangle \le \lambda_{\min}\left(B_{A_0}\right) \tr (A_-).\]
    So by \eqref{eq:Subeq-Reg-interior-hess-est-lower-eq5-alt-eq-eq1},
    \[\begin{aligned}
        \langle B_{A_0}, A_0 \rangle - \sqrt{n}C_1|B_{A_0}|\le&\ \langle B_{A_0}, A_- \rangle \le \lambda_{\min}\left(B_{A_0}\right) \tr (A_-) .
    \end{aligned}\]
    But by non-degeneracy of $F$ \eqref{eq:Subeq-Reg-non-degenerate-subequation} and Lemma~\ref{lemma:Subeq-Reg-convex-subequation-normals}, $B_{A_0}>0$, and in particular $\lambda_{\min}\left(B_{A_0} \right) >0$. And $\tr(A_-)$ is the sum of the non-positive eigenvalues of $A_-$, and is more negative than each eigenvalue of $A_-$. So,
    \[A = A_+ + A_- \ge A_- \ge \tr(A_-) I \ge -\frac{\max\left\{ |\langle B_{A_0}, A_0 \rangle|,\  \sqrt{n}|B_{A_0}| \right\}}{\lambda_{\min}\left(B_{A_0}\right)}(1+C_1)I. \tag*{\qedhere}\]
    \end{proof}

    \vspace{0.5em}\noindent\textit{Step 3.} \textit{Combining upper and lower bounds:}

    \vspace{0.5em}\noindent By \eqref{eq:Subeq-Reg-const-coeff-subequation-PDE}, $D^2u(x) \in \partial F \subset F$ for all $x\in {\Omega}$. So by \eqref{eq:Subeq-Reg-C2-max-prin-Step1-goal} and Claim~\ref{claim:Subeq-Reg-upper-bound-on-Hess-implies-lower-bound} applied to $D^2u$,
    \[|\langle D^2u\, e, e\rangle| \le C_2\left(1+\max_{\partial \Omega} |D^2u|\right), \qquad \forall\ e\in S^{n-1}.\]
    Hence, \eqref{eq:Subeq-Reg-interior-hess-est} follows by \eqref{eq:Subeq-Reg-interior-hess-est-spectral-norm}.
\end{proof}

\subsection{Tangential Hessian estimate}
\label{subsection:Subeq-Reg-tangential-hessian-estimate}

This section is concerned with an a priori estimate on the tangential Hessian components of a solution to \eqref{eq:Subeq-Reg-const-coeff-subequation-PDE}. These components are controlled in terms of the gradient estimate (Proposition~\ref{prop:Subeq-Reg-unif-grad-est}) and a formula \eqref{eq:Subeq-Reg-C2-TT-key-lemma-conclusion} for the tangential Hessian.

\begin{lemma}
\label{lemma:Subeq-Reg-C2-TT-estimate}
    For a bounded domain $\Omega \subset \RR^n$, a boundary condition $\varphi \in C^2(\partial \Omega)$, a subequation $F \subset \Sym^2(\RR^n)$, and a solution $u \in C^3(\overline{\Omega})$ to \eqref{eq:Subeq-Reg-const-coeff-subequation-PDE}, there exists $C = C(F, \Omega, \|\varphi\|_{C^3})>0$ such that
    \begin{equation}
        \label{eq:Subeq-Reg-tangential-hessian-estimate}
        \left| \langle D^2u(x) V,W\rangle \right| \le C|V||W|, \qquad \forall\ x\in \partial \Omega,\ V,W\in T_{x}\partial \Omega.
    \end{equation}
\end{lemma}

Since $u=\varphi$ on the boundary $\partial \Omega$, their intrinsic Hessians with respect to the Riemannian metric on the hypersurface $\partial \Omega$ are equal. The difference between the intrinsic Hessian  and the extrinsic Hessian is governed by the second fundamental form of $\partial \Omega$:
    \begin{equation}
        \label{eq:Subeq-Reg-second-fundamental-form}
        \langle \II_{\partial \Omega}(x) V, W \rangle := -\langle \partial_V \nu(x), W\rangle, \qquad \forall\ V,W\in T_x\partial \Omega,
    \end{equation}
where $\nu \in C^2(\overline{\Omega})$ denotes any $C^2$ extension of the inward unit normal \eqref{eq:Subeq-Reg-inward-normal-to-domain-bdry-notation}. Recall that $ \II_{\partial \Omega}(x): T_x\partial \Omega \rightarrow T_x\partial \Omega$ is well-defined, linear, symmetric (with respect to the induced inner product on $T_x\partial \Omega \subset \RR^n$), and independent of the chosen extension of $\nu$~\cite[Chapter~5, Exercise (8)]{Petersen-Riem-Geom-MR3469435}

Lemma~\ref{lemma:Subeq-Reg-C2-TT-estimate} follows  from the gradient estimate (Proposition~\ref{prop:Subeq-Reg-unif-grad-est}), and a general formula \eqref{eq:Subeq-Reg-C2-TT-key-lemma-conclusion} relating the intrinsic and extrinsic Hessian of $u|_{\partial \Omega}$ via the second fundamental form, recalled below without proof \cite[(5.1)]{caffarelli-nirenberg-spruck-dirichlet-problem-III-MR0806416}:

\begin{claim}
\label{claim:Subeq-Reg-C2-TT-key-claim}
    For $\Omega \subset \RR^n$ a bounded domain with $C^2$ boundary $\partial \Omega$, $u,v \in C^2( \overline{ \Omega})$ equal on the boundary, $u=v$ on $\partial \Omega$, and a boundary point $x\in \partial \Omega$ with tangent vectors $V,W\in T_x\partial \Omega$,
    \begin{equation}
        \label{eq:Subeq-Reg-C2-TT-key-lemma-conclusion}
        \langle D^2u(x) V,W\rangle + \partial_\nu u(x) \left\langle \II_{\partial \Omega}(x)V,W \right\rangle = \langle D^2v(x) V,W\rangle + \partial_\nu v(x) \left\langle \II_{\partial \Omega}(x) V,W \right\rangle,
    \end{equation}
    where $\II_{\partial \Omega}(x)$ \eqref{eq:Subeq-Reg-second-fundamental-form} denotes the second fundamental form of $\partial \Omega$ at $x$ defined in terms of the inward unit normal $\nu$ \eqref{eq:Subeq-Reg-inward-normal-to-domain-bdry-notation}.
\end{claim}
\begin{allproofs}
\begin{proof}
    The idea is to parametrize $\partial \Omega$ as a graph using Lemma~\ref{lemma:Subeq-Reg-unif-convex-boundary-coordinates}, and to differentiate $u-v$ twice along such a parametrization to obtain \eqref{eq:Subeq-Reg-C2-TT-key-lemma-conclusion}.
    
    Fix a point
    \[x_0\in \partial \Omega,\]
    and a principal coordinate system (Definition~\ref{def:Subeq-Reg-principal-coordinate-system}) on $\RR^n$ based at $x_0$. In particular,
    \begin{equation}
        \label{eq:Subeq-Reg-tangential-Hess-formula-eq2}
        x_0 = 0 \in \partial \Omega \quad \text{and}\quad \nu(0) = e_n=(0,\dots,0,1),
    \end{equation}
    and by Lemma~\ref{lemma:Subeq-Reg-unif-convex-boundary-coordinates} there are open sets $U'\subset \RR^{n-1}$ and $U\subset \RR^n$, and $\chi \in C^3(U')$ such that $\partial \Omega$ is locally the graph of $\chi$,
    \[(\id, \chi)(U') = U\cap \partial \Omega,\]
    and \eqref{eq:Subeq-Reg-principal-coordinate-graph-sigma} holds,
    \begin{equation}
        \label{eq:Subeq-Reg-tangential-Hess-formula-eq1}
        \chi(x') = \frac{1}{2}\sum_{i=1}^{n-1} \kappa_i |x_i|^2 + O(|x'|^3),\qquad \forall\ x'\in U'.
    \end{equation}
    where $x' = (x_1,\dots, x_{n-1})$. Also use subscripts to denote partial derivatives, for example
    \[u_i:=\partial_iu \quad \text{and} \quad u_{ij} := \partial_i\partial_j u, \qquad i,j=1,\dots, n.\]
    
    By assumption,
    \[u(x',\chi(x')) = v(x',\chi(x')),\qquad \forall\ x'\in U'.\]
    Differentiating this equation,
    \[u_i(x',\chi(x')) + u_n(x',\chi(x'))\chi_i(x') = v_i(x',\chi(x')) + v_n(x',\chi(x'))\chi_i(x'),\quad x'\in U',\ i=1,\dots, n-1. \]
    Differentiating again and evaluating at $x'=0$, since $D\chi(0) = 0$ by \eqref{eq:Subeq-Reg-tangential-Hess-formula-eq1},
    \begin{equation}
        \label{eq:Subeq-Reg-tangential-Hess-formula-eq3}
        u_{ij}(0) + u_{n}(0) \chi_{ij}(0) = v_{ij}(0) + v_{n}(0) \chi_{ij}(0), \qquad i,j=1,\dots, n-1.
    \end{equation}
    By \eqref{eq:Subeq-Reg-tangential-Hess-formula-eq2},
    \[T_{0}\partial \Omega = \operatorname{span}\{e_1,\dots, e_n\},\]
    and by \eqref{eq:Subeq-Reg-inward-normal-to-domain-bdry-notation-2}
    \[\partial_\nu u(0) = \partial_n u(0),\quad  \partial_\nu v(0) = \partial_n v(0).\]
    Hence, by \eqref{eq:Subeq-Reg-tangential-Hess-formula-eq3},
    \[\begin{aligned}
        &\langle D^2u(0)V,W\rangle + \partial_{\nu} u(0) \langle D^2 \chi(0) V,W\rangle \\
        =\ &\ \langle D^2v(0)V,W\rangle + \partial_{\nu} v(0) \langle D^2 \chi(0) V,W\rangle, \qquad \forall\ V,W\in T_{x_0} \partial \Omega.
    \end{aligned}\]
    But $D^2\chi(0) = \II_{\partial \Omega}(0)$ by \eqref{eq:Subeq-Reg-second-fundamental-form-as-hessian-of-sigma}. So substituting $x_0=0$,
    \[\begin{aligned}
        &\langle D^2 u(0) V,W\rangle + \partial_{\nu} u(x_0) \langle \II_{\partial \Omega}(x_0)V, W\rangle\\
        =\ &\langle D^2 v(0) V,W\rangle + \partial_{\nu} v(x_0) \langle \II_{\partial \Omega}(x_0)V, W\rangle, \qquad \forall\ V,W\in T_{x_0} \partial \Omega.
    \end{aligned}\]
    Since $x_0$ was arbitrary, \eqref{eq:Subeq-Reg-C2-TT-key-lemma-conclusion} follows.
\end{proof}
\end{allproofs}

\section{Mixed Hessian estimates}
\label{sec:Subeq-Reg-mixed-hessian-estimates}

This section controls the tangent-normal components of the Hessian $D^2u$ of a solution $u$ to \eqref{eq:Subeq-Reg-const-coeff-subequation-PDE}.

\begin{proposition}
    \label{prop:Subeq-Reg-C2-TN-prop}
    For a convex,  non-degenerate subequation $F\subset \Sym^2(\RR^n)$ \eqref{eq:Subeq-Reg-non-degenerate-subequation} with smooth boundary $\partial F$, that is invariant under a Lie group $G\subset O(n)$ acting transitively on $S^{n-1}$, a bounded domain $\Omega \subset \RR^n$ whose boundary $\partial \Omega$ is $C^3$ and strictly $\vec{F}$-convex (Definition~\ref{def:Subeq-Reg-HL-boundary-convexity}), and a boundary condition $\varphi \in C^3(\partial \Omega)$, there exists a constant $C = C(n, \Omega, F, G, \|\varphi\|_{C^3})$ such that any solution $u\in C^3(\overline{\Omega}) $ to \eqref{eq:Subeq-Reg-const-coeff-subequation-PDE} satisfies 
    \begin{equation}
        \label{eq:Subeq-Reg-C2-TN-conclusion}
        |\langle D^2u(x) V, \nu(x) \rangle | \le C, \qquad \forall\ x\in \partial \Omega,\ V\in T_x\partial \Omega,\ |V|=1,
    \end{equation}
    where $\nu:\partial \Omega\rightarrow \RR^n$ is the inward unit normal \eqref{eq:Subeq-Reg-inward-normal-to-domain-bdry-notation}.
\end{proposition}

This argument is a subequation version of an argument of Caffarelli--Nirenberg--Spruck~\cite[pp. 281--284]{caffarelli-nirenberg-spruck-dirichlet-problem-III-MR0806416}. They work in the $O(n)$-invariant case, and we generalize the argument to the $G$-invariant case.

The strategy to prove Proposition~\ref{prop:Subeq-Reg-C2-TN-prop} is similar to Step 2 of Lemma~\ref{lemma:Subeq-Reg-boundary-grad-est}. Fix a boundary point and a unit tangent vector at this boundary point,
\[x_0\in \partial\Omega , \qquad V\in T_{x_0} \partial \Omega, \quad |V|=1.\]
Differentiating the equation \eqref{eq:Subeq-Reg-signed-dist-fcal-PDE} along $V$, it follows that the directional derivative
\[\partial_V u(x) := \langle Du(x),V\rangle,\]
is a solution to the linearization $\Delta_{D\signdist^F(D^2u)}$ \eqref{eq:Subeq-Reg-linearization-of-dij-elliptic-pde} of \eqref{eq:Subeq-Reg-signed-dist-fcal-PDE} (see Step 4 below),
\begin{equation}
    \label{eq:Subeq-Reg-C2-TN-sketch-partialVu-solution}
    \Delta_{D\signdist^F(D^2u)} \partial_V u = 0.
\end{equation}
Applying Corollary~\ref{corollary:Subeq-Reg-guan-corollary} to $\pm \partial_V u $, we construct a subsolutions $w_{\pm}$ to \eqref{eq:Subeq-Reg-linearization-of-dij-elliptic-pde},
\[\begin{aligned}
    \Delta_{D\signdist^F(D^2u)} w_{+} \le 0 = \Delta_{D\signdist^F(D^2u)} \partial_V u  \le \Delta_{D\signdist^F(D^2u)} w_{-} &\qquad\text{on } \Omega,\\
    w_{+} \ge \partial_V u \ge w_{-} &\qquad \text{on } \partial \Omega,\\
     w_{+}(x_0) = \partial_V u(x_0) = w_{-}(x_0) = \partial_V\varphi(x_0).&
\end{aligned}\]
By the maximum principle for linear elliptic PDE, $w_{+} \ge \partial_V u\ge w_{-}$ on $\Omega$, and by Claim~\ref{claim:Subeq-Reg-normal-gradient-claim},
\[ \left| \partial_{\nu(x_0)} \partial_V u (x_0) \right| \le  \max\left\{ \|w_{+}\|_{C^1} , \|w_{-}\|_{C^1}  \right\} .\]
Unfortunately, the argument outlined above has a gap. Since  Corollary~\ref{corollary:Subeq-Reg-guan-corollary} is applied to $\pm \partial_Vu$, the $C^1$ norms $\|w_{\pm}\|_{C^1}$ depends on the $\|\partial_V u\|_{C^2}$, which has no a priori bound.

Note that, nowhere in this argument have we used that $V$ is in fact a tangent vector to $\partial \Omega$ at $x_0$. By the boundary condition in \eqref{eq:Subeq-Reg-signed-dist-fcal-PDE},
\[\partial_Vu(x_0) = \partial_V\varphi(x_0), \qquad \forall\ V\in T_{x_0} \partial \Omega.\]
This equality need not hold at any other point of $\partial \Omega$ since the constant vector field $V$ is only tangent at $x_0$. We could instead replace $V$ by a tangential vector field $\tau_0: \RR^n \rightarrow \RR^n$,
\begin{equation}
    \label{eq:Subeq-Reg-TN-tangential-vector-field-to-bdry}
    \tau_0(x)\in T_x\partial \Omega, \qquad \forall\ x\in \partial \Omega.
\end{equation}
Then $\partial_{\tau_0} u = \partial_{\tau_0} \varphi$ on $\partial \Omega$. But for such an arbitrary choice of $\tau_0$, $\partial_{\tau_0} u$ need not be a solution to the linearization. \eqref{eq:Subeq-Reg-C2-TN-sketch-partialVu-solution}

The idea of the proof is to balance these two observations. By Lemma~\ref{lemma:Subeq-Reg-symmetries-of-linearization}, if $\tau = \tau_{e,Q}$ \eqref{eq:Subeq-Reg-tau-VQ-vector-field-def} for some $e\in \RR^n$ and $Q\in \mathfrak{g}$, then $\partial_\tau u $ is a solution to the linearized equation \eqref{eq:Subeq-Reg-symmetries-of-linearization}. In general, such a vector field $\tau$ need not be tangential to $\partial \Omega$ \eqref{eq:Subeq-Reg-TN-tangential-vector-field-to-bdry}. However, if we choose $e = V-Qx_0$ for some $V\in T_{x_0} \partial \Omega$, i.e.,
\[\tau = \tau_{e,Q} = e + Qx = V + Q(x-x_0),\]
then $\tau(x_0) = V$ and $\partial_\tau u(x_0) = \partial_\tau \varphi(x_0)$. Similarly, choosing $Q\in \mathfrak{g}$ carefully, we can ensure that $\tau = \tau_{e,Q}$ is tangent to $\partial \Omega$ at $x_0$ \textit{up to first order} (see \eqref{eq:Subeq-Reg-C2-TN-step2-goal}). Then, $ \partial_{\tau} u(x)$ will agree with $\partial_\tau \varphi$ on $\partial \Omega$ up to a quadratic correction \eqref{eq:Subeq-Reg-C2-TN-new-Q-G-step-3-goal}. And Corollary~\ref{corollary:Subeq-Reg-guan-corollary} can be applied to this quadratic correction of $\partial_\tau \varphi$ instead, to complete the sketch from above.

The original argument of Caffarelli--Nirenberg--Spruck~\cite[(5.5)]{caffarelli-nirenberg-spruck-dirichlet-problem-III-MR0806416} constructed the vector field $\tau$ only for $O(n)$-invariant equations. The following lemma allows us to extend their argument to the $G$-invariant case (see Step 2 in \S\ref{subsection:Subeq-Reg-proof-of-mixed-Hessian-estimate}):
\begin{lemma}
    \label{lemma:Subeq-Reg-claim-about-lie-algebra-surjectivity}
    For a Lie group $G\subset O(n)$ acting transitively on $S^{n-1}$, with Lie algebra $\mathfrak{g} := T_IG$, there exists $C=C(G)>0$ such that
    \begin{equation}
    \label{eq:Subeq-Reg-claim-about-lie-algebra-surjectivity}
        \forall\ e\in S^{n-1},\ v\in T_eS^{n-1},\qquad \exists\ Q_v\in \mathfrak{g}  \quad \text{such that} \quad Q_ve = v \quad \text{and} \quad |Q_v|\le C|v|.
    \end{equation}
\end{lemma}
Lemma~\ref{lemma:Subeq-Reg-claim-about-lie-algebra-surjectivity} is used in Step 2 of the proof of Proposition~\ref{prop:Subeq-Reg-C2-TN-prop}.
\begin{proof}
    Consider the map $\Psi_e:G \rightarrow S^{n-1}$,
    \[\Psi_e(O) := Oe, \qquad \forall\ O\in G.\]
    Observe that $\Psi_e$ is smooth, surjective (by assumption) and equivariant under the transitive left multiplication action of $G \subset O(n)$. Computing its derivative at identity $I\in G$,
    \[d_I\Psi_e(Q) = \left. \frac{d}{dt} \right|_{t=0} e^{tQ}e = Qe, \qquad \forall\ Q\in \mathfrak{g} = T_IG.\]
    By the equivarience of rank~\cite[Theorem~7.25]{John-lee-smooth-mflds-MR2954043}, $\Psi_e$ is a submersion and in particular the derivative $d_I\Psi_e$ is surjective,
    \[\forall\ v\in T_eS^{n-1},\qquad \exists\ Q_v\in \mathfrak{g} := T_IG \quad \text{such that} \quad d_I\Psi_e(Q_v) = Q_ve = v.\]
    In other words, \eqref{eq:Subeq-Reg-claim-about-lie-algebra-surjectivity} holds. To find an estimate $C= C(G)$ on $Q_v$, note that we can make the explicit choice
    \begin{equation}
        \label{eq:Subeq-Reg-claim-about-lie-algebra-surjectivity-eq1}
        Q_v := \left( d_I\Psi_e|_{\left(\ker d_I\Psi_e\right)^\perp} \right)^{-1} v, \qquad \forall\ v\in T_eS^{n-1}.
    \end{equation}
    Finally, since $d\Psi_e$ is continuous is $e$ and $S^{n-1}$ is compact, we can choose $C(G)$ to be the maximum of the operator norms over $S^{n-1}$,
    \[{|Q_v|} \le \left( \max_{e\in S^{n-1}}\left\| \left( d_I\Psi_e|_{\left(\ker d_I\Psi_e\right)^\perp} \right)^{-1} \right\| \right) |v| =:C(G) |v|. \tag*{\qedhere}\]
\end{proof}

\subsection{Proof of Proposition~\ref{prop:Subeq-Reg-C2-TN-prop}}
\label{subsection:Subeq-Reg-proof-of-mixed-Hessian-estimate}

\begin{proof}
    Fix a boundary point and unit tangent vector,
    \begin{equation}
    \label{eq:Subeq-Reg-C2-TN-fixed-at-beginning-x0-V}
        x_0\in \partial \Omega \qquad \text{and} \qquad V\in T_{x_0} \partial \Omega,\ |V|=1.
    \end{equation}
    Steps 1--3, are concerned with constructing a vector field $\tau$ \eqref{eq:Subeq-Reg-invariant-subeq-TN-tau-formula-1-Q}, whose integral flow is drawn from the given Lie group $G\subset O(n)$, such that $\partial_\tau u$ equals $\partial_\tau \varphi$ modulo a quadratic correction \eqref{eq:Subeq-Reg-C2-TN-new-Q-G-step-3-goal}.
    
    \vspace{0.5em}\noindent\textit{Step 1:} \textit{For a $C^2$ domain $\Omega \subset \RR^n$ with inward unit normal $\nu$ \eqref{eq:Subeq-Reg-inward-normal-to-domain-bdry-notation}, $u,\varphi \in C^2(\overline{\Omega})$ such that $u=\varphi$ on $\partial \Omega$, and any vector field $\tau:\RR^n \rightarrow \RR^n$,}
    \begin{equation}
        \label{eq:Subeq-Reg-C2-TN-step1-goal}
        \partial_\tau (u - \varphi)  = \langle \tau, \nu \rangle\ \partial_\nu (u-\varphi) \qquad \text{on } \partial \Omega.
    \end{equation}

     \vspace{0.5em}\noindent Consider the tangential projection of $\tau$ along $\partial \Omega$,
     \[\tau_{\text{tan}}(x) := \tau(x) - \langle \tau(x),\nu(x) \rangle \nu(x),\qquad \forall\ x\in \partial \Omega.\]
     Since $\nu(x) \perp T_x\partial \Omega$,
     \[\tau_{\text{tan}}(x) \in T_x\partial \Omega, \qquad \forall\ x\in \partial \Omega.\]
     Since $u=\varphi$ on $\partial \Omega$. Hence,
     \[\langle D(u-\varphi)(x),\tau_{\text{tan}}(x) \rangle = 0,\qquad \forall\ x\in \partial \Omega. \]
     Hence,
     \[\begin{aligned}
         \langle D(u-\varphi),\tau \rangle =&\ \langle D(u-\varphi),\tau_{\text{tan}} \rangle + \langle \tau, \nu \rangle \langle D(u-\varphi), \nu \rangle\\
         =&\ 0+ \langle \tau, \nu \rangle \langle D(u-\varphi), \nu \rangle \qquad \text{on } \partial \Omega.
     \end{aligned}\]
     In other words, \eqref{eq:Subeq-Reg-C2-TN-step1-goal} holds.

     \vspace{0.5em}\noindent\textit{Step 2:} \textit{Recall \eqref{eq:Subeq-Reg-C2-TN-fixed-at-beginning-x0-V}; fix also $Q\in \mathfrak{g} = T_IG$, and set (cf. \eqref{eq:Subeq-Reg-tau-VQ-vector-field-def})}
    \begin{equation}
        \label{eq:Subeq-Reg-invariant-subeq-TN-tau-formula-1-Q}
        \tau(x) := V + Q(x-x_0),
    \end{equation}
    \textit{There exists $C_2 = C_2(G, \Omega)>0$ such that}
    \begin{equation}
        \label{eq:Subeq-Reg-C2-TN-step2-goal}
         \begin{aligned}
            |Q|\le&\ C_2,\\
            \left|\langle \tau(x), \nu(x)\rangle \right| \le&\ C_2|x-x_0|^2, \qquad \forall\ x\in \partial \Omega.
         \end{aligned}
    \end{equation}

    \vspace{0.5em}\noindent This  construction is a modification of that of Caffarelli--Nirenberg--Spruck~\cite[(5.5)]{caffarelli-nirenberg-spruck-dirichlet-problem-III-MR0806416}. By \eqref{eq:Subeq-Reg-invariant-subeq-TN-tau-formula-1-Q}, $\tau(x_0) = V \in T_{x_0}\partial \Omega$, and in particular, $\langle \tau(x_0), \nu(x_0) \rangle = 0$ for any choice of $Q$. The idea is to choose $Q \in \mathfrak{g}$ so that $\langle \tau, \nu\rangle \in C^2(\partial \Omega)$ has zero tangential derivative along $\partial \Omega$. Then, $\langle \tau, \nu\rangle $ would be bounded by the second order term in its Taylor expansion along $\partial \Omega$ to obtain the quadratic bound \eqref{eq:Subeq-Reg-C2-TN-step2-goal}.

    Choose an orthonormal basis for $\RR^n$ in which
    \begin{equation}
        \label{eq:Subeq-Reg-C2-TN-eq1.1}
        \nu(x_0) = e_n=(0,\dots,0,1).
    \end{equation}
    Then $T_{x_0}\partial \Omega$ is spanned by the remaining standard basis elements $e_1,\dots, e_{n-1}$. Differentiating the inner product $\langle \tau(x), \nu(x)\rangle$ along such tangent vectors at $x=x_0$,
    \begin{align}
        \left. \frac{\partial}{\partial x_k}\right|_{x=x_0} \langle \tau(x), \nu(x) \rangle =&\ \left\langle \partial_k \tau(x_0), \nu(x_0)  \right \rangle + \left\langle \tau(x_0), \partial_k\nu(x_0)  \right \rangle \nonumber \\
        =&\  \left\langle Qe_k, e_n  \right \rangle + \left\langle V, \partial_k\nu(x_0)  \right \rangle, \nonumber \\
        =&\  -\left\langle e_k, Qe_n  \right \rangle + \left\langle V, \partial_k\nu(x_0)  \right \rangle, \qquad k=1,\dots, n-1, \label{eq:Subeq-Reg-C2-TN-eq1.111}
    \end{align}
    where the last line follows because $Q^T=-Q$, because $Q \in \mathfrak{g} \subset \mathfrak{o}(n)$, and the Lie algebra $\mathfrak{o}(n)$ of $O(n)$ is the set of skew-symmetric matrices.

    Observe that the vector
    \[v_0 = \big( \left\langle V, \partial_1\nu(x_0) \right\rangle,\, \dots,\, \left\langle V, \partial_{n-1}\nu(x_0)  \right\rangle,\, 0 \big),\]
    is perpendicular to $e_n$, $\langle v,e_n\rangle = 0$. So by Lemma~\ref{lemma:Subeq-Reg-claim-about-lie-algebra-surjectivity}, there exists $Q\in \mathfrak{g}$ such that
    \begin{equation}
        \label{eq:Subeq-Reg-C2-TN-new-Q-new-def}
        Qe_n = v_0 = \big( \left\langle V, \partial_1\nu(x_0) \right\rangle,\, \dots,\, \left\langle V, \partial_{n-1}\nu(x_0)  \right\rangle,\, 0 \big),
    \end{equation}
    such that
    \begin{equation}
        \label{eq:Subeq-Reg-C2-TN-eq1.113-for-Q}
        |Q| \le C(G) |v_0| \le C(G)\sqrt{n-1}\|\nu\|_{C^1}   =:C_2' ,
    \end{equation}
    since $V$ is a unit vector \eqref{eq:Subeq-Reg-C2-TN-fixed-at-beginning-x0-V}. Since $\|\nu\|_{C^1}$ depends only on $\Omega$, $C_2' = C_2'(\Omega, G)$.
    So by \eqref{eq:Subeq-Reg-C2-TN-eq1.111}, for $Q\in \mathfrak{g}$ satisfying \eqref{eq:Subeq-Reg-C2-TN-new-Q-new-def}--\eqref{eq:Subeq-Reg-C2-TN-eq1.113-for-Q}, and $\tau = V+Q(x-x_0)$ \eqref{eq:Subeq-Reg-invariant-subeq-TN-tau-formula-1-Q},
    \begin{equation}
        \label{eq:Subeq-Reg-C2-TN-eq1.222}
        \left. \frac{\partial}{\partial x_k}\right|_{x=x_0} \langle \tau(x), \nu(x) \rangle  = 0, \qquad k =1,\dots, n-1.
    \end{equation}
    Thus $\langle \tau , \nu \rangle$ vanishes to first order at $x_0 $ along $ \partial \Omega$. It remains to convert this pointwise statement into the quadratic estimate in \eqref{eq:Subeq-Reg-C2-TN-step2-goal}. 

    \begin{figure}[t]
    \centering
    \begin{tikzpicture}[scale=3.5, domain=-1:1]
        \foreach \x in {-0.95,-0.90,...,1.00}
        {
            \draw[very thin, color=black!40!white] (\x,+0.5*\x*\x + 0.125*\x*\x*\x*\x ) -- (\x, 0.75);
        }
        \draw (0.2,-0.07) node {$x_0$};
        \draw (0.8,0.81) node {${\Omega}$};
        \draw[thick, <->] (-1,0) -- (1,0) node[anchor=north west] {$(x_1,\dots, x_{n-1})=x'$};
        \draw[thick, <->] (0, -0.1) -- (0, 0.85) node[above] {$x_n$}; 
        \draw[domain=-1:1, thick, smooth, variable=\x] plot ({\x}, {0.5*\x*\x + 0.125*\x*\x*\x*\x}) node[right] {$\operatorname{Gr}(\chi)$};
    \end{tikzpicture}
    \caption{In coordinates where $\nu(x_0) = e_n$ \eqref{eq:Subeq-Reg-C2-TN-eq1.1}, $\partial {\Omega}$ can be locally expressed as the graph of a $C^3$ function $\chi(x')$.}
    \label{fig:Subeq-Reg-partial-omega-is-locally-a-graph}
\end{figure}
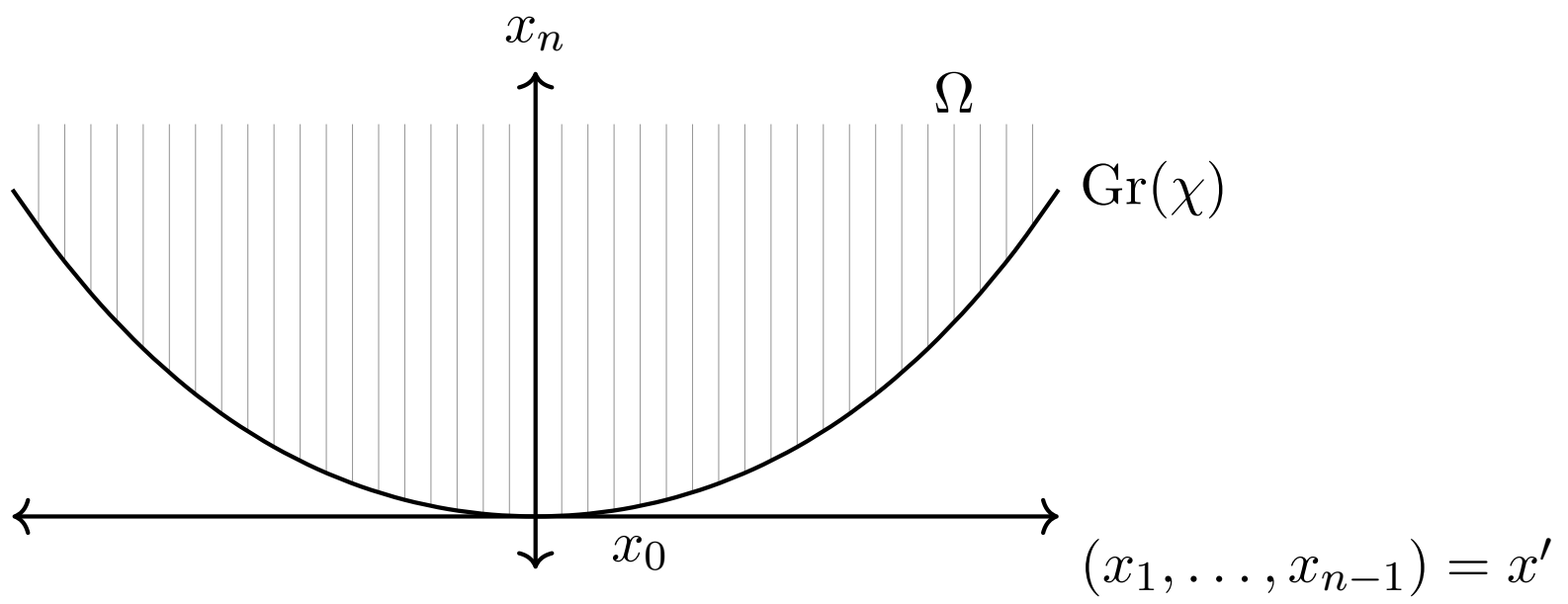
    
    Since $\partial \Omega$ is a $C^3$ manifold, by the choice of coordinates in \eqref{eq:Subeq-Reg-C2-TN-eq1.1}, there are open sets $U'\subset \RR^{n-1}$ and $U\subset \RR^n$, such that $\partial \Omega$ is the graph of $\chi \in C^3(U')$ over $U'$ Figure~\ref{fig:Subeq-Reg-partial-omega-is-locally-a-graph},
    \begin{equation}
        \label{eq:Subeq-Reg-C2-TN-eq1.3}
        (\id, \chi)(U') = U\cap \partial \Omega.
    \end{equation}
    By \eqref{eq:Subeq-Reg-invariant-subeq-TN-tau-formula-1-Q} and \eqref{eq:Subeq-Reg-C2-TN-eq1.222}, the function $\langle \tau, \nu \rangle \circ (\id,\chi)$ vanishes up to first order at $x_0' \in U'$, where $x_0'$ is the projection of $x_0$ onto the first $(n-1)$ components. Hence,
    \[\big| \left\langle \tau, \nu \right\rangle \circ (\id,\chi)(x')\big| \le C_2''|x' - x_0'|^2 \le C_2''\big|(x',\chi(x')) - (x_0', \chi(x_0'))\big|^2, \qquad \forall\ x'\in U',\]
    for some $C_2'' = C_2''\left(G, \Omega\right)$, since $\|\tau\|_{C^2}$ is controlled by $C_2'$ \eqref{eq:Subeq-Reg-C2-TN-eq1.113-for-Q} (by \eqref{eq:Subeq-Reg-invariant-subeq-TN-tau-formula-1-Q}), $\|\chi\|_{C^2}$ is controlled uniformly by $\Omega$, since $\Omega$ is bounded. Rephrasing the previous equation using \eqref{eq:Subeq-Reg-C2-TN-eq1.3},
    \begin{equation}
        \label{eq:Subeq-Reg-C2-TN-eq1.333}
        |\langle \tau(x), \nu(x) \rangle| \le C_2''|x-x_0|^2, \qquad \forall\ x\in U\cap \partial \Omega.
    \end{equation}
    Since $x_0 \in U$, $B(x_0 ,\delta) \subset U$ for some $\delta>0$, and
    \begin{equation}
        \label{eq:Subeq-Reg-C2-TN-eq1.344}
        |\langle \tau(x), \nu(x) \rangle| \le \|\langle \tau, \nu \rangle\|_{C^0(\partial \Omega)} \frac{|x-x_0|^2}{\delta^2}, \qquad \forall\ x\in \partial \Omega,\ |x-x_0|\ge \delta.
    \end{equation}
    Combining \eqref{eq:Subeq-Reg-C2-TN-eq1.113-for-Q}, \eqref{eq:Subeq-Reg-C2-TN-eq1.333}, and \eqref{eq:Subeq-Reg-C2-TN-eq1.344}, \eqref{eq:Subeq-Reg-C2-TN-step2-goal} holds for
    \[C_2 := \max\left\{ C_2',\ C_2'',\ \frac{\|\langle \tau, \nu \rangle\|_{C^0(\partial \Omega)}}{\delta^2} \right\}.\]
    
    \vspace{0.5em}\noindent\textit{Step 3:} \textit{Combining Steps 1--2, there exists $C_3 = C_3(n,F,G,\Omega, \|\varphi\|_{C^3})$ such that}
    \begin{equation}
        \label{eq:Subeq-Reg-C2-TN-new-Q-G-step-3-goal}
        -C_3|x-x_0|^2 + \partial_\tau \varphi(x) \le \partial_\tau u(x) \le  \partial_\tau \varphi(x) + C_3|x-x_0|^2, \qquad \forall\ x\in \partial \Omega.
    \end{equation}
    By Proposition~\ref{prop:Subeq-Reg-unif-grad-est}, there exists $C_3' = C_3'(F,\Omega, \|\varphi\|_{C^3})$ such that
    \[\|u-\varphi\|_{C^1} \le C_3'.\]
    So by Steps 1--2,
    \[\left|\partial_\tau (u-\varphi) \right| = \left|\langle \tau, \nu \rangle \right| \left| \partial_\nu(u-\varphi) \right| \le C_2|x-x_0|^2 C_3', \qquad \text{on } \partial \Omega. \]
    Thus, \eqref{eq:Subeq-Reg-C2-TN-new-Q-G-step-3-goal} holds.
    
    \vspace{0.5em}\noindent By \eqref{eq:Subeq-Reg-C2-TN-new-Q-G-step-3-goal}, $\partial_\tau u$ is bounded from both above and below on $\partial \Omega$, with equality at the chosen boundary point $x_0 \in \partial \Omega$. This bound is extended to the interior by applying Corollary~\ref{corollary:Subeq-Reg-guan-corollary}.

    \vspace{0.5em}\noindent\textit{Step 4:}\textit{There exists $w_-,w_+ \in C^2(\overline{\Omega})$ and $C_{4} = C_4(n, F, \Omega, G, \|\varphi\|_{C^3})$ such that}
    \begin{equation}
        \label{eq:Subeq-Reg-invariant-subeq-TN-step-5-new-goal}
        \begin{aligned}
            \|w_-\|_{C^1}, \|w_+\|_{C^1} \le&\ C_4,\\
        w_- \le \partial_\tau u \le&\ w_+ \qquad \text{on } \overline{\Omega},\\
        w_-(x_0) = \partial_\tau u(x_0) =&\ w_+(x_0).
        \end{aligned}
    \end{equation}
    
    \vspace{0.5em}\noindent By Proposition~\ref{prop:Subeq-Reg-existence-of-delta-subsol-and-supersol}~(i), there exists an admissible $\delta$-subsolution $\underline{u} \in C^2(\overline{\Omega})$ satisfying \eqref{eq:Subeq-Reg-strongly-admissible-subsolution}, and $C_4' = C_4'(F,\Omega, \|\varphi\|_{C^3})$ such that
    \[\|\underline{u}\|_{C^2} \le C_4'.\]
    So by Corollary~\ref{corollary:Subeq-Reg-guan-corollary} applied to $h = \pm \partial_\tau \varphi + C_3|x-x_0|^2$, there exists $w_+, w_- \in C^2(\overline{\Omega})$ and $C_4 = C_4(n,\Omega,F, C_4', \delta, C_2, \|\varphi\|_{C^3})$ such that by \eqref{eq:Subeq-Reg-C2-TN-new-Q-G-step-3-goal} (cf. \eqref{eq:Subeq-Reg-guan-corollary-assump-conclusion}),
    \begin{align}
        \|w_-\|_{C^1}, \|w_+\|_{C^1} \le C_4,& \label{eq:Subeq-Reg-C2-TN-new-Q-wpm-C1-norm-bound}\\
        - \Delta_{D\signdist^F(D^2u)} w_- \ge 0 = \Delta_{D\signdist^F(D^2u)} \partial_\tau u \ge \Delta_{D\signdist^F(D^2u)} w_+& \qquad \text{on } \Omega,\label{eq:Subeq-Reg-C2-TN-new-Q-wpm-Delta-ineq} \\
        -w_- \le \partial_\tau \varphi - C_2|x-x_0|^2 \le \partial_\tau u \le \partial_\tau \varphi + C_2|x-x_0|^2 \le  w_+& \qquad \text{on } \partial \Omega, \label{eq:Subeq-Reg-C2-TN-new-Q-wpm-bdry-ineq-1}\\
        -w_-(x_0) = \partial_\tau \varphi(x_0) = \partial_\tau u(x_0) =  w_+(x_0),& \label{eq:Subeq-Reg-C2-TN-new-Q-wpm-bdry-ineq}
    \end{align}
    where $ \Delta_{D\signdist^F(D^2u)} \partial_{\tau} u = 0$ in \eqref{eq:Subeq-Reg-C2-TN-new-Q-wpm-Delta-ineq} by Lemma~\ref{lemma:Subeq-Reg-symmetries-of-linearization} applied to $\tau = \tau_{V-Qx_0, Q}$ (see \eqref{eq:Subeq-Reg-tau-VQ-vector-field-def} and \eqref{eq:Subeq-Reg-invariant-subeq-TN-tau-formula-1-Q}), and where $\partial_\tau u(x_0) = \partial_\tau \varphi(x_0) $ in \eqref{eq:Subeq-Reg-C2-TN-new-Q-wpm-bdry-ineq} since $\tau(x_0) = V\in T_{x_0}\partial \Omega$ \eqref{eq:Subeq-Reg-invariant-subeq-TN-tau-formula-1-Q}.
    
    By \eqref{eq:Subeq-Reg-C2-TN-new-Q-wpm-C1-norm-bound}--\eqref{eq:Subeq-Reg-C2-TN-new-Q-wpm-bdry-ineq}, and the maximum principle for linear elliptic PDE \cite[Theorem~3.1]{gilbarg-trudinger-MR1814364} applied to $\Delta_{D\signdist^F(D^2u)}$ \eqref{eq:Subeq-Reg-linearization-of-dij-elliptic-pde} (which is elliptic by Claim~\ref{claim:Subeq-Reg-dij-is-elliptic-PDE}),
    \[\begin{aligned}
        -w_- \le \partial_\tau u \le&\  w_+ \qquad \text{on } \partial \overline{\Omega},\\
        -w_-(x_0) = \partial_\tau u(x_0) =&\  w_+(x_0).
    \end{aligned}\]
    Thus \eqref{eq:Subeq-Reg-invariant-subeq-TN-step-5-new-goal} holds.

    \vspace{0.5em}\noindent\textit{Step 5:} \textit{Equation \eqref{eq:Subeq-Reg-C2-TN-conclusion} holds:}

    \vspace{0.5em}\noindent By Claim~\ref{claim:Subeq-Reg-normal-gradient-claim} and \eqref{eq:Subeq-Reg-invariant-subeq-TN-step-5-new-goal},
    \begin{equation}
        \label{eq:Subeq-Reg-invariant-subeq-TN-step-12-goal}
        \left| \left\langle D \partial_\tau u (x_0), \ \nu(x_0) \right\rangle \right| \le \max\{\|\underline{w}\|_{C^1}, \|\overline{w}\|_{C^1}\} \le   C_4.
    \end{equation}
    Even though $\tau(x_0) = V$ \eqref{eq:Subeq-Reg-invariant-subeq-TN-tau-formula-1-Q}, observe that \eqref{eq:Subeq-Reg-invariant-subeq-TN-step-12-goal} is not equivalent to \eqref{eq:Subeq-Reg-C2-TN-conclusion} since the left hand side of \eqref{eq:Subeq-Reg-invariant-subeq-TN-step-12-goal} involves a derivative of $\tau$.
    
    Observe that
    \[D\partial_\tau u = D\langle \tau, Du\rangle = D^2u\, \tau + D\tau\, Du,\]
    and so,  
    \[\langle D \partial_\tau u ,\ \nu\rangle = \langle  D^2u\, \tau,\ \nu \rangle  + \langle D\tau\, Du,\ \nu \rangle, \qquad \text{on } \partial \Omega,\]
    where $\nu$ \eqref{eq:Subeq-Reg-inward-normal-to-domain-bdry-notation} is the inward normal. By Proposition~\ref{prop:Subeq-Reg-unif-grad-est},
    \[\|u\|_{C^1} \le C_1 = C_1( F,\Omega, \|\varphi\|_{C^3}).\]
    So by \eqref{eq:Subeq-Reg-invariant-subeq-TN-step-12-goal},
    \begin{align*}
        |\langle  D^2u(x_0)\, \tau(x_0),\ \nu(x_0) \rangle|\le&\ \left|\langle D \partial_\tau u (x_0) ,\ \nu(x_0) \rangle\right| +  |\langle D\tau(x_0)\, Du(x_0),\ \nu(x_0) \rangle | \\
        \le&\ C_4 + C_1\|\tau\|_{C^1}\|\nu \|_{C^0} =: C_5, \qquad \text{on } \partial \Omega.
    \end{align*}
    Note that $\|\tau\|_{C^1}$ depends only on $G$ and $\Omega$ \eqref{eq:Subeq-Reg-C2-TN-step2-goal}, and $\|\nu\|_{C^1}$ depends only on $\Omega$ \eqref{eq:Subeq-Reg-inward-normal-to-domain-bdry-notation}. So $C_5 = C_5(n, \Omega, F, \|\varphi\|_{C^3})$. And by \eqref{eq:Subeq-Reg-invariant-subeq-TN-tau-formula-1-Q},
    \[|\langle  D^2u(x_0)\, \tau(x_0),\ \nu(x_0) \rangle|  \le C_4 + C_1\|\tau\|_{C^1}\|\nu \|_{C^0}.\]
    Since $x_0 = 0\in \partial \Omega$ and $V \in T_{x_0} \partial \Omega$, $|V|=1$ \eqref{eq:Subeq-Reg-C2-TN-fixed-at-beginning-x0-V}, are arbitrary, and since $\tau(x_0) = V$ \eqref{eq:Subeq-Reg-invariant-subeq-TN-tau-formula-1-Q},
    \begin{equation*}
        |D^2u(x_0)V, \nu(x_0)\rangle| \le  C_5, \qquad \forall\ x_0\in \partial \Omega,\ V\in T_{x_0} \partial \Omega,
    \end{equation*}
    i.e, \eqref{eq:Subeq-Reg-C2-TN-conclusion} holds for $C:=C_5$.
\end{proof}

\section{Normal Hessian estimate via the conormal subequation}
\label{sec:Subeq-Reg-normal-hessian-estimates}

This section estimates the purely normal component of the Hessian $D^2u$ of a solution $u$ to \eqref{eq:Subeq-Reg-const-coeff-subequation-PDE}, and hence completes the proof of Theorem~\ref{thm:Subeq-Reg-a-priori-estimates}.

\subsection{Strategy for normal Hessian estimate}

Consider a boundary point $x_0\in \partial \Omega$. By Proposition~\ref{prop:Subeq-Reg-unif-grad-est}, Lemma~\ref{lemma:Subeq-Reg-C2-TT-estimate}, and Proposition~\ref{prop:Subeq-Reg-C2-TN-prop}, all the components of the Hessian $D^2u(x_0)$ are uniformly estimated except for the normal Hessian component. Furthermore, the normal Hessian component is easily lower bounded by convexity of $F$ (see Lemma~\ref{lemma:Subeq-Reg-C2-NN-lower-bound-lin-alg}). So a bound on $D^2u(x_0)$ is reduced to an upper bound on its normal component. 

To obtain a Hessian estimate on $u$, we only need to estimate $D^2u(x_0)$ at a point $x_0\partial \Omega$ maximizing the norm, $x_0 = \arg\max_{\partial \Omega}|D^2u|$. But by Proposition~\ref{prop:Subeq-Reg-C2-NN-trudinger-prop}, upper bounds of the form $|D^2u(x)|\le C$ are equivalent to a lower bound of the form $\signdist^{F_+}\circ \pi_+(D^2u(x) ) \ge \delta >0$, where $\signdist^{F_+}$ is the signed-distance to the boundary of the conormal subequation $\partial F_+$ \eqref{eq:Subeq-Reg-F+-projected-subequation}. So, it suffices to estimate $D^2u(x)$ at the point (Corollary~\ref{corollary:Subeq-Reg-C2-NN-trudinger-prop})
\begin{equation}
    \label{eq:Subeq-Reg-Trudingers-point}
    x_{\min} \in  \underset{\substack{x\in \partial \Omega,\ O_x\in G\\ O_xe_n = \nu(x)}}{{\arg\min}} \signdist^{F_+} \circ \pi_+ \left( O_x^T D^2u(x) O_x \right) ,
\end{equation}
where $\signdist^{F_+}$ is the signed-distance \eqref{eq:Subeq-Reg-signed-dist-def} to the boundary $\partial F_+$ of the conormal subequation \eqref{eq:Subeq-Reg-F+-projected-subequation} (which is a subequation by Lemma~\ref{lemma:Subeq-Reg-projection-F+-is-subeq}), and $\nu:\partial \Omega \rightarrow \RR^n$ \eqref{eq:Subeq-Reg-inward-normal-to-domain-bdry-notation} is the inward unit normal (see Remark~\ref{remark:Subeq-Reg-why-is-xmin-well-defined}).

Finally, convexity of $F_+$ is exploited to construct a barrier to $\partial_\nu u$ at $x_{\min}$. Corollary~\ref{corollary:Subeq-Reg-guan-corollary} and Claim~\ref{claim:Subeq-Reg-normal-gradient-claim} applied to this barrier yield an upper bound on the normal derivative of $\partial_\nu u$, i.e., on the normal Hessian component, similar to Step~4 in the proof of Proposition~\ref{prop:Subeq-Reg-C2-TN-prop}.

\begin{remark}[Why is $x_{\min}$ \eqref{eq:Subeq-Reg-Trudingers-point} well-defined?]~
\label{remark:Subeq-Reg-why-is-xmin-well-defined}
    \begin{itemize}
        \item $x_{\min}$ is not uniquely determined by \eqref{eq:Subeq-Reg-Trudingers-point}, but this notation is used to represent any minimizer. Furthermore, if $F_+ = \Sym^2(\RR^{n-1})$, then $\signdist^{F_+} \equiv \infty$ \eqref{eq:Subeq-Reg-signdist-F+-inf-or-convex}, and $x_{\min}$ is taken to be an arbitrary point in $\partial \Omega$.

        \item For each $x\in \partial \Omega$, there exists $O_x\in G$ such that $O_x e_n = \nu(x)$ by transitivity of $G$. Then, $O_xe_1,\dots, O_xe_{n-1}$ span the tangent space $T_x\partial \Omega$, and $\pi_+(O_x^T D^2u(x) O_x )$ is the tangential Hessian of $D^2u$ with respect to an orthonormal frame chosen from $G$.
        
        The reason to use such an identification, is because the projection $\pi_+$ \eqref{eq:Subeq-Reg-projection-onto-the-principal-minor} is chosen to be fixed and independent of $x\in \partial \Omega$. We could instead leave the Hessian $D^2u(x)$ fixed, and vary the projection $\pi_+$ as a $G$-adapted function on $\partial \Omega$.
        
        \item By Lemma~\ref{lemma:Subeq-Reg-sign-dist-projection-invariance}, $\signdist^{F_+}\circ \pi_+(O_x D^2u(x) O_x)$ is independent of the matrix $O_x\in G$ satisfying $O_xe_n = \nu(x)$. So the expression inside the $\arg\max$ in \eqref{eq:Subeq-Reg-Trudingers-point} only depends on $x\in \partial \Omega$.
    \end{itemize}
\end{remark}

\begin{lemma}
\label{lemma:Subeq-Reg-sign-dist-projection-invariance}
    For a $G$-invariant subequation $F\subset \Sym^2(\RR^n)$, and $O, O' \in G$ such that $O'^TO e_n = e_n$,
    \begin{equation}
        \label{eq:Subeq-Reg-Trud-point-xmin-indep-of-Ox}
    \end{equation}
    \[\signdist^{F_+}\circ \pi_+(O^TAO) = \signdist^{F_+}(O'^TAO'), \qquad \forall\ A\in \Sym^2(\RR^n).\]
\end{lemma}
\begin{proof}
    Since $O', O \in G$ and $O '^TO e_n = e_n$,
        \[O '^TO  = \begin{pmatrix}
            O '' & 0 \\
            0& 1
        \end{pmatrix} \in G, \qquad \text{for some } O '' \in O(n-1).\]
        By $G$-invariance of $F$,
        \[O ^T O ' FO '^TO = F.\]
        Applying $\pi_+$ \eqref{eq:Subeq-Reg-projection-onto-the-principal-minor} on both sides and taking closures \eqref{eq:Subeq-Reg-F+-projected-subequation},
        \[O ''^T F_+ O '' = F_+.\]
        By Lemma~\ref{lemma:Subeq-Reg-signed-distance-lemma}~(iv),
        \[\signdist^{F_+}(O ''^T A_+ O '')
        = \signdist^{F_+}(A_+),
        \qquad \forall\ A_+\in\Sym^2(\RR^{n-1}).\]
        But observe that,
        \[\pi_+(O ^T   A O )
        = O ''^T \pi_+(O '^T   A O ')O '',\]
        conclude that
        \[\signdist^{F_+}\circ \pi_+(O ^T   A O )
        = \signdist^{F_+} \big( O''^T \pi_+ (O '^T   A O ') O'' \big)= \signdist^{F_+}\circ \pi_+ (O '^T   A O ').\]
\end{proof}

\begin{remark}[Comparison to Trudinger]
    Trudinger's approach to the normal Hessian estimate for equations of the form \eqref{eq:Subeq-Reg-intro-purely-second-order-PDE}, was to define 
    \begin{equation}
        \label{eq:Subeq-Reg-Trudinger-asymp-operator}
        \fcal^\infty(B) := \lim_{R\rightarrow \infty} \fcal\begin{pmatrix}
        B& 0 \\
        0 & R
    \end{pmatrix}.
    \end{equation}
    By \eqref{eq:Subeq-Reg-formula-for-sign-dist-to-F+}, $\signdist^{F_+}$ coincides with an asymptotic operator when $\fcal = \signdist^F$. Trudinger observed that a $C^2$ estimate on $D^2u(x_0)$, at a point where $\nu(x_0) = e_n$ \eqref{eq:Subeq-Reg-inward-normal-to-domain-bdry-notation}, is equivalent to a positive lower bound $\fcal^\infty(\pi_+(D^2u(x_0)))$.
    
    Proposition~\ref{prop:Subeq-Reg-C2-NN-trudinger-prop} reveals that this observation is really about the geometry of the conormal subequation $F_+ = \overline{\pi_+(F)}$ \eqref{eq:Subeq-Reg-F+-projected-subequation}, and how far the tangential projection $\pi_+(D^2u(x_0))$ lies from the boundary $\partial F_+$.
\end{remark}

\subsection{Sufficient to estimate at \texorpdfstring{$x_{\min}$}{x\_min}}

Proposition~\ref{prop:Subeq-Reg-C2-NN-trudinger-prop} reduces the Hessian estimate to a Hessian estimate at $x_{\min}$ \eqref{eq:Subeq-Reg-Trudingers-point}:

\begin{corollary}
\label{corollary:Subeq-Reg-C2-NN-trudinger-prop}
Under the assumptions of Theorem~\ref{thm:Subeq-Reg-a-priori-estimates}, with $\signdist^{F_+}$ being the signed-distance to $\partial F_+$ \eqref{eq:Subeq-Reg-F+-projected-subequation}, if $x\in \partial \Omega$ and any $O_x\in G$ satisfies $O_x e_n = \nu(x)$, then the following hold: 
\begin{enumerate}[label={(\roman*)}]
    \item If $\signdist^{F_+} \circ\pi_+(O_xD^2u(x)O_x^T) \ge \delta>0$ then there exists $C = C(n,\Omega, F, \|\varphi\|_{C^3}, \delta)>0$ such that
    \[|D^2u(x)|\le C.\]
    \item If $|D^2u(x)|\le C$, then there exists $\delta = \delta(F,C)>0$ such that
    \[\signdist^{F_+} \circ\pi_+(O_xD^2u(x)O_x^T) \ge \delta.\]
    \item If $|D^2u(x_{\min})|\le C_1$ (see \eqref{eq:Subeq-Reg-Trudingers-point}), then there exists $C = C(n,\Omega, F, \|\varphi\|_{C^3}, C_1)>0$ such that
    \[|D^2u(x)|\le C, \qquad \forall\ x\in \partial \Omega.\]
\end{enumerate}
\end{corollary}
\begin{proof}
    By Proposition~\ref{prop:Subeq-Reg-unif-grad-est}, Lemma~\ref{lemma:Subeq-Reg-C2-TT-estimate}, and Proposition~\ref{prop:Subeq-Reg-C2-TN-prop}, all components of $D^2u(x)$, $x\in \partial \Omega$, are uniformly estimated, except for the normal component. In particular, the matrix
    \begin{equation}
    \label{eq:Subeq-Reg-trudinger-cor-matrix-A-def}
        A := O_x^TD^2u(x)O_x ,
    \end{equation}
    satisfies \eqref{eq:Subeq-Reg-C2-NN-trudinger-prop-tan-mix-bound}, and Proposition~\ref{prop:Subeq-Reg-C2-NN-trudinger-prop} applies. Note that $A\in \partial F$ since $D^2u(x)\in \partial F$ \eqref{eq:Subeq-Reg-const-coeff-subequation-PDE} and since $F$ is $G$-invariant.
    
    By transitivity of $G\subset O(n)$, pick $O_x\in G$ such that $O_x e_n = \nu(x)$. Then, $O_xe_1,\dots, O_xe_{n-1}$ span $T_x\partial \Omega$. By \eqref{eq:Subeq-Reg-tangential-hessian-estimate} and \eqref{eq:Subeq-Reg-C2-TN-conclusion}, there exists a constant $C_2 =C_2(n,\Omega, F, \|\varphi\|_{^3})>0$ such that
    \[|\langle D^2u(x) O_xe_i, O_xe_j \rangle| <C_2,\qquad i,j=1,\dots,n,\ (i,j)\neq (n,n). \]
    Then Corollary~\ref{corollary:Subeq-Reg-C2-NN-trudinger-prop}~(i)--(ii) follow by applying Proposition~\ref{prop:Subeq-Reg-C2-NN-trudinger-prop} to $A = O_x^TD^2u(x)O_x$ \eqref{eq:Subeq-Reg-trudinger-cor-matrix-A-def}. 

    Finally, if $|D^2u(x_{\min})|\le C_1$, then by \eqref{eq:Subeq-Reg-Trudingers-point} and part (ii), there exists $\delta_1 = \delta_1(F,C_1)>0$ such that
    \[\signdist^{F_+} \circ\pi_+(O_xD^2u(x)O_x^T) \ge \signdist^{F_+} \circ \pi_+ (O_{x_{\min}}D^2u(x_{\min})O_{x_{\min}}^T) \ge \delta_1, \qquad \forall\ x\in \partial \Omega,\]
    where $O_{x}, O_{x_{\min}} \in G$ are any two matrices such that $O_xe_n =\nu(x)$ and $O_{x_{\min}}e_n = \nu(x_{\min})$. Hence, Corollary~\ref{corollary:Subeq-Reg-C2-NN-trudinger-prop}~(iii) follows by applying part (i) at each $x\in \partial \Omega$. 
\end{proof}

\subsection{Normal Hessian estimate at \texorpdfstring{$x_{\min}$}{x\_min}}
\label{subsection:Subeq-Reg-C2-NN-trudingers-point}

This section is concerned with an upper bound on the normal component of the Hessian $D^2u(x_{\min})$ at a boundary point $x_{\min} \in \partial \Omega$ \eqref{eq:Subeq-Reg-Trudingers-point}. By \eqref{eq:Subeq-Reg-signdist-F+-inf-or-convex}, there are two cases to consider---(i) $F_+ = \Sym^2(\RR^{n-1})$ and $\signdist^{F_+}\equiv \infty$ and (ii) $\signdist^{F_+}$ is finite and concave. In the first case, a Hessian estimate follows easily from Proposition~\ref{prop:Subeq-Reg-C2-NN-trudinger-prop} (Step 1). Hence, only the second case needs an additional argument.

\begin{proposition}
    \label{prop:Subeq-Reg-C2-NN-upper-bound-at-ymin}
    Under the assumptions of Theorem~\ref{thm:Subeq-Reg-a-priori-estimates}, there exists $C = C(n,\Omega, F,\|\varphi\|_{C^4})$ such that
    \begin{equation}
        \label{eq:Subeq-Reg-C2-NN-upper-bound-at-ymin-conclusion}
        \langle D^2u(x_{\min})\, \nu(x_{\min}),\nu(x_{\min}) \rangle \le C,
    \end{equation}
    where $x_{\min}$ is the point \eqref{eq:Subeq-Reg-Trudingers-point} and $\nu:\partial \Omega \rightarrow \RR^n$ is the inward unit normal \eqref{eq:Subeq-Reg-inward-normal-to-domain-bdry-notation}.
\end{proposition}

This argument is inspired by the normal Hessian estimate of Trudinger~\cite[pp. 157--158]{Trudinger1995-DHE}.

\begin{proof}
    \noindent\textit{Step 1.} \textit{If $F = \Sym^2(\RR^{n-1})$, then \eqref{eq:Subeq-Reg-C2-NN-upper-bound-at-ymin-conclusion} holds:}

    \vspace{0.5em}\noindent By transitivity of $G$ on $S^{n-1}$, there exists $O_{0}\in G$ such that $O_{0} e_n = \nu({x_{\min}})$ \eqref{eq:Subeq-Reg-inward-normal-to-domain-bdry-notation}. In particular, $O_{0}e_1,\dots, O_{0} e_{n-1}$ span $T_x\partial \Omega$, and by Lemma~\ref{lemma:Subeq-Reg-C2-TT-estimate} and Proposition~\ref{prop:Subeq-Reg-C2-TN-prop},
    \[\left| \langle D^2u(x_{\min})\, O_0 e_i,\ O_0e_j \rangle \right| \le C_0 = C_0(n,F,\Omega, G,\|\varphi\|_{C^3}).\]
    But since $F= \Sym^2(\RR^{n-1})$, by \eqref{eq:Subeq-Reg-signdist-F+-inf-or-convex},
    \[ \signdist^{F_+} \circ \pi_+\left( O_0^T D^2u(x_{\min}) O_0 \right)  = \infty \ge 1. \]
    So by Proposition~\ref{prop:Subeq-Reg-C2-NN-trudinger-prop}, there exists $C = C(1,n,F,C_0)$ such that
    \[ |O_0^TD^2u(x_{\min})O_0| = |D^2u(x_{\min})| \le C. \]
    This completes the proof of Step 1.
    
    \vspace{0.5em}\noindent By Step 1 and \eqref{eq:Subeq-Reg-signdist-F+-inf-or-convex} assume for the rest of the argument that
    \begin{equation}
        \label{eq:Subeq-Reg-NN-C2-new-F+-finite-concave}
        F_+ \subsetneq \Sym^2(\RR^{n-1}) \qquad \text{and} \qquad  \signdist^{F_+} \quad\text{is finite and concave.}
    \end{equation}

    \vspace{0.5em}\noindent\textit{Step 2.} \textit{There exists a neighborhood $U\subset \partial \Omega$ of $x_{\min}$ and a $C^4$ map $O: U \rightarrow G \subset O(n)$, such that}
    \begin{equation}
        \label{eq:Subeq-Reg-C2-NN-trud-prop-new-step-1-goal-Ox}
        O_x e_n = \nu(x), \qquad \forall\ x\in U.
    \end{equation}
    
    \vspace{0.5em}\noindent As seen in the proof of Lemma~\ref{lemma:Subeq-Reg-claim-about-lie-algebra-surjectivity}, the map $\Psi:G \rightarrow S^{n-1}$,
    \[\Psi(O) := Oe_n, \qquad \forall\ O\in G,\]
    is a submersion by surjectivity of $\Psi$, and equivarience of rank \cite[Theorem~7.25]{John-lee-smooth-mflds-MR2954043}. Since $\Psi(I) = \nu(x_{\min})$, by the implicit function theorem, there exists a neighborhood $V\subset S^{n-1}$ of $\nu(x_{\min})$ and a smooth map $\psi:V \rightarrow G$, such that
    \begin{equation}
        \label{eq:Subeq-Reg-C2-NN-trud-prop-new-step-1-psi-def-prop}
        \begin{aligned}
            \Psi\circ \psi(v) = \psi(v)e_n =&\ v, \qquad \forall\ v\in V \subset S^{n-1}.
        \end{aligned}
    \end{equation}
    Consider the open set $U = \nu^{-1} (V)$, and the map $O: U \rightarrow G \subset O(n)$,
    \[O_x := \psi(\nu(x)).\]
    Then, by \eqref{eq:Subeq-Reg-C2-NN-trud-prop-new-step-1-psi-def-prop},
    \[\begin{aligned}
            O_x e_n= \psi(\nu(x))e_n =&\ \nu(x), \qquad \forall\ x\in U = \nu^{-1}(V).
        \end{aligned}\]
    This completes the proof of Step 2.

    \vspace{0.5em}\noindent\textit{Notation:} For $U \subset \partial \Omega$ and $O: U \rightarrow G$ as in Step 1, consider
    \begin{equation}
        \label{eq:Subeq-Reg-C2-NN-tangential-hessian-def-new-new}
        \begin{aligned}
            [D^2u]_{O_x} :=&\ \pi_+\left( O_x^T D^2u(x) O_x \right) \\
            =&\  \left( \left\langle D^2u(x)\, O_xe_i,\, O_xe_j  \right\rangle \right)_{i,j=1}^{n-1} \in \Sym^2(\RR^{n-1}), \qquad \forall\ x\in U.
        \end{aligned}
    \end{equation}
    Since $O_x e_n = \nu(x)$, it follows that $O_xe_1,\dots, O_xe_{n-1}$ span the tangent space $T_x\partial \Omega$. So $[D^2u]_{O_x}$ is a matrix representation of the tangential minor of the Hessian $D^2u(x)$ for $x\in U$. By \eqref{eq:Subeq-Reg-tangential-hessian-estimate} and \eqref{eq:Subeq-Reg-C2-TN-conclusion}, there exists $C_1 = C_1(n,F,\Omega, \|\varphi\|_{C^3})>0$ such that
    \begin{equation}
        \label{eq:Subeq-Reg-C2-NN-est-at-xmin-tan-mix}
        \left| \langle D^2u(x)\, O_{x} e_i,\, O_xe_j \rangle \right| \le C_1,\qquad \forall\ x\in U,\ i,j=1,\dots, n,\ (i,j)\neq (n,n).
    \end{equation}
    In particular, there is $C_1' = C_1'(C_1,n)$ such that
    \begin{equation}
        \label{eq:Subeq-Reg-C2-NN-est-at-xmin-only-tan-part}
        \left| [D^2u]_{O_x} \right| \le C_1',\qquad \forall\ x\in U.
    \end{equation}
    Analogous to \eqref{eq:Subeq-Reg-C2-NN-tangential-hessian-def-new-new}, also define the $G$-adapted tangential Hessian $[D^2{v}]_{O_x}$, for any $v\in C^2(\overline{\Omega})$.
    
    Similarly, use the local orthonormal frame $O_xe_1,\dots, O_xe_{n-1} \in T_x\partial \Omega$, $x\in U$, to identify the second fundamental form $\II_{\partial \Omega}$ \eqref{eq:Subeq-Reg-second-fundamental-form} as an $(n-1)$ dimensional symmetric matrix,
    \begin{equation}
        \label{eq:Subeq-Reg-C2-NN-tangential-hessian-notation-for-arb-A}
        [\II_{\partial \Omega}]_{O}: U \rightarrow \Sym^2(\RR^{n-1}),\qquad [\II_{\partial \Omega}]_{O_x} :=    \big( \left\langle \II_{\partial \Omega}(x)\, O_xe_i,\ O_xe_j \right\rangle \big)_{i,j=1}^{n-1}.
    \end{equation}

    \vspace{0.5em}\noindent\textit{Step 3.} \textit{Any subdifferential $B_1 \in \Sym^2(\RR^{n-1})$ to $\signdist^{F_+}$ at $[D^2u]_{O_x}$ \eqref{eq:Subeq-Reg-C2-NN-tangential-hessian-def-new-new} satisfies}
    \begin{equation}
        \label{eq:Subeq-Reg-C2-NN-bound-on-subdiff-B1}
        B_1 \ge 0,\qquad |B_1|<C_1'', \qquad  \tr(B_1) \ge \delta_1,
    \end{equation}
    \textit{for some constants $C_1'', \delta_1>0$ depending only on $C_1'$ \eqref{eq:Subeq-Reg-C2-NN-est-at-xmin-only-tan-part} and $F$:}
    
    \vspace{0.5em}\noindent By concavity of $\signdist^{F_+} $ \eqref{eq:Subeq-Reg-NN-C2-new-F+-finite-concave}, it admits a subdifferential at $[D^2u]_{O_{x_{\min}}}$ \eqref{eq:Subeq-Reg-C2-NN-tangential-hessian-def-new-new}, say
    \[B_1 \in \partial \signdist^{F_+}\left( [D^2u]_{O_{x_{\min}}} \right), \]
    such that
    \begin{equation}
        \label{eq:Subeq-Reg-C2-NN-fcal-tilde-concave-subdifferential}
        \signdist^{F_+}(B) \le \signdist^{F_+}\left( [D^2u]_{O_{x_{\min}}}\right) + \left\langle B_1,\ B - [D^2u]_{O_{x_{\min}}}\right\rangle,\qquad \forall\ B\in \Sym^2(\RR^{n-1}).
    \end{equation}
    Letting $B = P+ [D^2u|_{T}] (x_{\min})$ in \eqref{eq:Subeq-Reg-C2-NN-fcal-tilde-concave-subdifferential}, by \eqref{eq:Subeq-Reg-basic-monotonicity-for-sign-dist}
    \[\langle B_1, P \rangle  \ge  \signdist^{F_+}\left( P+ [D^2u]_{O_{x_{\min}}}\right) - \signdist^{F_+}\left( [D^2u]_{O_{x_{\min}}}\right) \ge 0, \qquad \forall\ P\ge 0.\]
    In other words,
    \begin{equation}
        \label{eq:Subeq-Reg-C2-NN-bound-on-subdiff-B1-part-0}
        B_1\ge 0.
    \end{equation}
    Furthermore, $|[D^2u]_{O_{x_{\min}}}|\le C_1'$ by \eqref{eq:Subeq-Reg-C2-NN-est-at-xmin-only-tan-part} and in particular $[D^2u]_{O_{x_{\min}}}$ lies in a bounded set. So the subdifferential $B_1$ to $\signdist^{F_+}$ at $[D^2u]_{O_{x_{\min}}}$ is bounded in terms of $C_1'$ and a $C^{0,1}$ norm of $\signdist^{F_+}$ (which depends only on $F$),
    \begin{equation}
        \label{eq:Subeq-Reg-C2-NN-bound-on-subdiff-B1-part-1}
        |B_1|<C_1'',
    \end{equation}
    for some $C_1'' = C_1''(C_1', F)$. It remains to prove the trace lower bound in \eqref{eq:Subeq-Reg-C2-NN-bound-on-subdiff-B1}.
    
    By \eqref{eq:Subeq-Reg-formula-for-sign-dist-to-F+} and monotonicity of $\signdist^F$ \eqref{eq:Subeq-Reg-basic-monotonicity-for-sign-dist}, $\signdist^{F_+}(tI_{n-1}) = \signdist^{F_+}\circ \pi_+(tI) \ge \signdist^F(tI)$, where $I_{n-1}$ denotes the $(n-1)$-dimensional identity matrix. So by \eqref{eq:Subeq-Reg-C2-NN-est-at-xmin-only-tan-part}, \eqref{eq:Subeq-Reg-C2-NN-fcal-tilde-concave-subdifferential}, and \eqref{eq:Subeq-Reg-C2-NN-bound-on-subdiff-B1-part-1},
    \[\begin{aligned}
        \langle B_1, tI_{n-1}\rangle \ge&\  \signdist^{F_+}(t_{I_{n-1}}) - \signdist^{F_+} ([D^2u]_{O_{x_{\min}}}) +  \langle B_1, [D^2u]_{O_{x_{\min}}}\rangle\\
        \ge&\ \signdist^{F}(tI) - \signdist^{F_+} ([D^2u]_{O_{x_{\min}}}) +  \langle B_1, [D^2u]_{O_{x_{\min}}}\rangle \\
        \ge&\ \signdist^{F}(tI) - \max_{|B|\le C_1'}\signdist^{F_+}(B) - C_1'C_1'',
    \end{aligned}\]
    Applying the linear growth of $\signdist^{F}$ \eqref{eq:Subeq-Reg-signed-distance-lemma-linear-growth},
    \[\begin{aligned}
        \tr(B_1) = \langle B_1, I_{n-1}\rangle \ge&\  \frac{\signdist^{F}(tI)}{t} - \frac{\max_{|B|\le C_1'}\signdist^{F_+}(B) + C_1'C_1''}{t} \\
        \ge&\ 1 - \frac{M+\max_{|B|\le C_1'}\signdist^{F_+}(B) + C_1'C_1''}{t}.
    \end{aligned}\]
    So \eqref{eq:Subeq-Reg-C2-NN-bound-on-subdiff-B1} follows, after picking $t$ sufficiently large, and combining with \eqref{eq:Subeq-Reg-C2-NN-bound-on-subdiff-B1-part-0}--\eqref{eq:Subeq-Reg-C2-NN-bound-on-subdiff-B1-part-1}.

    \vspace{0.5em}\noindent\textit{Step 4.} \textit{There exists a constant $C_2 = C_2(\delta_1, n, F, \Omega, \|\varphi\|_{C^3}) >0$ such that either}
    \[|D^2u(x_{\min})|\le C_2,\]
    \textit{(and in particular \eqref{eq:Subeq-Reg-C2-NN-upper-bound-at-ymin-conclusion} holds), or for any subdifferential $B_1 \in \Sym^2(\RR^{n-1})$ to $\signdist^{F_+}$ at $[D^2u]_{O_{x_{\min}}}$ \eqref{eq:Subeq-Reg-C2-NN-tangential-hessian-notation-for-arb-A} as in Step 3, there is a small neighborhood $U_0 \subset U$ depending only on $C_1''$ \eqref{eq:Subeq-Reg-C2-NN-bound-on-subdiff-B1} and $\Omega$ such that}
    \begin{equation}
        \label{eq:Subeq-Reg-C2-NN-est-at-xmin-step3.5-goal}
        \left\langle B_1, [\II_{\partial  \Omega}]_{O_x}  \right\rangle \ge C_2^{-1} >0, \qquad \forall\ x\in U_0,
    \end{equation}
    \textit{where $[\II_{\partial \Omega}]_{O_x}$ is the matrix representation \eqref{eq:Subeq-Reg-C2-NN-tangential-hessian-notation-for-arb-A} of the second fundamental form \eqref{eq:Subeq-Reg-second-fundamental-form}:}
    
    \vspace{0.5em}\noindent Since $\partial \Omega$ is strictly $\vec{F}$-convex, by Proposition~\ref{prop:Subeq-Reg-existence-of-delta-subsol-and-supersol}~(i), for any $\delta>0$, there exists an admissible $\delta$-subsolution $\underline{u} \in C^2(\overline{\Omega})$, satisfying \eqref{eq:Subeq-Reg-strongly-admissible-subsolution}
    \[\begin{aligned}
        D^2\underline{u} - \delta I \in F & \qquad \text{on } \Omega,\\
        \underline{u} = \varphi& \qquad \text{on } \partial \Omega,
    \end{aligned}\]
    with a $C^2$ bound
    \begin{equation}
    \label{eq:Subeq-Reg-C2-NN-at-xmin-usubsol-C2-bound}
        \|\underline{u}\|_{C^2(\overline{\Omega})} \le C_{sub} = C_{sub}(n,\Omega, F, \|\varphi\|_{C^2}, \delta).
    \end{equation}
    By \eqref{eq:Subeq-Reg-signed-dist-def} and \eqref{eq:Subeq-Reg-strongly-admissible-subsolution}, observe that
    \[\signdist^F(D^2\underline{u} - \delta I) \ge 0.\]
    Since $O_x\in G$ \eqref{eq:Subeq-Reg-C2-NN-trud-prop-new-step-1-goal-Ox}, by $G$-invariance of $\signdist^F$ (Lemma~\ref{lemma:Subeq-Reg-signed-distance-lemma}~(iv)),
    \[\signdist^F(O_x^T D^2u(x)O_x - \delta I) \ge 0, \qquad \forall\ x\in U.\]
    So by monotonicity of $\signdist^F$ \eqref{eq:Subeq-Reg-basic-monotonicity-for-sign-dist}, and \eqref{eq:Subeq-Reg-formula-for-sign-dist-to-F+} applied to $O_x^T D^2\underline{u}(x)O_x$, using the notation $[D^2\underline{u}]_{O_x} = \pi_+(O_x^T D^2\underline{u}(x)O_x)$ \eqref{eq:Subeq-Reg-C2-NN-tangential-hessian-def-new-new},
    \begin{equation}
    \label{eq:Subeq-Reg-adm-subsol-sign-dist-infty-of-subsol-delta-non-neg}
        \signdist^{F_+}([D^2\underline{u}]_{O_x} -\delta I_{n-1}) \ge 0, \qquad \forall\ x\in U,
    \end{equation}
    where $I_{n-1}$ is the $(n-1)$-dimensional identity matrix, and $U\subset \partial \Omega$ is the open set from Steps 2. By \eqref{eq:Subeq-Reg-C2-NN-bound-on-subdiff-B1}--\eqref{eq:Subeq-Reg-C2-NN-fcal-tilde-concave-subdifferential} applied for $B = [D^2\underline{u}]_{O_{x_{\min}}} - \delta I$, and \eqref{eq:Subeq-Reg-adm-subsol-sign-dist-infty-of-subsol-delta-non-neg},
    \begin{align}
        \signdist^{F_+} \left([D^2u]_{O_{x_{\min}}} \right)\ge&\ \signdist^{F_+}\left([D^2\underline{u}]_{O_{x_{\min}}} - \delta I_{n-1} \right)  \nonumber \\
        &\ - \left\langle B_1,\ [D^2\underline{u}]_{O_{x_{\min}}} - \delta I_{n-1} - [D^2u]_{O_{x_{\min}}} \right\rangle  \nonumber \\
        \ge&\ 0-\left\langle B_1,\ [D^2\underline{u}]_{O_{x_{\min}}}- [D^2u]_{O_{x_{\min}}} \right\rangle + \delta \left\langle B_1, I_{n-1} \right\rangle  \nonumber \\
        \ge&\ -\left\langle B_1,\ [D^2\underline{u}]_{O_{x_{\min}}}- [D^2u]_{O_{x_{\min}}} \right\rangle + \delta \tr(B_1) \nonumber \\
        \ge&\ -\left\langle B_1,\ [D^2\underline{u}]_{O_{x_{\min}}}- [D^2u]_{O_{x_{\min}}} \right\rangle + \delta \delta_1. \label{eq:Subeq-reg-admis-subsol-sign-dist-infty-subdiff-ineq-1}
    \end{align}
    Suppose that
    \[\left\langle B_1,\ [D^2\underline{u}]_{O_{x_{\min}}}- [D^2u]_{O_{x_{\min}}} \right\rangle \le \delta \delta_1.\]
    Then by \eqref{eq:Subeq-reg-admis-subsol-sign-dist-infty-subdiff-ineq-1},
    \[\begin{aligned}
        \signdist^{F_+} ([D^2u]_{O_{x_{\min}}}) \ge&\ - \left\langle B_1,\ [D^2\underline{u}]_{O_{x_{\min}}}- [D^2u]_{O_{x_{\min}}} \right\rangle + \delta\delta_1 \\
        \ge&\ \delta\delta_1/2 >0.
    \end{aligned}\]
    So by Proposition~\ref{prop:Subeq-Reg-C2-NN-trudinger-prop}, exists $C_2' = C_2'(n, F, C_1 ,\delta, \delta_1)>0$ \eqref{eq:Subeq-Reg-C2-NN-est-at-xmin-tan-mix} such that
    \begin{equation}
        \label{eq:Subeq-Reg-adm-subsol-3-case-a-final}
        |D^2u(x_{\min})| \le C_2'.
    \end{equation}
    In particular, \eqref{eq:Subeq-Reg-C2-NN-upper-bound-at-ymin-conclusion} holds if $\left\langle B_1,\ [D^2\underline{u}]_{O_{x_{\min}}}- [D^2u]_{O_{x_{\min}}} \right\rangle \le \delta \delta_1.$.

    On the other hand, suppose
    \begin{equation}
        \label{eq:Subeq-Reg-C2-NN-est-at-xmin-new-case-b-assump}
        \left\langle B_1,\ [D^2\underline{u}]_{O_{x_{\min}}}- [D^2u]_{O_{x_{\min}}} \right\rangle \ge \delta\delta_1/2.
    \end{equation}
    By \eqref{eq:Subeq-Reg-const-coeff-subequation-PDE} and \eqref{eq:Subeq-Reg-strongly-admissible-subsolution}, $u=\underline{u}$ on $\partial \Omega$, and so by Claim~\ref{claim:Subeq-Reg-C2-TT-key-claim},
    \begin{equation}
        \label{eq:Subeq-reg-admis-subsol-sign-dist-infty-subdiff-ineq-3.5}
        [D^2\underline{u}]_{O_{x_{\min}}} - [D^2u]_{O_{x_{\min}}} = \left( \partial_\nu u-\partial_\nu \underline{u}\right)(x_{\min}) [\II_{\partial \Omega}]_{O_{x_{\min}}},
    \end{equation}
    where $\partial_\nu$ \eqref{eq:Subeq-Reg-inward-normal-to-domain-bdry-notation-2} is the derivative along the inward unit normal $\nu$ \eqref{eq:Subeq-Reg-inward-normal-to-domain-bdry-notation} to $\partial \Omega$, and $[\II_{\partial \Omega} ]_{O_{x_{\min}}}$ is the matrix representation \eqref{eq:Subeq-Reg-C2-NN-tangential-hessian-notation-for-arb-A}.
    
    By \eqref{eq:Subeq-Reg-const-coeff-subequation-PDE}, \eqref{eq:Subeq-Reg-strongly-admissible-subsolution}, and the comparison principle~\cite[Theorem~6.4]{HarLaw-dirichletduality},
    \[\begin{aligned}
        \underline{u}\le u \qquad &\text{on } \overline{\Omega},\\
        \underline{u}= u \qquad &\text{on } \partial \Omega.
    \end{aligned}\]
    In particular, by Claim~\ref{claim:Subeq-Reg-normal-gradient-claim},
    \begin{equation}
        \label{eq:Subeq-Reg-NN-new-Step3-case-b-nu-eq}
        0\le \partial_\nu u(x_{\min}) -\partial_\nu \underline{u}(x_{\min}) \le C_2'',
    \end{equation}
    where the upper bound $C_2'' = C_2''(n,\Omega, F, \|\varphi\|_{C^3}, \delta)$ comes from Proposition~\ref{prop:Subeq-Reg-unif-grad-est} and \eqref{eq:Subeq-Reg-C2-NN-at-xmin-usubsol-C2-bound}. So, by \eqref{eq:Subeq-reg-admis-subsol-sign-dist-infty-subdiff-ineq-3.5} and the assumption \eqref{eq:Subeq-Reg-C2-NN-est-at-xmin-new-case-b-assump},
    \begin{align}
        \delta \delta_1/3 \le&\ \left\langle B_1,\ [D^2\underline{u}]_{O_{x_{\min}}}- [D^2u]_{O_{x_{\min}}} \right\rangle \nonumber \\
         =&\ \left( \partial_\nu u(x_{\min})-\partial_\nu \underline{u}(x_{\min})\right) \left\langle B_1,\ [\II_{\partial \Omega}]_{O_{x_{\min}}} \right\rangle. \label{eq:Subeq-Reg-adm-subsol-3-case-b-eq2}
    \end{align}
    Observe that $\partial_\nu u(x_{\min}) -\partial_\nu \underline{u}(x_{\min})$ is non-zero by \eqref{eq:Subeq-Reg-adm-subsol-3-case-b-eq2} and so positive by \eqref{eq:Subeq-Reg-NN-new-Step3-case-b-nu-eq}. So by \eqref{eq:Subeq-Reg-NN-new-Step3-case-b-nu-eq} and \eqref{eq:Subeq-Reg-adm-subsol-3-case-b-eq2},
    \[\left\langle B_1,\ [\II_{\partial \Omega}]_{O_{x_{\min}}} \right\rangle \ge \frac{\delta \delta_1}{3C_2''} >0.\]
    By \eqref{eq:Subeq-Reg-C2-NN-tangential-hessian-notation-for-arb-A}, $[\II_{\partial \Omega}]_{O_x}$ is continuous on $U$, and by \eqref{eq:Subeq-Reg-C2-NN-est-at-xmin-step3.5-goal}, $|B_1|\le C_1''$. So, there exists a neighborhood $U_0 \subset U$ of $x_{\min}$, where $U\subset \partial \Omega$ is from Step 2, such that
    \begin{equation}
        \label{eq:Subeq-Reg-adm-subsol-3-case-b-eq3-final}
        \left\langle B_1,\ [\II_{\partial \Omega}]_{O_x} \right\rangle \ge \frac{\delta \delta_1}{6C_2''} >0, \qquad \forall\ x\in U_0.
    \end{equation}
    Furthermore, $U_0$ depends only on $C_1''$, $\II_{\partial \Omega}$ and $x\mapsto O_x$ \eqref{eq:Subeq-Reg-C2-NN-trud-prop-new-step-1-goal-Ox}. Finally, by \eqref{eq:Subeq-Reg-adm-subsol-3-case-a-final} and \eqref{eq:Subeq-Reg-adm-subsol-3-case-b-eq3-final}, Step 4 follows for $C_2 := \max\{C_2',\, {6C_2''} /{\delta \delta_1}\}$.
    
    \vspace{0.5em}\noindent By Step 4, without loss of generality, assume \eqref{eq:Subeq-Reg-C2-NN-est-at-xmin-step3.5-goal} holds for the rest of the argument.

    \vspace{0.5em}\noindent\textit{Step 5.} \textit{For the neighborhood $U_0\subset \partial \Omega$ of $x_{\min}$ from Step 4 on which \eqref{eq:Subeq-Reg-C2-NN-est-at-xmin-step3.5-goal} holds, there exits a smooth function $h_0\in C^2(U_0)$ whose $C^2$ norm depends only on $n$, $\Omega$, $F$, and $\|\varphi\|_{C^4}$ such that} 
    \begin{equation}
        \label{eq:Subeq-Reg-C2-NN-est-at-xmin-step3-goal}
        \begin{aligned}
            \partial_\nu u\le&\ h_0\qquad \text{on } U_0\subset \partial \Omega,\\
            \partial_\nu u(x_{\min}) =&\ h_0(x_{\min}).
        \end{aligned}
    \end{equation}
    \vspace{0.5em}\noindent This is the key step and uses concavity of $\signdist^{F_+}$ as well as the fact that $x_{\min}$ is a minimizer.
    
    \vspace{0.5em}\noindent Since $[D^2u]_{O_x} = \pi_+(O_x^T D^2u(x) O_x)$ \eqref{eq:Subeq-Reg-C2-NN-tangential-hessian-def-new-new}, and $\signdist^{F_+}\circ \pi_+(O_x^T D^2u(x) O_x)$ is independent of the choice of $O_x \in G$ such that $O_xe_n = \nu(x)$ \eqref{eq:Subeq-Reg-Trud-point-xmin-indep-of-Ox},
    \[\signdist^{F_+} \left( [D^2u]_{O_{x_{\min}}}\right) \le \signdist^{F_+} \left( [D^2u]_{O_x}\right), \qquad \forall\ x\in U_0,\]
    with equality at $x=x_{\min}$ \eqref{eq:Subeq-Reg-Trudingers-point}.
    Coupled with \eqref{eq:Subeq-Reg-C2-NN-fcal-tilde-concave-subdifferential},
    \begin{align}
        0 \le&\ \signdist^{F_+} \left( [D^2u]_{O_x}\right) -  \signdist^{F_+} \left( [D^2u]_{O_{x_{\min}}}\right) \nonumber \\
        \le&\ \left\langle B_1,\  [D^2u]_{O_x} -  [D^2u]_{O_{x_{\min}}} \right\rangle, \qquad \forall\ x\in U_0 \subset \partial \Omega, \label{eq:Subeq-Reg-C2-NN-est-at-xmin-concavity-eq1}
    \end{align}
    with \textit{both} inequalities achieving equality at $x=x_{\min}$. It is crucial here that both inequalities achieve equality at $x_{\min}$. This eliminates the non-linear term $\signdist^{F_+} \left( [D^2u]_{O_x}\right)$ from \eqref{eq:Subeq-Reg-C2-NN-est-at-xmin-concavity-eq1} (see the discussion following \eqref{eq:Subeq-Reg-C2-NN-est-at-xmin-final-step3-eq2}).

    Since $u = \varphi$ on $\partial \Omega$, and $[D^2u]_{O_x}$ \eqref{eq:Subeq-Reg-C2-NN-tangential-hessian-def-new-new} is the tangential Hessian of $u$, by \eqref{eq:Subeq-Reg-C2-TT-key-lemma-conclusion},
    \[[D^2u]_{O_x} + \partial_{\nu} u(x) \, [\II_{\partial \Omega}]_{O_x} = [D^2\varphi]_{O_x} + \partial_{\nu} \varphi(x) \, [\II_{\partial \Omega}]_{O_x}, \qquad \forall\ x\in U_0,\]
    using the notation \eqref{eq:Subeq-Reg-C2-NN-tangential-hessian-notation-for-arb-A}, where $\partial_\nu u$ is the directional derivative along $\nu$ \eqref{eq:Subeq-Reg-inward-normal-to-domain-bdry-notation-2}. In particular,
    \[\begin{aligned}
        [D^2u]_{O_{x_{\min}}} =&\ [D^2\varphi|_T](x_{\min}) + \partial_{\nu} \varphi(x_{\min}) \, [\II_{\partial \Omega}]_{O_{x_{\min}}} - \partial_{\nu} u(x_{\min}) \, [\II_{\partial \Omega}]_{O_{x_{\min}}},\\
        [D^2u]_{O_x} =&\ [D^2\varphi]_{O_x} + \partial_{\nu} \varphi(x) \, [\II_{\partial \Omega}]_{O_x} - \partial_{\nu} u(x) \, [\II_{\partial \Omega}]_{O_x}, \qquad \forall\ x\in U_0.\\
    \end{aligned}\]
    The difference of this pair of equations allows us to eliminate the second order term $[D^2u]_{O_x}$ from \eqref{eq:Subeq-Reg-C2-NN-est-at-xmin-concavity-eq1} leaving us with an estimate on $\partial_\nu u$ over the boundary $\partial \Omega$,
    \begin{align}
        0\le&\ \left\langle B_1,\  [D^2\varphi|_{T}](x) -  [D^2\varphi|_{T}](x_{\min}) \right\rangle - \partial_\nu \varphi(x_{\min})\, \langle B_1,\ [\II_{\partial \Omega}(x_{\min}]) \rangle \nonumber \\
        &\ + \partial_\nu u(x_{\min})\, \langle B_1,\ [\II_{\partial \Omega}]_{O_{x_{\min}}} \rangle +  \partial_\nu \varphi(x)\, \langle B_1,\ [\II_{\partial \Omega}]_{O_x} \rangle\nonumber \\
        &\ - \partial_\nu u(x)\, \langle B_1,\ [\II_{\partial \Omega}]_{O_x} \rangle, \qquad \forall\ x\in U_0, \label{eq:Subeq-Reg-C2-NN-est-at-xmin-final-step3-eq2}
    \end{align}
    with equality at $x=x_{\min}$ (because both inequalities in \eqref{eq:Subeq-Reg-C2-NN-est-at-xmin-concavity-eq1} achieve equality at $x_{\min}$)
    
    The special choice of the point $x_{\min}$ \eqref{eq:Subeq-Reg-Trudingers-point} eliminates the bad term $ \signdist^{F_+} \left( [D^2u]_{O_x}\right)$ from \eqref{eq:Subeq-Reg-C2-NN-est-at-xmin-concavity-eq1}. So the ``worst'' term remaining in \eqref{eq:Subeq-Reg-C2-NN-est-at-xmin-final-step3-eq2} is the last one, since it is the only term that is neither constant nor a function of $\varphi$ and its derivatives. So \eqref{eq:Subeq-Reg-C2-NN-est-at-xmin-final-step3-eq2} can control $\partial_\nu u$ as long as the coefficient of $\partial_\nu u$ is bounded from below. But Step 4, specifically \eqref{eq:Subeq-Reg-C2-NN-est-at-xmin-step3.5-goal}, provides this lower bound on the coefficient $\langle B_1,\ [\II_{\partial \Omega}]_{O_x}$.
    
    Rearranging \eqref{eq:Subeq-Reg-C2-NN-est-at-xmin-final-step3-eq2} and applying \eqref{eq:Subeq-Reg-C2-NN-est-at-xmin-step3.5-goal},
    \begin{align}
        \partial_\nu u(x) \le&\ \frac{1}{\langle B_1,\ [\II_{\partial \Omega}]_{O_x} \rangle} \Big[ \left\langle B_1,\ [ D^2\varphi|_{T}](x) -  [D^2\varphi|_{T}](x_{\min}) \right\rangle + \partial_\nu \varphi(x)\, \langle B_1,\ [\II_{\partial \Omega}]_{O_x} \rangle   \nonumber \\
        &\   - \partial_\nu \varphi(x_{\min})\, \langle B_1,\ [\II_{\partial \Omega}]_{O_{x_{\min}}} \rangle + \partial_\nu u(x_{\min})\, \langle B_1,\ [\II_{\partial \Omega}]_{O_{x_{\min}}} \rangle \Big]   \nonumber \\
        \le&\ \frac{1}{\langle B_1,\ \II_{\partial \Omega}(x) \rangle} \Big[ \left\langle B_1,\  [D^2\varphi|_{T}](x) -  [D^2\varphi|_{T}](x_{\min}) \right\rangle + \partial_\nu \varphi(x)\, \langle B_1,\ [\II_{\partial \Omega}]_{O_x} \rangle    \nonumber   \\
        &\   - \partial_\nu \varphi(x_{\min})\, \langle B_1,\ [\II_{\partial \Omega}]_{O_{x_{\min}}} \rangle + \partial_\nu u(x_{\min})\, \langle B_1,\ [\II_{\partial \Omega}]_{O_{x_{\min}}} \rangle \Big]   \nonumber \\
        =:&\ h_0(x), \qquad \forall\ x\in U_0 \subset \partial \Omega, \label{eq:Subeq-Reg-C2-NN-est-at-xmin-final-step3-eq5}
    \end{align}
    with equality at $x=x_{\min}$,
    \[\partial_\nu u(x_{\min}) = h_0(x_{\min}).\]
    Furthermore, by Proposition~\ref{prop:Subeq-Reg-unif-grad-est}, there exits $C_3' = C_3'(n,\Omega, F, \|\varphi\|_{C^3}) >0$ such that
    \begin{equation}
        \label{eq:Subeq-Reg-C2-NN-est-at-xmin-final-step3-eq4}
        \|u\|_{C^1} \le C_3'.
    \end{equation}
    Therefore, since $x_{\min}$ is a constant in \eqref{eq:Subeq-Reg-C2-NN-est-at-xmin-final-step3-eq5},
    \begin{equation}
        \label{eq:Subeq-Reg-C2-NN-est-at-xmin-final-step3-eq7}
        \|h_0\|_{C^2(U_0)}\le C_3
    \end{equation}
    for a constant $C_3$ depending on $C_3'$ \eqref{eq:Subeq-Reg-C2-NN-est-at-xmin-final-step3-eq4}, $\|\varphi\|_{C^4}$, $C_1'$ \eqref{eq:Subeq-Reg-C2-NN-est-at-xmin-only-tan-part}, $C_1''$ \eqref{eq:Subeq-Reg-C2-NN-bound-on-subdiff-B1}, and $C_2$ \eqref{eq:Subeq-Reg-C2-NN-est-at-xmin-step3.5-goal}. Furthermore, note that since $h_0$ involves the second fundamental form \eqref{eq:Subeq-Reg-second-fundamental-form}, we require $\partial \Omega$ to have $C^4$ boundary to get a $C^2$ bound on $h_0$.

    

    \vspace{0.5em}\noindent\textit{Step 6.} \textit{There exists a smooth function $h\in C^2(\partial{\Omega})$ whose $C^2$ norm depends only on $n$, $\Omega$, $F$, and $\|\varphi\|_{C^4}$ such that if $e:= \nu(x_{\min})$, then} 
    \begin{equation}
        \label{eq:Subeq-Reg-C2-NN-est-at-xmin-step4-goal}
        \begin{aligned}
            \partial_e u\le&\ h\qquad \text{on } \partial \Omega,\\
            \partial_e u(x_{\min}) =&\ h(x_{\min}).
        \end{aligned}
    \end{equation}

    \vspace{0.5em}\noindent From Step 5, recall $h_0\in C^2(U_0)$ satisfying \eqref{eq:Subeq-Reg-C2-NN-est-at-xmin-step3-goal}. Since $e = \nu(x_{\min}) $, in particular,
    \[\partial_e u(x_{\min}) = \partial_\nu u(x_{\min}) = h_0(x_{\min}).\]
    The strategy is to perturb this equation for $x$ near $x_{\min}$, to  replace the vector field $\nu$ with the constant vector field $e$ in \eqref{eq:Subeq-Reg-C2-NN-est-at-xmin-step3-goal}. By \eqref{eq:Subeq-Reg-C2-TN-step1-goal} applied to the constant vector field $\tau = e$,
    \begin{equation}
        \label{eq:Subeq-Reg-C2-NN-est-at-xmin-final-step3-eq10}
        \partial_e u - \partial_e\varphi = \langle e, \nu \rangle (\partial_\nu u - \partial_\nu \varphi).
    \end{equation}
    Furthermore, since $e = \nu(x_{\min})$, $\langle e, \nu(x_{\min})\rangle = 1$. So there is a neighborhood $U_1$ of $x_{\min}$,
    \[U_1 \subset U_0 \subset \partial \Omega,\]
    depending only on $U_0$ and $\Omega$, such that
    \begin{equation}
        \label{eq:Subeq-Reg-C2-NN-est-at-xmin-final-step3-eq9}
        \langle e, \nu(x)\rangle \ge 1/2, \qquad \forall\ x\in U_1.
    \end{equation}
    Since $\langle e, \nu(x)\rangle$ is non-negative on $U_1$, substituting \eqref{eq:Subeq-Reg-C2-NN-est-at-xmin-final-step3-eq10} into \eqref{eq:Subeq-Reg-C2-NN-est-at-xmin-final-step3-eq5},
    \begin{align}
        \partial_e u(x) \le&\ \partial_e\varphi (x) + \langle e, \nu(x) \rangle \partial_\nu u(x) - \langle e, \nu(x) \rangle\partial_\nu \varphi(x) \nonumber\\
        \le&\ \partial_e\varphi (x) + \langle e, \nu(x) \rangle h_0(x) - \langle e, \nu(x) \rangle   \partial_\nu \varphi(x)  \nonumber\\
        =:&\ h_1(x), \qquad \forall\ x\in U_1, \label{eq:Subeq-Reg-C2-NN-est-at-xmin-final-step3-eq5-new}
    \end{align}
    and
    \[\partial_e u(x_{\min}) =\partial_\nu u(x_{\min}) = h_0(x_{\min}) =  h_1(x_{\min}) .\]
    By \eqref{eq:Subeq-Reg-C2-NN-est-at-xmin-final-step3-eq5-new}, $h_1$ depends only on $h_0$, $\Omega$ and $\varphi$, and by \eqref{eq:Subeq-Reg-C2-NN-est-at-xmin-final-step3-eq7} there exists $C_4= C_4(C_3, \|\varphi\|_{C^3}, \Omega)$ such that
    \begin{equation}
        \label{eq:Subeq-Reg-C2-NN-est-at-xmin-final-step3-eq7-new}
        \|h_1\|_{C^2(U_1)}\le C_4
    \end{equation}

    By \eqref{eq:Subeq-Reg-C2-NN-est-at-xmin-final-step3-eq5-new}, \eqref{eq:Subeq-Reg-C2-NN-est-at-xmin-step4-goal} holds, but only in a small neighborhood $U_1$ of $x_{\min}$. We extend \eqref{eq:Subeq-Reg-C2-NN-est-at-xmin-final-step3-eq5-new} beyond $U_1$ using a cutoff function.
    
    Let $\eta \in C^\infty(\partial \Omega)$ be a smooth non-negative cut-off function that is 1 on a neighborhood $U_2 \subset U_1$ of $x_{\min}$, with support in $U_1$, 
    \begin{equation}
        \label{eq:Subeq-Reg-C2-NN-est-at-xmin-final-step3-eq6}
        \eta \ge 0, \qquad \eta \equiv 1 \quad \text{on } U_2 \subset U_1, \qquad \text{and}\qquad \eta \equiv 0 \quad \text{on } \partial \Omega \setminus U_1.
    \end{equation}
    Note that $\eta$ depends only on $\Omega$ and $U_1$, which in turn depends only $U_0$ from Step 4. And by Step 4, $U_0$ depends only on $C_1'$ \eqref{eq:Subeq-Reg-C2-NN-bound-on-subdiff-B1} and $\Omega$. Hence,
    \begin{equation}
        \label{eq:Subeq-Reg-C2-NN-at-xmin-new-step5-eq-eta-bound}
        \|\eta\|_{C^2(\partial \Omega)} \le C_4' = C_4'(C_1', \Omega)
    \end{equation}
    Consider the function $h:\partial \Omega \rightarrow \RR$,
    \[  h(x):= h_1(x) \eta(x) + \left( 1-\eta(x)\right)\|u\|_{C^1(\overline{\Omega})} ,\qquad \forall\ x\in \partial \Omega.\]
    Note that $\|h\|_{C^2}$ is estimated in terms of $\|\eta\|_{C^2}$ \eqref{eq:Subeq-Reg-C2-NN-at-xmin-new-step5-eq-eta-bound}, $\|h_1\|_{C^2(U_1)}$ \eqref{eq:Subeq-Reg-C2-NN-est-at-xmin-final-step3-eq7-new}, and $\|u\|_{C^1}$ \eqref{eq:Subeq-Reg-C2-NN-est-at-xmin-final-step3-eq4}. And, by \eqref{eq:Subeq-Reg-C2-NN-est-at-xmin-final-step3-eq5-new} and \eqref{eq:Subeq-Reg-C2-NN-est-at-xmin-final-step3-eq6},
    \[\begin{aligned}
        \partial_e u(x) \le&\ h(x), \qquad \forall\ x\in \partial \Omega,\\
        \partial_e u(x_{\min}) =&\  h_1 (x_{\min}) = h(x_{\min}). 
    \end{aligned}\]
    Thus, $h$ satisfies \eqref{eq:Subeq-Reg-C2-NN-est-at-xmin-step4-goal}.

    \vspace{0.5em}\noindent\textit{Step 7.} \textit{\eqref{eq:Subeq-Reg-C2-NN-upper-bound-at-ymin-conclusion} follows from \eqref{eq:Subeq-Reg-C2-NN-est-at-xmin-step4-goal}, Corollary~\ref{corollary:Subeq-Reg-guan-corollary}, and Claim~\ref{claim:Subeq-Reg-normal-gradient-claim}:}
    
    \vspace{0.5em}\noindent Since $\partial \Omega$ is strictly $\vec{F}$-convex, by Proposition~\ref{prop:Subeq-Reg-existence-of-delta-subsol-and-supersol}, \eqref{eq:Subeq-Reg-const-coeff-subequation-PDE} admits an admissible $\delta$-subsolution $\underline{u} \in C^2(\overline{\Omega})$ satisfying \eqref{eq:Subeq-Reg-strongly-admissible-subsolution}. Hence, by  Corollary~\ref{corollary:Subeq-Reg-guan-corollary} applied to the function $h$ from Step 6 at $x_{\min} \in \partial \Omega$, there exists $w\in C^2(\overline{\Omega})$ and $C_5 = C_5(n,\Omega, F,\|\underline{u}\|_{C^2}, \delta, \|h\|_{C^2})>0$ satisfying \eqref{eq:Subeq-Reg-guan-corollary-assump-conclusion},
    \begin{equation}
        \label{eq:Subeq-Reg-C2-NN-est-at-xmin-new-step6-eq3}
        \begin{aligned}
        h\le&\ w \qquad \text{on } \partial \Omega,\\
        h (x_{\min})=&\ w  (x_{\min}),\\
        \Delta_{D\signdist^F(D^2u)} w \le&\ 0,\\
        \|w\|_{C^2(\overline{\Omega})} \le&\ C_5.
    \end{aligned}
    \end{equation}
    Combined with \eqref{eq:Subeq-Reg-C2-NN-est-at-xmin-step4-goal},
    \begin{equation}
        \label{eq:Subeq-Reg-C2-NN-est-at-xmin-Step4-eq1}
        \begin{aligned}
        \partial_e u\le h\le&\ w \qquad \text{on } \partial \Omega,\\
        \partial_e u(x_{\min})=  h (x_{\min})=&\ w  (x_{\min}).
    \end{aligned}
    \end{equation}
    But differentiating \eqref{eq:Subeq-Reg-const-coeff-subequation-PDE}, by \eqref{eq:Subeq-Reg-linearization-of-dij-elliptic-pde} and \eqref{eq:Subeq-Reg-C2-NN-est-at-xmin-new-step6-eq3},
    \[\partial_e(\signdist^{F}(D^2u)) = \Delta_{D \signdist^{F}(D^2u)} \partial_e u = 0 \ge \Delta_{D \signdist^{F}(D^2u)} w.\]
    By the maximum principle for linear elliptic PDE \cite[Theorem~3.1]{gilbarg-trudinger-MR1814364},
    \[\begin{aligned}
        \partial_e u \le&\ w \qquad \text{on } \overline{\Omega},\\
        \partial_n u(x_{\min})  =&\ w (x_{\min}).
    \end{aligned}\]
    Since $\nu(x_{\min}) = e$, by Claim~\ref{claim:Subeq-Reg-normal-gradient-claim} applied at $x_{\min}$ and \eqref{eq:Subeq-Reg-C2-NN-est-at-xmin-new-step6-eq3},
    \[\partial_e \partial_e u(x_{\min}) \le \partial_e w(x_{\min}) \le \|w\|_{C^1} \le C_5.\]
    In other words, \eqref{eq:Subeq-Reg-C2-NN-upper-bound-at-ymin-conclusion} holds.
\end{proof}

\subsection{Proof of Theorem~\ref{thm:Subeq-Reg-a-priori-estimates}}
\label{section:Subeq-Reg-proof-of-main-a-priori-theorem}

    By \eqref{eq:Subeq-Reg-tangential-hessian-estimate} and \eqref{eq:Subeq-Reg-C2-TN-conclusion}, there is a constant $C_1 = C_1(n, \Omega, F, \|\varphi\|_{C^3})$ bounding the tangential and mixed Hessian components,
    \begin{equation}
        \label{eq:Subeq-Reg-C2-NN-tan-mix-hess-estimate}
        |\langle D^2u(x) V, W \rangle| \le C_1,\qquad \forall\ V,W\in S^{n-1},\ V\perp \nu(x),\ x\in \partial \Omega.
    \end{equation}
    
    \vspace{0.5em}\noindent\textit{Step 1.} \textit{Lower bound:} By \eqref{eq:Subeq-Reg-C2-NN-tan-mix-hess-estimate} and Lemma~\ref{lemma:Subeq-Reg-C2-NN-lower-bound-lin-alg},
    \begin{equation}
        \label{eq:Subeq-Reg-C2-NN-main-lower-bound}
        \langle D^2u(x) \nu(x), \nu(x)\rangle  \ge -C_2 ,\qquad \forall\  x\in \partial \Omega,
    \end{equation}
    for some $C_2 = C_2(n,F, C_1) > 0$.

    \vspace{0.5em}\noindent\textit{Step 2:} \textit{Hessian estimate at $x_{\min}$:} By Proposition~\ref{prop:Subeq-Reg-C2-NN-upper-bound-at-ymin},
    \[\langle D^2u(x_{\min}) \nu(x_{\min}), \nu(x_{\min})\rangle  \le C_3,\]
    and combined with \eqref{eq:Subeq-Reg-C2-NN-tan-mix-hess-estimate}--\eqref{eq:Subeq-Reg-C2-NN-main-lower-bound},
    \begin{equation}
    \label{eq:Subeq-Reg-C2-NN-upper-bound-at-ymin}
        |D^2u(x_{\min})|\le C_4,
    \end{equation}
    where $C_3,C_4>0$ depend only on $n, \Omega, F, \|\varphi\|_{C^4}$.
    
    \vspace{0.5em}\noindent\textit{Step 3:} \textit{Hessian estimate on $\partial \Omega$:} By Corollary~\ref{corollary:Subeq-Reg-C2-NN-trudinger-prop}~(iii),
    \[ |D^2u(x)|\le C_5, \qquad \forall\ x\in \partial \Omega,\]
    where $C_5 = C_5(C_1, n, \Omega, F, \|\varphi \|_{C^4}, C_4)$.

\section{Higher regularity and the continuity method}
\label{sec:Subeq-Reg-higher-regularity}

This section records how standard higher order regularity theory can be adapted to subequations \eqref{eq:Subeq-Reg-const-coeff-subequation-PDE}. The idea is to apply higher regularity to the signed-distance $\signdist^F$ and the PDE \eqref{eq:Subeq-Reg-signed-dist-fcal-PDE}. Traditionally, a priori $C^2$ estimates on a solution $u$ to \eqref{eq:Subeq-Reg-const-coeff-subequation-PDE} (Theorem~\ref{thm:Subeq-Reg-a-priori-estimates}) yield $C^{2,\alpha}$ estimates via the Evans--Krylov theorem. The Evans--Krylov theorem is usually stated for globally smooth and concave operators, while $\signdist^F$ is smooth only in a neighborhood of $\partial F$. So, we state a version that applies to $\signdist^F$.

The linearization of $\Delta_{D\signdist^{F}(D^2u)}$ \eqref{eq:Subeq-Reg-linearization-of-dij-elliptic-pde} is uniformly elliptic (Lemma~\ref{lemma:Subeq-Reg-sign-dist-is-unif-elliptic}) by Theorem~\ref{thm:Subeq-Reg-a-priori-estimates}. Hence, the Schauder estimates can be used to upgrade $C^{2,\alpha}$ regularity of $u$ to higher order regularity.

\subsection{Evans--Krylov theorem}
\label{subsection:Subeq-Reg-evans-krylov-theorem}

A solution $u$ to \eqref{eq:Subeq-Reg-const-coeff-subequation-PDE} with an a priori bounded Hessian also admits a H\"older estimate on its Hessian by an argument due to Evans~\cite{Evans-C2alpha-est-1982-MR649348} and Krylov~\cite{Krylov-C2alpha-domain-bdry-MR688919}. The following version can be found in Han~\cite[Theorem~5.5.1]{Qing-Han-MR3468839}:

\begin{theorem}[Evans--Krylov Theorem]
    \label{thm:Subeq-Reg-Evans-Krylov-theorem}
    For a bounded domain $\Omega \subset \RR^n$ with $C^3$ boundary, continuous function $\fcal :\Sym^2(\RR^n) \rightarrow \RR$, and a solution $u \in C^4(\overline{\Omega})$ to
    \[\begin{aligned}
        \fcal(D^2u) = f \qquad &\text{on } \Omega,\\
        u = \varphi  \qquad &\text{on } \partial \Omega,
    \end{aligned}\]
    for $f\in C^2(\overline{\Omega})$ and $\varphi \in C^3(\partial{\Omega})$, if $\fcal$ is
    \begin{enumerate}[label={(\roman*)}]
        \item $C^3$ in a neighborhood of $D^2u(\overline{\Omega})$,
        \item uniformly elliptic relative to $u$, i.e., there exists a constant $\Lambda>0$ such that
        \begin{equation}
            \label{eq:Subeq-Reg-Evans-Krylov-assump-unif-elliptic}
            \Lambda^{-1} |\xi|^2 \le \left\langle D\fcal(D^2u(x)) \xi, \, \xi \right\rangle \le \Lambda |\xi|^2, \qquad \forall\ x\in \overline{\Omega},\ \xi \in \RR^n, 
        \end{equation}
        \item concave relative to $u$, i.e.,
        \begin{equation}
            \label{eq:Subeq-Reg-Evans-Krylov-assump-concave}
            D^2\fcal(D^2u(x))( A, A) \le 0, \qquad \forall\ x\in \overline{\Omega},\ A\in \Sym^2(\RR^n),
        \end{equation}
        where the Hessian $D^2\fcal$ is regarded as a symmetric 2-tensor on $\Sym^2(\RR^n)$,
    \end{enumerate}
    then there exists constants $\alpha = \alpha(\Lambda, \Omega) \in (0,1)$ and $C = C(\Lambda, \Omega) >0$ such that
    \[\|D^2u\|_{C^{\alpha}(\overline{\Omega})} \le C \left( \|u\|_{C^{2}(\overline{\Omega})} + \|f\|_{C^{2}(\overline{\Omega})} + \|\varphi\|_{C^{3}(\partial {\Omega})} \right). \]
\end{theorem}


The signed-distance $\signdist^{F}$ \eqref{eq:Subeq-Reg-signed-dist-def} is smooth and concave on a neighborhood of $\partial F$ (Lemma~\ref{lemma:Subeq-Reg-signed-distance-lemma}~(iii) and (v)), hence $\fcal:= \signdist^{F}$ is concave relative to $u$ \eqref{eq:Subeq-Reg-Evans-Krylov-assump-concave}. Furthermore, the $C^2$ estimate from Theorem~\ref{thm:Subeq-Reg-a-priori-estimates} ensures that $\signdist^{F}$ is also uniformly elliptic relative to $u$ \eqref{eq:Subeq-Reg-Evans-Krylov-assump-unif-elliptic}:

\begin{lemma}[Uniform ellipticity]
\label{lemma:Subeq-Reg-sign-dist-is-unif-elliptic}
    Under the assumptions of Theorem~\ref{thm:Subeq-Reg-a-priori-estimates}, there exists a constant $\Lambda = \Lambda(n,F,\Omega, \|\varphi\|_{C^4})>0$, such that the signed-distance $\signdist^{F}$ \eqref{eq:Subeq-Reg-signed-dist-def} is uniformly elliptic relative to any solution $u$ to \eqref{eq:Subeq-Reg-const-coeff-subequation-PDE} (cf. \eqref{eq:Subeq-Reg-Evans-Krylov-assump-unif-elliptic}),
    \begin{equation}
        \label{eq:Subeq-Reg-higher-reg-unif-ellipticity-of-sign-dist}
        \Lambda^{-1} |\xi|^2 \le \left\langle D \signdist^{F} (D^2u(x)) \xi, \, \xi \right\rangle \le \Lambda |\xi|^2, \qquad \forall\ x\in \overline{\Omega},\ \xi \in \RR^n, 
    \end{equation}
\end{lemma}
\begin{proof}
    By Theorem~\ref{thm:Subeq-Reg-a-priori-estimates}, there exists a constant $C_1 = C_1(n, F, \Omega, \|\varphi\|_{C^4})$ such that
    \begin{equation}
        \label{eq:Subeq-Reg-Evans-Krylov-cor-eq1}
        \|u\|_{C^2(\overline{\Omega})} <C_1.
    \end{equation}
    By Lemma~\ref{lemma:Subeq-Reg-signed-distance-lemma}~(iii), $\signdist^{F}$ is smooth and concave in a neighborhood $U$ of $\partial F$. Shrinking $U$ if necessary, assume that $\signdist$ is smooth on the closure $\overline{U}$. By Lemma~\ref{lemma:Subeq-Reg-subequation-is-elliptic}~(iii) and \eqref{eq:Subeq-Reg-sign-dist-has-pos-grad},
    \[D\signdist^{F}(A)>0, \qquad \forall\ A\in \overline{U}.\]
    Let $\Lambda_1$ denote the maximum spectral norm \eqref{eq:Subeq-Reg-interior-hess-est-spectral-norm} of $D\signdist^{F}$ over the compact set of matrices in $\overline{U}$ with norm bounded by $C_1$ \eqref{eq:Subeq-Reg-Evans-Krylov-cor-eq1},
    \[\Lambda_1:= \max \left\{ \frac{\langle  D\signdist^{F}(A)\xi, \xi \rangle }{\langle \xi, \xi \rangle}: A\in \overline{U} ,\quad |A|\le C_1, \quad \xi \in \RR^n \setminus \{0\} \right\}, \]
    and let $\Lambda_2^{-1}$ be the smallest eigenvalue of $D\signdist^{F}$ over the same compact set,
    \[\Lambda_2 := \max\left\{ \frac{\langle \xi, \xi \rangle }{\langle  D\signdist^{F}(A)\xi, \xi \rangle}: A\in \overline{U} ,\quad |A|\le C_1, \quad \xi \in \RR^n \setminus \{0\} \right\} .\]
    By \eqref{eq:Subeq-Reg-const-coeff-subequation-PDE} and \eqref{eq:Subeq-Reg-Evans-Krylov-cor-eq1}, $D^2u(x)$ is a candidate in the max defining both $\Lambda_1$ and $\Lambda_2$. Hence,
    \[\Lambda_2^{-1}  |\xi|^2 \le \langle  D\signdist^{F}(D^2u(x))\xi, \xi \rangle \le \Lambda_1 |\xi|^2. \]
    Thus \eqref{eq:Subeq-Reg-higher-reg-unif-ellipticity-of-sign-dist} holds for $\Lambda = \max\{\Lambda_1,\Lambda_2\}$.
\end{proof}

\begin{corollary}
    \label{corollary:Subeq-Reg-Evans-Krylov-theorem}
    Under the assumptions of Corollary~\ref{corollary:Subeq-Reg-existence-thm}, there exists $\alpha \in (0,1)$ and $C>0$ depending only on $n$, $\Omega$, $F$, and $\|\varphi\|_{C^4(\partial \Omega)}$ such that any solution $u \in C^4(\overline{\Omega})$ to \eqref{eq:Subeq-Reg-const-coeff-subequation-PDE} satisfies
    \begin{equation}
        \label{eq:Subeq-Reg-Evans-Krylov-conclusion}
        \|D^2u\|_{C^{\alpha}(\overline{\Omega})} \le C .
    \end{equation}
\end{corollary}

\subsection{Schauder estimates}
\label{subsection:Subeq-Reg-schauder-estimates}

The a priori $C^{2,\alpha}$ estimates of Corollary~\ref{corollary:Subeq-Reg-Evans-Krylov-theorem} are improved to higher order regularity via the Schauder estimates for linear elliptic PDE \cite[Theorems~6.2, 6.6]{gilbarg-trudinger-MR1814364}. If $u$ is a solution to \eqref{eq:Subeq-Reg-const-coeff-subequation-PDE}, differentiating \eqref{eq:Subeq-Reg-signed-dist-fcal-PDE}, its derivative $Du$ satisfies the linear elliptic PDE \eqref{eq:Subeq-Reg-linearization-of-dij-elliptic-pde}
\begin{equation}
    \label{eq:Subeq-Reg-higher-schauder-derivative-linearization-soln}
    \Delta_{D\signdist^{F}(D^2u)} Du = 0 \qquad \text{on } \Omega.
\end{equation}
This equation is uniformly elliptic by Lemma~\ref{lemma:Subeq-Reg-sign-dist-is-unif-elliptic}. The argument is indicated in Caffarelli--Cabr\'e~\cite[Chapter 9, Remark~1]{caffarelli-cabre-MR1351007}, and we record a version for the signed-distance $\signdist^F$ \eqref{eq:Subeq-Reg-signed-dist-def}: 

\begin{proposition}
    \label{prop:Subeq-Reg-higher-reg-boundary-schauder}
    Under the assumptions of Corollary~\ref{corollary:Subeq-Reg-existence-thm}, for $k\ge 2$, $\alpha \in (0,1)$ from Corollary~\ref{corollary:Subeq-Reg-Evans-Krylov-theorem}, there exists $C = C(n,k,\alpha, F, \Omega, \|\varphi\|_{C^{k,\alpha}}, \|\varphi\|_{C^4})>0$ such that any solution $u\in C^4 \cap C^{k,\alpha}(\overline{\Omega}) $ satisfies
    \[\|u\|_{C^{k,\alpha}(\overline{\Omega})} \le C.\]
\end{proposition}

\begin{proof}[Sketch of proof]
    The interior Schauder estimate~\cite[Theorems~6.2]{gilbarg-trudinger-MR1814364} applied to \eqref{eq:Subeq-Reg-higher-schauder-derivative-linearization-soln} can be bootstrapped to yield interior higher order regularity of $u$.
    
    Since $Du|_{\partial \Omega}$ is not controlled by the boundary condition of \eqref{eq:Subeq-Reg-const-coeff-subequation-PDE}, the boundary Schauder estimates~\cite[Theorems~6.6]{gilbarg-trudinger-MR1814364} do not directly yield higher order estimates up to the boundary. However boundary Schauder estimates can be applied to the tangential components of $Du|_{\partial \Omega}$, which in particular yield higher regularity for the tangential and mixed components of the Hessian $D^2u(x_0)$, at any boundary point $x_0$. This can be seen for example by flattening the boundary of the domain $\partial \Omega$ in a neighborhood of $x_0$ (See \cite[Chapter 9, Remark 1]{caffarelli-cabre-MR1351007}).
    
    Finally, the normal component of the Hessian $D^2u(x_0)$ can be estimated by applying an implicit function theorem argument to the PDE \eqref{eq:Subeq-Reg-signed-dist-fcal-PDE} using uniform ellipticity \eqref{eq:Subeq-Reg-higher-reg-unif-ellipticity-of-sign-dist} of $\signdist^{F}$.
\end{proof}

\begin{allproofs}
\begin{proof}
    Argue by induction on $k\ge 2$. For the base case $k=2$, by Corollary~\ref{corollary:Subeq-Reg-Evans-Krylov-theorem}, there exists $C_1>0$ depending only on $n$, $\Omega$, $F$ and $\|\varphi\|_{C^4}$ such that
    \begin{equation}
        \label{eq:Subeq-Reg-higher-reg-bdry-induction-base-case}
        \|u\|_{C^{2,\alpha}(\overline{\Omega})} \le C_1 .
    \end{equation}
    
    Assume by induction that there exists $C_2>0$ depending on $k-1$, $\alpha$, $n$, $F$, $\Omega$, and $\|\varphi\|_{C^4(\partial \Omega)}$ and $\|\varphi\|_{C^{k-1,\alpha}(\partial \Omega)}$ such that
    \begin{equation}
        \label{eq:Subeq-Reg-higher-reg-boundary-induction-hypothesis}
        \|u\|_{C^{k-1,\alpha}(\overline{\Omega})} \le C_2 .
    \end{equation}
    Differentiating \eqref{eq:Subeq-Reg-signed-dist-fcal-PDE},
    \begin{equation}
        \label{eq:Subeq-Reg-higher-reg-boundary-schauder-eq1}
        \Delta_{D\signdist^{F} (D^2u)} Du = 0 \qquad \text{on } \Omega,
    \end{equation}
    where $\Delta_{D\signdist^{F}(D^2u)}$ \eqref{eq:Subeq-Reg-linearization-of-dij-elliptic-pde} acts component wise on the vector $Du$. By Lemma~\ref{lemma:Subeq-Reg-sign-dist-is-unif-elliptic}, $\Delta_{D\signdist^{F}(D^2u)}$ is uniformly  elliptic with parameter $\Lambda$ \eqref{eq:Subeq-Reg-higher-reg-unif-ellipticity-of-sign-dist},i.e.,
    \[\Lambda^{-1} |\xi|^2 \le \left\langle D\signdist^{F}(D^2u) \xi,\ \xi  \right\rangle \le \Lambda |\xi|^2, \qquad \forall\ \xi \in \RR^n.\]
    Furthermore, $\Delta_{D\signdist^{F}(D^2u)}$ has $C^{k-3,\alpha}$ coefficients by \eqref{eq:Subeq-Reg-higher-reg-boundary-induction-hypothesis}. If we had an a priori $C^{k-1,\alpha}$ estimate on $Du|_{\partial \Omega}$, then a $C^{k-1,\alpha}(\overline{\Omega})$ estimate on $Du$ would follow from the boundary Schauder estimates \cite[Theorem~6.6]{gilbarg-trudinger-MR1814364}. Unfortunately the induction hypothesis \eqref{eq:Subeq-Reg-higher-reg-boundary-induction-hypothesis} only yields a priori $C^{k-2,\alpha}$ estimates on $Du|_{\partial \Omega}$. The strategy is the estimate the tangential part of $Du|_{\partial \Omega}$ using the boundary condition $\varphi$ and get an $\alpha$-H\"older estimates on an $(k-1)^{\text{th}}$ order derivative of $u$ along any tangential direction. Then the normal component is estimated using the equation \eqref{eq:Subeq-Reg-signed-dist-fcal-PDE}, uniform ellipticity \eqref{eq:Subeq-Reg-higher-reg-unif-ellipticity-of-sign-dist} and the implicit function theorem.

    \vspace{0.5em}\noindent\textit{Step 1.} \textit{For a tangential vector field $\tau \in C^{k}( \Omega, \RR^n)$,}
    \begin{equation}
        \label{eq:Subeq-Reg-high-reg-bdry-schauder-step1-assump}
        \tau(x) \in T_x\partial \Omega,\qquad \forall\ x\in \partial \Omega,
    \end{equation}
    \textit{there exists a constant $C_3>0$ depending on $C_1$, $\Lambda$, $n$, $F$, $k$, $\alpha$, $\|\varphi\|_{C^{k,\alpha}}$ and $\|\tau\|_{C^{k}}$ such that}
    \begin{equation}
        \label{eq:Subeq-Reg-high-reg-bdry-schauder-step1-goal}
        \left\| \langle Du, \tau\rangle \right\|_{C^{k-1,\alpha}(\overline{\Omega})} \le C_3.
    \end{equation}

    \vspace{0.5em}\noindent Let $\tau = (\tau_1,\dots, \tau_n)$ and use the notation
    \[\partial_\tau u:= \langle Du, \tau\rangle = \sum_{l=1}^n \tau_l \partial_l u.\]
    By \eqref{eq:Subeq-Reg-higher-reg-boundary-schauder-eq1}, it turns out $\partial_\tau u$ satisfies a linear PDE with second order part $\Delta_{D\signdist^{F}(D^2u)}$ \eqref{eq:Subeq-Reg-linearization-of-dij-elliptic-pde}. Observe that
    \[\begin{aligned}
        \Delta_{D\signdist^{F}(D^2u)} \partial_\tau u =&\ \sum_{i,j=1}^n \signdist^{F}_{ij}(D^2u) \partial_i \partial_j \partial_{\tau}  u \\
        =&\ \sum_{i,j,l=1}^n \signdist^{F}_{ij}(D^2u) \partial_i \partial_j\left( \tau_l \partial_{l}  u \right) \\
        =&\ \sum_{i,j,l=1}^n \signdist^{F}_{ij}(D^2u) \left( \tau_l \partial_i \partial_j \partial_l  u + \partial_i\tau_l \partial_j \partial_l  u + \partial_j\tau_l \partial_i \partial_l  u + \partial_i\partial_j \tau_l \partial_l  u \right) \\
        =&\ \sum_{l=1}^n \tau_l\Delta_{D\signdist^{F}(D^2u)} \partial_l  u + \sum_{i,j,l=1}^n \signdist^{F}_{ij}(D^2u) \left( \partial_i\tau_l \partial_j u + \partial_j\tau_l \partial_i u + u \partial_i\partial_j \tau_l \right) \\
        =&\ 0 + \sum_{i,j,l=1}^n \signdist^{F}_{ij}(D^2u) \left( \partial_i\tau_l \partial_j u + \partial_j\tau_l \partial_i u + u \partial_i\partial_j \tau_l \right),
    \end{aligned}\]
    where the last line follows from \eqref{eq:Subeq-Reg-higher-reg-boundary-schauder-eq1}. In other words, $ v= \partial_\tau u$ is a solution to the linear PDE
    \begin{equation}
        \label{eq:Subeq-Reg-high-reg-bdry-schauder-eq2}
        \Delta_{D\signdist^{F}(D^2u)}\, v = \sum_{i,j,l=1}^n \signdist^{F}_{ij}(D^2u) \left( \partial_i\tau_l \partial_j u + \partial_j\tau_l \partial_i u + u \partial_i\partial_j \tau_l \right).
    \end{equation}
    By the induction hypothesis \eqref{eq:Subeq-Reg-higher-reg-boundary-induction-hypothesis} and the assumption $\tau \in C^{k}(\overline{\Omega})$, the coefficients of the linear PDE \eqref{eq:Subeq-Reg-high-reg-bdry-schauder-eq2} involve at most two derivatives of $u$ or of $\tau$. Hence their $C^{k-3,\alpha}(\overline{\Omega})$-norm is bounded uniformly in terms of $C_2$, $F$, and $\|\tau\|_{C^{k}}$. Furthermore, the ellipticity of \eqref{eq:Subeq-Reg-high-reg-bdry-schauder-eq2} is bounded in terms of $\Lambda$ \eqref{eq:Subeq-Reg-higher-reg-unif-ellipticity-of-sign-dist}.

    In addition to solving the equation \eqref{eq:Subeq-Reg-high-reg-bdry-schauder-eq2}, by \eqref{eq:Subeq-Reg-const-coeff-subequation-PDE}, and the assumption that $\tau$ is tangential to $\partial \Omega$ \eqref{eq:Subeq-Reg-high-reg-bdry-schauder-step1-assump}, $v=\partial_\tau u$ satisfies the boundary condition
    \begin{equation}
        \label{eq:Subeq-Reg-high-reg-bdry-schauder-eq3}
        v = \partial_\tau u = \partial_\tau \varphi \qquad \text{on } \partial \Omega.
    \end{equation}
    By \eqref{eq:Subeq-Reg-high-reg-bdry-schauder-eq2}--\eqref{eq:Subeq-Reg-high-reg-bdry-schauder-eq3} and the boundary Schauder estimate \cite[Theorem~6.6]{gilbarg-trudinger-MR1814364}, there is a constant $C_3>0$ depending on a $C^{k-3,\alpha}$ estimate on the coefficients of \eqref{eq:Subeq-Reg-high-reg-bdry-schauder-eq2} (which depend in turn on $C_2$, $F$, and $\|\tau\|_{C^{k}}$), the ellipticity $\Lambda$ \eqref{eq:Subeq-Reg-higher-reg-unif-ellipticity-of-sign-dist}, as well as $n$, $k$ and $\alpha$ such that
    \[\|\partial_\tau u\|_{C^{k-1,\alpha}(\overline{\Omega}))} \le C_3.\]

    \vspace{0.5em}\noindent\textit{Step 2.} \textit{There exists a constant $C_4>0$ depending on $C_2$ \eqref{eq:Subeq-Reg-higher-reg-boundary-induction-hypothesis}, $\Lambda$ \eqref{eq:Subeq-Reg-higher-reg-unif-ellipticity-of-sign-dist}, $n$, $F$, $k$, $\alpha$, $\Omega$ and $\|\varphi\|_{C^{k,\alpha}}$ such that, each boundary point $x_0\in \partial \Omega$ admits a neighborhood $U_{x_0}$ such that}
    \begin{equation}
        \label{eq:Subeq-Reg-high-reg-bdry-schauder-step2-goal}
        \|D^2u\|_{C^{k-2,\alpha}(\overline{\Omega} \cap U_{x_0})} \le C_4.
    \end{equation}

    \vspace{0.5em}\noindent Fix a boundary point $x_0\in \partial \Omega$ and two small bounded open sets $U_1,U_2\subset \RR^n$ such that
    \[x_0 \in U_1\subset \overline{U_1}\subset  U_2.\]
    Shrinking $U_1$ and $U_2$ if necessary, construct $n$ vector fields $\tau_1,\dots, \tau_n: \overline{\Omega} \rightarrow \RR^n$ that are orthonormal on $U_1\cap \overline{\Omega}$, supported in $U_2\cap \overline{\Omega}$, and such that the first $n-1$ are tangent to $T_x\partial \Omega$, i.e.,
    \[\begin{aligned}
        \left\langle \tau_i(x),\tau_j(x) \right\rangle = \delta_{ij},& \qquad \forall\ x\in U_1\cap \overline{\Omega}, \\
        \tau_i(x) = 0,&\qquad \forall\ x\in \overline{\Omega} \setminus U_2,\ i=1,\dots, n, \\
        \tau_i(x) \in T_x\partial \Omega,&\qquad  \forall\ x\in \partial \Omega,\quad i=1,\dots, n-1.
    \end{aligned}\]
    By the first and third condition, the last vector equals the unit inward normal $\nu$ to $\partial \Omega$ \eqref{eq:Subeq-Reg-inward-normal-to-domain-bdry-notation} over $U_1$,
    \[\tau_n(x) =\nu(x), \qquad \forall\ x\in U_1\cap \partial \Omega.\]
    Consider the matrix of the Hessian $D^2u$ in the orthonormal basis $\tau_1,\dots, \tau_n$ over $U_1$,
    \[(u_{\tau_i,\tau_j}) := \left( \langle D^2u\, \tau_i,\ \tau_j \rangle  \right)_{i,j=1}^n = \begin{pmatrix}
        \tau_1^T \\ \vdots \\ \tau_n^T
    \end{pmatrix} D^2u \begin{pmatrix}
        \tau_1 & \dots & \tau_n
    \end{pmatrix}.\]
    By \eqref{eq:Subeq-Reg-signed-dist-fcal-PDE} and $O(n)$-invariance of $\signdist^{F}$ \eqref{eq:Subeq-Reg-sign-dist-is-orthogonal-invariant},
    \[\signdist^{F}(u_{\tau_i\tau_j}(x)) = 0, \qquad \forall\ x\in \overline{\Omega} \cap U_1.\]
    Furthermore, since $\tau_1,\dots, \tau_{n-1}$ satisfy \eqref{eq:Subeq-Reg-high-reg-bdry-schauder-step1-assump}, so by Step 1 \eqref{eq:Subeq-Reg-high-reg-bdry-schauder-step1-goal} there exists a constant $C_3>0$ depending only on $C_1$, $\Lambda$, $n$, $F$, $k$, $\alpha$, $\|\varphi\|_{C^{k,\alpha}}$ and $\|\tau_1\|_{C^{k}}, \dots, \|\tau_{n-1}\|_{C^k}$ such that
    \[\|\partial_{\tau_i}u\|_{C^{k-1,\alpha}} \le C_3, \qquad i=1,\dots, n-1.\]
    In particular, taking another derivative,
    \begin{equation}
        \label{eq:Subeq-Reg-high-reg-bdry-schauder-step-2-implicit-eq1}
        \|u_{\tau_i\tau_j}\|_{C^{k-2,\alpha}} = \|\langle D^2u \, \tau_i,\ \tau_j \rangle \|_{C^{k-2,\alpha}} \le C_3', \qquad i,j=1,\dots, n,\ (i,j) \neq (n,n),
    \end{equation}
    where $C_3'$ depends on $C_3$ and $\tau_1,\dots, \tau_n$. By uniform ellipticity \eqref{eq:Subeq-Reg-higher-reg-unif-ellipticity-of-sign-dist}, the derivative of $\signdist^{F}$ along the $(n,n)$-component is uniformly positive at $x_0$, 
    \[\langle D\signdist^{F}(u_{\tau_i,\tau_j}(x_0)), e_ne_n^T\rangle \ge \Lambda^{-1}>0. \]
    So, by the implicit function theorem, there exists a neighborhood $U_3$ of $x_0$ and a smooth function $h$ such that
    \[u_{\tau_n\tau_n}(x)  = h\left( u_{\tau_i\tau_j}(x):i,j=1,\dots, n,\ (i,j)\neq (n,n)\right), \qquad \forall\ x\in \overline{\Omega} \cap U_3,\]
    with derivatives of $h$ bounded in terms of $\Lambda$. Hence, by \eqref{eq:Subeq-Reg-high-reg-bdry-schauder-step-2-implicit-eq1}, there exists $C_3''>0$ depending on $C_3'$ and $\Lambda$ such that
    \begin{equation}
        \label{eq:Subeq-Reg-high-reg-bdry-schauder-step-2-implicit-eq2}
        \|u_{\tau_n\tau_n}\|_{C^{k-2,\alpha}( \overline{\Omega} \cap U_3)} \le C_3''.
    \end{equation}
    Then, \eqref{eq:Subeq-Reg-high-reg-bdry-schauder-step2-goal} follows from \eqref{eq:Subeq-Reg-high-reg-bdry-schauder-step-2-implicit-eq1} and \eqref{eq:Subeq-Reg-high-reg-bdry-schauder-step-2-implicit-eq2}.

    \vspace{0.5em}\noindent\textit{Step 3.} \textit{There exists a constant $C_5>0$ depending on $C_2$ \eqref{eq:Subeq-Reg-higher-reg-boundary-induction-hypothesis}, $\Lambda$ \eqref{eq:Subeq-Reg-higher-reg-unif-ellipticity-of-sign-dist}, $n$, $F$, $k$, $\alpha$, $\Omega$ and $\|\varphi\|_{C^{k,\alpha}}$ such that}
    \begin{equation}
        \label{eq:Subeq-Reg-high-reg-bdry-schauder-step3-goal}
        \|u\|_{C^{k,\alpha}(\overline{\Omega})} \le C_5.
    \end{equation}

    \vspace{0.5em}\noindent Cover $\partial \Omega$ by the neighborhoods $U_{x_0}$ from Step 2 as $x_0 $ varies in $ \partial \Omega$. By compactness of $\partial{\Omega}$, there exists finitely many points $x_1,\dots, x_{N} \in \partial \Omega$ and corresponding neighborhoods $U_{x_1},\dots, U_{x_N}$. By \eqref{eq:Subeq-Reg-high-reg-bdry-schauder-step2-goal},
    \begin{equation}
        \label{eq:Subeq-Reg-high-reg-bdry-schauder-step3-eq1}
        \|D^2u\|_{C^{k-2,\alpha}(\overline{\Omega}\cap U_{x_i})} \le C_4,\qquad, i=1,\dots, N. 
    \end{equation}
    Now pick a smooth cutoff function $\eta$ that is supported in the interior of the domain $\Omega$, and identically 1 outside the union of the neighborhoods $U_{x_1},\dots, U_{x_N}$, i.e.,
    \[\eta\in C^\infty(\Omega), \qquad \eta\ge 0,\qquad \eta(x) = 1, \quad \forall\ x\in \Omega \setminus \bigcup_{i=1}^N U_{x_i},  \qquad \operatorname{supp} \eta := \overline{\{\eta \neq 0\}} \subset \Omega.\]
    And consider the vector fields
    \[\tau_i := \eta\, e_i, \qquad i=1,\dots, n-1,\]
    where $e_1,\dots, e_n$ are the standard basis vectors. Then, all of the $\tau_i$ satisfy \eqref{eq:Subeq-Reg-high-reg-bdry-schauder-step1-assump}, and by Step 1 \eqref{eq:Subeq-Reg-high-reg-bdry-schauder-step1-goal} there exists a constant $C_4'>0$ such that,
    \begin{equation}
        \label{eq:Subeq-Reg-high-reg-bdry-schauder-step3-eq2}
        \|\partial_i u\|_{C^{k-1,\alpha} \left(\Omega \setminus \bigcup_{i=1}^N U_{x_i} \right)} \le C_4', \qquad i=1,\dots, n. 
    \end{equation}
    Combining the induction hypothesis \eqref{eq:Subeq-Reg-higher-reg-boundary-induction-hypothesis} with \eqref{eq:Subeq-Reg-high-reg-bdry-schauder-step3-eq1}--\eqref{eq:Subeq-Reg-high-reg-bdry-schauder-step3-eq2}, \eqref{eq:Subeq-Reg-high-reg-bdry-schauder-step3-goal} follows. In other words, the induction step of the argument is completed. Hence, Proposition~\ref{prop:Subeq-Reg-higher-reg-boundary-schauder} follows by induction.
\end{proof}    
\end{allproofs}


\subsection{Continuity method}
\label{section:continuity-method}

Corollary~\ref{corollary:Subeq-Reg-existence-thm} is proved in this section via the continuity method.

\begin{proof}[Proof of Corollary~\ref{corollary:Subeq-Reg-existence-thm}]
    Fix a matrix $A\in \partial F$ and consider the quadratic
    \begin{equation}
        \label{eq:Subeq-Reg-existence-QA-def}
        Q_A(x) : = \frac{1}{2} \langle Ax, x\rangle, \qquad \forall\ x\in \overline{\Omega}.
    \end{equation}
    Let
    \[\varphi_t := (1-t) Q_A + t \varphi, \qquad t\in[0,1],\]
    and consider the family of Dirichlet problems for $t\in [0,1]$,
    \begin{equation}
        \label{eq:Subeq-Reg-existence-cty-path-1}
        \begin{aligned}
            D^2u_t \in \partial F\quad &\text{on } \Omega,\\
            u_t = \varphi_t \quad &\text{on } \partial \Omega.
        \end{aligned}
    \end{equation}
    First, the path \eqref{eq:Subeq-Reg-existence-cty-path-1} can be rephrased in terms of the signed-distance $\signdist^{F}$ \eqref{eq:Subeq-Reg-signed-dist-def},
    \begin{equation}
        \label{eq:Subeq-Reg-existence-cty-path-2}
        \begin{aligned}
            \signdist^{F}(D^2u_t) =&\ 0\quad \text{on } \Omega,\\
            u_t =&\ \varphi_t \quad \text{on } \partial \Omega.
        \end{aligned}
    \end{equation}
    Furthermore, under the change of variables $v_t = u_t - \varphi_t$, \eqref{eq:Subeq-Reg-existence-cty-path-2} is equivalent
    \begin{equation}
        \label{eq:Subeq-Reg-existence-cty-path-3}
        \begin{aligned}
            \signdist^{F}(D^2v_t + \varphi_t) =&\ 0\quad \text{on } \Omega,\\
            v_t =&\ 0 \quad \text{on } \partial \Omega.
        \end{aligned}
    \end{equation}

    Observe that $\|\varphi_t\|_{C^k}$ is uniformly bounded in terms of $\|Q_A\|_{C^k}$ and $\| \varphi\|_{C^k}$, for any $k\ge 0$. So by Theorem~\ref{thm:Subeq-Reg-a-priori-estimates}, Corollary~\ref{corollary:Subeq-Reg-Evans-Krylov-theorem} and Proposition~\ref{prop:Subeq-Reg-higher-reg-boundary-schauder}, there exist $\alpha \in(0,1)$ and $C_1>0$ depending only on $n$, $F$, $\Omega$, $\|\varphi\|_{C^5}$ and $\|Q_A\|_{C^5}$, such that for any solution $u_t$ to \eqref{eq:Subeq-Reg-existence-cty-path-1},
    \begin{equation}
        \label{eq:Subeq-Reg-cty-method-unif-a-priori-est}
        \|u_t\|_{C^{4,\alpha}(\overline{\Omega})} \le C_1, \qquad t\in [0,1].
    \end{equation}
    Consider
    \begin{equation}
        \label{eq:Subeq-Reg-cty-ical-def}
        \ical := \left\{ t\in [0,1]: \text{exists} \quad u_t \in C^{4,\alpha}(\overline{\Omega})  \quad \text{satisfying \eqref{eq:Subeq-Reg-existence-cty-path-1}} \right\},
    \end{equation}
    The goal is to prove that $1\in \ical$.
    
    \vspace{0.5em}\noindent\textit{Non-empty:} By \eqref{eq:Subeq-Reg-existence-QA-def}, $D^2Q_A = A \in \partial F$ and $Q_A = \varphi_0$. Hence, $u_0 = Q_A$ solves the \eqref{eq:Subeq-Reg-existence-cty-path-1} and so $0\in \ical \neq \varnothing$.

    \vspace{0.5em}\textit{$\ical$ is open in $[0,1]$:} Let $t_0 \in \ical$. By \eqref{eq:Subeq-Reg-cty-ical-def}, \eqref{eq:Subeq-Reg-existence-cty-path-1} has solution $u_{t_0} \in C^{4,\alpha}(\overline{\Omega})$ at $t=t_0$. Then
    \[v_{t_0} := u_{t_0} + \varphi_{t_0}\]
    is a solution to \eqref{eq:Subeq-Reg-existence-cty-path-3} at $t=t_0$. Consider the Banach space $C^{4,\alpha}_0(\overline{\Omega})$ of $C^{4,\alpha}$ functions vanishing on $\partial \Omega$, and let
    \[\fcal: [0,1]\times C^{4,\alpha}_0(\overline{\Omega}) \rightarrow C^{2,\alpha}(\overline{\Omega}),\qquad \fcal(t,v) := \signdist^{F}(D^2v + D^2\varphi_t).\]
   Since $\signdist^{F}$ \eqref{eq:Subeq-Reg-signed-dist-def} is smooth on a neighborhood of $\partial F$, $\fcal$ is only well-defined on a neighborhood of the given pair $(t_0,v_{t_0})$. The derivative of $\fcal$ at $(t_0,v_{t_0})$ is
    \[\begin{aligned}
        d_{(t_0,v_{t_0})}\fcal(t,v) :=&\ \left. \frac{d}{ds} \right|_{s=0} \fcal(t_0+st, v_{t_0}+sv) \\
        =&\ \left. \frac{d}{ds} \right|_{s=0} \signdist^{F}(D^2v_{t_0}+sD^2v + D^2\varphi_{t_0+st}) \\
        =&\ \left\langle D\signdist^{F}(D^2v_{t_0}- D^2\varphi_{t_0}),\ D^2v + D^2\left. \frac{d}{ds} \right|_{s=0} \varphi_{t_0+st} \right\rangle, \\
        =&\ \left\langle D\signdist^{F}(D^2u_{t_0}),\ D^2v + D^2\left. \frac{d}{ds} \right|_{s=0} \varphi_{t_0+st} \right\rangle \\
        =&\ \Delta_{ D\signdist^{F}(D^2u_{t_0})} \left(v+\left. \frac{d}{ds} \right|_{s=0} \varphi_{t_0+st} \right),
    \end{aligned} \]
    in terms of the linearization \eqref{eq:Subeq-Reg-linearization-of-dij-elliptic-pde}. Since $\Delta_{ D\signdist^{F}(D^2u_{t_0})}$ is linear and strictly elliptic (Claim~\ref{claim:Subeq-Reg-dij-is-elliptic-PDE}), the Dirichlet problem
    \[\begin{aligned}
        \Delta_{ D\signdist^{F}(D^2u_{t_0})} w =  \psi& \qquad \text{on } \Omega,\\
         w = 0& \qquad \text{on } \partial \Omega,
    \end{aligned}\]
    has a unique solution $v\in C^{4,\alpha}(\overline{\Omega})$ for any $\psi \in C^{2,\alpha}(\overline{\Omega})$. In other words, the partial derivative of $\fcal(t,v)$ along $v$,
    \[w\mapsto d_{(t_0,v_{t_0})}\fcal (0, w) = \Delta_{ D\signdist^{F}(D^2u_{t_0})} w,\qquad w\in C^{4,\alpha}_0 (\overline{\Omega}), \]
    is an isomorphism from $C^{4,\alpha}_0 (\overline{\Omega})$ to $C^{2,\alpha}(\overline{\Omega})$. By the Banach space valued implicit function theorem, there exists $\epsilon>0$ a smooth function
    \[v:(t_0-\epsilon, t_0 + \epsilon)\cap [0,1] \rightarrow C^{4,\alpha}_0(\overline{\Omega}),\]
    such that
    \[\fcal(t,v(t)) = \signdist^{F}(D^2v(t) + D^2\varphi_t) = 0, \qquad \forall\ t\in (t_0-\epsilon, t_0+\epsilon)\cap [0,1],\]
    i.e., $v(t)$ satisfies \eqref{eq:Subeq-Reg-existence-cty-path-3}. Then,
    \[u_t:= v_t+ \varphi_t \quad \text{solves}\quad \eqref{eq:Subeq-Reg-existence-cty-path-1}, \qquad t\in(t_0-\epsilon, t_0+\epsilon) \cap [0,1].\]
    By \eqref{eq:Subeq-Reg-cty-ical-def}, $(t_0-\epsilon, t_0+\epsilon)\cap [0,1] \subset \ical$, and so $\ical$ is open.

    \vspace{0.5em}\noindent\textit{Step 3:} \textit{$\ical$ is closed:} Consider a sequence $\ical \ni t_k \rightarrow t_\infty \in [0,1]$. By \eqref{eq:Subeq-Reg-cty-ical-def},
    \[\text{exists} \quad u_{t_k} \in C^{4,\alpha}(\overline{\Omega}) \quad \text{solving}\quad \eqref{eq:Subeq-Reg-existence-cty-path-1}\quad \text{at}\quad t=t_k,\qquad k=1,2,\dots.\]
    By \eqref{eq:Subeq-Reg-cty-method-unif-a-priori-est}, $u_{t_k}$ is uniformly bounded in $C^{4,\alpha} (\overline{\Omega})$. So by the Arzela--Ascoli theorem, replacing the sequence $u_{t_k}$ by a subsequence if necessary,
    \[u_{t_\infty} := \lim_{k\rightarrow \infty} u_{t_{k}} \in C^{4,\alpha} (\overline{\Omega}). \]
    Since $u_{t_k} \rightarrow u_{t_\infty}$ in $C^2$, taking a limit of \eqref{eq:Subeq-Reg-existence-cty-path-1} along $t_k$, it follows that $u_{t_\infty} \in C^{4,\alpha}(\overline{\Omega})$ is a solution to \eqref{eq:Subeq-Reg-existence-cty-path-1} at $t=t_\infty$. Thus $t_\infty \in \ical$, and thus $\ical$ is closed.

    \vspace{0.5em}\noindent Since $\ical$ is a non-empty, open and closed subset of $[0,1]$, it follows that $\ical = [0,1]$, and in particular $1\in \ical$. In other words, by \eqref{eq:Subeq-Reg-cty-ical-def}, \eqref{eq:Subeq-Reg-const-coeff-subequation-PDE} admits a solution in $C^{4,\alpha}(\overline{\Omega})$.
\end{proof}

\section{Applications}

In the spirit of Harvey--Lawson's masterful exposition \cite[\S10]{HarLaw-dirichletduality}, we list
a selection of equations to which Theorem~\ref{thm:Subeq-Reg-a-priori-estimates} and Corollary~\ref{corollary:Subeq-Reg-existence-thm} apply, and ``...readers are encouraged
to look at examples close to their interests and bear them in mind while reading other 
parts of the paper.'' \cite[p. 430]{HarLaw-dirichletduality}

\subsection{Monge--Amp\`ere and special Lagrangian equations}

Theorem~\ref{thm:Subeq-Reg-a-priori-estimates} applies to the top branch of the inhomogeneous real, complex, and quaternionic Monge--Amp\`ere equations, corresponding respectively to the subequations
\begin{equation}
    \label{eq:Subeq-Reg-MA-CMA-QMA-subequations}
    \begin{aligned}
    F_{\operatorname{MA}} =&\ \left\{ A\in \Sym^2(\RR^n): A\ge 0,\quad \det A \ge 1 \right\} \\
    F_{\operatorname{CMA}} =&\ \left\{ A\in \Sym^2(\RR^{2m}): A^{\Herm}\ge 0,\quad \det A^{\Herm} \ge 1 \right\}\\
    F_{\operatorname{QMA}} =&\ \left\{ A\in \Sym^2(\RR^{4m}): A^{\QHerm} \ge 0,\quad \det A^{\QHerm} \ge 1 \right\},
\end{aligned}
\end{equation}
where $A^{\Herm}$ is the \textit{Hermitian part},
\begin{equation}
    \label{eq:Subeq-Reg-hermitian-part-of-a-matrix}
    A^{\Herm} := \frac{A- J AJ}{2},
\end{equation}
and $A^{\QHerm}$ is the \textit{quaternionic Hermitian part}
\begin{equation}
    \label{eq:Subeq-Reg-quaternionic-hermitian-part-of-a-matrix}
    A^{\QHerm} := \frac{A-IAI-JAJ-KAK}{4}.
\end{equation}

$F_{\operatorname{MA}}$, $F_{\operatorname{CMA}}$ and $F_{\operatorname{QMA}}$ are all smooth and non-degenerate subequations by Lemma~\ref{lemma:Subeq-Reg-subequation-is-elliptic}~(ii), since $\det$ is smooth and has positive definite derivative, $D \det >0$ on $\intr \pcal$ \eqref{eq:Subeq-Reg-positive-semi-def-subeq}. Furthermore, they are convex, since $\log \det$ is concave, and are respectively $O(n)$-, $U(m)$- and $Sp(m)$-invariant. Hence, Theorem~\ref{thm:Subeq-Reg-a-priori-estimates} and Corollary~\ref{corollary:Subeq-Reg-existence-thm} apply to yield existence of smooth solutions to the Dirichlet boundary value problem on a smooth domain $\Omega$ that is uniformly convex, uniformly pseudoconvex, or uniformly quaternionic pseudoconvex respectively.

Similarly, the subequations corresponding to the special Lagrangian equation are
\begin{equation}
    \label{eq:Subeq-Reg-subequations-for-SL}
    F_c:= \left\{A\in \Sym^2(\RR^n): f(\lambda(A)) := \sum_{i=1}^n \tan^{-1} \lambda_i(A) \ge c\right\}, \qquad c\in \left( - \frac{n\pi}{2}, \frac{n\pi}{2}\right),
\end{equation}
where $\lambda(A) = (\lambda_1(A),\dots, \lambda_n(A))$ is the vector eigenvalues of $A$, and we use the odd branch $\tan^{-1}:\RR \rightarrow \left( -\frac{\pi}{2}, \frac{\pi}{2} \right)$. First, $F_c$ \eqref{eq:Subeq-Reg-subequations-for-SL} are smooth and non-degenerate by Lemma~\ref{lemma:Subeq-Reg-subequation-is-elliptic}~(ii) since
\[\frac{\partial f}{\partial \lambda_i} = \frac{1}{1+\lambda_i^2}>0, \qquad i=1,\dots, n.\]
Furthermore, $F_c$ is $O(n)$-invariant, and convex for $c\in [\frac{n-2}{2}\pi, \frac{n}{2}\pi)$ \cite[Lemma~8.1]{caffarelli-nirenberg-spruck-dirichlet-problem-III-MR0806416} \cite[Lemma~2.1]{yuan-global-solutions-MR2199179}. Hence, Theorem~\ref{thm:Subeq-Reg-a-priori-estimates} and Corollary~\ref{corollary:Subeq-Reg-existence-thm} yield existence of smooth solutions on strictly $\vec{F}_c$-convex domains $\Omega$ (see Harvey--Lawson~\cite[Theorem~4.6]{Har-Law-pseudoconvex-MR4179860}).

Analogous to \eqref{eq:Subeq-Reg-MA-CMA-QMA-subequations}, one can also define the complex and quaternionic special Lagrangian subequations, by replacing $\lambda(A)$ in \eqref{eq:Subeq-Reg-subequations-for-SL} by either the eigenvalues of the Hermitian part, $\lambda(A^{\Herm})$ when $n=2m$ \eqref{eq:Subeq-Reg-hermitian-part-of-a-matrix}, or the quaternionic Hermitian part $\lambda(A^{\QHerm})$ when $n=4m$ \eqref{eq:Subeq-Reg-quaternionic-hermitian-part-of-a-matrix}. These complex and quaternionic analogs are also similarly smooth, non-degenerate and convex for $c\in [\frac{n-2}{2}\pi, \frac{n}{2}\pi)$, and respectively $U(m)$- and $Sp(m)$-invariant. The complex case was studied by Collins--Picard--Wu~\cite[Theorem~1.2]{CPW-MR3659638}, who proved existence of solutions to the BVP \eqref{eq:Subeq-Reg-const-coeff-subequation-PDE}. Theorem~\ref{thm:Subeq-Reg-a-priori-estimates} and Corollary~\ref{corollary:Subeq-Reg-existence-thm} yield existence of smooth solutions to the BVP \eqref{eq:Subeq-Reg-const-coeff-subequation-PDE} in both cases.

The following examples illustrate why the hypotheses of convexity and non-degeneracy are needed in Theorem~\ref{thm:Subeq-Reg-a-priori-estimates} (see also Remark~\ref{remark:Subeq-Reg-why-are-results-optimal}):

\begin{figure}[t]
    \centering
    \begin{subfigure}[t]{0.45\textwidth}
        \centering
        \begin{tikzpicture}[scale=3]
            \fill[white!80!black] (1,0) -- (0,0) --(0,1) -- (1,1);
            \draw[thin, <->, white!70!black] (-0.3,0) -- (1,0) node[below] {$\lambda_1$};
            \draw[thin, <->, white!60!black] (0,-0.3) -- (0,1) node[left] {$\lambda_2$};
            \draw (1,0) -- (0,0) -- (0,1);
            \draw (0.3,0.3) node {$\pcal$};
        \end{tikzpicture}
        \caption{}
    \end{subfigure}
    \hfill
    \begin{subfigure}[t]{0.45\textwidth}
        \centering
        \begin{tikzpicture}[scale=3]
            \fill[white!80!black] (1,0) -- (0.3,0) .. controls (0.1,0) and (0,0.1) .. (0,0.3) -- (0,1) -- (1,1);
            \draw[thin, <->, white!70!black] (-0.3,0) -- (1,0) node[below] {$\lambda_1$};
            \draw[thin, <->, white!60!black] (0,-0.3) -- (0,1) node[left] {$\lambda_2$};
            \draw (1,0) -- (0.3,0) .. controls (0.1,0) and (0,0.1) .. (0,0.3) -- (0,1);
            \draw (0.3,0.3) node {$\pcal'$};
        \end{tikzpicture}
        \caption{}
    \end{subfigure}
    \caption{(a)The subequation $\pcal \subset \Sym^2(\RR^2)$ \eqref{eq:Subeq-Reg-positive-semi-def-subeq}, plotted in the space of eigenvalues. (b) The subequation $\pcal'$ obtained by smoothing $\pcal$ near the origin admits a weak solution $u\notin C^{1,1}(\overline{\Omega})$ \eqref{eq:Subeq-Reg-CNS-counter-example-to-HRMA}.}
    \label{fig:Subeq-Reg-smoothing-of-pcal-near-the-origin}
\end{figure}
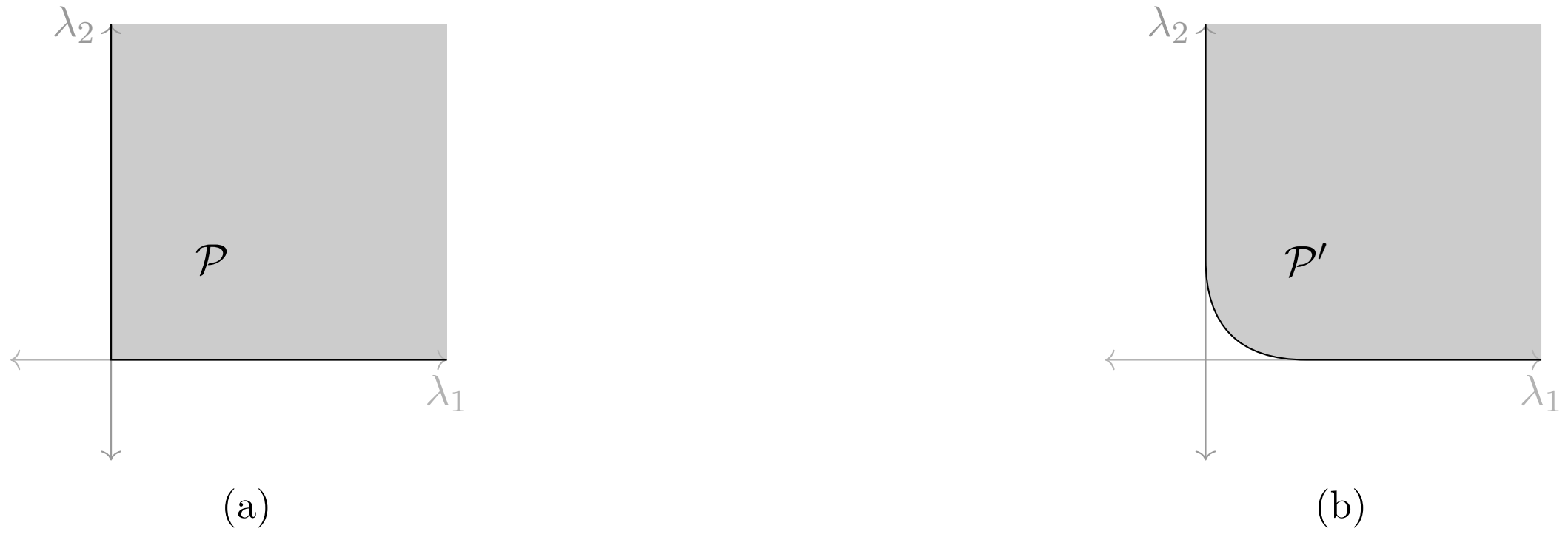

\begin{example}~
\label{example:Subeq-Reg-examples-for-optimality}
\begin{enumerate}[label={(\roman*)}]
    \item (Non-convex) When $|c|< (n-2)\frac{\pi}{2}$, $F_c$ \eqref{eq:Subeq-Reg-subequations-for-SL} is smooth, non-degenerate and $O(n)$-invariant but is \textit{non-convex}. The $F_c$-BVP \eqref{eq:Subeq-Reg-const-coeff-subequation-PDE} admits weak solutions that are not $C^2$, due to the work of Nadirashvili--Vl\u{a}du\c{t}~\cite{nadirashvili-Vladut-MR2683755}, with generalizations by Wang--Yuan, and Mooney--Savin~\cite{wang-yuan-singular-solutions-MR3117304, Connor-Savin-MR4836007}.


    \item (Degenerate) The BVP \eqref{eq:Subeq-Reg-intro-purely-second-order-PDE} for the homogeneous real Monge--Amp\`ere equation (HRMA) is the PDE
    \begin{equation}
        \label{eq:Subeq-Reg-homogeneous-RMA-intro}
        \begin{aligned}
        u\quad \text{convex},\qquad &\\
        \det(D^2u) = 0 \qquad& \text{on } \Omega,\\
        u = \varphi \qquad& \text{on } \partial \Omega,
    \end{aligned}
    \end{equation}
    and is equivalent to the BVP \eqref{eq:Subeq-Reg-const-coeff-subequation-PDE} for $F = \pcal$ \eqref{eq:Subeq-Reg-positive-semi-def-subeq}. It turns out \eqref{eq:Subeq-Reg-homogeneous-RMA-intro} admits non-classical solutions. The original idea for non-classical solutions to the HRMA goes back to Hartman--Nirenberg~\cite{Hartman-Nirenberg-MR126812} who used the method of characteristics to construct solutions that are not $C^2$. This idea can also be generalized to the (degenerate) homogeneous complex Monge--Amp\'ere equation~\cite{Bedford-Fornaess-MR517774} and  the (degenerate)  quaternionic Monge--Amp\`ere equation~\cite{Alesker-quaternionic-MA-MR1967025}. 

    Note that $\pcal$ is convex and $O(n)$-invariant but is neither non-degenerate nor smooth. In fact, $\pcal$ can be smoothed out while retaining the counter examples (see Figure~\ref{fig:Subeq-Reg-smoothing-of-pcal-near-the-origin}) \cite[Ex.~2]{CNS-degenerate-RMA-MR864651}:
    
    In dimension $n=2$, denote the unit disc by $\DD = \{x^2 + y^2 < 1\} \subset \RR^2$, and observe that,
    \begin{equation}
        \label{eq:Subeq-Reg-CNS-counter-example-to-HRMA}
        u(x, y) := (1+y)^{2-\epsilon}, \qquad \forall\ (x,y) \in \overline{\DD},
    \end{equation}
    solves \eqref{eq:Subeq-Reg-homogeneous-RMA-intro} with $\varphi = u|_{\partial \DD}$. Near $(0,-1) \in \partial \DD$, the boundary $\partial \DD$ can be expressed as the graph $y(x) = -\sqrt{1-x^2} \approx -1+ \frac{1}{2}x^2$, and
    \[\varphi = u|_{\partial \DD}(x) = (1+y(x))^{2-\epsilon} \approx \left(\frac{x^2}{2}\right)^{2-\epsilon}  = \frac{x^{4-2\epsilon}}{2^{2-\epsilon}} \in C^{3,1-2\epsilon}(\partial \DD). \]
    But by \eqref{eq:Subeq-Reg-CNS-counter-example-to-HRMA} $u \notin C^{1,1}(\overline{\DD})$.

    Furthermore, since $D^2u$ stays outside a small neighborhood of the origin in $\Sym^2(\RR^2)$, $\pcal$ can be smoothed out near the origin to obtain $\pcal'$ (Figure~\ref{fig:Subeq-Reg-smoothing-of-pcal-near-the-origin}), which is smooth, convex and $O(n)$-invariant, but not non-degenerate such that $u$ is also a solution to the $\pcal'$-BVP.
\end{enumerate}
\end{example}

\subsection{Equations determined by pairs of subequations}
\label{section:Subeq-Reg-equations-determined-by-pairs-of-subequations}

Harvey--Lawson~\cite{Har-Law-PDE-Krylov-MR4127162} generalized the Dirichlet problem \eqref{eq:Subeq-Reg-const-coeff-subequation-PDE} to
\begin{equation}
    \label{eq:Subeq-Reg-two-subequations-const-coeff-PDE-krylov-point-of-view}
    \begin{aligned}
        D^2u \in E_{F_1, F_2} := F_1\cap -\widetilde{F_2} \qquad &\text{on } \Omega,\\
        u = \varphi \qquad &\text{on } \partial \Omega,
    \end{aligned}
\end{equation}
where $F_1, F_2 \subset \Sym^2(\RR^n)$ are both subequations. The order of the two subequations matters, and $E_{F_2,F_1} $ is called \textit{the mirror} to the equation $E_{F_1, F_2}$. When $F_1=F_2 = F$, then the equation $E_{F_1, F_2}$ is equal to its mirror, and \eqref{eq:Subeq-Reg-two-subequations-const-coeff-PDE-krylov-point-of-view} reduces to \eqref{eq:Subeq-Reg-const-coeff-subequation-PDE}. In general, existence and uniqueness of solutions to \eqref{eq:Subeq-Reg-two-subequations-const-coeff-PDE-krylov-point-of-view} depends on whether the equation $E_{F_1, F_2}$ or its mirror $E_{F_2,F_1}$ has interior:

\begin{theorem}
\label{thm:Subeq-Reg-pair-of-subequations-harlaw-existence}
    {\emph{\cite[Theorem~2.6]{Har-Law-PDE-Krylov-MR4127162}}} For subequations $F_1, F_2\subset \Sym^2(\RR^n)$, $\varphi \in C^0(\partial \Omega)$ and $\Omega \subset \RR^n$ bounded and uniformly convex,
    \begin{enumerate}[label={(\roman*)}]
        \item (uniqueness) there is at most one weak solution $u\in C^0(\overline{D})$ to \eqref{eq:Subeq-Reg-two-subequations-const-coeff-PDE-krylov-point-of-view} $\quad\iff\quad$
        \[\intr E_{F_1, F_2} = \varnothing\quad\iff\quad E_{F_1, F_2} = \partial F_1 \cap \partial F_2 \quad\iff\quad F_1\subset F_2,\]
        \item (existence) there exists a weak solution $u\in C^0(\overline{D})$ to \eqref{eq:Subeq-Reg-two-subequations-const-coeff-PDE-krylov-point-of-view} $\quad\iff\quad$ $\intr E_{F_2,F_1} = \varnothing$, i.e., the mirror has empty interior.
    \end{enumerate}
\end{theorem}

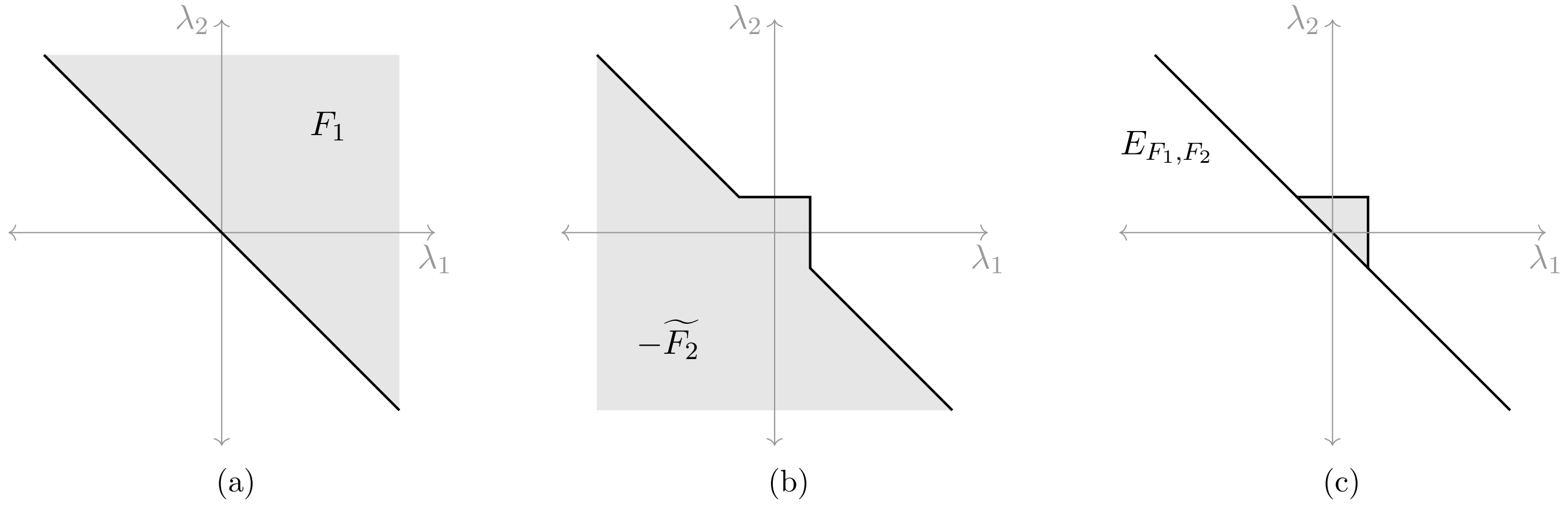
\begin{figure}[t]
    \centering
    \begin{subfigure}[t]{0.3\textwidth}
    \centering
        \begin{tikzpicture}[scale=0.8]
            \fill[white!90!black] (-2.5,2.5) -- (-0.5,0.5) -- (0.5,-0.5) -- (2.5,-2.5) -- (2.5,2.5) -- (-2.5,2.5);
            \draw[<->, white!60!black, thin] (-3,0) -- (3,0) node[below] {$\lambda_1$};
            \draw[<->, white!60!black, thin] (0,-3) -- (0,3) node[left] {$\lambda_2$};
            \draw[thick] (-2.5,2.5) -- (-0.5,0.5) -- (0.5,-0.5) -- (2.5,-2.5);
            \draw (1.5,1.5) node {$F_1$};
        \end{tikzpicture}
    \caption{}
    \end{subfigure}
    \hfill
    \begin{subfigure}[t]{0.3\textwidth}
    \centering
        \begin{tikzpicture}[scale=0.8]
            \fill[white!90!black] (-2.5,2.5) -- (-0.5,0.5) -- (0.5,0.5) -- (0.5,-0.5) -- (2.5,-2.5) -- (-2.5,-2.5) -- (-2.5,2.5);
            \draw[<->, white!60!black, thin] (-3,0) -- (3,0) node[below] {$\lambda_1$};
            \draw[<->, white!60!black, thin] (0,-3) -- (0,3) node[left] {$\lambda_2$};
            \draw[thick] (-2.5,2.5) -- (-0.5,0.5) -- (0.5,0.5) -- (0.5,-0.5) -- (2.5,-2.5);
            \draw (-1.5,-1.5) node {$-\widetilde{F_2}$};
        \end{tikzpicture}
    \caption{}
    \end{subfigure}
    \hfill
    \begin{subfigure}[t]{0.3\textwidth}
    \centering
        \begin{tikzpicture}[scale=0.8]
            \filldraw[thick, fill=white!90!black] (-0.5,0.5) -- (0.5,0.5) -- (0.5,-0.5) -- (-0.5,0.5);
            \draw[<->, white!60!black, thin] (-3,0) -- (3,0) node[below] {$\lambda_1$};
            \draw[<->, white!60!black, thin] (0,-3) -- (0,3) node[left] {$\lambda_2$};
            \draw[thick] (-2.5,2.5) -- (-0.5,0.5);
            \draw[thick] (0.5,-0.5) -- (2.5,-2.5);
            \draw (-1.5,1.6) node[anchor=north east] {$E_{F_1, F_2}$};
        \end{tikzpicture}
    \caption{}
    \end{subfigure}
    \caption{An example where Theorem~\ref{thm:Subeq-Reg-Har-Law-generalized-equations-regularity-part-2} applies plotted in the space of eigenvalues: (a) $F_1 = \{\tr A\ge 0\}$, (b) $-\widetilde{F_2} = \{ \tr A\le 0 \ \text{or}\ A\le I\}$ (c) $E_{F_1, F_2} = F_1\cap -\widetilde{F_2}$.}
    \label{fig:Subeq-Reg-existence-for-pair-of-subequations-case-1}
\end{figure}

Suppose there exists a strong solution $u\in C^4(\overline{\Omega})$ to \eqref{eq:Subeq-Reg-two-subequations-const-coeff-PDE-krylov-point-of-view}. We could try and obtain a priori estimates as done in \S\S\ref{sec:Subeq-Reg-unif-grad-estimates}--\ref{sec:Subeq-Reg-normal-hessian-estimates}. The difficulty is that, in general when $\intr E_{F_1,F_2} \neq \varnothing$ (i.e., in the absence of uniqueness), the maximum principle and Hessian estimates might fail at points $x_0\in \overline{\Omega}$ such that $D^2u(x_0) \in \intr E_{F_1, F_2}$. However, we can prove a priori estimates assuming the interior of the equation, $\intr E_{F_1,F_2}$, is bounded:

\begin{theorem}
    \label{thm:Subeq-Reg-Har-Law-generalized-equations-regularity-part-2}
    Suppose $\subset O(n)$ is a Lie group acting transitively on $S^{n-1}$, $F_1 \subset \Sym^2(\RR^n)$ is a convex, non-degenerate \eqref{eq:Subeq-Reg-non-degenerate-subequation}, $G$-invariant subequations with smooth boundary $\partial F_1$, $F_2\subset \Sym^2(\RR^n)$ is any subequation, $\varphi \in C^4(\partial \Omega)$ is a boundary condition, and $\Omega \subset \RR^n$ is a bounded uniformly convex domain with $C^4$ boundary $\partial \Omega$, and finally, suppose there exist $R >0$ such that $\intr E_{F_1,F_2}$ lies in the ball of radius $R$,
    \begin{equation}
        \label{eq:Subeq-Reg-assumption-on-pairs-interior-for-estimate}
        \intr E_{F_1, F_2}  \subset \{A\in \Sym^2(\RR^n): |A|\le R\}.
    \end{equation}
    Then any solution $u\in C^4( \overline{\Omega})$ to \eqref{eq:Subeq-Reg-two-subequations-const-coeff-PDE-krylov-point-of-view}, satisfies
    \[\|u\|_{C^2(\overline{\Omega})} \le C,\]
    for a constant $C = C(n,F_1, F_2,  \Omega, \|\varphi\|_{C^4}, R)$.
\end{theorem}
\begin{proof}[Sketch of proof]
    \noindent\textit{Step 1.} \textit{$C^0$ estimate:} A uniformly convex domain $\Omega$ is $\vec{F}$-convex for any subequation $F$ (Proposition~\ref{prop:Subeq-Reg-unif-convex-is-F-vec-for-all-F}). By Proposition~\ref{prop:Subeq-Reg-existence-of-delta-subsol-and-supersol}~(i) applied to the subequation $F_1\cap F_2$ and to $\widetilde{F_1}$ respectively, there exists an admissible $\delta$-subsolution to $F_1\cap F_2$ and an admissible $\delta$-supersolution to $F_1$,
    \[\begin{aligned}
        D^2 \underline{u} - \delta I \in F_1\cap F_2,\quad D^2\overline{u} + \delta I \notin F_1,\\
        \underline{u} = \overline{u} =\varphi \quad \text{on } \partial \Omega.
    \end{aligned}\]
    Since a solution $u$ to \eqref{eq:Subeq-Reg-two-subequations-const-coeff-PDE-krylov-point-of-view} is both an $F_1$-subsolution and a $F_2$-supersolution. So by the comparison principle~\cite[Theorem~6.4]{HarLaw-dirichletduality}, applied first to the $F_2$-supersolution $u$ and the $F_2$-subsolution $\underline{u}$, and second to the $F_1$-subsolution $u$ and $F_1$-supersolution $\overline{u}$,
    \begin{equation}
    \label{eq:Subeq-Reg-pairs-of-subeq-general-step1}
        \begin{aligned}
        \underline{u} \le&\ u \le  \overline{u}  \qquad \text{on } \Omega, \\
        \underline{u} =&\ u = \overline{u} = \varphi \quad \text{on } \partial \Omega.
    \end{aligned}
    \end{equation}
    In particular, a $C^0$ estimate holds.
    
    \vspace{0.5em}\noindent\textit{Step 2.} \textit{Interior gradient estimate:} The idea is to break the domain into two regions,
    \begin{equation}
        \label{eq:Subeq-Reg-estimates-for-pairs-two-domains}
        \Omega_1:= \{x \in \Omega: |D^2u(x)|<2R\} \qquad \text{and} \qquad \Omega_2:= \{x \in \Omega: |D^2u(x)|>2R\}.
    \end{equation}
    By the assumption \eqref{eq:Subeq-Reg-assumption-on-pairs-interior-for-estimate},
    \[D^2u(x)\notin \overline{\intr E_{F_1, F_2}}, \qquad \forall\ x\in \overline{\Omega_2}.\]
    Since $E_{F_1,F_2} \setminus \overline{\intr E_{F_1, F_2}} \subset \partial F_1 \cap \partial F_2$
    \[D^2u(x) \in \partial F_1\cap \partial F_2, \qquad \forall\ x\in \overline{\Omega_2}.\]
    In particular, $D^2 u \in \partial F_1$ on $\Omega_1$. So by the gradient maximum principle Lemma~\ref{lemma:Subeq-Reg-interior-grad-est},
    \begin{equation}
        \label{eq:Subeq-Reg-estimates-for-pairs-Step-2-eq1}
        \max_{\overline{\Omega_2}} |Du| \le C(n) \max_{\partial \Omega_2} |Du|.
    \end{equation}
    Moreover, by \eqref{eq:Subeq-Reg-assumption-on-pairs-interior-for-estimate} and \eqref{eq:Subeq-Reg-estimates-for-pairs-two-domains}, $D^2u$ is bounded by the constant $2R$ on $\Omega_1$. So by the fundamental theorem of calculus applied to $Du$ along line segments connecting a point of $\Omega_1$ to its boundary, since $\Omega$ is bounded, 
    \begin{equation}
        \label{eq:Subeq-Reg-estimates-for-pairs-Step-2-eq2}
        \max_{\overline{\Omega_1}}|Du| \le \max_{\partial \Omega_1}|Du| + 2R\operatorname{diam}(\Omega).
    \end{equation}
    By \eqref{eq:Subeq-Reg-estimates-for-pairs-Step-2-eq1}--\eqref{eq:Subeq-Reg-estimates-for-pairs-Step-2-eq2} the interior gradient estimate is reduced to a gradient estimate on
    \[\partial \Omega_1 \cup \partial \Omega_2 \subset \partial \Omega \cup \{x\in \Omega: |D^2u(x)| = 2R\}.\]
    To reduce to gradient estimate to the boundary $\partial \Omega$, we must deal with the region $|D^2u| = 2R$. Let
    \[v_e := \partial_e u + 3R\langle \,\cdot\,,e \rangle, \qquad \forall\ e\in S^{n-1}.\]
    Observe that
    \[ \partial_e v_e(x) = \partial_e\partial_e u(x) + 3R \ge -2R+3R =R>0, \qquad \text{if } x\in \Omega,\ |D^2u(x)|=2R,\]
    since $\langle A e,e\rangle \ge -2R$ when $|A|=2R$. Conclude that $v_e$ doesn't attain an interior maximum or minimum the set $\{x\in \Omega: |D^2u(x)|=2R\}$. Hence,
    \begin{align}
        |\partial_e u(x)| \le&\ \max_{e\in S^{n-1}} \left[\max_{x\in \overline{\Omega}} \langle x,e\rangle + \max_{x\in \partial \Omega}[\partial_eu(x) + \langle x,e\rangle ] \right] \nonumber\\
        \le&\ C(\Omega)(1+ \max_{\partial \Omega}|Du|), \qquad \forall\ e\in S^{n-1},\ x\in \Omega,\ |D^2u(x)| = 2R. \label{eq:Subeq-Reg-estimates-for-pairs-Step-2-eq3}
    \end{align}
    Combining \eqref{eq:Subeq-Reg-estimates-for-pairs-Step-2-eq1}--\eqref{eq:Subeq-Reg-estimates-for-pairs-Step-2-eq3},
    \[\max_{\overline{\Omega}}|Du| \le C(n,R,\Omega) (1+ \max_{\partial \Omega}|Du|).\]

    \vspace{0.5em}\noindent\textit{Step 3.} \textit{Boundary gradient estimate:} The tangential components of the gradient $Du$ are bounded by the boundary condition $u|_{\partial \Omega} =\varphi |_{\partial \Omega}$. And the normal component is bounded by $\max\{ \|\underline{u}\|_{C^1}, \|\overline{u}\|_{C^1}\}$ by \eqref{eq:Subeq-Reg-pairs-of-subeq-general-step1} and Claim~\ref{claim:Subeq-Reg-normal-gradient-claim}.
    
    \vspace{0.5em}\noindent\textit{Step 4.} \textit{Hessian estimate:} 
    A $C^2$ estimate trivially holds on $\overline{\Omega_1}$ \eqref{eq:Subeq-Reg-estimates-for-pairs-two-domains}. By the Hessian maximum principle applied to $F_1$ on on ${\Omega_2}$,
    \[\max_{\overline{\Omega_2}}|D^2u| \le C(1+\max_{\partial \Omega_2}|D^2u|).\]
    But note that a point $x\in \partial \Omega_2$ either satisfies $|D^2u(x)|=2R$ or that $x\in \partial \Omega$. Hence,
    \[\max_{\overline{\Omega_2}}|D^2u| \le C(1+\max_{\partial \Omega_2}|D^2u|)\le C(1+ 2R+\max_{\partial \Omega}|D^2u|).\]
    Finally, the Hessian estimate on $\partial \Omega$ follows the same argument presented in this article. The boundary Hessian estimate needs to only be proved at boundary points $x\in \partial \Omega$ where $|D^2u(x)|> 2R$. Note that $D^2u \in \partial F_1$ near any such $x$ and $\underline{u}$ \eqref{eq:Subeq-Reg-pairs-of-subeq-general-step1} is an admissible $\delta$-subsolution to $F_1$. So essentially the same argument from \S\ref{subsection:Subeq-Reg-tangential-hessian-estimate} and \S\S\ref{sec:Subeq-Reg-mixed-hessian-estimates}--\ref{sec:Subeq-Reg-normal-hessian-estimates} can be applied to $F_1$ to obtain estimates at such boundary point $x$.
    
    Note that, in \S\ref{sec:Subeq-Reg-normal-hessian-estimates}, if the point $x_{\min}$ (see \eqref{eq:Subeq-Reg-Trudingers-point}) lies on the set $\{x\in \partial \Omega: |D^2u(x)|\le 2R\}$, then we directly have a Hessian estimate at $x_{\min}$. So, we may assume in the whole argument that $|D^2u(x_{\min})| >2R$, in which case $D^2u(x) \in \partial F_1$ for $x$ near $x_{\min}$, and the argument of \S\ref{subsection:Subeq-Reg-C2-NN-trudingers-point} applies.
\end{proof}

\subsection{G\r{a}rding--Dirichlet polynomials}
\label{subsection:Subeq-Reg-Garding-Dirichlet-polynomials}

A special class of polynomials generalizing the determinant are the G\r{a}rding--Dirichlet (GD) polynomials, whose associated subequations generalize the Monge--Amp\`ere \eqref{eq:Subeq-Reg-MA-CMA-QMA-subequations} \cite{Har-Law-det-majorization-MR4589710, HL-duke-paper-MR4965192}.

A homogeneous polynomial $\fcal:\Sym^2(\RR^n)\rightarrow \RR$  of degree $N$ with
\[\fcal(I)>0\]
is \textit{$I$-G\r{a}rding hyperbolic} if,
\[\RR\ni t\mapsto \fcal(A+tI) \quad \text{has $N$ real roots }t_1,\dots, t_N,  \quad \forall\ A\in \Sym^2(\RR^n),\]
where $I$ is the identity matrix. Its G\r{a}rding cone is 
\[\Gamma_{\fcal} := \text{the connected component of } \left\{A \in \Sym^2(\RR^n): \fcal(A) >0 \right\} \text{ containing I}. \]
$\fcal$ is called \textit{G\r{a}rding--Dirichlet (GD)} if in addition,
\[ \intr \pcal = \left\{ A\in \Sym^2(\RR^n): A >0 \right\} \subset \Gamma_{\fcal}.\]
In this case, the closure $\overline{\Gamma}_{\fcal}$ is a subequation \cite[(3.1)]{Har-Law-det-majorization-MR4589710}. 
Finally, $\fcal$ is said to be \textit{complete} if
\[\overline{\Gamma}_{\fcal} \cap -\overline{\Gamma}_{\fcal}  \cap \pcal = \{0\}. \]

Harvey--Lawson studied complete G\r{a}rding--Dirichlet polynomials $\fcal$, and proved a determinant majorization inequality \cite[(1.1)]{HL-duke-paper-MR4965192}. They were inspired by Guo--Phong--Tong~\cite{Guo-Phong-Tong-MR4593734}, who showed earlier that for $U(m)$-invariant partial differential operators on the complex Hessian satisfying such a determinant majorization inequality, any solution admits an a priori $C^0$ estimate.

Our results assert existence of smooth solutions to the BVP \eqref{eq:Subeq-Reg-intro-purely-second-order-PDE}
\begin{equation}
    \label{eq:Subeq-Reg-BVP-for-GD-polynomials}
    \begin{aligned}
    \fcal(D^2u) = 1 \qquad& \text{on } \Omega,\\
    u = \varphi \qquad& \text{on } \partial \Omega,
\end{aligned}
\end{equation}
when $\fcal$ is a complete GD polynomial on $\Sym^2(\RR^n)$ invariant under a group $G$ acting transitively on $S^{n-1}$. Indeed in this case,
\begin{equation}
    \label{eq:Subeq-Reg-the-subequation-for-GD-PDE}
    F_{\fcal} :=  \left\{A \in \Gamma_{\fcal}: \fcal(A) \ge 1 \right\}
\end{equation}
is a smooth, convex, non-degenerate, and $G$-invariant subequation \cite[(3.2), Proposition~3.4]{Har-Law-det-majorization-MR4589710}. Hence, Theorem~\ref{thm:Subeq-Reg-a-priori-estimates} and Corollary~\ref{corollary:Subeq-Reg-existence-thm} applied to $F_{\fcal}$ to yield existence:

\begin{theorem}
\label{thm:Subeq-Reg-existence-for-GD-polynomials}
    For a complete GD polynomial $\fcal \in C^4( \Sym^2(\RR^n))$ that is invariant under a group $G$ acting transitively on $S^{n-1}$, 
    for $\varphi \in C^5(\overline{\Omega})$ and $\Omega \subset \RR^n$ bounded with $C^4$ and strictly $\vec{F_{\fcal}}$-convex boundary \eqref{eq:Subeq-Reg-the-subequation-for-GD-PDE}, the BVP~\eqref{eq:Subeq-Reg-BVP-for-GD-polynomials} admits a solution $u\in C^{4,\alpha}(\overline{\Omega})$.
\end{theorem}

\subsection{Complex and quaternionic G\r{a}rding--Dirichlet polynomials}
\label{subsec:Subeq-Reg-GD-thm-applied-to-CMA-and_QMA}

    Theorem~\ref{thm:Subeq-Reg-existence-for-GD-polynomials} also captures the BVP for invariant GD polynomials on the space of either Hermitian symmetric matrices $\Herm^2(\CC^m)$ or quaternionic Hermitian matrices $\QHerm^2(\HH^m)$. We illustrate the basic idea in the case of $\Herm^2(\CC^m)$.

    Recall that $\Herm^2(\CC^m)$ injects naturally as a subspace of $\Sym^2(\RR^{2m})$. Consider a GD polynomial $\fcal_{\CC}: \Herm^2(\CC^m)\rightarrow \RR$. Define
    \[\fcal_{\RR}:\Sym^2(\RR^{2m}) \rightarrow \RR, \qquad \fcal_{\RR}(A) := \fcal_{\CC}(A^{\Herm}).\]
    where $A^{\Herm}$ is the Hermitian part \eqref{eq:Subeq-Reg-hermitian-part-of-a-matrix}. Since $A^{\Herm}$ is a linear function in the entries of $A$, $\fcal_{\RR}$ is a real GD polynomial. Furthermore, since
    \[[D^2u]^{\Herm} = 2 \dbard u,\]
    in particular,
    \[\fcal_{\RR}(D^2u) = \fcal_{\CC}(\dbard u).\]
    And when $\fcal_{\CC}$ is $U(m)$-invariant, so is $\fcal_{\RR}$. Hence Theorem~\ref{thm:Subeq-Reg-existence-for-GD-polynomials} can be applied to $\fcal_{\RR}$ to solve the BVP
    \[\begin{aligned}
        \fcal_{\CC}(\dbard u) =&\ 1 \qquad \text{on } \Omega,\\
        u =&\ \varphi \qquad \text{on } \partial  \Omega.
    \end{aligned}\]
     The simplest example of a complex GD polynomial is the complex determinant $\fcal_{\CC} = \det: \Herm^2(\CC^m) \rightarrow \RR$, and so Theorem~\ref{thm:Subeq-Reg-existence-for-GD-polynomials} applies to the complex Monge--Amp\`ere equation.

     By similar arguments, Theorem~\ref{thm:Subeq-Reg-existence-for-GD-polynomials} also applies to the BVP~\eqref{eq:Subeq-Reg-BVP-for-GD-polynomials} for quaternionic GD polynomials $\fcal_{\HH}: \QHerm^2(\HH^m) \rightarrow \RR$, such as the quaternionic/Moore determinant $\det_{\HH}: \QHerm^2(\HH^n) \rightarrow \RR$ and the corresponding quaternionic Monge--Amp\`ere equation.

\printbibliography
\noindent\textsc{Department of Mathematics, University of Maryland}\\
\noindent\href{mailto:pvasanth@umd.edu}{\tt pvasanth@umd.edu}, \href{mailto:yanir@alum.mit.edu}{\tt yanir@alum.mit.edu}

\end{document}